\documentclass[a4paper,11pt]{article}

\usepackage[USenglish]{babel}
\usepackage[T1]{fontenc}
\usepackage[utf8]{inputenc}

\usepackage{graphicx,subfigure,tikz,pgfplots,epstopdf}

\usepackage{amsmath,amssymb,amsfonts,amsthm,array,geometry,stmaryrd,booktabs,paralist,todonotes,mathtools,pdflscape,siunitx,placeins}

\usepackage{ifpdf}
\ifpdf
\usepackage[pdftex,colorlinks=true,linkcolor=blue,citecolor=green,urlcolor=blue,bookmarks]{hyperref}
\else
\usepackage[dvipdfm,colorlinks=true,linkcolor=blue,citecolor=green,urlcolor=blue,bookmarks]{hyperref}
\fi 

\numberwithin{equation}{section}

\newcommand{\dd}{\textup{d}}
\newcommand{\bb}[1]{\boldsymbol{#1}}
\newcommand{\dbl}{\left\llbracket}
\newcommand{\dbr}{\right\rrbracket}

\newtheorem{thm}{Theorem}[section]

\newtheorem{rem}[thm]{Remark}

\begin{document}
    %
    \title{Exact and numerical solutions of the Riemann problem for a conservative model of compressible two-phase flows}
    \author{Ferdinand Thein\footnotemark[1], Evgeniy Romenski\footnotemark[2] \footnotemark[3], and Michael Dumbser\footnotemark[3]}
    \date{} 
    
    \maketitle
    \begin{abstract}
        In this work we study the solution of the Riemann problem for the barotropic version of the conservative symmetric hyperbolic and thermodynamically compatible (SHTC) two-phase flow model introduced in \cite{Romenski2007,Romenski2009}. All characteristic fields are carefully studied and explicit expressions are derived for the Riemann invariants and the Rankine-Hugoniot conditions.  
        Due to the presence of multiple characteristics in the system under consideration, non-standard wave phenomena can occur.
        Therefore we briefly review admissibility conditions for discontinuities and then discuss possible wave interactions.
        In particular we will show that overlapping rarefaction waves are possible and moreover we may have shocks that lie inside a rarefaction wave. 
        In contrast to nonconservative two phase flow models, such as the Baer-Nunziato system, we can use the advantage of the conservative form of the model under consideration.
        Furthermore, we show the relation between the considered conservative SHTC system and the corresponding barotropic version of the nonconservative Baer-Nunziato model.   
        Additionally, we derive the reduced four equation Kapila system for the case of instantaneous relaxation, which is the common limit system of both, the conservative SHTC model and the non-conservative Baer-Nunziato model. 
        Finally, we compare exact solutions of the Riemann problem with numerical results obtained for the conservative two-phase flow model under consideration, for the non-conservative Baer-Nunziato system and for the Kapila limit. 
        The examples underline the previous analysis of the different wave phenomena, as well as differences and similarities of the three systems.
    \end{abstract}
    \renewcommand{\thefootnote}{\fnsymbol{footnote}}
    \footnotetext[1]{Institute for Analysis and Numerics, Otto-von-Guericke University Magdeburg,\\ PSF 4120, D-39016 Magdeburg, Germany.
    \href{mailto:ferdinand.thein@ovgu.de}{\textit{ferdinand.thein@ovgu.de}}}
    \footnotetext[2]{Sobolev Institute of Mathematics, Novosibirsk, Russia
    \href{mailto:evrom@math.nsc.ru}{\textit{evrom@math.nsc.ru}}}
    \footnotetext[3]{Department of Civil, Environmental and Mechanical Engineering, University of Trento,\\ Via Mesiano, 77 - I-38123 Trento, Italy.
    \href{mailto:michael.dumbser@unitn.it}{\textit{michael.dumbser@unitn.it}}}
    \renewcommand{\thefootnote}{\arabic{footnote}}



    \section{Introduction}\label{sec:introduction}
%

The aim of the paper is to construct exact solutions for the Riemann problem of the one-dimensional barotropic version of the conservative Symmetric Hyperbolic and Thermodynamically Compatible (SHTC) model of compressible two-phase flows introduced in \cite{Romenski2007,Romenski2009}.
The results obtained can be useful for a qualitative understanding of the physical processes occurring in two-phase flows, for a comparative analysis of various models, and also for testing numerical methods. 

Note that the existence of a well-developed theory of exact solutions of the Riemann problem ensured the success in the development of modern shock-capturing numerical methods \cite{godunov} for solving the Euler equations of compressible gasdynamics, which are of fundamental importance in science and engineering, in particular aerospace engineering and astrophysics. See \cite{Toro2009} for an exhaustive overview of shock capturing schemes based on the exact or approximate solution of the Riemann problem.   
And even at present, the known exact solutions of the Euler equations are successfully used to test new numerical methods for solving hyperbolic systems of equations. 
Therefore, in our opinion, the construction of exact solutions for the equations of two-phase flows will influence the formation of a common point of view in the field of multiphase flow modeling and the development of new numerical methods for solving problems related to this area. 

In contrast to single-phase gas dynamics, there is still no universally accepted model of compressible multiphase flows, even for the case of only two phases, see, for example \cite{Saurel2017dilute}.
The generally accepted approach to design a two-phase flow model is based on the assumption that a mixture is a system of two interacting single phase continua, see, for example, \cite {Ishii2010}.
The most widely used PDE system is the Baer-Nunziato model \cite{Baer1986}, various modifications of which have been applied by many authors to study different types of flows, including flows with phase transitions and chemical reactions, see, for example recent papers \cite{Chiapolino2020,Furfaro2020,Saurel2017} and references therein. 

In this paper, we will study a model developed on the basis of the theory of SHTC systems, which was developed in \cite {Godrom2003,Godunov1996,Godunov:1995a,Rom1998,Romenski2001,Peshkov2018}.
This theory connects the local well-posedness of the governing partial differential equations of continuum physics (symmetric hyperbolicity in the sense of Friedrichs) with the fulfillment of the laws of thermodynamics (the law of conservation of energy and the law of increasing entropy). 
The general master system of SHTC equations can be derived from an underlying variational principle \cite{Peshkov2018}, and the study of numerous models of continuum mechanics has shown that their governing differential equations belong to the SHTC class of PDEs. 
This theory can be used to formulate new, well-posed models of processes in complex media, including unified models for the description of viscous Newtonian and non-Newtonian fluids and nonlinear elasto-plastic solids at large deformations  \cite{Peshkov2016,Dumbser2016,DPRZ2017,GPRNonNewtonian}, its extension to general relativity \cite{Romenski2020PTRSA}, SHTC models of rupture dynamics \cite{Tavelli2020,Gabriel2021PTRSA} and flows in deformable porous media  \cite{Romenski2020,Romenski2021AIP}.

The object of study in this paper is the SHTC model of compressible two-phase flow, the general equations of which are formulated and discussed in \cite{Romenski2004,Romenski2007,Romenski2009,Romenski2016}.
In the model under consideration, the two-phase medium is assumed to be a single continuum, the properties of which take into account the features of the two-phase flow. 
This means that the element of the medium is described by the average field of phase velocities and pressures for the mixture, but the flow of phases through this element is allowed. 
We consider the simplest model of a one-dimensional flow of a barotropic mixture. In this case, the SHTC equations can be written in a fully conservative form and therefore allow a direct formulation of discontinuous solutions. 
Nevertheless, the construction of exact solutions for the Riemann problem turns out to be a rather difficult task.

It is necessary to note that the governing equations of the Baer-Nunziato model \cite{Baer1986,Andrianov2003} in general can not be completely written in a conservative form, even in the one-dimensional case. This creates
difficulties in the definition of discontinuous solutions of the shock-wave type. But it turns out that the barotropic Baer-Nunziato system can be presented in a conservative form for one 
particular choice of interfacial velocity and interfacial pressure. And this conservative system is exactly the same as the barotropic case of the SHTC model studied in this paper.  
We discuss the similarities and differences between the different models by rewriting the equations of the SHTC model in the form of the Baer-Nunziato model for two interacting single phase continua.  
We also consider the reduced one-dimensional limit system obtained in the stiff pressure and velocity relaxation limit i.e. the resulting single-velocity single-pressure approximation of the SHTC and Baer-Nunziato models for the barotropic flow and show that they are identical and reduce to the well-known Kapila model \cite{Kapila2001}. 

The conservative form of the SHTC model of compressible two-phase flows has advantages not only for the construction of exact solutions of the Riemann problem, but also when using advanced shock capturing numerical methods. For the numerical simulations shown in this paper we will therefore rely on classical second order high resolution shock capturing TVD finite volume schemes as presented in \cite{Toro2009}. 

Degenerate behaviour in hyperbolic systems is a topic of constant interest. However, one has to be careful identifying the source of degeneracy.
Non-strictly hyperbolic systems, i.e.\ systems with multiple eigenvalues, naturally arise in multi dimensions \cite{Dafermos2016}.
Specific cases of non-strictly hyperbolic systems where studied in \cite{Keyfitz1980,Barkve1989,Li2015} and recently by Freist\"uhler and Pellhammer \cite{Freistuehler2020}.
There mostly systems of two equations where studied and the coincidence of eigenvalues often occurred in points where also the character of the related field changed from genuine nonlinear to
linearly degenerated. In particular Freist\"uhler showed in \cite{Freistuhler1991} that coinciding eigenvalues in the presence of a discontinuity must belong linearly degenerated fields.
Considering the Euler equations or the present system such a situation is related to the vanishing of the fundamental derivative $\mathcal{G}$.
This may happen for certain equations of state, but in general does not necessarily imply multiple eigenvalues.
A prominent example is the system of Euler equations and we highly recommend \cite{Wendroff1972,Menikoff1989,Mueller2006}.
Another reason for degeneracies may be the loss of an eigenvector as described in \cite{Sever2002,Choudhury2013}.
Systems with missing eigenvectors are often called \textit{weakly hyperbolic} and systems with coinciding eigenvalues are sometimes also called \emph{hyperbolic degenerate}. 
For further reading and a more detailed survey we recommend \cite{Chen2010}.
Coinciding eigenvalues may lead to difficult situations, but as long as there is a full set of eigenvectors which span the complete space the system is still diagonisable.
For the construction of \textit{complete Riemann solvers} a full set of eigenvectors is of great importance, see e.g. \cite{Roe81,OS82,DOT}.  
However, we want to emphasize that the consequences depend crucially on the system under consideration.
In some cases the numerical methods break down when applied to weakly hyperbolic systems, see e.g. \cite{SHTCSurfaceTension}, whereas in other situations a numerical treatment of weakly hyperbolic systems 
is still possible, see e.g. \cite{Matern2016,Hantke2020}.

%

The rest of the paper is organized as follows: in Section \ref{sec:model} we present the SHTC system under consideration in this paper and study its eigenstructure. In Section \ref{sec:degen_adm_cond} we briefly summarize degeneracies and admissibility conditions. The wave relations, i.e. the Riemann invariants and the Rankine-Hugoniot conditions of the model are presented in Section \ref{sec:waverel}, while the possible wave configurations are shown in Section \ref{sec:wave_config}. The relation of the conservative SHTC system with the non-conservative Baer-Nunziato model and the common Kapila limit are shown in Section \ref{sec:relmod}. Some examples of exact solutions and corresponding numerical results are presented in Section \ref{sec:num_res}. The paper closes with some concluding remarks and an outlook to future work in Section \ref{sec:conclusion}.

    \section{Conservative barotropic SHTC model for compressible two-phase flows}\label{sec:model}
\subsection{Multi-dimensional case}
The PDE system for compressible barotropic two-phase flows was discussed in Romenski et al.\ \cite{Romenski2007,Romenski2009}.
Written in terms of the specific energy $E = E(\alpha_1,c_1,\rho,w^k)$ it reads
\begin{subequations}\label{eqn.HPRFF}
    \begin{align}
        \frac{\partial \rho \alpha_1}{\partial t} + \frac{\partial \rho \alpha_1 u^k }{\partial x_k} &= \xi_1,\label{eqn.alphaFF}\\
        \frac{\partial \rho c_1}{\partial t} + \frac{\partial (\rho c_1 u^k+\rho E_{w_k})}{\partial x_k} &=  \xi_2,\label{eqn.contiFF1}\\
        \frac{\partial \rho}{\partial t} + \frac{\partial \rho u^k}{\partial x_k} &= \xi_3,\label{eqn.contiFF}\\
        \frac{\partial \rho u^i}{\partial t} + \frac{\partial (\rho u^i u^k + p \delta_{ik} + \rho w^iE_{w^k} )}{\partial x_k} & = \xi_4,\label{eqn.momentumFF}\\
        \frac{\partial w^k}{\partial t} + \frac{\partial(w^lu^l+E_{c_1})}{\partial x_k} + u^l\left(\frac{\partial w^k}{\partial x_l} - \frac{\partial w^l}{\partial x_k}\right)
        &= \xi_5.\label{eqn.relvel}
        %
    \end{align}
\end{subequations}
Here, $\alpha_1$ is the volume fraction of the first phase, which is connected with the volume fraction of the second phase $\alpha_2$ by the saturation law $\alpha_1+\alpha_2 = 1$,
$\rho$ is the mixture mass density, which is connected with the phase mass densities $\rho_1,\rho_2$ by the relation $\rho = \alpha_1\rho_1 + \alpha_2\rho_2$.
The phase mass fractions are defined as $c_1=\alpha_1 \rho_1/\rho,\, c_2=\alpha_2 \rho_2/\rho$ and it is easy to see that $c_1 + c_2 = 1$.
The mixture velocity is given by $u^i = c_1u_1^i + c_2u_2^i$ and $w^i=u_1^i - u_2^i$ is the relative phase velocity.
The equations describe the balance law for the volume fraction, the balance law for the mass fraction, the conservation of total mass, the total momentum balance law
and the balance for the relative velocity. 
The phase interaction is present via algebraic source terms in \eqref{eqn.alphaFF} and \eqref{eqn.relvel}, which are proportional to thermodynamic forces. 
These source terms are phase pressure relaxation to the common value
\begin{align}
    \xi_1 = -\rho{\phi}/{\theta_1} = -\rho{E_{\alpha_1}}/{\theta_1}\label{pressure_relax}
\end{align}
and interfacial friction
\begin{align}
    \xi_5 = -{\lambda_{k} }/{\theta_2} = -{E_{w^k}}/{\theta_2}.\label{velocity_relax}
\end{align}
The coefficients $\theta_1,\theta_2$ characterize the rate of pressure and velocity relaxation and can depend on parameters of state.
Due to mass and momentum conservation throughout this work we assume $\xi_2 = \xi_3 = \xi_4 = 0$.
\subsection{Discussion of the mixture equation of state}
We now want to further specify the derivatives of the generalized energy $E$. 
Due to the relation for the mass fractions and the concentrations we write $\alpha \equiv \alpha_1$, $\alpha_2 = 1 - \alpha$ and $c \equiv c_1$, $c_2 = 1 - c$, when appropriate.
The mixture equation of state (EOS) is defined as the sum of the mass averaged phase equations of state and the kinematic energy of relative motion
\begin{align}
    E(\alpha, c, \rho, w_1, w_2, w_3) &= e(\alpha,c,\rho) + c_1c_2\frac{w_iw^i}{2} = e(\alpha,c,\rho) + c(1 - c)\frac{w_iw^i}{2},\label{eq:mix_eos}\\
    e(\alpha,c,\rho) &= c_1e_1(\rho_1) + c_2e_2(\rho_2) = ce_1\left(\frac{c\rho}{\alpha}\right) + (1 - c)e_2\left(\frac{(1-c)\rho}{1 - \alpha}\right).\label{eq:mix_int_e}
\end{align}
where $e_i(\rho_i)$ is the specific internal energy of the $i$-th phase and is assumed to be known.
Special attention has to be paid when the derivatives of the internal energies are needed. 
\begin{rem}[isentropic vs.\ isothermal]
    When we want to calculate the derivative of the internal energy with respect to the density in the barotropic case,
    we have to specify whether we are considering an \emph{isentropic} or an \emph{isothermal} process.
    If we keep the entropy constant we have the well known derivative
    \begin{align}
        \left(\frac{\partial e}{\partial\rho}\right)_s = \frac{p}{\rho^2}.\label{e_deriv:iso_S}
    \end{align}
    However, if the temperature is held constant the derivative is given by
    \begin{align}
        \left(\frac{\partial e}{\partial\rho}\right)_T = \frac{p}{\rho^2} + T\left(\frac{\partial s}{\partial \rho}\right)_T.\label{e_deriv:iso_T}
    \end{align}
\end{rem}
Note that the relations for the barotropic case were not given in previous work and in particular the isothermal case is no straightforward simplification of the general case. The detailed calculations are given in the Appendix \ref{sec:app} and we will just summarize the results needed here. We now introduce the \emph{mixture pressure}
\begin{align}
    p = \alpha_1p_1 + \alpha_2p_2,\label{def:mix_p}
\end{align}
the \emph{specific enthalpy} of phase $i$
\begin{align}
    h_i(\rho_i) = e_i(\rho_i) + \frac{p_i(\rho_i)}{\rho_i}\label{def:enthalpy}
\end{align}
and the \emph{specific Gibbs energy} of phase $i$
\begin{align}
    g_i(\rho_i) = e_i(\rho_i) - Ts_i(\rho_i) + \frac{p_i(\rho_i)}{\rho_i}.\label{def:gibbs_energy}
\end{align}
The Gibbs energy, or sometimes \emph{free enthalpy}, may be generalized to the \emph{chemical potential} $\mu$ when more substances are involved, see \cite{Landau1987}.
For the mixture EOS we have the following derivatives
\begin{align}
    \frac{\partial E}{\partial\alpha} = \frac{\partial e}{\partial\alpha},\quad \frac{\partial E}{\partial c} = \frac{\partial e}{\partial c} + (1 - 2c)\frac{w_iw^i}{2},\quad
    \frac{\partial E}{\partial\rho} = \frac{\partial e}{\partial\rho}\quad\text{and}\quad
    \frac{\partial E}{\partial w_i} = c(1 - c)w_i.\label{eqn:mixE_deriv}
\end{align}
The derivatives for the internal energy $e(\alpha,c,\rho) = c_1e_1(\rho_1) + c_2e_2(\rho_2)$ of the mixture are 
\begin{equation}
	 \frac{\partial e}{\partial\alpha} = \frac{p_2 - p_1}{\rho}, \qquad 
	 \frac{\partial e}{\partial c}     = h_1(\rho_1) - h_2(\rho_2), \qquad 
	 \frac{\partial e}{\partial\rho}   = \frac{p}{\rho^2}
	 \label{eqn.der.isentropic} 
\end{equation}
for the isentropic case and 
\begin{equation}
	\frac{\partial e}{\partial\alpha} = \frac{p_2 - p_1}{\rho} + T(c_1s_1 - c_2s_2), \qquad 
	\frac{\partial e}{\partial c}     =  g_1(\rho_1) - g_2(\rho_2), \qquad 
	\frac{\partial e}{\partial\rho}   = \frac{p}{\rho^2} - \frac{T(c_1s_1 - c_2s_2)}{\rho} 
	\label{eqn.der.isothermal} 
\end{equation}
for the isothermal case.

%

\subsection{One-dimensional model}
In this paper we focus on the one dimensional case and thus the equations simplify to
\begin{subequations}\label{sys:gpr_1d_v1}
    \begin{align} 
        \frac{\partial \alpha\rho}{\partial t} + \frac{\partial \alpha\rho u}{\partial x} &= \xi_1,\label{eq:alpha_balance_v1}\\
        \frac{\partial \rho c}{\partial t} + \frac{\partial (\rho c u + \rho E_w)}{\partial x} &= \xi_2,\label{eq:partial_mass_balance_v1}\\
        \frac{\partial\rho}{\partial t} + \frac{\partial\rho u}{\partial x} &= \xi_3,\label{eq:mix_mass_balance_v1}\\
        \frac{\partial\rho u}{\partial t} + \frac{\partial\left(\rho u^2 + p + \rho wE_w\right)}{\partial x} &= \xi_4,\label{eq:mix_mom_balance_v1}\\
        \frac{\partial w}{\partial t} + \frac{\partial(wu + E_c)}{\partial x} &= \xi_5.\label{eq:rel_velocity_balance_v1}
    \end{align} 
\end{subequations}
%
Note, that the curl term in equation \eqref{eqn.relvel} vanishes. Applying the obtained results and relations and introducing the notation 
\begin{align*}
	\Psi_i(\rho_i) = \begin{cases}
		h_i(\rho_i),\, \qquad \text{isentropic}\\
		g_i(\rho_i),\, \qquad \text{isothermal}
	\end{cases}
\end{align*}
the system can be rewritten in the following form:

\begin{subequations}\label{sys:gpr_1d_v2}
    \begin{align} 
        \frac{\partial\alpha_1\rho}{\partial t} + \frac{\partial\alpha_1\rho u}{\partial x} &= \xi_1,\label{eq:alpha_balance}\\
        \frac{\partial\alpha_1\rho_1}{\partial t} + \frac{\partial\alpha_1\rho_1 u_1}{\partial x} &= \xi_2,\label{eq:partial_mass_balance}\\
        \frac{\partial\rho}{\partial t} + \frac{\partial\rho u}{\partial x} &= \xi_3,\label{eq:mix_mass_balance}\\
        \frac{\partial(\alpha_1\rho_1 u_1 + \alpha_2\rho_2 u_2)}{\partial t}
        + \frac{\partial\left(\alpha_1\rho_1 u_1^2 + \alpha_2\rho_2 u_2^2 + \alpha_1 p_1(\rho_1) + \alpha_2 p_2(\rho_2)\right)}{\partial x} &= \xi_4,\label{eq:mix_mom_balance}\\
        \frac{\partial w}{\partial t} + \dfrac{\partial}{\partial x} \left(\dfrac{1}{2}u_1^2 - \dfrac{1}{2}u_2^2 + \Psi_1(\rho_1) - \Psi_2(\rho_2)\right) &= \xi_5,\label{eq:rel_velocity_balance}
    \end{align} 
\end{subequations}
For the isentropic case the total energy inequality, which serves as mathematical entropy inequality, reads
\begin{align}
    \sum_{i=1}^2\frac{\partial\alpha_i\rho_i\left(e_i + \frac{1}{2}u_i^2\right)}{\partial t}
        + \frac{\partial\alpha_i\rho_i u_i\left(h_i + \frac{1}{2}u_i^2\right)}{\partial x}\leq 0, \label{ineq:isoS_energy}
\end{align}
while for the isothermal case we have the inequality
\begin{align}
    \sum_{i=1}^2\frac{\partial\alpha_i\rho_i\left(e_i - Ts_i + \frac{1}{2}u_i^2\right)}{\partial t}
        + \frac{\partial\alpha_i\rho_i u_i\left(g_i + \frac{1}{2}u_i^2\right)}{\partial x} \leq 0.\label{ineq:isoT_energy}
\end{align}
Inequality \eqref{ineq:isoT_energy} is the mathematical formulation of the physical statement that the free energy of a system under consideration is minimized in an isothermal process. 
Moreover, this formulation is a straightforward generalization of the inequality obtained for the isothermal Euler  equations, see \cite{Dafermos2016,Ruggeri2021,Thein2018}.  
We want to note that it is beneficial to derive the isothermal model from the more general model including the thermal impulse. More detailed information are given in the Appendix.
\subsection{Conservative formulation}
The system (\ref{eq:alpha_balance}) - (\ref{eq:rel_velocity_balance}) can be written in the conservative form given by 
\begin{align*}
	\frac{\partial}{\partial t}\bb{W} + \frac{\partial}{\partial x}\bb{F}(\bb{W}) = \bb{\Xi}, 
\end{align*}
where the vector of conserved quantities reads 
\begin{align*}
    \bb{W} = (w_1,w_2,w_3,w_4,w_5)^T \equiv \left(\alpha_1\rho,\alpha_1\rho_1,\rho,\alpha_1\rho_1 u_1 + \alpha_2\rho_2 u_2,u_1-u_2\right)^T.
\end{align*}
Using the equations obtained so far we can write the conservative flux as follows

\begin{align}
    \bb{F}(\bb{W}) = \begin{pmatrix} w_1\dfrac{w_4}{w_3}\\[10pt] w_2\dfrac{(w_3 - w_2)w_5 + w_4}{w_3}\\[10pt] w_4\\[10pt]
    w_2\left(\dfrac{(w_3 - w_2)w_5 + w_4}{w_3}\right)^2 + (w_3 - w_2)\left(\dfrac{w_4 - w_2w_5}{w_3}\right)^2 + \dfrac{w_1}{w_3}p_1(\bb{W}) + \dfrac{w_3-w_1}{w_3}p_2(\bb{W})\\[12pt]
    \dfrac{1}{2}w_5\left(2\dfrac{(w_3 - w_2)w_5 + w_4}{w_3} - w_5\right) + \Psi_1(\bb{W}) - \Psi_2(\bb{W})
    \end{pmatrix}.\label{eq:cons_flux_vec}
\end{align}
\subsection{Primitive formulation}
We want to reformulate the barotropic system (\ref{sys:gpr_1d_v2}) in terms of the primitive variables $\alpha_1, \rho_1, \rho_2, u_1$ and $u_2$.
The other quantities are then obtained using the relations 
\begin{equation*}
    \alpha_2 = 1 - \alpha_1, \quad 
    \rho = \alpha_1\rho_1 + (1 - \alpha_1)\rho_2, \quad 
    \rho u = \alpha_1\rho_1u_1 + (1 - \alpha_1)\rho_2u_2, \quad 
    w = u_1 - u_2.
\end{equation*}
Here we assume $\alpha_1 \in (0,1)$ and $\rho_1,\rho_2 > 0$, i.e.\ we exclude vacuum states.
We further introduce the speed of sound $a_i$ of phase $i$ given by the following relation
\begin{align}
    a_i^2 = \begin{dcases}
        \rho_i\left(\dfrac{\partial h_i}{\partial\rho_i}\right)_s,\, \qquad \text{isentropic}, \\
        \rho_i\left(\dfrac{\partial g_i}{\partial\rho_i}\right)_T,\, \qquad \text{isothermal}. 
    \end{dcases}\label{def:sound_speed}
\end{align}
For a thermodynamically consistent equation of state the speed of sound is well defined and we do not have to differ between the two cases for the mathematical considerations. 
With $\bb{W} \equiv (\alpha_1,\rho_1,\rho_2,u_1,u_2)$ the Jacobian of this system is given by
\begin{align}
    \bb{A}(\bb{W}) = \begin{pmatrix}
        u                                   & 0                     & 0                     & 0         & 0\\[8pt]
        \dfrac{\rho_1}{\alpha_1}(u_1 - u)   & u_1                   & 0                     & \rho_1    & 0\\[8pt]
        \dfrac{\rho_2}{\alpha_2}(u - u_2)   & 0                     & u_2                   & 0         & \rho_2\\[8pt]
        \dfrac{p_1 - p_2}{\rho}             & \dfrac{a_1^2}{\rho_1} & 0                     & u_1       & 0\\[8pt]
        \dfrac{p_1 - p_2}{\rho}             & 0                     & \dfrac{a_2^2}{\rho_2} & 0         & u_2
    \end{pmatrix}\label{eq:jacobian_prim_var}
\end{align}
and the source terms can be transformed using the following matrix
\begin{align}
    \bb{B}(\bb{W}) = \begin{pmatrix}
        \dfrac{1}{\rho} & 0 & -\dfrac{\alpha_1}{\rho} & 0 & 0\\[8pt]
        -\dfrac{\rho_1}{\alpha_1\rho} & \dfrac{1}{\alpha_1} & \dfrac{\rho_1}{\rho} & 0 & 0\\[8pt]
        \dfrac{\rho_2}{\alpha_2\rho} & -\dfrac{1}{\alpha_2} & \dfrac{c_1}{\alpha_2} + c_2 & 0 & 0\\[8pt]
        0 & \dfrac{u_2 - u_1}{\rho} & -\dfrac{u_2}{\rho} & \dfrac{1}{\rho} & c_2\\[8pt]
        0 & \dfrac{u_2 - u_1}{\rho} & -\dfrac{u_2}{\rho} & \dfrac{1}{\rho} & -c_1
    \end{pmatrix}.\label{eq:sourcemat_prim_var}
\end{align}
Now we can write the system in the following compact form
\begin{align*}
    \partial_t\bb{W} + \bb{A}(\bb{W})\partial_x\bb{W} = \bb{B}(\bb{W})\bb{\Xi}.
\end{align*}
The eigenvalues can be computed as 
\begin{align}
    \lambda_{1\pm} = u_1 \pm a_1,\quad\lambda_C = u,\quad\lambda_{2\pm} = u_2 \pm a_2\label{sys_eigenvalues}
\end{align}
and we have (up to scaling) the following right eigenvectors
\begin{align}
    \bb{R}_{1\pm} = \begin{pmatrix} 0 \\ 1 \\ 0 \\ \pm\dfrac{a_1}{\rho_1} \\ 0\end{pmatrix},\quad
    \bb{R}_C = \begin{pmatrix} \varepsilon_1\varepsilon_2\\ \delta_1\varepsilon_2\\ \delta_2\varepsilon_1\\ \hphantom{-}(u-u_1)\varepsilon_2\gamma_1\\ -(u-u_2)\varepsilon_1\gamma_2\end{pmatrix},
    \quad
    \bb{R}_{2\pm} = \begin{pmatrix} 0 \\ 0 \\ 1 \\ 0 \\ \pm\dfrac{a_2}{\rho_2}\end{pmatrix}.\label{sys_eigenvectors_prim}
\end{align}
Here we introduced the following abbreviations
\begin{align*}
    \delta_1 &= \frac{p_1 - p_2}{\rho} - \frac{(u-u_1)^2}{\alpha_1},&\quad \delta_2 &= \frac{p_1 - p_2}{\rho} + \frac{(u-u_2)^2}{\alpha_2},\\
    \varepsilon_1 &= \frac{(u - u_1)^2 - a_1^2}{\rho_1},&\quad \varepsilon_2 &= \frac{(u - u_2)^2 - a_2^2}{\rho_2},\\
    \gamma_1 &= \frac{\alpha_1(p_1 - p_2) - \rho a_1^2}{\alpha_1\rho_1\rho},&\quad\gamma_2 &= -\frac{\alpha_2(p_1 - p_2) + \rho a_2^2}{\alpha_2\rho_2\rho}.
\end{align*}
The subscripts $1,2$ refer to the corresponding phases and the pair $(\lambda_C, \bb{R}_C)$ takes a special role as we will see in a moment.
We further want to investigate the fields and see whether they are genuine nonlinear or linearly degenerated.
The gradients of the eigenvalues with respect to the given variables $\bb{W} = (\alpha_1,\rho_1,\rho_2,u_1,u_2)$ are given by
\begin{align}
    \nabla_{\bb{W}}\lambda_{1\pm} = \begin{pmatrix} 0\\ \pm\dfrac{\partial a_1}{\partial\rho_1}\\ 0\\ 1\\ 0\end{pmatrix},\quad
    \nabla_{\bb{W}}\lambda_{C} = \begin{pmatrix} \dfrac{\rho_1\rho_2}{\rho^2}(u_1 - u_2)\\ \dfrac{\alpha_1c_2(u_1 - u_2)}{\rho}\\ \dfrac{\alpha_2c_1(u_2 - u_1)}{\rho}\\ c_1\\ c_2\end{pmatrix},
    \quad
    \nabla_{\bb{W}}\lambda_{2\pm} = \begin{pmatrix} 0\\ 0\\ \pm\dfrac{\partial a_2}{\partial\rho_2}\\ 0\\ 1\end{pmatrix}.\label{grad_eig_val}
\end{align}
We immediately obtain
\begin{align*}
    \nabla_{\bb{W}}\lambda_{1\pm}\cdot\bb{R}_{1\pm} &= \pm\frac{\partial a_1}{\partial\rho_1} \pm \frac{a_1}{\rho_1}
    = \pm\frac{1}{\rho_1}\dfrac{\partial(\rho_1 a_1)}{\partial\rho_1} = \pm\frac{a_1}{\rho_1}\mathcal{G}_1,\\
    %
    %
    \nabla_{\bb{W}}\lambda_{2\pm}\cdot\bb{R}_{2\pm} &= \pm\frac{\partial a_2}{\partial\rho_2} \pm \frac{a_2}{\rho_2}
    = \pm\frac{1}{\rho_2}\dfrac{\partial(\rho_2 a_2)}{\partial\rho_2} = \pm\frac{a_2}{\rho_2}\mathcal{G}_2.
\end{align*}
Here we have introduced the \emph{fundamental derivative} $\mathcal{G}$ and the following relation holds
\begin{align*}
    \mathcal{G} = 1 + \frac{\rho}{a}\frac{\partial a}{\partial\rho}\quad\Leftrightarrow\quad\frac{1}{\rho}\frac{\partial(\rho a)}{\partial\rho} = \frac{a}{\rho}\mathcal{G}.
\end{align*}
For more details on $\mathcal{G}$ and its crucial influence on the flow we recommend Menikoff and Plohr \cite{Menikoff1989} and M\"uller and Voss \cite{Mueller2006}.
Throughout this work we assume $\mathcal{G} > 0$ and thus the fields $1\pm$ and $2\pm$ are genuine nonlinear. Indeed for an ideal gas we have
\begin{align*}
    \mathcal{G} =
    \begin{dcases}
        \frac{\gamma + 1}{2},\quad &\text{isentropic with } \gamma > 1\\
        1,\quad &\text{isothermal}
    \end{dcases}
\end{align*}
We further have the remarkable property that
\begin{align*}
    \langle\bb{R}_{1+},\bb{R}_{1-}\rangle\perp\langle\bb{R}_{2+},\bb{R}_{2-}\rangle\quad
    \text{and}\quad\nabla_{\bb{W}}\lambda_{1\pm}\in\langle\bb{R}_{1+},\bb{R}_{1-}\rangle,\quad\nabla_{\bb{W}}\lambda_{2\pm}\in\langle\bb{R}_{2+},\bb{R}_{2-}\rangle.
\end{align*}
Thus the eigenvectors $\bb{R}_{1\pm}$ and $\bb{R}_{2\pm}$ span a four dimensional hyperplane in the state space where $\alpha$ is constant.
It remains to discuss the field $C$. Therefore we use
\begin{align*}
    u - u_1 &= (c_1 - 1)u_1 + c_2u_2 = -c_2(u_1 - u_2) = -c_2w,\\
    u - u_2 &= c_1u_1 + (c_2 - 1)u_2 = \hphantom{-}c_1(u_1 - u_2) = \hphantom{-}c_1w.
\end{align*}
We obtain
\begin{align*}
    \nabla_{\bb{W}}\lambda_C\cdot\bb{R}_C
    &= \begin{pmatrix} \dfrac{\rho_1\rho_2}{\rho^2}(u_1 - u_2)\\ \dfrac{\alpha_1c_2(u_1 - u_2)}{\rho}\\ \dfrac{\alpha_2c_1(u_2 - u_1)}{\rho}\\ c_1\\ c_2\end{pmatrix}
    \cdot\begin{pmatrix} \varepsilon_1\varepsilon_2\\ \delta_1\varepsilon_2\\ \delta_2\varepsilon_1\\ \hphantom{-}(u-u_1)\varepsilon_2\gamma_1\\ -(u-u_2)\varepsilon_1\gamma_2\end{pmatrix}\\
    &= \frac{w}{\rho^2}\left(\rho_1\rho_2\varepsilon_1\varepsilon_2 + \alpha_1\alpha_2\rho_2\delta_1\varepsilon_2 - \alpha_1\alpha_2\rho_1\delta_2\varepsilon_1 
    - \alpha_1\rho_1\alpha_2\rho_2\varepsilon_2\gamma_1 - \alpha_1\rho_1\alpha_2\rho_2\varepsilon_1\gamma_2\right)\\
    &= \frac{w}{\rho^2}\left(\left((c_2w)^2 - a_1^2\right)\left((c_1w)^2 - a_2^2\right)
    + \alpha_1\alpha_2\left((c_1w)^2 - a_2^2\right)\left(\frac{p_1 - p_2}{\rho} - \frac{(c_2w)^2}{\alpha_1}\right)\right.\\
    &- \alpha_1\alpha_2\left((c_2w)^2 - a_1^2\right)\left(\frac{p_1 - p_2}{\rho} + \frac{(c_1w)^2}{\alpha_2}\right)
    - \alpha_2\left((c_1w)^2 - a_2^2\right)\frac{\alpha_1(p_1 - p_2) - \rho a_1^2}{\rho}\\
    &\left.+ \alpha_1\left((c_2w)^2 - a_1^2\right)\frac{\alpha_2(p_1 - p_2) + \rho a_2^2}{\rho}\right)\\
    &= \frac{w}{\rho^2}\left(\left((c_2w)^2 - a_1^2\right)\left((c_1w)^2 - a_2^2\right)
    - \alpha_2\left((c_1w)^2 - a_2^2\right)(c_2w)^2\right.\\
    &\left.- \alpha_1\left((c_2w)^2 - a_1^2\right)(c_1w)^2 + \alpha_2\left((c_1w)^2 - a_2^2\right)a_1^2
    + \alpha_1\left((c_2w)^2 - a_1^2\right)a_2^2\right)\\
    &= 0.
\end{align*}
Thus this field is linearly degenerated and hence discontinuities associated to this field are contact waves.
In view of the above results it is clear that each character of the present fields is independent of the flow.
Note that up to now we cannot exclude situations where eigenvalues coincide, i.e.\ have a multiplicity larger than one.
Different possible phenomena related to these special situations will be discussed in Section \ref{sec:wave_config}.

    \section{Degeneracies of the system and admissibility conditions}\label{sec:degen_adm_cond}
By construction, the system under consideration is symmetric hyperbolic in the sense of Friedrichs \cite{Friedrichs1954}.
The initial value problem (Cauchy problem) for such a system is well-posed locally in time \cite{Dafermos2016,Serre2007}. But the question of the solvability of the Cauchy problem in the large for
symmetric hyperbolic systems is still a challenging problem and, for example, nontrivial interesting phenomena (such as the resonance effect) can arise.
As mentioned in the previous Section we now will have to study under which conditions this system is hyperbolic in the sense that we have real eigenvalues and a full set of eigenvectors, or, more precisely, under which conditions certain degeneracies may occur. 
Furthermore another crucial point is, as mentioned before, that up to now the order of the waves is not clear, i.e.\ the order of the eigenvalues.
Thus we will briefly review results concerning the admissibility conditions for discontinuities which will play an important role throughout this work.
\subsection{Coincidence of Eigenvalues, Parabolic Degeneracy and Hyperbolic Resonance}
We now want to investigate situations where the eigenvectors may become linearly dependent.
Considering the eigenvectors (\ref{sys_eigenvectors_prim}) it is obvious that the eigenvectors $\bb{R}_{i\pm}$ become linearly dependent in each phase iff $a_i = 0$.
Since we assume a strictly positive speed of sound in each phase this is excluded and the four eigenvectors $\bb{R}_{i\pm}$ remain linearly independent.
For the field $(\lambda_C,\mathbf{R}_C)$ we can state for a fixed field $(\lambda_{i\pm},\mathbf{R}_{i\pm})$ that
\begin{align}
    \mathbf{R}_C \in \langle \mathbf{R}_{i\pm}\rangle \quad\Leftrightarrow\quad \lambda_{i\pm} = \lambda_C.\label{cond:absent_contact}
    %
\end{align}
This can be seen by direct calculations which are given in the Appendix \ref{sec:app} for convenience.
We want to discuss the consequences and interpretation of this case later on.
We note that in the extreme situation that three eigenvalues coincide, which is the maximum for reasonable EOS, the eigenvector $\bb{R}_C$ becomes the null vector. 
This situation can only happen when $\lambda_C$ coincides with one eigenvalue of each phase and thus we have 
\begin{align*}
    (u - u_1)^2 = a_1^2\quad\text{and}\quad (u - u_2)^2 = a_2^2.
\end{align*}
\subsection{Admissibility Conditions for Discontinuities Revisited}
In the previous Section \ref{sec:model} we presented the eigenvalues (\ref{sys_eigenvalues}) and eigenvectors (\ref{sys_eigenvectors_prim}) of the system under consideration (\ref{sys:gpr_1d_v2}).
As noted before different situations may occur due to coinciding eigenvalues. Thus it is important to review suitable criteria to single out admissible solutions.
In particular we are interested in the case when discontinuities occur. For a given and fixed state $\bb{W}$ we have for the eigenvalues
\begin{align}
    \lambda^{(1)}(\bb{W}) \leq \dots \leq \lambda^{(p)}(\bb{W})\label{non_str_eigval_order}
\end{align}
with $p \leq n$ and $n$ being the total number of the possible eigenvalues. In our case we have $n = 5$ and further
\begin{align*}
    \left\{\lambda^{(1)}(\bb{W}),\dots,\lambda^{(p)}(\bb{W})\right\} \subseteq \left\{\lambda_{1\pm}(\bb{W}),\lambda_C(\bb{W}),\lambda_{2\pm}(\bb{W})\right\}.
\end{align*}
With this we allow for situations where the order of the eigenvalues changes and that eigenvalues may coincide.
We start with the original work by Lax \cite{Lax1984} and assume $p = n$. In particular we consider a strictly hyperbolic system.
Further we consider the states left and right of a discontinuity denoted by $\bb{W}_L$ and $\bb{W}_R$.
According to Lax a $k$-shock with the speed $S$ satisfies the condition
\begin{subequations}\label{Lax_cond}
    \begin{align}
        \lambda^{(k)}(\bb{W}_L) &> S > \lambda^{(k-1)}(\bb{W}_L),\label{Lax_cond1}\\
        \lambda^{(k+1)}(\bb{W}_R) &> S > \lambda^{(k)}(\bb{W}_R).\label{Lax_cond2}
    \end{align}
\end{subequations}
From these inequalities we deduce that $n - (k - 1)$ characteristics impinge from the left of the discontinuity and $k$ from the right.
Following the presentation given in Dafermos \cite{Dafermos2016} the situation can be generalized as follows. We consider the eigenvalues for the left and right state
\begin{align*}
    &\lambda^{(1)}(\bb{W}_L) \leq \dots \leq \lambda^{(i-1)}(\bb{W}_L) < S < \lambda^{(i)}(\bb{W}_L) \leq \dots \leq \lambda^{(n)}(\bb{W}_L),\\
    &\lambda^{(1)}(\bb{W}_R) \leq \dots \leq \lambda^{(j)}(\bb{W}_R) < S < \lambda^{(j+1)}(\bb{W}_R) \leq \dots \leq \lambda^{(n)}(\bb{W}_R),
\end{align*}
with the agreement of $\lambda^{(0)}(\bb{W}_L) = -\infty$ and $\lambda^{(n+1)}(\bb{W}_R) = \infty$.
If the inequalities are satisfied with $i = j$ we obtain the Lax condition (\ref{Lax_cond}) given above and the shock is called \emph{compressive}.
In the case of $i < j$ the shock is called \emph{overcompressive} and for $i > j$ the shock is called \emph{undercompressive}.
Now we want to deal with the non strict situation, i.e.\ there is no well-defined ordering of the eigenvalues and they may coincide.
In Keyfitz and Kranzer \cite{Keyfitz1980} a generalization to non-strictly hyperbolic systems is given as follows
\begin{enumerate}[1)]
    \item $n+1$ characteristics enter the shock and $n-1$ leave it or
    \item $n-1$ characteristics enter and leave the shock whereas the remaining two are tangent to the shock and belong to linearly degenerated fields.
\end{enumerate}
When an $i$-shock has $n-1$ leaving characteristics and the corresponding eigenvectors fulfill
\begin{align}
    \det\left[\bb{R}^{(1)}(\bb{W}_L),\dots,\bb{R}^{(i-1)}(\bb{W}_L),\bb{W}_R - \bb{W}_L,\bb{R}^{(i+1)}(\bb{W}_R),\dots,\bb{R}^{(n)}(\bb{W}_R)\right] \neq 0\label{evol_cond_det}
\end{align}
then the shock is called \emph{evolutionary}, see again \cite{Dafermos2016} or the book of Kulikovskii et al.\ \cite{Kulikovskii2001}.
We want to end this brief review of admissibility conditions with the results presented in \cite{Andrianov2003}.
There the results obtained in the references given above are basically collected and generalized to incorporate different situations that are of interest for non-strictly hyperbolic systems.
Let us consider the situation (\ref{non_str_eigval_order}).
Having $p$ eigenvalues we conclude that we have $p$ unknowns on each side of the discontinuity and additionally the shock speed $S$, i.e.\ $N = 2p + 1$ unknowns.
These $N$ unknowns can be determined as follows. Let us assume we have $m$ relations across the discontinuity, e.g.\ the jump conditions.
Second, the unknowns should be determined by the flow using the characteristics. Let the number of incoming characteristics be $i$, the number of outgoing characteristics $o$ and the number of
coinciding characteristics be $c$. The incoming and coinciding characteristics are determined by the past and thus provide further information.
Hence in order to determine the unknowns we demand
\begin{align}
    N = i + c + m.\label{evo_cond}
\end{align}
In \cite{Andrianov2003} this is called \emph{evolutionarity condition} and a discontinuity is called \emph{evolutionary} if the condition (\ref{evo_cond}) holds (in agreement with the results above).
This implies that a discontinuity is evolutionary iff $o = m - 1$, see \cite{Andrianov2003}. A closely related concept was introduced by Freist\"uhler in \cite{Freistuehler1992a}.
The idea is quite similar as again a linearised problem is investigated.
The eigenvalues are distinguished, using there relation to the speed of the discontinuity, into slow and fast characteristics.
In this sense outgoing characteristics are slow characteristics on the left side of the discontinuity and fast ones on the right side. For the incoming characteristics the situation is reversed.
Coinciding characteristics move at the characteristic speed of the discontinuity.

    \section{Wave relations}\label{sec:waverel}
We now want to obtain the relations that are valid across the different types of waves.
In particular we derive the Riemann invariants which are constant across rarefaction waves and the Rankine Hugoniot jump conditions across discontinuities.
\subsection{Rarefaction waves}\label{subsec:rarefaction}
Given a smooth solution we may apply a nonlinear transformation of the variables aiming to simplify the system using another appropriate choice of variables.
Such a special set of variables is given by the \emph{Riemann invariants}, cf.\ \cite{Dafermos2016,Evans1998,Smoller2012}.
The Riemann invariants for $\bb{R}_{1\pm}$ and $\bb{R}_{2\pm}$ can be calculated using
\begin{align*}
    \frac{\dd\alpha_1}{\dd s} &= 0,\;\frac{\dd\rho_\mu}{\dd s} = \left(\pm\frac{a_\mu}{\rho_\mu}\mathcal{G}_\mu\right)^{-1},\;\frac{\dd\rho_\nu}{\dd s} = 0,\;
    \frac{\dd u_\mu}{\dd s} = \pm\left(\pm\frac{a_\mu}{\rho_\mu}\mathcal{G}_\mu\right)^{-1}\frac{a_\mu}{\rho_\mu}\quad\text{and}\quad\frac{\dd u_\nu}{\dd s} = 0.
\end{align*}
Here $\mu,\nu \in \{1,2\},\mu\neq\nu$ where $\mu = 1$ for $\bb{R}_{1\pm}$ and $\mu = 2$ for $\bb{R}_{2\pm}$.
The relations for $\rho_\mu$ and $u_\mu$ can be combined to
\begin{align*}
    \frac{\dd u_\mu}{\dd\rho_\mu} = \pm\frac{a_\mu}{\rho_\mu}.
\end{align*}
Thus we obtain the following invariants
\begin{align}
    \alpha_1 = const.,\quad \mathcal{R}_{\pm} = u_\mu \pm \int \frac{a_\mu}{\rho_\mu}\,\dd\rho_\mu,\quad\rho_\nu = const.\quad\text{and}\quad u_\nu = const.\label{riem_inv_raref}
\end{align}
Furthermore, the slope inside a left rarefaction wave is given by
\begin{align}
    \frac{\dd x}{\dd t} = \frac{x}{t} = \lambda_{\mu-} = u_\mu - a_\mu\label{left_raref_slope}
\end{align}
and hence we obtain that the solution inside the rarefaction fan is given by
\begin{align}
    u_\mu = \frac{x}{t} + a_\mu\quad\text{and}\quad F(\rho) = u_\mu - u_{\mu,L} + \int_{\rho_{\mu,L}}^{\rho} \frac{a_\mu}{\rho_\mu}\,\dd\rho_\mu = 0.\label{left_raref_fan}
\end{align}
Here, $\rho$ is obtained as the root of $F(\rho)$. Similar we obtain the results for right rarefaction waves corresponding to $\lambda_{\mu+}$
\begin{align}
    %
    \frac{\dd x}{\dd t} = \frac{x}{t} = \lambda_{\mu+} = u_\mu + a_\mu,\;
    u_\mu &= \frac{x}{t} - a_\mu\quad\text{and}\quad F(\rho) = u_{\mu,R} - u_\mu - \int_\rho^{\rho_{\mu,R}} \frac{a_\mu}{\rho_\mu}\,\dd\rho_\mu = 0.\label{right_raref}
\end{align}
Although the Riemann invariants (\ref{riem_inv_raref}) state that $\rho_\nu$ and $u_\nu$ remain constant, we prefer the notion that these quantities are not affected by the rarefaction wave.
As we see later on they might change due to other waves.
\subsection{Shock waves}\label{subsec:shock}
In the presence of discontinuities, such as shock waves, the situation is different. However, since the present system is conservative, corresponding Rankine Hugoniot jump conditions
\begin{align*}
    \dbl\mathbf{F}(\bb{W})\dbr = S\dbl \bb{W}\dbr.
\end{align*}
have to hold at a shock with speed $S$, see for example \cite{Dafermos2016,LeFloch2002}.
In the present case we obtain the following conditions across discontinuities
\begin{subequations}\label{rhc_v1}
    \begin{align}
        \dbl\alpha_1\rho u\dbr &= S\dbl\alpha_1\rho\dbr,\label{jc:alpha_v1}\\
        \dbl\alpha_1\rho_1 u_1\dbr &= S\dbl\alpha_1\rho_1\dbr,\label{jc:Q1_v1}\\
        \dbl\rho u\dbr &= S\dbl\rho\dbr,\label{jc:Q_v1}\\
        \dbl \alpha_1\rho_1u_1^2 + \alpha_2\rho_2u_2^2 + \alpha_1p_1 + \alpha_2p_2\dbr &= S\dbl\alpha_1\rho_1u_1 + \alpha_2\rho_2u_2\dbr,\label{jc:momentum_v1}\\
        \dbl \frac{1}{2}\left(u_1^2 - u_2^2\right) + \Psi_1 - \Psi_2\dbr &= S\dbl u_1 - u_2\dbr.\label{jc:w_v1}
    \end{align}
\end{subequations}
The jump conditions (\ref{jc:alpha_v1})-(\ref{jc:Q_v1}) can be reformulated to
\begin{align}
    \dbl\alpha_1\rho(u - S)\dbr &= 0,\label{jc:alpha_v2}\\
    \dbl\alpha_1\rho_1(u_1 - S)\dbr &= 0,\label{jc:Q1_v2}\\
    \dbl\rho(u - S)\dbr &= 0.\label{jc:Q_v2}
\end{align}
Using the third jump condition (\ref{jc:Q_v2}) we obtain for the first equation (\ref{jc:alpha_v2})
\begin{align*}
    \rho(u - S)\dbl \alpha_1\dbr &= 0.
\end{align*}
We introduce the abbreviations
\begin{align*}
    Q = -\rho(u - S),\quad Q_1 = -\rho_1(u_1 - S)\quad\text{and}\quad Q_2 = -\rho_2(u_2 - S)
\end{align*}
for the mixture mass flux and the mass fluxes of the phases, respectively. In particular we have
\begin{align}
    Q &= -\rho(u - S) = -\left(\rho u - \rho S\right) = -\left(\alpha_1\rho_1u_1 + \alpha_2\rho_2u_2 - (\alpha_1\rho_1 + \alpha_2\rho_2)S\right)\notag\\
    &= -\left(\alpha_1\rho_1(u_1 - S) + \alpha_2\rho_2(u_2 - S)\right) = \alpha_1Q_1 + \alpha_2Q_2.\label{eq:mass_flux}
\end{align}
Thus we can derive a jump condition for the partial mass flux of the second phase,
using the continuity of the mixture mass flux (\ref{jc:Q_v2}) and the continuity mass flux of the first phase (\ref{jc:Q1_v2}), i.e.
\begin{align}
    0 = \dbl Q\dbr \stackrel{(\ref{eq:mass_flux})}{=} \dbl\alpha_1Q_1 + \alpha_2Q_2\dbr \stackrel{(\ref{jc:Q1_v2})}{=} \dbl\alpha_2Q_2\dbr.\label{jc:Q2_v1}
\end{align}
Using the jump conditions for the partial mass fluxes (\ref{jc:Q1_v2}) and (\ref{jc:Q2_v1}) we obtain for the fourth jump condition (\ref{jc:momentum_v1})
\begin{align}
    &\phantom{=\ } \dbl \alpha_1\rho_1u_1^2 + \alpha_2\rho_2u_2^2 + \alpha_1p_1 + \alpha_2p_2\dbr = S\dbl\alpha_1\rho_1u_1 + \alpha_2\rho_2u_2\dbr\notag\\
    \Leftrightarrow\quad
    0 &= \dbl \alpha_1\rho_1u_1(u_1 - S) + \alpha_2\rho_2u_2(u_2 - S) + \alpha_1p_1 + \alpha_2p_2\dbr\notag\\
    &= \alpha_1\rho_1(u_1 - S)\dbl u_1\dbr + \dbl\alpha_1 p_1\dbr + \alpha_2\rho_2(u_2 - S)\dbl u_2\dbr + \dbl\alpha_2 p_2\dbr\notag\\
    &= -\alpha_1Q_1\dbl u_1\dbr + \dbl\alpha_1 p_1\dbr - \alpha_2Q_2\dbl u_2\dbr + \dbl\alpha_2 p_2\dbr.\label{jc:momentum_v2}
\end{align}
The fifth jump condition (\ref{jc:w_v1}) can be reformulated as follows
\begin{align}
    &\dbl \frac{1}{2}\left(u_1^2 - u_2^2\right) + \Psi_1 - \Psi_2\dbr = S\dbl u_1 - u_2\dbr\notag\\
    \Leftrightarrow\quad0 = &\dbl \frac{1}{2}\left(u_1^2 - u_2^2\right) - S(u_1 - u_2) + \Psi_1 - \Psi_2\dbr\notag\\
    = &\dbl\frac{1}{2}\left(u_1 - u_2\right)(u_1 + u_2 - 2S) + \Psi_1 - \Psi_2\dbr\notag\\
    = &\dbl \frac{1}{2}\left(u_1 - u_2\right)(u_1 - S + u_2 - S) + \Psi_1 - \Psi_2\dbr\notag\\
    = &\dbl\frac{1}{2}\left(u_1 - u_2\right)(u_1 - S) + \Psi_1\dbr + \dbl \frac{1}{2}\left(u_1 - u_2\right)(u_2 - S) - \Psi_2\dbr\notag\\
    = &\dbl\frac{1}{2}(u_1 - S)^2 + \Psi_1\dbr +\frac{1}{2}\dbl(S - u_2)(u_1 - S)\dbr\notag\\
    + &\dbl-\frac{1}{2}(u_2 - S)^2 - \Psi_2\dbr +\frac{1}{2}\dbl(u_1 - S)(u_2 - S)\dbr\notag\\
    = &\dbl\frac{1}{2}(u_1 - S)^2 + \Psi_1\dbr - \dbl\frac{1}{2}(u_2 - S)^2 + \Psi_2\dbr.\label{jc:w_v2}
\end{align}
Summarizing we have the following system of jump conditions
\begin{subequations}\label{rhc_v2}
    \begin{align}
        \rho(u - S)\dbl \alpha_1\dbr &= 0,\label{jc:alpha}\\
        \dbl\alpha_1\rho_1(u_1 - S)\dbr &= 0,\label{jc:Q1}\\
        \dbl\rho(u - S)\dbr &= 0,\label{jc:Q}\\
        -\alpha_1Q_1\dbl u_1\dbr + \dbl\alpha_1 p_1\dbr - \alpha_2Q_2\dbl u_2\dbr + \dbl\alpha_2 p_2\dbr &= 0,\label{jc:momentum}\\
        \dbl\frac{1}{2}(u_1 - S)^2 + \Psi_1\dbr - \dbl\frac{1}{2}(u_2 - S)^2 + \Psi_2\dbr &= 0.\label{jc:w}
    \end{align}
\end{subequations}
Note that (\ref{jc:Q}) may be replaced by (\ref{jc:Q2_v1}).
Further, equations (\ref{jc:momentum}) and (\ref{jc:w}) may also be reformulated to
\begin{align}
    -Q\dbl u\dbr + \dbl\rho c_1c_2w^2\dbr + \dbl p\dbr &= 0,\label{jc:momentum_v3}\\
    \frac{1}{2}\dbl\left(\frac{Q_1}{\rho_1}\right)^2\dbr - \frac{1}{2}\dbl\left(\frac{Q_2}{\rho_2}\right)^2\dbr + \dbl \Psi_1 - \Psi_2\dbr &= 0.\label{jc:w_v3}
\end{align}
To give the complete picture we briefly summarize the jump conditions according to the system (\ref{sys:gpr_1d_v1}) in terms of the mixture EOS
\begin{subequations}\label{rhc_v3}
    \begin{align}
        \rho(u - S)\dbl \alpha_1\dbr &= 0,\\
        \dbl\rho c_1(u - S)\dbr + \dbl \rho E_w\dbr &= 0,\\
        \dbl\rho(u - S)\dbr &= 0,\\
        \rho(u - S)\dbl u\dbr + \dbl p + \rho wE_w\dbr &= 0,\\
        \dbl w(u - S)\dbr + \dbl E_c\dbr &= 0.
    \end{align}
\end{subequations}
We also want to emphasize the analogy to the jump conditions of the Euler equations, cf.\ \cite{Dafermos2016,Toro2009}.
Equation (\ref{jc:momentum}) is basically a volume fraction weighted combination of an individual momentum jump condition for each phase as it appears in the Euler equations.
Moreover, a single jump term in (\ref{jc:w}) agrees with the jump bracket of the energy equation for the Euler equations.
\subsection{Exploiting the jump conditions for a Lax shock}\label{subsec:lax_shock}
In the following we assume that we have a Lax-shock, i.e.\ with \emph{no} tangential eigenvalues. In particular this implies $u \neq S$ with $S$ being the shock speed.
Indeed, as we will see this cannot happen for shock waves.
From equation (\ref{jc:alpha}) we thus have that $\alpha_1$ is continuous across a shock corresponding to the eigenvalues $\lambda_{1\pm}$ and $\lambda_{2\pm}$.
Up to now we have made no further assumption on $\alpha_1$. With the aim to further exploit and simplify the jump conditions (\ref{jc:alpha}) - (\ref{jc:w}) we assume $\alpha_1 \in (0,1)$ from now on.
First we have, due to the continuity of $\alpha_1$ (\ref{jc:alpha}), that also the mass fluxes $Q_1$ and $Q_2$ are continuous, see (\ref{jc:Q1}) and (\ref{jc:Q2_v1}).
Hence the partial mass fluxes may be written outside the jump brackets in (\ref{jc:w_v3}), i.e.
\begin{align}
    \frac{Q_1^2}{2}\dbl\frac{1}{\rho_1^2}\dbr - \frac{Q_2^2}{2}\dbl\frac{1}{\rho_2^2}\dbr + \dbl \Psi_1 - \Psi_2\dbr = 0.\label{jc:w_v4}
\end{align}
Further we can use the continuity of the partial mass fluxes and write with $\mu\in\{1,2\}$
\begin{align}
    u_{\mu,R} - S = -\frac{Q_\mu}{\rho_{\mu,R}},\quad u_{\mu,L} - S = -\frac{Q_\mu}{\rho_{\mu,L}}
    \quad\Rightarrow\quad \dbl u_\mu\dbr = -Q_\mu\dbl\frac{1}{\rho_\mu}\dbr.\label{jc:velocity}
\end{align}
Using (\ref{jc:velocity}) together with the momentum jump condition (\ref{jc:momentum_v2}) we get 
\begin{align}
    0 = \alpha_1\left(Q_1^2\dbl \frac{1}{\rho_1}\dbr + \dbl p_1\dbr\right) + \alpha_2\left(Q_2^2\dbl \frac{1}{\rho_2}\dbr + \dbl p_2\dbr\right).\label{jc:momentum_v4}
\end{align}
Thus equations (\ref{jc:w_v4}) and (\ref{jc:momentum_v4}) form a linear system for the partial mass fluxes
\begin{align}
    \underbrace{
    \begin{pmatrix*}[r]
        \alpha_1\dbl \dfrac{1}{\rho_1}\dbr       & \alpha_2\dbl \dfrac{1}{\rho_2}\dbr\\[10pt]
        \dfrac{1}{2}\dbl\dfrac{1}{\rho_1^2}\dbr   & -\dfrac{1}{2}\dbl\dfrac{1}{\rho_2^2}\dbr
    \end{pmatrix*}}_{=: \bb{M}}
    \cdot
    \begin{pmatrix} Q_1^2\\[10pt] Q_2^2\end{pmatrix}
    =
    -\begin{pmatrix} \alpha_1\dbl p_1\dbr + \alpha_2\dbl p_2\dbr\\[10pt] \dbl \Psi_1 - \Psi_2\dbr\end{pmatrix}.\label{lin_sys_Q_square}
\end{align}
The inverse is given by
\begin{align}
    \bb{M}^{-1} &= \frac{1}{\det(\bb{M})}\begin{pmatrix*}[r]
        -\dfrac{1}{2}\dbl\dfrac{1}{\rho_2^2}\dbr    & -\alpha_2\dbl \dfrac{1}{\rho_2}\dbr\\[10pt]
        -\dfrac{1}{2}\dbl\dfrac{1}{\rho_1^2}\dbr    & \alpha_1\dbl \dfrac{1}{\rho_1}\dbr
    \end{pmatrix*},\label{shock_sys_mat_mix}\\
    \det(\bb{M}) &= -\frac{1}{2}\left(\alpha_1\dbl \dfrac{1}{\rho_1}\dbr\dbl\dfrac{1}{\rho_2^2}\dbr + \alpha_2\dbl\dfrac{1}{\rho_1^2}\dbr\dbl \dfrac{1}{\rho_2}\dbr\right).\label{detshock_sys_mat_mix}
\end{align}
Note that across a shock the phase densities necessarily jump and thus the determinant is not equal to zero, see Appendix \ref{sec:app}.
Hence we obtain for the partial mass fluxes
\begin{align}
    Q^2_1 &= \frac{1}{\det(\bb{M})}\left(\dfrac{1}{2}\dbl\dfrac{1}{\rho_2^2}\dbr\left(\alpha_1\dbl p_1\dbr + \alpha_2\dbl p_2\dbr\right)
    + \alpha_2\dbl \dfrac{1}{\rho_2}\dbr\dbl \Psi_1 - \Psi_2\dbr\right),\label{sol_Q1_square}\\
    Q^2_2 &= \frac{1}{\det(\bb{M})}\left(\dfrac{1}{2}\dbl\dfrac{1}{\rho_1^2}\dbr\left(\alpha_1\dbl p_1\dbr + \alpha_2\dbl p_2\dbr\right)
    + \alpha_1\dbl \dfrac{1}{\rho_1}\dbr\dbl \Psi_2 - \Psi_1\dbr\right).\label{sol_Q2_square}
\end{align}
Furthermore, these partial mass fluxes are also related to each other using
\begin{align*}
    Q_1 = -\rho_1(u_1 - S)\;\wedge\;Q_2 = -\rho_2(u_2 - S)\quad\Rightarrow\quad u_1 - u_2 = -\frac{Q_1}{\rho_1} + \frac{Q_2}{\rho_2}.
\end{align*}
Thus one may express one partial mass flux by the other.
Up to now we only have the squares of the mass flux. The correct sign of the corresponding root is given by the Lax condition for the present shock.
The Lax criterion states for a shock in phase $\mu$ that
\begin{align*}
    \lambda_{\mu\pm}(\bb{W}_L) > S > \lambda_{\mu\pm}(\bb{W}_R)
\end{align*}
where $\bb{W}_L$ and $\bb{W}_R$ are the left and right states adjacent to the shock.
For a left shock in phase $\mu \in \{1,2\}$ we have
\begin{align*}
    &u_{\mu,L} - a_{\mu,L} > S > u_{\mu,R} - a_{\mu,R}\\
    \Leftrightarrow\quad \phantom{-}&u_{\mu,L} - S > a_{\mu,L}\;\wedge\; a_{\mu,R} > u_{\mu,R} - S\\
    \Leftrightarrow\quad -&\rho_{\mu,L}(u_{\mu,L} - S) < -\rho_{\mu,L}a_{\mu,L}\;\wedge\; -\rho_{\mu,R}(u_{\mu,R} - S) > -\rho_{\mu,R}a_{\mu,R}\\
    \Leftrightarrow\quad -&\rho_{\mu,R}a_{\mu,R} < Q_\mu < -\rho_{\mu,L}a_{\mu,L} < 0.
\end{align*}
Similar we obtain for a right shock
\begin{align*}
    &u_{\mu,L} + a_{\mu,L} > S > u_{\mu,R} + a_{\mu,R}\\
    \Leftrightarrow\quad \phantom{-}&u_{\mu,L} - S > -a_{\mu,L}\;\wedge\; -a_{\mu,R} > u_{\mu,R} - S\\
    \Leftrightarrow\quad -&\rho_{\mu,L}(u_{\mu,L} - S) < \rho_{\mu,L}a_{\mu,L}\;\wedge\; -\rho_{\mu,R}(u_{\mu,R} - S) > \rho_{\mu,R}a_{\mu,R}\\
    \Leftrightarrow\quad \phantom{-}&\rho_{\mu,L}a_{\mu,L} > Q_\mu > \rho_{\mu,R}a_{\mu,R} > 0.
\end{align*}
Thus the partial mass flux has a strict sign.
Hence across a $\lambda_{1\pm}$-shock the sign of $Q_1$ is determined and for a $\lambda_{2\pm}$-shock the sign of $Q_2$ is given, respectively.
Therefore we choose the following
\begin{align}
    \begin{dcases}
        Q_2 &= \rho_2\left(\frac{Q_1}{\rho_1} + w\right),\quad \lambda_{1\pm} - \text{Shock},\\
        Q_1 &= \rho_1\left(\frac{Q_2}{\rho_2} - w\right),\quad \lambda_{2\pm} - \text{Shock}.
    \end{dcases}\label{rel_mass_fluxes_mix}
\end{align}
The values for the densities and $w$ may assumed to be given, i.e.\ the values on one side of the shock.
Once we have obtained the partial mass flux we can eliminate a further unknown using the velocity jump condition (\ref{jc:velocity})
\begin{align*}
    \dbl u_\mu\dbr = -Q_\mu\dbl \frac{1}{\rho_\mu}\dbr,\;\mu\in\{1,2\}.
\end{align*}
\subsection{Entropy inequality}
Given the (mathematical) entropy inequalities (\ref{ineq:isoS_energy}) and (\ref{ineq:isoT_energy}) we have in the presence of a discontinuity
\begin{align}
    0 \geq \sum_{i=1}^2 -S\dbl \alpha_i\rho_i\left(e_i + \frac{1}{2}u_i^2\right) \dbr + \dbl \alpha_i\rho_i u_i\left(h_i + \frac{1}{2}u_i^2\right)\dbr\label{ineq_disc:isoS_energy_v1}
\end{align}
for the isentropic case and
\begin{align}
    0 \geq \sum_{i=1}^2 -S\dbl \alpha_i\rho_i\left(e_i - Ts_i + \frac{1}{2}u_i^2\right)\dbr + \dbl \alpha_i\rho_i u_i\left(g_i + \frac{1}{2}u_i^2\right) \dbr.\label{ineq_disc:isoT_energy_v1}
\end{align}
for the isothermal case, respectively.
With algebraic manipulations using the continuity of the partial mass fluxes (\ref{jc:Q1_v2}), (\ref{jc:Q2_v1}) and the jump condition for the momentum (\ref{jc:momentum}) we obtain for both cases
\begin{align}
    0 \geq -\sum_{i=1}^2 \alpha_iQ_i\dbl \Psi_i + \frac{1}{2}\left(u_i - S\right)^2 \dbr.\label{ineq_disc:baro_energy_v1}
\end{align}
Note that the bracket terms correspond to the terms in the jump condition for the relative velocity (\ref{jc:w}).
Thus we can replace the term corresponding to one phase by the other if (\ref{jc:w}) holds and obtain
\begin{align}
    0 \geq -Q\dbl \Psi_1 + \frac{1}{2}\left(u_1 - S\right)^2 \dbr\quad\text{or}\quad 0 \geq -Q\dbl \Psi_2 + \frac{1}{2}\left(u_2 - S\right)^2 \dbr.\label{ineq_disc:baro_energy_v2}
\end{align}
In particular the inequality for the isothermal case is perfectly analogous to the one used in \cite{Hantke2019a}.
\subsection{Contact wave}\label{subsec:contact}
The wave corresponding to the linearly degenerated field $(\lambda_C,\bb{R}_C)$ is a contact wave and because of
\begin{align*}
    0 = \nabla_{\bb{W}}\lambda_C\cdot\bb{R}_C = \nabla_{\bb{W}}u\cdot\bb{R}_C
\end{align*}
it is immediately clear that we have the Riemann invariant
\begin{align}
    u = const.\label{riem_inv_contact}
\end{align}
The characteristic condition for the contact wave is $\lambda_C(\bb{W}_L) = S = \lambda_C(\bb{W}_R)$ which further gives for the velocity of the contact $S = u$.
Thus we have $Q = -\rho(u - S) = 0$ and hence $\dbl \alpha_1\dbr$ may be non-zero.
Indeed, since $\alpha_1$ is continuous across the other waves the jump of $\alpha_1$ is given by the initial data.
From the continuity of the partial mass fluxes we obtain 
\begin{align}
    0 &= Q = \alpha_1Q_1 + \alpha_2Q_2\notag\\
    0 &= \dbl \alpha_1Q_1\dbr = -\dbl \alpha_1\rho_1(u_1 - u)\dbr = -\dbl\alpha_1\rho_1c_2w\dbr = -\dbl\rho c_1c_2 w\dbr.\label{jc_c:Q_v1}
\end{align}
For the mixture momentum at the contact we obtain using (\ref{jc:momentum_v3}) and (\ref{jc_c:Q_v1})
\begin{align}
    0 &= -Q\dbl u\dbr + \dbl\rho c_1c_2w^2\dbr + \dbl p\dbr = \dbl\rho c_1c_2w^2\dbr + \dbl p\dbr = \rho c_1c_2w\dbl w\dbr + \dbl p\dbr.\label{jc_c:momentum_v1}
\end{align}
From the jump condition for the relative velocity (\ref{jc:w}) we get 
\begin{align}
    0 &= \dbl\frac{1}{2}(u_1 - S)^2 + \Psi_1\dbr - \dbl\frac{1}{2}(u_2 - S)^2 + \Psi_2\dbr\notag\\
    &= \frac{1}{2}\dbl (c_2w)^2 - (c_1w)^2\dbr + \dbl \Psi_1 - \Psi_2\dbr\notag\\
    &= \frac{1}{2}\dbl \left(c_2^2 - c_1^2\right)w^2\dbr + \dbl \Psi_1 - \Psi_2\dbr\notag\\
    &= \frac{1}{2}\dbl \left(c_2 - c_1\right)w^2\dbr + \dbl \Psi_1 - \Psi_2\dbr.\label{jc_c:w_v1}
\end{align}
Thus we have the following equations at the contact
\begin{subequations}\label{rhc_contact_v1}
    \begin{align}
        \dbl u\dbr &= 0,\label{jc_c:u}\\
        \dbl \rho c_1c_2w\dbr &= 0,\label{jc_c:Q}\\
        \rho c_1c_2w\dbl w\dbr + \dbl p\dbr &= 0,\label{jc_c:momentum}\\
        \frac{1}{2}\dbl \left(c_2 - c_1\right)w^2\dbr + \dbl \Psi_1 - \Psi_2\dbr &= 0.\label{jc_c:w}
    \end{align}
\end{subequations}
Note that $\rho c_1c_2 w^2$ is a dynamic mixture pressure related to the relative velocity. 
By introducing the generalized total mixture pressure $\bar{p} = \rho c_1c_2 w^2 + p$ as the sum of the dynamic pressure and the static mixture pressure $p$ we can rewrite (\ref{jc_c:momentum}) as
\begin{align}
    0 = \rho c_1c_2w\dbl w\dbr + \dbl p\dbr = \dbl\rho c_1c_2w^2 + p\dbr = \dbl\bar{p}\dbr.\label{jc_c:momentum_v2}
\end{align}
%
We again give the jump conditions in terms of the mixture EOS for completeness:  
\begin{subequations}\label{rhc_contact_v2}
    \begin{align}
        \dbl u\dbr &= 0,\\
        \dbl \rho E_w\dbr &= 0,\\
        \dbl p + \rho wE_w\dbr &= 0,\\
        \dbl E_c\dbr &= 0.
    \end{align}
\end{subequations}
Finally, it is easy to see from (\ref{ineq_disc:baro_energy_v2}) and $u = S$ that the entropy inequality is fulfilled with the right side being identically zero.

    \section{Wave configurations and relations}\label{sec:wave_config}
In this section we want to study particular wave configurations.
The crucial point is, as mentioned before, that up to now the order of the waves is not clear, i.e.\ the order of the eigenvalues.
In particular eigenvalues may coincide in certain points.
Thus the reviewed results concerning the admissibility conditions for discontinuities will play a crucial role throughout this section.
\subsection{Contact}
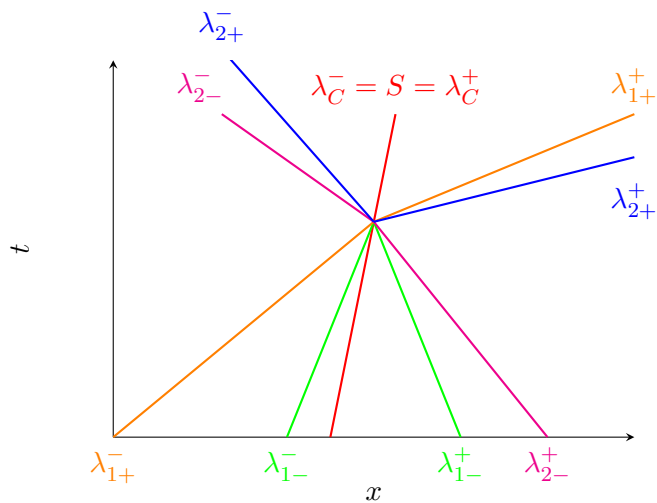
\begin{figure}[h!]
    \center
    \begin{tikzpicture}
        \begin{axis}[xlabel={$x$},axis x line=bottom,ylabel={$t$},axis y line=left,ymin=0,ymax=3.5,domain=-5:5,axis equal image,ytick=\empty,xtick=\empty,clip mode=individual]
            %
            \addplot[mark=none,smooth,red,thick,solid,domain=0:3/5] (\x,{5*\x});
            \node[red,above] at (axis cs:3/5,3) {$\lambda_C^- = S = \lambda_C^+$};
            \addplot[mark=none,smooth,green,thick,solid,domain=-2/5:2/5] (\x,{5*\x/2 + 1});
            \node[green,below] at (axis cs:-2/5,0) {$\lambda_{1-}^-$};
            \addplot[mark=none,smooth,green,thick,solid,domain=2/5:6/5] (\x,{-5*\x/2 + 3});
            \node[green,below] at (axis cs:6/5,0) {$\lambda_{1-}^+$};
            \addplot[mark=none,smooth,orange,thick,solid,domain=-2:2/5] (\x,{5*\x/6 + 5/3});
            \node[orange,below] at (axis cs:-2,0) {$\lambda_{1+}^-$};
            \addplot[mark=none,smooth,orange,thick,solid,domain= 2/5:14/5] (\x,{5*\x/12 + 11/6});
            \node[orange,above] at (axis cs:14/5,3) {$\lambda_{1+}^+$};
            \addplot[mark=none,smooth,magenta,thick,solid,domain=-1:2/5] (\x,{-5*\x/7 + 16/7});
            \node[magenta,above] at (axis cs:-1.2,3) {$\lambda_{2-}^-$};
            \addplot[mark=none,smooth,magenta,thick,solid,domain= 2/5:2] (\x,{-5*\x/4 + 5/2});
            \node[magenta,below] at (axis cs:2,0) {$\lambda_{2-}^+$};
            \addplot[mark=none,smooth,blue,thick,solid,domain= 2/5:14/5] (\x,{3*\x/12 + 19/10});
            \node[blue,above] at (axis cs:14/5,29/15) {$\lambda_{2+}^+$};
            \addplot[mark=none,smooth,blue,thick,solid,domain=-1:2/5] (\x,{-8*\x/7 + 86/35});
            \node[blue,above] at (axis cs:-1,25/7) {$\lambda_{2+}^-$};
        \end{axis}
    \end{tikzpicture}
    \caption{Shock with $u = S$ (example).}
    \label{fig:shock_contact}
\end{figure}

We first want to discuss phenomena related to the linearly degenerates field $(\lambda_C,\mathbf{R}_C)$, i.e.\ the contact wave.
Obviously we can determine the position of the contact wave due to the fact that $\lambda_C = u$ is the convex combination of the individual phase velocities, i.e.
\begin{align*}
    &\min\{u_1,u_2\} \leq u = c_1u_1 + c_2u_2 \leq \max\{u_1,u_2\}\\
    \Leftrightarrow\quad&\min\{\lambda_{1-},\lambda_{2-}\} < \lambda_C < \max\{\lambda_{1+},\lambda_{2+}\}.
\end{align*}
Let us now consider the situation that we have an isolated shock (w.l.o.g.) corresponding to $\lambda_{\mu-}$ moving with speed $S$.
A priori it is not obvious at all whether we may encounter the situation that $u = S$. In the following we will show that this is not possible.
In this situation we have a discontinuity with the tangential eigenvalues $\lambda_C^- = \lambda_C^+ = S$.
Here a superscript $-$ refers to the state $\mathbf{W}^-$ left of the discontinuity and a superscript $+$ refers to the state $\mathbf{W}^+$ right of the discontinuity, respectively.
According to the results obtained by Keyfitz et al.\ \cite{Keyfitz1980} cited above this corresponds to the second situation.
Thus we would further need $n-1$ incoming and $n-1$ outgoing characteristics, with $n=5$ in our case.
If we consider the results given in \cite{Andrianov2003} we obtain the same results with $m = 5$ jump conditions, $c=2$ tangential eigenvalues, $N = 2n + 1 = 2\cdot 5 + 1 = 11$ unknowns
and thus $i = 4$ incoming and $o = m - 1 = 4$ outgoing characteristics. Hence it is clear that for an evolutionary discontinuity we must not have further coinciding eigenvalues.
Let us picture the situation more precisely and assume w.l.o.g. that we have a shock corresponding to $\lambda_{1-}$.
Thus we would have
\begin{align}
    \lambda_{1-}^- > S > \lambda_{1-}^+\quad\text{and}\quad \lambda_C^- = \lambda_C^+ = S.\label{cond_shock_cont}
\end{align}
We can conclude immediately that $\lambda_{1+}^- > S$ and $S > \lambda_{2-}^-$. Using the continuity of $\alpha_1Q_1$ it follows that $\lambda_{1+}^+ > S$ and hence $S > \lambda_{2-}^+$.
Since we now already have four ingoing characteristics it follows that $\lambda_{2+}^-$ and $\lambda_{2+}^+$ must be outgoing characteristics. The situation
\begin{align*}
    \mathcal{I} &= \{\lambda_{1-}^-,\lambda_{1-}^+,\lambda_{1+}^-,\lambda_{2-}^+\},\\
    \mathcal{C} &= \{\lambda_C^-,\lambda_C^+\},\\
    \mathcal{O} &= \{\lambda_{1+}^+,\lambda_{2-}^-,\lambda_{2+}^-,\lambda_{2+}^+\}
\end{align*}
is depicted as an example in Figure \ref{fig:shock_contact}.
So clearly this situation is excluded, since the information related to the eigenvalue $\lambda_{2+}$ is coming out of the discontinuity.
Another short argument would be to say that we are in a strictly hyperbolic situation on both sides of the discontinuity and hence only a classical wave is allowed.
Note that this is not by far that obvious for systems with multiple linearly degenerated fields.
Thus we can state that in our system an evolutionary discontinuity with $u = S$ is a contact and vice versa.
Even more important is the statement that due to $u\neq S$ for a shock we always have the continuity of $\alpha$ across the shock.\\
This also gives another view on the results obtained above in Section \ref{sec:degen_adm_cond} and in particular the situation of coinciding eigenvectors described by (\ref{cond:absent_contact}).
In this particular situation $\alpha$ will not even jump across $\mathbf{R}_C$ and thus remains constant in the complete fan.
Hence we can exclude this situation by simply prescribing different values for $\alpha$ initially.
Or in other words this situation may only occur when $\alpha_L = \alpha_R$ holds for the initial states of the Riemann problem.
We therefore could interpret the case (\ref{cond:absent_contact}) as a consequence of a redundant $\alpha$ equation.
With $\alpha$ constant everywhere we can reformulate the system as
\begin{align*}
    &\frac{\partial}{\partial t}\alpha = 0,\\
    &\frac{\partial}{\partial t}\tilde{\bb{W}} + \frac{\partial}{\partial x}\tilde{\bb{F}}(\alpha,\tilde{\bb{W}}) = 0\\
    \text{with}\quad &\tilde{\bb{W}} = (w_2,w_3,w_4,w_5)^T\\
    \text{and}\quad &\tilde{\bb{F}}(\alpha,\tilde{\bb{W}}) = (F_2(\alpha,\tilde{\bb{W}}),F_3(\alpha,\tilde{\bb{W}}),F_4(\alpha,\tilde{\bb{W}}),F_5(\alpha,\tilde{\bb{W}})))^T.
\end{align*}
The situation $\lambda_{i\pm} = \lambda_C$ then corresponds to the case of \emph{hyperbolic resonance} discussed by Isaacson and Temple \cite{Isaacson1992}.
In the literature you also find the phrasing \emph{parabolic degeneracy} or \emph{weak hyperbolicity} for missing eigenvectors, but we think resonance is the term best suited here.
\subsubsection{Contact inside Rarefaction}
\begin{figure}[h!]
    \center
    \begin{tikzpicture}
        \begin{axis}[xlabel={$x$},axis x line=bottom,ylabel={$t$},axis y line=center,ymin=0,ymax=3.5,domain=-5:5,axis equal image,ytick=\empty,xtick=\empty,clip mode=individual]
            \fill[blue!20!white,opacity = 0.3] (axis cs: 0,0) -- (axis cs: -3,3) -- (axis cs: -3/5,3) -- cycle;
            \fill[blue!20!white,opacity = 0.5] (axis cs: 0,0) -- (axis cs: 6/5,3) -- (axis cs: 18/5,3) -- cycle;
            \addplot[mark=none,smooth,blue,thick,solid,domain=-3:0] (\x,{-\x});
            \node[black,above] at (axis cs:-3,3) {$\lambda_{\mu-}^{(L)}$};
            \addplot[mark=none,smooth,blue,thick,solid,domain=-12/5:0] (\x,{-5*\x/4});
            \addplot[mark=none,smooth,blue,thick,solid,domain=-9/5:0] (\x,{-5*\x/3});
            \addplot[mark=none,smooth,blue,thick,solid,domain=-6/5:0] (\x,{-2.5*\x});
            \addplot[mark=none,smooth,blue,thick,solid,domain=-3/5:0] (\x,{-5*\x});
            \node[blue,above] at (axis cs:-1/2,3) {$\lambda_{\mu-}^-$};
            \addplot[mark=none,smooth,red,thick,solid,domain=0:3/5] (\x,{5*\x});
            \node[red,above] at (axis cs:3/5,3) {$S$};
            \addplot[mark=none,smooth,blue,thick,solid,domain=0:6/5] (\x,{2.5*\x});
            \node[blue,above] at (axis cs:1.2,3) {$\lambda_{\mu-}^+$};
            \addplot[mark=none,smooth,blue,thick,solid,domain=0:9/5] (\x,{5*\x/3});
            \addplot[mark=none,smooth,blue,thick,solid,domain=0:12/5] (\x,{5*\x/4});
            \addplot[mark=none,smooth,blue,thick,solid,domain=0:3] (\x,{\x});
            \addplot[mark=none,smooth,blue,thick,solid,domain=0:18/5] (\x,{5*\x/6});
            \node[black,above] at (axis cs:18/5,3) {$\lambda_{\mu-}^{(R)}$};
            \addplot[mark=none,smooth,green,thick,solid,domain=-2/5:2/5] (\x,{5*\x/2 + 1});
            \node[green,below] at (axis cs:-2/5,0) {$\lambda_{\nu+}^-$};
            \addplot[mark=none,smooth,green,thick,solid,domain=2/5:6/5] (\x,{-5*\x/2 + 3});
            \node[green,below] at (axis cs:6/5,0) {$\lambda_{\nu+}^+$};
            \addplot[mark=none,smooth,orange,thick,solid,domain=-2:2/5] (\x,{5*\x/6 + 5/3});
            \node[orange,below] at (axis cs:-2,0) {$\lambda_{\mu+}^-$};
            \addplot[mark=none,smooth,orange,thick,solid,domain= 2/5:14/5] (\x,{5*\x/12 + 11/6});
            \node[orange,above] at (axis cs:14/5,3) {$\lambda_{\mu+}^+$};
            \addplot[mark=none,smooth,magenta,thick,solid,domain=-1:2/5] (\x,{-5*\x/7 + 16/7});
            \node[magenta,above] at (axis cs:-1.2,3) {$\lambda_{\nu-}^-$};
            \addplot[mark=none,smooth,magenta,thick,solid,domain= 2/5:2] (\x,{-5*\x/4 + 5/2});
            \node[magenta,below] at (axis cs:2,0) {$\lambda_{\nu-}^+$};
        \end{axis}
    \end{tikzpicture}
    \caption{Contact inside a rarefaction (example).}
    \label{fig:contact_in_rarefaction_L}
\end{figure}
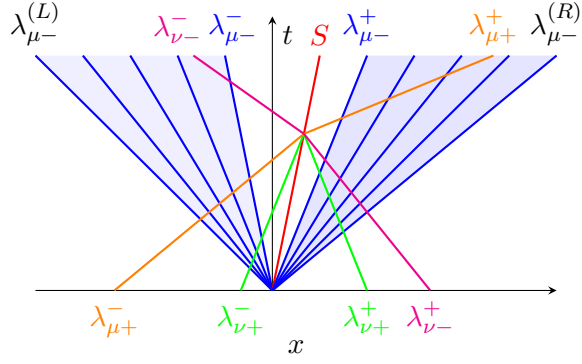

Let us assume that a contact lies inside a rarefaction wave, see e.g. the sketch shown in Figure \ref{fig:contact_in_rarefaction_L}.
It cannot be attached to one side (or both) of the rarefaction, because then we would have the situation of hyperbolic resonance discussed above (i.e. coincidence of eigenvalues and eigenvectors).
Thus the contact will tear the rarefaction wave into two parts. Assume we have a rarefaction wave corresponding to $\lambda_{\mu-}$ and thus we consider the following situation
\begin{align*}
    \lambda_{\mu-}^- < \lambda_C^- = S = \lambda_C^+ < \lambda_{\mu-}^+.
\end{align*}
According to the results obtained above concerning the admissibility of discontinuities we have $c = 2$ coinciding characteristics and thus need $i = o = 4$ characteristics going in and out.
From the given relation for the eigenvalues $\lambda_{\mu-}$ and $\lambda_C$ we can directly conclude $\lambda_{\mu+}^+ > S$.
We can then conclude that that $\lambda_{\nu-}^+ < S$.
Using the inequalities for $\lambda_{\mu-}$ we see that $Q_\mu < 0$ across the contact. Since $0 = Q = \alpha_\mu Q_\mu + \alpha_\nu Q_\nu$ we thus have $Q_\nu > 0$.
Using $Q_\mu < 0$ we yield $\lambda_{\mu+}^- > S$. From $Q_\nu > 0$ we yield $\lambda_{\nu-}^- < S$ and hence finally $\lambda_{\nu+}^+ < S$.
Altogether we therefore obtain
\begin{align*}
    \mathcal{I} = \{\lambda_{\mu+}^-,\lambda_{\nu-}^+,\lambda_{\nu+}^-,\lambda_{\nu+}^+\},\quad
    \mathcal{C} = \{\lambda_C^-,\lambda_C^+\}\quad\text{and}\quad
    \mathcal{O} = \{\lambda_{\mu-}^-,\lambda_{\mu-}^+,\lambda_{\mu+}^+,\lambda_{\nu-}^-\}.
\end{align*}
Such a discontinuity violates the admissibility criteria given above due to the fact,
that the characteristics $\lambda_{\mu-}^-$ and $\lambda_{\mu-}^+$ do not contribute any information to the discontinuity.
Note that as mentioned before the phrasing \emph{outgoing characteristic} seems to be not quite suited here for these two eigenvalues.
However, the results remain valid and one could use the notation of slow and fast characteristics with respect to the side of the discontinuity, cf.\ \cite{Freistuehler1992a}.
In this sense the situation is analogue to that of a contact coinciding with a shock.
\subsection{Overlapping Rarefaction Waves}
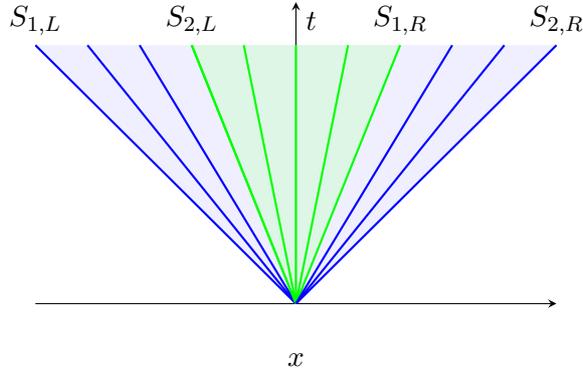
\begin{figure}[h!]
    \center
    \begin{tikzpicture}
        \begin{axis}[xlabel={$x$},axis x line=bottom,ylabel={$t$},axis y line=center,ymin=0,ymax=3.5,domain=-5:5,axis equal image,ytick=\empty,xtick=\empty,clip mode=individual]
            \fill[blue!20!white,opacity = 0.3] (axis cs: 0,0) -- (axis cs: -3,3) -- (axis cs: 3,3) -- cycle;
            \fill[green!20!white,opacity = 0.5] (axis cs: 0,0) -- (axis cs: -6/5,3) -- (axis cs: 6/5,3) -- cycle;
            \addplot[mark=none,smooth,green,thick,solid,domain=-6/5:0] (\x,{-2.5*\x});
            \addplot[mark=none,smooth,blue,thick,solid,domain=-9/5:0] (\x,{-5*\x/3});
            \addplot[mark=none,smooth,blue,thick,solid,domain=-12/5:0] (\x,{-5*\x/4});
            \addplot[mark=none,smooth,blue,thick,solid,domain=-3:0] (\x,{-\x});
            \node[black,above] at (axis cs:-3,3) {$S_{1,L}$};
            \addplot[mark=none,smooth,green,thick,solid,domain=-3/5:0] (\x,{-5*\x});
            \node[black,above] at (axis cs:-6/5,3) {$S_{2,L}$};
            \addplot[mark=none,smooth,green,thick,solid,domain=-6/5:0] (\x,{-2.5*\x});
            \addplot[mark=none,smooth,green,thick,solid] coordinates {(0, 0) (0,3)};
            \addplot[mark=none,smooth,green,thick,solid,domain=0:3/5] (\x,{5*\x});
            \addplot[mark=none,smooth,green,thick,solid,domain=0:6/5] (\x,{2.5*\x});
            \node[black,above] at (axis cs:1.2,3) {$S_{1,R}$};
            \addplot[mark=none,smooth,blue,thick,solid,domain=0:9/5] (\x,{5*\x/3});
            \addplot[mark=none,smooth,blue,thick,solid,domain=0:12/5] (\x,{5*\x/4});
            \addplot[mark=none,smooth,blue,thick,solid,domain=0:3] (\x,{\x});
            \node[black,above] at (axis cs:3,3) {$S_{2,R}$};
        \end{axis}
    \end{tikzpicture}
    \caption{Overlapping rarefaction waves (example).}
    \label{fig:overlap_rarefaction}
\end{figure}

A possible wave configuration might be two rarefaction waves that overlap, see Figure \ref{fig:overlap_rarefaction}.
In this situation one rarefaction wave belongs to phase one and the other to phase two, respectively.
The rarefaction fans are cones given by
\begin{align*}
    \mathcal{C}_1 &= \left\{(t,x)| 0 < t, S_{1,L}t \leq x \leq S_{1,R}t \right\}\\
    \mathcal{C}_2 &= \left\{(t,x)| 0 < t, S_{2,L}t \leq x \leq S_{2,R}t \right\}\\
    \mathcal{C}^\ast &= \mathcal{C}_1\cap\mathcal{C}_2
\end{align*}
In the case of an empty intersection we have the classical situation and the solution is obtained using the corresponding eigenvector.
If however, $\mathcal{C}^\ast\neq\emptyset$ we have
\begin{align*}
    \mathcal{C}^\ast = \left\{(t,x)| 0 < t, \max\{S_{1,L},S_{2,L}\}t \leq x \leq \min\{S_{1,R},S_{2,R}\}t\right\}
\end{align*}
Due to the special structure of the eigenvectors the obtained invariants remain unchanged and the formulas (\ref{left_raref_fan}) and (\ref{right_raref}) stay valid.
Indeed if we write down the (homogeneous) system using the primitive variables $\bb{W} = (\alpha_1,\rho_1,\rho_2,u_1,u_2)$ we have the Jacobian (\ref{eq:jacobian_prim_var}).
For a rarefaction wave $\alpha_1$ remains constant and thus the system simplifies to
\begin{subequations}\label{sys:gpr_1d_v3}
    \begin{align}
        \frac{\partial\rho_1}{\partial t} + \frac{\partial\rho_1 u_1}{\partial x} &= 0,\label{eq:partial_mass1_balance_v2}\\
        \frac{\partial\rho_1u_1}{\partial t} + \frac{\partial\left(\rho_1 u_1^2 + p_1\right)}{\partial x} &= 0,\label{eq:partial_mom1_balance}\\
        \frac{\partial\rho_2}{\partial t} + \frac{\partial\rho_2 u_2}{\partial x} &= 0,\label{eq:partial_mass2_balance}\\
        \frac{\partial\rho_2u_2}{\partial t} + \frac{\partial\left(\rho_2 u_2^2 + p_2\right)}{\partial x} &= 0.\label{eq:partial_mom2_balance}
    \end{align}
\end{subequations}
Hence the system decouples into two barotropic Euler systems for each phase and the rarefaction waves can be obtained individually.
The solution is then obtained as the superposition of the individual solutions.
\subsection{Shock interacting with a Rarefaction Wave}
\begin{figure}[h!]
    \center
    \begin{tikzpicture}
        \begin{axis}[xlabel={$x$},axis x line=bottom,ylabel={$t$},axis y line=center,ymin=0,ymax=3.5,domain=-5:5,axis equal image,ytick=\empty,xtick=\empty]
            %
            \fill[blue!20!white,opacity = 0.3] (axis cs: 0,0) -- (axis cs: -0.3,3) -- (axis cs: -3,3) -- cycle;
            \addplot[mark=none,smooth,blue,thick,solid,domain=-0.3:0] (\x,{-10*\x});
            \addplot[mark=none,smooth,blue,thick,solid,domain=-3/5:0] (\x,{-5*\x});
            \addplot[mark=none,smooth,blue,thick,solid,domain=-9/5:0] (\x,{-5*\x/3});
            \addplot[mark=none,smooth,blue,thick,solid,domain=-12/5:0] (\x,{-5*\x/4});
            \addplot[mark=none,smooth,blue,thick,solid,domain=-3:0] (\x,{-\x});
            \addplot[mark=none,smooth,red,thick,solid,domain=-6/5:0] (\x,{-2.5*\x});
            \node[black,above] at (axis cs:-1.2,3) {$\lambda_{\nu-}^- = S = \lambda_{\nu-}^+$};
        \end{axis}
    \end{tikzpicture}
    \caption{Shock inside a rarefaction wave - Case (i).}
    \label{fig:shock_in_raref_left_i}
\end{figure}
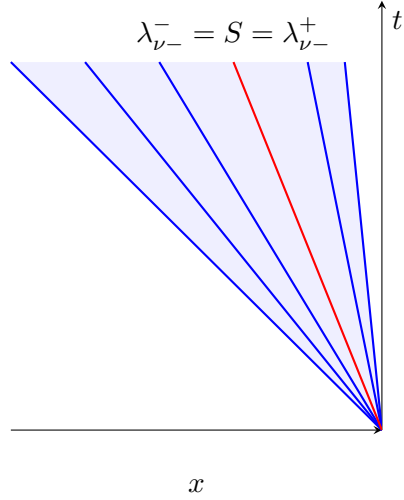

Another situation that can occur is a shock which lies inside a rarefaction fan, see Figure \ref{fig:shock_in_raref_left_i}.
%
%
%
Let us assume for the moment that we have a left shock corresponding to phase $\mu\in\{1,2\}$ and a left rarefaction corresponding to phase $\nu\in\{1,2\},\,\mu\neq\nu$.
Clearly, the phase $\mu$ is only affected by the shock wave due to the structure of the eigenvector corresponding to the rarefaction wave.
For the phase $\nu$ the situation is more complicated. We have for the eigenvalue $\lambda_{\nu-} = \xi := x/t$ corresponding to the rarefaction wave
\begin{align*}
    \lambda_{\nu-}(\bb{W}_L) \leq \xi \leq \lambda_{\nu-}(\bb{W_R})
\end{align*}
where $\bb{W}_L$ and $\bb{W}_R$ denote the states left and right of the rarefaction wave.
Let us denote the shock speed with $S$ and quantities left of the shock are denoted with a superscript $-$ and a $+$ when they are on the right, respectively.
There are four cases which are possible in this situation
\begin{enumerate}[(i)]
    \item $\lambda_{\nu-}^- = S = \lambda_{\nu-}^+$,
    \item $\lambda_{\nu-}^- < S < \lambda_{\nu-}^+$,
    \item $\lambda_{\nu-}^- = S < \lambda_{\nu-}^+$,
    \item $\lambda_{\nu-}^- < S = \lambda_{\nu-}^+$.
\end{enumerate}
Further we demand the Lax condition $\lambda_{\mu-}^- > S > \lambda_{\mu-}^+$.\\
\newline
\textbf{\underline{Case (i):}} The first case can be excluded since it implies linear degeneracy of the field $(\lambda_{\nu-},\bb{R}_{\nu-})$, see \cite{Keyfitz1980,Freistuhler1991}.
This is obviously not the case as long as we have $\mathcal{G}_\nu \neq 0$, which we may assume for our EOS.
Further discussion of the fundamental derivative can be found in \cite{Menikoff1989,Mueller2006}.\\
\newline
\textbf{\underline{Case (ii):}} Considering the characteristics in the second case we obviously have
\begin{align*}
    \lambda_{\nu-}^- < S < \lambda_{\mu-}^-\quad\text{and}\quad \lambda_{\mu-}^+ < S < \lambda_{\nu-}^+.
\end{align*}
Thus we directly obtain
\begin{align*}
    S < \lambda_{\nu-}^+ < u_\nu^+ < \lambda_{\nu+}^+\quad\text{and}\quad S < \lambda_{\mu-}^- < u_\nu^- < \lambda_{\mu+}^-.
\end{align*}
Due to the continuity of the mass fluxes we further yield
\begin{align*}
    u_\nu^+ > S\quad\Leftrightarrow\quad 0 > Q_\nu^+ = Q_\nu^- \quad\Leftrightarrow\quad u_\nu^- > S\quad\Rightarrow\quad \lambda_{\nu+}^- > S,\\
    u_\mu^- > S\quad\Leftrightarrow\quad 0 > Q_\mu^- = Q_\mu^+ \quad\Leftrightarrow\quad u_\mu^+ > S\quad\Rightarrow\quad \lambda_{\mu+}^+ > S.
\end{align*}
Since $u$ is a convex combination of the phase velocities we also conclude $\lambda_C^- > S$ and $\lambda_C^+ > S$.
Summarizing we have the following situation
\begin{align*}
    \mathcal{I} = \{\lambda_{\mu-}^-,\lambda_{\mu-}^+,\lambda_{\mu+}^-,\lambda_{\nu+}^-,\lambda_C^-\},\quad
    \mathcal{C} = \emptyset,\quad
    \mathcal{O} = \{\lambda_{\nu-}^-,\lambda_{\nu-}^+,\lambda_{\nu+}^+,\lambda_{\mu+}^+,\lambda_C^+\}.
\end{align*}
Therefore this situation is not admissible.
According to the admissibility criteria we have $N = 11$ unknowns $m = 5$ equations and thus we would need $i = 6$ incoming and $o = 4$ outgoing characteristics.
\begin{figure}[h!]
    \center
    \begin{tikzpicture}
        \begin{axis}[xlabel={$x$},axis x line=bottom,ylabel={$t$},axis y line=center,ymin=0,ymax=3.5,domain=-5:5,axis equal image,ytick=\empty,xtick=\empty]
            %
            \fill[blue!20!white,opacity = 0.3] (axis cs: 0,0) -- (axis cs: -9/5,3) -- (axis cs: -3,3) -- cycle;
            \fill[blue!20!white,opacity = 0.3] (axis cs: 0,0) -- (axis cs: -3/5,3) -- (axis cs: -0.3,3) -- cycle;
            \addplot[mark=none,smooth,blue,thick,solid,domain=-0.3:0] (\x,{-10*\x});
            \addplot[mark=none,smooth,blue,thick,solid,domain=-3/5:0] (\x,{-5*\x});
            \node[black,above] at (axis cs:-3/5,3) {$\lambda_{\nu-}^+$};
            \addplot[mark=none,smooth,blue,thick,solid,domain=-9/5:0] (\x,{-5*\x/3});
            \node[black,above] at (axis cs:-9/5,3) {$\lambda_{\nu-}^-$};
            \addplot[mark=none,smooth,blue,thick,solid,domain=-12/5:0] (\x,{-5*\x/4});
            \addplot[mark=none,smooth,blue,thick,solid,domain=-3:0] (\x,{-\x});
            \addplot[mark=none,smooth,red,thick,solid,domain=-6/5:0] (\x,{-2.5*\x});
            \node[black,above] at (axis cs:-1.2,3) {$S$};
        \end{axis}
    \end{tikzpicture}
    \caption{Shock inside a rarefaction wave - Case (ii).}
    \label{fig:shock_in_raref_left_ii}
\end{figure}
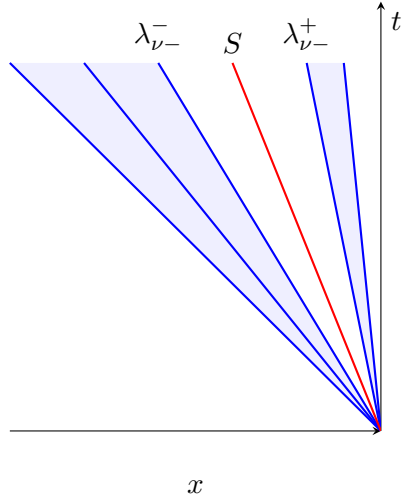
\newline
\textbf{\underline{Case (iii):}} Now, for the third case we have
\begin{align*}
    \lambda_{\nu-}^- &< S < \lambda_{\mu-}^- < \lambda_{\mu+}^-\\
    \text{and}\quad \lambda_{\mu-}^+ &< S = \lambda_{\nu-}^+ < \lambda_{\nu+}^+.
\end{align*}
%
We have for the mass flux of phase $\mu$ that $Q_\mu < 0$. Thus we have $u_\mu^+ > S$ and hence $\lambda_{\mu+}^+ > S$.
For the mass flux of phase $\nu$ we obtain
\begin{align*}
    &\lambda_{\nu-}^- = u_\nu^- - a_\nu^- < S = \lambda_{\nu-}^+ = u_\nu^+ - a_\nu^+\\
    \Leftrightarrow\quad \phantom{-}&u_\nu^- - S < a_\nu^-\;\wedge\; a_\nu^+ = u_\nu^+ - S\\
    \Leftrightarrow\quad -&\rho_\nu^-(u_\nu^- - S) > -\rho_\nu^-a_\nu^-\;\wedge\; -\rho_\nu^+(u_\nu^+ - S) = -\rho_\nu^+a_\nu^+\\
    \Leftrightarrow\quad -&\rho_\nu^-a_\nu^- < Q_\nu = -\rho_\nu^+a_\nu^+ < 0.
\end{align*}
Thus we also have $Q < 0$ and hence
\begin{align*}
    0 &> Q = -\rho^-(u^- - S)\quad\Leftrightarrow\quad \lambda_C^- = u^- > S\\
    0 &> Q = -\rho^+(u^+ - S)\quad\Leftrightarrow\quad \lambda_C^+ = u^+ > S.
\end{align*}
We therefore obtain that $\lambda_C^-$ is an ingoing characteristic and $\lambda_C^+$ an outgoing characteristic.
In the situation under consideration we have as before $N = 11$ unknowns $m = 5$ equations.
Additionally we have $c = 1$ coinciding wave and thus we need $i = 5$ incoming and $o = 4$ outgoing characteristics. This implies that $\lambda_{\nu+}^-$ is also an ingoing characteristic.
Summing up we have the following situation
\begin{align*}
    \mathcal{I} = \{\lambda_{\mu-}^-,\lambda_{\mu-}^+,\lambda_{\mu+}^-,\lambda_{\nu+}^-,\lambda_C^-\},\quad
    \mathcal{C} = \{\lambda_{\nu-}^+\}\quad\text{and}\quad
    \mathcal{O} = \{\lambda_{\mu+}^+,\lambda_{\nu-}^-,\lambda_{\nu+}^+,\lambda_C^+\}.
\end{align*}
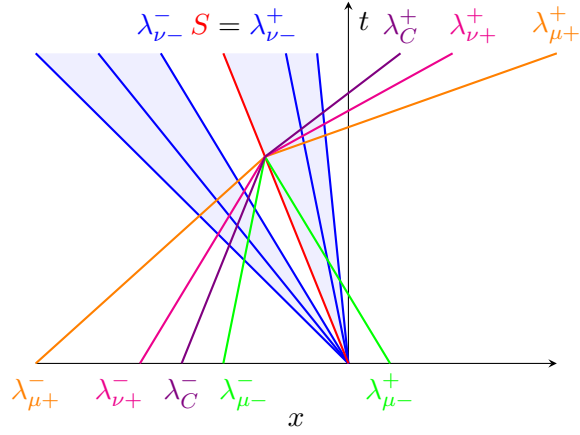
\begin{figure}[h!]
    \center
    %
    \begin{tikzpicture}
        \begin{axis}[xlabel={$x$},axis x line=bottom,ylabel={$t$},axis y line=center,ymin=0,ymax=3.5,domain=-5:5,axis equal image,ytick=\empty,xtick=\empty,clip mode=individual]
            %
            \fill[blue!20!white,opacity = 0.3] (axis cs: 0,0) -- (axis cs: -1.2,3) -- (axis cs: -0.3,3) -- cycle;
            \fill[blue!20!white,opacity = 0.3] (axis cs: 0,0) -- (axis cs: -9/5,3) -- (axis cs: -3,3) -- cycle;
            \addplot[mark=none,smooth,blue,thick,solid,domain=-0.3:0] (\x,{-10*\x});
            \addplot[mark=none,smooth,blue,thick,solid,domain=-3/5:0] (\x,{-5*\x});
            \addplot[mark=none,smooth,blue,thick,solid,domain=-9/5:0] (\x,{-5*\x/3});
            \node[black,above] at (axis cs:-9/5,3) {$\textcolor{blue}{\lambda_{\nu-}^-}$};
            \addplot[mark=none,smooth,blue,thick,solid,domain=-12/5:0] (\x,{-5*\x/4});
            \addplot[mark=none,smooth,blue,thick,solid,domain=-3:0] (\x,{-\x});
            \addplot[mark=none,smooth,red,thick,solid,domain=-6/5:0] (\x,{-2.5*\x});
            \node[black,above] at (axis cs:-1,3) {$\textcolor{red}{S} = \textcolor{blue}{\lambda_{\nu-}^+}$};
            \addplot[mark=none,smooth,green,thick,solid,domain=-4/5:2/5] (\x,{-5*\x/3 + 2/3});
            \node[green,below] at (axis cs:2/5,0) {$\lambda_{\mu-}^+$};
            \addplot[mark=none,smooth,green,thick,solid,domain=-6/5:-4/5] (\x,{5*\x + 6});
            \node[green,below] at (axis cs:-1,0) {$\lambda_{\mu-}^-$};
            \addplot[mark=none,smooth,orange,thick,solid,domain=-4/5:2] (\x,{5*\x/14 + 16/7});
            \node[orange,above] at (axis cs:2,3) {$\lambda_{\mu+}^+$};
            \addplot[mark=none,smooth,orange,thick,solid,domain=-3:-4/5] (\x,{10*\x/11 + 30/11});
            \node[orange,below] at (axis cs:-3,0) {$\lambda_{\mu+}^-$};
            \addplot[mark=none,smooth,magenta,thick,solid,domain=-4/5:1] (\x,{5*\x/9 + 22/9});
            \node[magenta,above] at (axis cs:1.2,3) {$\lambda_{\nu+}^+$};
            \addplot[mark=none,smooth,magenta,thick,solid,domain=-2:-4/5] (\x,{5*\x/3 + 10/3});
            \node[magenta,below] at (axis cs:-2.2,0) {$\lambda_{\nu+}^-$};
            \addplot[mark=none,smooth,violet,thick,solid,domain=-4/5:1/2] (\x,{10*\x/13 + 34/13});
            \node[violet,above] at (axis cs:1/2,3) {$\lambda_C^+$};
            \addplot[mark=none,smooth,violet,thick,solid,domain= -8/5:-4/5] (\x,{5*\x/2 + 4});
            \node[violet,below] at (axis cs:-8/5,0) {$\lambda_C^-$};
        \end{axis}
    \end{tikzpicture}
    \caption{Shock inside a rarefaction wave - Case (iii).}
    \label{fig:shock_in_raref_left_iii}
\end{figure}
\newline
\textbf{\underline{Case (iv):}} Finally, for the fourth case we have
\begin{align*}
    \lambda_{\nu-}^- = S < \lambda_{\mu-}^-\quad\text{and}\quad \lambda_{\mu-}^+ < S < \lambda_{\nu-}^+.
\end{align*}
Using similar arguments as for the third case we obtain the following situation
\begin{align*}
    \mathcal{I} = \{\lambda_{\mu-}^-,\lambda_{\mu-}^+,\lambda_{\mu+}^-,\lambda_{\nu+}^-,\lambda_C^-\},\quad
    \mathcal{C} = \{\lambda_{\nu-}^-\}\quad\text{and}\quad
    \mathcal{O} = \{\lambda_{\mu+}^+,\lambda_{\nu-}^+,\lambda_{\nu+}^+,\lambda_C^+\}.
\end{align*}
\begin{figure}[h!]
    \center
    %
    \begin{tikzpicture}
        \begin{axis}[xlabel={$x$},axis x line=bottom,ylabel={$t$},axis y line=center,ymin=0,ymax=3.5,domain=-5:5,axis equal image,ytick=\empty,xtick=\empty,clip mode=individual]
            %
            \fill[blue!20!white,opacity = 0.3] (axis cs: 0,0) -- (axis cs: -1.2,3) -- (axis cs: -3,3) -- cycle;
            \fill[blue!20!white,opacity = 0.3] (axis cs: 0,0) -- (axis cs: -3/5,3) -- (axis cs: -0.3,3) -- cycle;
            \addplot[mark=none,smooth,blue,thick,solid,domain=-0.3:0] (\x,{-10*\x});
            \addplot[mark=none,smooth,blue,thick,solid,domain=-3/5:0] (\x,{-5*\x});
            \node[black,above] at (axis cs:-1/2,3) {$\textcolor{blue}{\lambda_{\nu-}^+}$};
            \addplot[mark=none,smooth,blue,thick,solid,domain=-9/5:0] (\x,{-5*\x/3});
            \addplot[mark=none,smooth,blue,thick,solid,domain=-12/5:0] (\x,{-5*\x/4});
            \addplot[mark=none,smooth,blue,thick,solid,domain=-3:0] (\x,{-\x});
            \addplot[mark=none,smooth,red,thick,solid,domain=-6/5:0] (\x,{-2.5*\x});
            \node[black,above] at (axis cs:-1.4,3) {$\textcolor{blue}{\lambda_{\nu-}^-} = \textcolor{red}{S}$};
            \addplot[mark=none,smooth,green,thick,solid,domain=-4/5:2/5] (\x,{-5*\x/3 + 2/3});
            \node[green,below] at (axis cs:2/5,0) {$\lambda_{\mu-}^+$};
            \addplot[mark=none,smooth,green,thick,solid,domain=-6/5:-4/5] (\x,{5*\x + 6});
            \node[green,below] at (axis cs:-1,0) {$\lambda_{\mu-}^-$};
            \addplot[mark=none,smooth,orange,thick,solid,domain=-4/5:2] (\x,{5*\x/14 + 16/7});
            \node[orange,above] at (axis cs:2,3) {$\lambda_{\mu+}^+$};
            \addplot[mark=none,smooth,orange,thick,solid,domain=-3:-4/5] (\x,{10*\x/11 + 30/11});
            \node[orange,below] at (axis cs:-3,0) {$\lambda_{\mu+}^-$};
            \addplot[mark=none,smooth,magenta,thick,solid,domain=-4/5:1] (\x,{5*\x/9 + 22/9});
            \node[magenta,above] at (axis cs:1.2,3) {$\lambda_{\nu+}^+$};
            \addplot[mark=none,smooth,magenta,thick,solid,domain=-2:-4/5] (\x,{5*\x/3 + 10/3});
            \node[magenta,below] at (axis cs:-2.2,0) {$\lambda_{\nu+}^-$};
            \addplot[mark=none,smooth,violet,thick,solid,domain=-4/5:1/2] (\x,{10*\x/13 + 34/13});
            \node[violet,above] at (axis cs:1/2,3) {$\lambda_C^+$};
            \addplot[mark=none,smooth,violet,thick,solid,domain= -8/5:-4/5] (\x,{5*\x/2 + 4});
            \node[violet,below] at (axis cs:-8/5,0) {$\lambda_C^-$};
        \end{axis}
    \end{tikzpicture}
    \caption{Shock inside a rarefaction wave - Case (iv).}
    \label{fig:shock_in_raref_left_iv}
\end{figure}
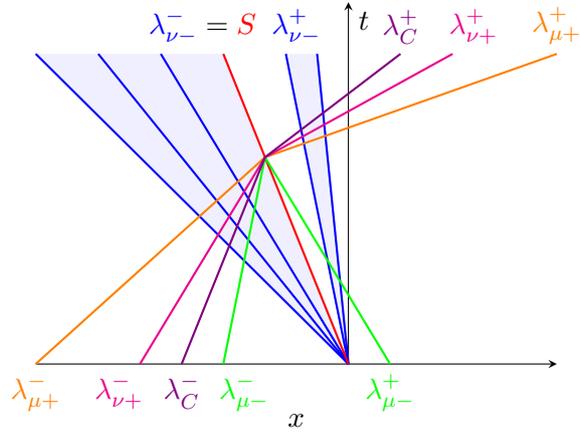
Now we have two possible configuration which seem to be allowed. A priori it is not obvious whether one of these cases can be ruled out or if both may occur.
Thus we will make use of the energy inequality (\ref{ineq_disc:baro_energy_v1}) to investigate both cases.\\
\newline
\underline{\textbf{Energy Inequality Case (iii):}} As in (\ref{ineq_disc:baro_energy_v2}) we use the entropy inequality for the phase $\nu$ which defines the rarefaction, i.e.
\begin{align*}
    0 \geq -Q\dbl \Psi_\nu + \frac{1}{2}\left(u_\nu - S\right)^2 \dbr.
\end{align*}
From the characteristics $\lambda_C^-,\lambda_C^+$ we directly conclude $Q < 0$ and hence we expect
\begin{align}
    0 \geq \dbl \Psi_\nu + \frac{1}{2}\left(u_\nu - S\right)^2 \dbr\label{adm_cond_shock_in_raref}
\end{align}
for this configuration in order to be admissible. We rewrite the kinetic energy in terms of the mass flux and assume the right state to be given.
Thus it is possible to write the jump bracket as a function in the density $\rho_\nu^-$ of the left side of the discontinuity, i.e.
\begin{align*}
    0 &\geq \dbl \Psi_\nu + \frac{1}{2}\left(u_\nu - S\right)^2 \dbr = \Psi_\nu^+ - \Psi_\nu(\rho_\nu^-) + \frac{Q_\nu^2}{2}\left(\frac{1}{(\rho_\nu^+)^2} - \frac{1}{(\rho_\nu^-)^2}\right),\\
    \text{with}\quad Q_\nu &= -\rho_\nu^+a_\nu(\rho_\nu^+).
\end{align*}
Investigating this function gives the following
\begin{align*}
    f(\rho) &:= \Psi_\nu^+ - \Psi_\nu(\rho) + \frac{Q_\nu^2}{2}\left(\frac{1}{(\rho_\nu^+)^2} - \frac{1}{\rho^2}\right),\quad f(\rho_\nu^+) = 0,\\
    f'(\rho) &= -\frac{a_\nu(\rho)^2}{\rho} + \frac{Q_\nu^2}{\rho^3},\quad f'(\rho_\nu^+) = 0,\\
    f''(\rho) &= -2\frac{a_\nu(\rho)^2}{\rho^2}\mathcal{G}_\nu + 3\frac{a_\nu(\rho)^2}{\rho^2} - 3\frac{Q_\nu^2}{\rho^4},\quad f''(\rho_\nu^+) < 0.
\end{align*}
In order to get further insight we need to discuss the function $g(\rho) = \rho a_\nu(\rho)$ which is sometimes called \emph{Lagrangian wave speed} \cite{Menikoff1989}.
More precisely, we already discussed it investigating the characteristic fields and we have for phase $\nu$
\begin{align}
    g(\rho) = \rho a_\nu(\rho) > 0\quad\text{and}\quad g'(\rho) = a_\nu\mathcal{G}_\nu > 0,\;\forall \rho > 0.\label{rel:lagrangewavespeed}
\end{align}
Thus we can conclude for all $0 < \rho < \rho_\nu^+$
\begin{align*}
    f'(\rho) &= -\frac{a_\nu(\rho)^2}{\rho} + \frac{Q_\nu^2}{\rho^3} = \frac{1}{\rho^3}\left(Q_\nu^2 - \rho^2a_\nu(\rho)^2\right) > 0,\\
    \text{and}\quad f''(\rho) &= -2\frac{a_\nu(\rho)^2}{\rho^2}\mathcal{G}_\nu - \frac{3}{\rho^4}\left(Q_\nu^ 2 - \rho^2a_\nu(\rho)^2\right) < 0.
\end{align*}
Due to this monotonicity behaviour, $\rho_\nu^- < \rho_\nu^+$ and since $f(\rho_\nu^+) = 0$ we can conclude
\begin{align*}
    \dbl \Psi_\nu + \frac{1}{2}\left(u_\nu - S\right)^2 \dbr = f(\rho_\nu^-) < 0.
\end{align*}
Hence configuration (iii) respects the mathematical entropy inequality.\\
\newline
\underline{\textbf{Energy Inequality Case (iv):}} We use a similar argumentation as before.
Now we assume the left state to be given and write the jump bracket as a function in the density $\rho_\nu^+$ of the right side of the discontinuity, i.e.
\begin{align*}
    0 &\geq \dbl \Psi_\nu + \frac{1}{2}\left(u_\nu - S\right)^2 \dbr = \Psi_\nu(\rho_\nu^+) - \Psi_\nu^- + \frac{Q_\nu^2}{2}\left(\frac{1}{(\rho_\nu^+)^2} - \frac{1}{(\rho_\nu^-)^2}\right),\\
    \text{with}\quad Q_\nu &= -\rho_\nu^-a_\nu(\rho_\nu^-).
\end{align*}
Investigating this function gives the following
\begin{align*}
    f(\rho) &:= \Psi_\nu(\rho_\nu^+) - \Psi_\nu^- + \frac{Q_\nu^2}{2}\left(\frac{1}{\rho^2} - \frac{1}{(\rho_\nu^-)^2}\right),\quad f(\rho_\nu^-) = 0,\\
    f'(\rho) &= \frac{a_\nu(\rho)^2}{\rho} - \frac{Q_\nu^2}{\rho^3},\quad f'(\rho_\nu^-) = 0,\\
    f''(\rho) &= 2\frac{a_\nu(\rho)^2}{\rho^2}\mathcal{G}_\nu - 3\frac{a_\nu(\rho)^2}{\rho^2} + 3\frac{Q_\nu^2}{\rho^4},\quad f''(\rho_\nu^-) > 0.
\end{align*}
Again we use the monotonicity of the Lagrangian wave speed (\ref{rel:lagrangewavespeed}) to conclude for all $0 < \rho_\nu^- < \rho$
\begin{align*}
    f'(\rho) &= \frac{a_\nu(\rho)^2}{\rho} - \frac{Q_\nu^2}{\rho^3} = \frac{1}{\rho^3}\left(\rho^2a_\nu(\rho)^2 - Q_\nu^2\right) > 0,\\
    %
\end{align*}
Due to this monotonicity behaviour, $\rho_\nu^- < \rho_\nu^+$ and since $f(\rho_\nu^-) = 0$ we can conclude
\begin{align*}
    \dbl \Psi_\nu + \frac{1}{2}\left(u_\nu - S\right)^2 \dbr = f(\rho_\nu^+) > 0.
\end{align*}
Hence configuration (iv) violates the energy inequality and finally we can say that only case (iii) is possible.
The treatment for the case considering $\lambda_{\mu+}$ and $\lambda_{\nu+}$ is completely analogous showing that the shock and the rarefaction characteristic coincide now on the left side.\\
\newline
If we reconsider the arguments above it becomes clear that the side on which the rarefaction characteristic coincides with the discontinuity is basically defined by the mixture mass flux $Q$.
Thus we can extend the above results to the cases of a $\lambda_{\mu+}$ shock interacting with a $\lambda_{\nu-}$ rarefaction or the analogous case with $\lambda_{\mu-}$ and $\lambda_{\nu+}$.
For $Q > 0$ in the first case we exemplary have the situation
\begin{align*}
    \mathcal{I} = \{\lambda_{\mu-}^+,\lambda_{\mu+}^-,\lambda_{\mu+}^+,\lambda_{\nu+}^-,\lambda_C^+\},\quad
    \mathcal{C} = \{\lambda_{\nu-}^-\}\quad\text{and}\quad
    \mathcal{O} = \{\lambda_{\nu-}^+,\lambda_{\nu+}^+,\lambda_{\mu-}^-,\lambda_C^-\}.
\end{align*}
and for $Q < 0$
\begin{align*}
    \mathcal{I} = \{\lambda_{\mu-}^+,\lambda_{\mu+}^-,\lambda_{\mu+}^+,\lambda_{\nu+}^-,\lambda_C^+\},\quad
    \mathcal{C} = \{\lambda_{\nu-}^+\}\quad\text{and}\quad
    \mathcal{O} = \{\lambda_{\nu-}^-,\lambda_{\nu+}^+,\lambda_{\mu-}^-,\lambda_C^-\}.
\end{align*}
The other cases can be discussed as before. In particular the case $Q=0$ is not admissible since it would contradict the genuine non-linearity of the field $\lambda_\nu$.
\subsection{Shock Resonance}
Due to the already mentioned result by Freist\"uhler \cite{Freistuhler1991} and the results above we can exclude multiple eigenvalues near a shock discontinuity.
However, we have to discuss the situation that we have a shock in each phase and both move at the same speed $S$.
Hence we first consider the situation
\begin{align*}
    \lambda_{\mu-}^- > S > \lambda_{\mu-}^+\quad\text{and}\quad\lambda_{\nu-}^- > S > \lambda_{\nu-}^+.
\end{align*}
We directly verify $u_\mu^- > S$, $u_\nu^- > S$ and hence $u^- > S$. Due to the continuity of $Q$ we also have $u^+ > S$.
Furthermore we have $\lambda_{\mu+}^- > S$ and $\lambda_{\nu+}^- > S$ and this gives in total seven incoming characteristics, i.e.
\begin{align*}
    \mathcal{I} = \{\lambda_{\mu-}^-,\lambda_{\mu-}^+,\lambda_{\nu-}^-,\lambda_{\nu-}^+,\lambda_{\mu+}^-,\lambda_{\nu+}^-,\lambda_C^-\}.
\end{align*}
There is no coinciding characteristic and we have three outgoing characteristics
\begin{align*}
    \mathcal{O} = \{\lambda_{\mu+}^+,\lambda_{\nu+}^+,\lambda_C^+\}.
\end{align*}
Thus such a configuration violates the admissibility conditions. Let us close with the case
\begin{align*}
    \lambda_{\mu-}^- > S > \lambda_{\mu-}^+\quad\text{and}\quad\lambda_{\nu+}^- > S > \lambda_{\nu+}^+.
\end{align*}
We directly verify $u_\mu^- > S$, $S > u_\nu^+$ and hence we have
\begin{align*}
    \lambda_{\mu+}^- > S,\, 0 > Q_\mu\quad\text{and}\quad S > \lambda_{\nu-}^+,\,Q_\nu > 0.
\end{align*}
From the continuity of the partial mass fluxes we conclude $u_\mu^+ > S$ and $S > u_\nu^-$ and hence we already have the following six incoming and two outgoing characteristics
\begin{align*}
    \mathcal{I} = \{\lambda_{\mu-}^-,\lambda_{\mu-}^+,\lambda_{\nu+}^-,\lambda_{\nu+}^+,\lambda_{\mu+}^-,\lambda_{\nu-}^+\},\quad\mathcal{O} = \{\lambda_{\mu+}^+,\lambda_{\nu-}^-\}.
\end{align*}
Thus the situation of $\lambda_C^- = \lambda_C^+ = S$ can be excluded.
But due to the continuity of the mass flux $Q$ this implies that for $Q < 0$ or $Q > 0$ we have seven incoming and three outgoing characteristics and thus this situation is also not admissible.
With this we have discussed most wave patterns that can appear in a Riemann problem for the model under consideration, including the interaction of two characteristic families.  
    \section{Related Models} \label{sec:relmod}
It is further interesting to discuss how other well established models are related to the studied system at hand.
In particular we discuss Kapila's limit of the barotropic SHTC model and the relation to the Baer-Nunziato model in the following.
\subsection{Instantaneous Relaxation Limit of the Barotropic SHTC Model}\label{sec:kapila_limit}
As it is noted above there are two relaxation terms in the system (\ref{eqn.HPRFF}).
One for the pressure relaxation (\ref{pressure_relax}) and one for the velocity relaxation (\ref{velocity_relax}), respectively.
Since in real processes these relaxation processes can be quite fast, it is useful to obtain asymptotic limits of the solution for small values of the relaxation times.
In the paper \cite{Murrone2005} the simplified equations for instantaneous pressure and relative velocity relaxations are accurately derived by the asymptotic analysis as a relaxation limit of the general Baer-Nunziato two-pressure two-velocity model for two-phase compressible fluid flows.
In this section, based on the results obtained in \cite{Murrone2005} we derive a reduced SHTC system of PDEs for the instantaneous velocity and pressure relaxation.
Note that it seems intuitive that the time scale of pressure relaxation is less than the time scale of velocity relaxation,
because the physical mechanism of pressure relaxation is the pressure waves propagation.
Thus, one can consider a possibility to study the instantaneous relaxation limit separately for the pressure and then for the velocity.
If to consider only the pressure relaxation limit then we arrive at a single pressure two-velocity system of governing equations.
It turns out that this PDE system is not hyperbolic.
If we want to deal with hyperbolic equations we should consider instantaneous velocity relaxation together with the instantaneous pressure relaxation.

We will not repeat the rigorous asymptotic analysis as done in the \cite{Murrone2005}, but will rely on the assumptions that follow from the rigorous theory.
So, for the sake of simplicity we start with  the instantaneous velocity relaxation.
Let us consider equation (\ref{eqn.relvel})
\begin{align}
    \frac{\partial w^k}{\partial t} + \frac{\partial(w^lu^l+E_{c_1})}{\partial x_k} + u^l\left(\frac{\partial w^k}{\partial x_l} - \frac{\partial w^l}{\partial x_k}\right)
    =-\dfrac{ E_{w^k}}{\theta_2}=-\dfrac{ c_1c_2{w^k}}{\theta_2},\label{eqn.relvelnew}
\end{align}
An asymptotic expansion of the relative velocity $w^k=w^k_0+\theta_2w^k_1+...$ for small $\theta_2$ gives us
\begin{align}
    w^k_0=0, \quad w^k_1=\frac{\partial E_{c_1}}{\partial x_k}.
\end{align}
For our purpose we only need to account for the zeroth order term of the expansion $w^k_0 = 0$.
This gives us a single velocity approximation of the model, whereas the second term of the expansion $w^k_1 = \dfrac{\partial E_{c_1}}{\partial x_k}$ gives us a phase diffusion Fick's law.

Now we substitute $w^k=0$ into the system (\ref{eqn.HPRFF}) and remove the equation for the relative velocity. This gives a single velocity model for two-phase flows
\begin{subequations}\label{eqn.HPRFFsv}
    \begin{align}
        \frac{\partial \rho u^i}{\partial t} + \frac{\partial(\rho u^i u^k + p \delta_{ik})}{\partial x_k} &= 0,\label{eqn.momentumFFsv}\\
        \frac{\partial \rho}{\partial t} + \frac{\partial \rho u^k}{\partial x_k} &= 0,\label{eqn.contiFFsv}\\
        \frac{\partial \rho c_1}{\partial t} + \frac{\partial \rho c_1 u^k}{\partial x_k} &= 0,\label{eqn.contiFF1sv}\\
        \frac{\partial \rho \alpha_1}{\partial t} + \frac{\partial \rho \alpha_1 u^k }{\partial x_k} &= -\frac{\rho \phi}{\theta_1},\label{eqn.alphaFFsv}
    \end{align}
\end{subequations}
where $u^i$ is the single velocity of the flow.

The above system (\ref{eqn.HPRFFsv}) is equivalent to
\begin{subequations}\label{eqn.HPRFFsv1}
    \begin{align}
        \frac{\partial \rho u^i}{\partial t} + \frac{\partial (\rho u^i u^k + (\alpha_1p_1 + \alpha_2p_2) \delta_{ik})}{\partial x_k} &= 0, \label{eqn.momentumFFsv1}\\
        \frac{\partial \alpha_1\rho_1}{\partial t} + \frac{\partial \alpha_1\rho_1 u^k}{\partial x_k} &= 0,\label{eqn.contiFF1sv1}\\
        \frac{\partial \alpha_2\rho_2}{\partial t} + \frac{\partial \alpha_2\rho_2 u^k}{\partial x_k} &= 0,\label{eqn.contiFF1sv2}\\
        \frac{\partial \alpha_1}{\partial t} + u^k\frac{\partial\alpha_1}{\partial x_k} &= \frac{p_1-p_2}{\theta_1}.\label{eqn.alphaFFsv1}
    \end{align}
\end{subequations}
The latter is obtained with the use of the definitions $c_1 = \alpha_1\rho_1/\rho, \rho=\alpha_1\rho_1+\alpha_2\rho_2$.

Let us now consider the instantaneous pressure relaxation limit for \eqref{eqn.HPRFFsv1} assuming $\theta_1 \rightarrow 0$.
Note that $\theta_1 \rightarrow 0$ gives us $p_1=p_2$ as a zeroth order approximation, but it is not correct to simply put $p_1=p_2=P$ into the equations, because this would give us
\begin{align}
    \frac{\partial \alpha_1}{\partial t} + u^k\frac{\partial\alpha_1}{\partial x_k} = 0.
\end{align}
This equation for $\alpha_1$ means that the volume fraction does not change along the trajectory, although the pressure can change.
But if the phase pressures change, then $\alpha_1$ must change, too, due to the different phase compressibility coefficients.

The correct way to derive equations for the instantaneous pressure relaxation limit is to account for the following consequences of $p_1=p_2=P$: 
\begin{align}\label{equal.p}
    \dd p_1 = \frac{K_1}{\rho_1}\dd\rho_1 = \dd p_2 = \frac{K_2}{\rho_2}\dd \rho_2,\quad
    \frac{K_1}{\rho_1}\frac{\partial\rho_1}{\partial t} = \frac{K_2}{\rho_2}\frac{\partial\rho_2}{\partial t},\quad
    \frac{K_1}{\rho_1}\frac{\partial\rho_1}{\partial x_k} = \frac{K_2}{\rho_2}\frac{\partial\rho_2}{\partial x_k},
\end{align}
where $K_i=\rho_i a_i^2$ is the phase bulk modulus, $a_i$ is the phase speed of sound.
Now we use the above equation \eqref{equal.p}, the phase mass conservation equations \eqref{eqn.contiFF1sv1} and \eqref{eqn.contiFF1sv2} to derive the following equation for the volume fraction
\begin{align}
    \frac{\partial \alpha_1}{\partial t} + u^k\frac{\partial\alpha_1}{\partial x_k} + \frac{\alpha_1 \alpha_2(K_1 - K_2)}{\alpha_1K_2 + \alpha_2K_1}\frac{\partial u^k}{\partial x_k} = 0.
\end{align}
Thus, in case of instantaneous velocity and pressure relaxation we arrive at the following system
\begin{subequations}\label{eqn.HPRFFK}
    \begin{align}
        \frac{\partial \rho u^i}{\partial t} + \frac{\partial(\rho u^i u^k + P\delta_{ik})}{\partial x_k} &= 0,\label{eqn.momentumFFK}\\
        \frac{\partial \alpha_1\rho_1}{\partial t} + \frac{\partial \alpha_1\rho_1 u^k}{\partial x_k} &= 0,\label{eqn.contiFF1K}\\
        \frac{\partial \alpha_2\rho_2}{\partial t} + \frac{\partial \alpha_2\rho_2 u^k}{\partial x_k} &= 0,\label{eqn.contiFF2K}\\
        \frac{\partial \alpha_1}{\partial t} + u^k\frac{\partial\alpha_1}{\partial x_k}
        + \frac{\alpha_1 \alpha_2(K_1 - K_2)}{\alpha_1K_2 + \alpha_2K_1}\frac{\partial u^k}{\partial x_k} &=0. \label{eqn.alphaFFK}
    \end{align}
\end{subequations}
The above equations are exactly the same as is in the Kapila model.

\subsection{Comparison with the Baer-Nunziato Model}
It is of further interest to compare the system (\ref{sys:gpr_1d_v2}) with the one dimensional barotropic Baer-Nunziato model given by the following equations
\begin{subequations}\label{sys:bn_1d}
    \begin{align}
        \frac{\partial\alpha_1}{\partial t} + u_I\frac{\partial \alpha_1}{\partial x} &= \zeta_1,\label{eq:alpha_balance_bn}\\
        \frac{\partial\alpha_1\rho_1}{\partial t} + \frac{\partial\alpha_1\rho_1 u_1}{\partial x} &= \zeta_2,\label{eq:partial_rho1_balance_bn}\\
        \frac{\partial\alpha_2\rho_2}{\partial t} + \frac{\partial\alpha_2\rho_2 u_2}{\partial x} &= \zeta_3,\label{eq:partial_rho2_balance_bn}\\
        \frac{\partial\alpha_1\rho_1 u_1}{\partial t}
        + \frac{\partial\left(\alpha_1\rho_1 u_1^2 + \alpha_1 p_1(\rho_1)\right)}{\partial x} - p_I\frac{\partial\alpha_1}{\partial x} &= \zeta_4,\label{eq:partial_mom1_balance_bn}\\
        \frac{\partial\alpha_2\rho_2 u_2}{\partial t}
        + \frac{\partial\left(\alpha_2\rho_2 u_2^2 + \alpha_2 p_2(\rho_2)\right)}{\partial x} - p_I\frac{\partial\alpha_2}{\partial x} &= \zeta_5,\label{eq:partial_mom2_balance_bn}
    \end{align}
\end{subequations}
It was already mentioned in \cite{Romenski2004} and \cite{Romenski2007} that for smooth solutions the systems can be reformulated into each other in one space dimension.
However, we want to recall this transformation with a slightly different purpose.
In the following we want to show the equivalence of these two systems (\ref{sys:gpr_1d_v2}) and (\ref{sys:bn_1d}) for smooth solutions
and in particular that there is no freedom of choice for the interface quantities $u_I$, $p_I$ and the sources in this case.
This is of special interest since there are several possible choices for the interface quantities in the context of these Baer-Nunziato type models, see e.g.\ \cite{Andrianov2003}.
\begin{enumerate}[(i)]
    \item For the choice $\zeta_2 = \xi_2$, equation (\ref{eq:partial_rho1_balance_bn}) and equation (\ref{eq:partial_mass_balance}) are equal.
    \item The sum of (\ref{eq:partial_rho1_balance_bn}) and (\ref{eq:partial_rho2_balance_bn}) gives the mixture mass balance (\ref{eq:mix_mass_balance}) with $\xi_3 = \zeta_2 + \zeta_3$.
    \item Subtracting (\ref{eq:partial_mass_balance}) from (\ref{eq:mix_mass_balance}) we obtain (\ref{eq:partial_rho2_balance_bn}) with $\zeta_3 = \xi_3 - \xi_2$.
    \item The balance equations (\ref{eq:alpha_balance}) and (\ref{eq:alpha_balance_bn}) for the volume fraction are equivalent for the choices $u_I = u$ and $\zeta_1 = (\xi_1 - \alpha_1\xi_3)/\rho$.
    \item With $\alpha_2 = 1 - \alpha_1$ and $\xi_4 = \zeta_4 + \zeta_5$ the sum of (\ref{eq:partial_mom1_balance_bn}) and (\ref{eq:partial_mom2_balance_bn})
          gives the mixture momentum balance (\ref{eq:mix_mom_balance}).
\end{enumerate}
To show how the barotropic SHTC model (\ref{sys:gpr_1d_v2}) can be derived from the Baer-Nunziato model (\ref{sys:bn_1d}),
we need to derive the balance for the relative velocity (\ref{eq:rel_velocity_balance}).
Therefore we want to reformulate the partial momentum balances (\ref{eq:partial_mom1_balance_bn})
\begin{align}
    \zeta_4 &= \frac{\partial\alpha_1\rho_1 u_1}{\partial t} + \frac{\partial\left(\alpha_1\rho_1 u_1^2 + \alpha_1 p_1(\rho_1)\right)}{\partial x} - p_I\frac{\partial\alpha_1}{\partial x}\notag\\
    &= u_1\left(\frac{\partial\alpha_1\rho_1}{\partial t} + \frac{\partial\alpha_1\rho_1 u_1}{\partial x}\right)
    + \alpha_1\rho_1\left(\frac{\partial u_1}{\partial t} + \frac{1}{2}\frac{\partial u_1^2}{\partial x}\right)\notag
    + \frac{\partial\alpha_1p_1(\rho_1)}{\partial x} - p_I\frac{\partial\alpha_1}{\partial x}\notag\\
    \Leftrightarrow\quad
    \frac{\zeta_4 - u_1\zeta_2}{\alpha_1\rho_1} &= \frac{\partial u_1}{\partial t} + \frac{1}{2}\frac{\partial u_1^2}{\partial x}
    + \frac{1}{\alpha_1\rho_1}\frac{\partial\alpha_1p_1(\rho_1)}{\partial x} - \frac{p_I}{\alpha_1\rho_1}\frac{\partial \alpha_1}{\partial x}\label{eq:partial_mom1_balance_bn_v2}
\end{align}
and (\ref{eq:partial_mom2_balance_bn})
\begin{align}
    \frac{\zeta_5 - u_2\zeta_3}{\alpha_2\rho_2} &= \frac{\partial u_2}{\partial t} + \frac{1}{2}\frac{\partial u_2^2}{\partial x}
    + \frac{1}{\alpha_2\rho_2}\frac{\partial\alpha_2p_2(\rho_2)}{\partial x} - \frac{p_I}{\alpha_2\rho_2}\frac{\partial\alpha_2}{\partial x}.\label{eq:partial_mom2_balance_bn_v2}
\end{align}
Subtracting equation (\ref{eq:partial_mom2_balance_bn_v2}) from equation (\ref{eq:partial_mom1_balance_bn_v2}) gives
\begin{align*}
    &\frac{\partial(u_1 - u_2)}{\partial t} + \frac{1}{2}\frac{\partial\left(u_1^2 - u_2^2\right)}{\partial x}
    + \frac{1}{\alpha_1\rho_1}\frac{\partial \alpha_1p_1(\rho_1)}{\partial x} - \frac{1}{\alpha_2\rho_2}\frac{\partial \alpha_2p_2(\rho_2)}{\partial x}
    - \frac{p_I}{\alpha_1\rho_1}\frac{\partial \alpha_1}{\partial x} + \frac{p_I}{\alpha_2\rho_2}\frac{\partial\alpha_2}{\partial x}\\
    = &\frac{\partial (u_1 - u_2)}{\partial t} + \frac{1}{2}\frac{\partial\left(u_1^2 - u_2^2\right)}{\partial x}
    + \left(\frac{p_1}{\alpha_1\rho_1} + \frac{p_2}{\alpha_2\rho_2}\right)\frac{\partial\alpha_1}{\partial x}\\
    + &\frac{a_1^2}{\rho_1}\frac{\partial\rho_1}{\partial x} - \frac{a_2^2}{\rho_2}\frac{\partial\rho_2}{\partial x}
    - p_I\left(\frac{1}{\alpha_1\rho_1} + \frac{1}{\alpha_2\rho_2}\right)\frac{\partial\alpha_1}{\partial x}
    = \frac{\zeta_4 - u_1\zeta_2}{\alpha_1\rho_1} - \frac{\zeta_5 - u_2\zeta_3}{\alpha_2\rho_2}.
\end{align*}
Using relation (\ref{def:sound_speed}) and $p_I = (\alpha_2\rho_2p_1 + \alpha_1\rho_1p_2)/\rho$ we obtain
\begin{align*}
    \frac{\partial(u_1 - u_2)}{\partial t} + \dfrac{\partial\left(\dfrac{1}{2}u_1^2 - \dfrac{1}{2}u_2^2 + \Psi_1(\rho_1) - \Psi_2(\rho_2)\right)}{\partial x}
    = \underbrace{\frac{\zeta_4 - u_1\zeta_2}{\alpha_1\rho_1} - \frac{\zeta_5 - u_2\zeta_3}{\alpha_2\rho_2}}_{=\xi_5}
\end{align*}
Thus we have shown that the SHTC model can be derived from the Baer-Nunziato model with the following choices
\begin{align*}
    u_I &= u,\quad p_I = \frac{\alpha_2\rho_2p_1 + \alpha_1\rho_1p_2}{\rho},\quad\bb{\Xi} = \bb{B}\bb{\zeta}\\
    \text{with}\quad\bb{B} &= \begin{pmatrix}
        \rho & \alpha_1                     & \alpha_1                    & 0                         & 0                           \\
        0    & 1                            & 0                           & 0                         & 0                           \\
        0    & 1                            & 1                           & 0                         & 0                           \\
        0    & 0                            & 0                           & 1                         & 1                           \\
        0    & -\dfrac{u_1}{\alpha_1\rho_1} & \dfrac{u_2}{\alpha_2\rho_2} & \dfrac{1}{\alpha_1\rho_1} & -\dfrac{1}{\alpha_2\rho_2}
    \end{pmatrix}
\end{align*}
It remains to show how the partial momentum balances can be derived using the SHTC model. We start with the balance for the relative velocity (\ref{eq:rel_velocity_balance})
\begin{align}
    \xi_5 &= \frac{\partial(u_1 - u_2)}{\partial t} + \dfrac{\partial\left(\dfrac{1}{2}u_1^2 - \dfrac{1}{2}u_2^2 + \Psi_1(\rho_1) - \Psi_2(\rho_2)\right)}{\partial x}\notag\\
    &= \frac{\alpha_1\rho_1}{\alpha_1\rho_1}\left(\frac{\partial u_1}{\partial t} + \dfrac{\partial\left(\dfrac{1}{2}u_1^2 + \Psi_1(\rho_1)\right)}{\partial x}\right)
    -  \frac{\alpha_2\rho_2}{\alpha_2\rho_2}\left(\frac{\partial u_2}{\partial t} + \dfrac{\partial \left(\dfrac{1}{2}u_2^2 + \Psi_2(\rho_2)\right)}{\partial x}\right)\notag\\
    %
    &\Leftrightarrow\quad \alpha_1\rho_1\alpha_2\rho_2\xi_5 + \alpha_2\rho_2u_1\xi_2 - \alpha_1\rho_1u_2\zeta_3\notag\\
    &= \alpha_2\rho_2\left(\frac{\partial\alpha_1\rho_1u_1}{\partial t} + \frac{\partial\left(\alpha_1\rho_1u_1^2 + \alpha_1p_1(\rho_1)\right)}{\partial x}
    - p_1\frac{\partial\alpha_1}{\partial x}\right)\notag\\
    &- \alpha_1\rho_1\left(\frac{\partial\alpha_2\rho_2u_2}{\partial t} + \frac{\partial\left(\alpha_2\rho_2u_2^2 + \alpha_2p_2(\rho_2)\right)}{\partial x}
    - p_2\frac{\partial\alpha_2}{\partial x}\right)\label{eq:rel_velocity_balance_v2}
\end{align}
Now we multiply the mixture momentum balance (\ref{eq:mix_mom_balance}) with $\alpha_1\rho_1$ and add it to the previous equation (\ref{eq:rel_velocity_balance_v2}) to obtain
\begin{align*}
    %
    &\frac{\partial\alpha_1\rho_1u_1}{\partial t} + \frac{\partial\left(\alpha_1\rho_1u_1^2 + \alpha_1p_1(\rho_1)\right)}{\partial x}
    - \frac{\alpha_2\rho_2p_1 + \alpha_1\rho_1p_2}{\rho}\frac{\partial\alpha_1}{\partial x}\\
    = &\frac{1}{\rho}\left(\alpha_1\rho_1\alpha_2\rho_2\xi_5 + \alpha_2\rho_2u_1\xi_2 - \alpha_1\rho_1u_2\zeta_3 + \alpha_1\rho_1\xi_4\right)
\end{align*}
With the same choice for the interface pressure as before and the corresponding source term we have obtained the partial momentum balance for the first phase.
The second balance can be obtained in a similar way (multiply with $\alpha_2\rho_2$ and subtract)
\begin{align*}
    &\frac{\partial\alpha_2\rho_2u_2}{\partial t} + \frac{\partial\left(\alpha_2\rho_2u_2^2 + \alpha_2p_2(\rho_2)\right)}{\partial x}
    - \frac{\alpha_2\rho_2p_1 + \alpha_1\rho_1p_2}{\rho}\frac{\partial\alpha_2}{\partial x}\\
    = -&\frac{1}{\rho}\left(\alpha_1\rho_1\alpha_2\rho_2\xi_5 + \alpha_2\rho_2u_1\xi_2 - \alpha_1\rho_1u_2\zeta_3 - \alpha_2\rho_2\xi_4\right)
\end{align*}
Summarizing the relations we have
\begin{align*}
    u_I &= u,\quad p_I = \frac{\alpha_2\rho_2p_1 + \alpha_1\rho_1p_2}{\rho},\quad\bb{\zeta} = \bb{C}\bb{\Xi}\\
    \text{with}\quad\bb{C} &= \begin{pmatrix}
        \dfrac{1}{\rho} & 0 &-\dfrac{\alpha_1}{\rho} & 0 & 0       \\
        0 & 1                   & 0         & 0     & 0            \\
        0 & -1                  & 1         & 0     & 0            \\
        0 & c_2u_1 + c_1u_2     & -c_1u_2   & c_1   & c_1c_2\rho   \\
        0 & -(c_2u_1 + c_1u_2)  & c_1u_2    & c_2   & -c_1c_2\rho
    \end{pmatrix}
\end{align*}
One immediately verifies $\bb{B}\bb{C} = \bb{I}$ and thus the equivalence in the smooth case is proven.
Following from this equivalence, we verify that the systems share the same Jacobian (\ref{eq:jacobian_prim_var})
and thus the same eigenvalues (\ref{sys_eigenvalues}) and eigenvectors (\ref{sys_eigenvectors_prim}).
However, there are crucial differences between the Baer-Nunziato system (\ref{sys:bn_1d}) and the SHTC system (\ref{sys:gpr_1d_v2}).
Apart from the obvious distinction that system (\ref{sys:gpr_1d_v2}) can be written in conservative form, the most remarkable difference are the jump conditions for discontinuities.
For system (\ref{sys:bn_1d}) the equation for $\alpha$ is not in conservative form and hence it is in general not possible to write down jump conditions, cf.\ \cite{Murrone2005}.
However, by the argument that $\alpha$ stays constant across the shock the equations are decouple into Euler systems for each phase with corresponding jump conditions.
In contrast to system (\ref{sys:gpr_1d_v2}) where the phases remain coupled across discontinuities and a shock in one phase also affects the other.
Note that the particular choice for $u_I$ and $p_I$ made above is also mentioned as a consistent choice in a work by H\'{e}rard \cite{Herard2016}.
In the isentropic single temperature case it also coincides with the choice presented in Coquel et al. \cite{Coquel2002}.

    \section{Numerical Results}\label{sec:num_res}
%

For the SHTC model, all numerical results shown in this section have been obtained with a classical second order MUSCL-Hancock scheme, see \cite{Toro2009} for details. 
For the non-conservative Baer-Nunziato system, which we solve for comparison at the end of this section, we employ a second-order path-conservative version of the MUSCL Hancock scheme, see e.g. \cite{Castro2006,Pares2006,OsherNC,PCOsher,USFORCE2}. 
In all cases we use the simple Rusanov-flux (local Lax-Friedrichs flux) as approximate Riemann solver.

\subsection{Exact vs. numerical solution of the homogeneous system}
The exact solution for the examples shown in this section is obtained as follows: 
\begin{enumerate}[1)]
    \item Assume that the eigenvalues $\lambda_{i-}$ are located left of $\lambda_C$ and the eigenvalues $\lambda_{i+}$ on the right, respectively.
    \item Prescribe the state left of the contact, $\alpha_1$ on the right and then solve for remaining quantities on the right.
    \item Choose eigenvalue next to the contact on the left side. Prescribe wave type (shock/rarefaction) and choose wave speed (head speed for a rarefaction).
    \item Repeat with remaining eigenvalue on the left side while the admissibility is constantly checked and thus certain waves might be excluded already.
    \item Repeat the afore mentioned steps for the right side.
    \item Sample the complete solution.
\end{enumerate}
\subsubsection{Shock in rarefaction}
\begin{figure}[h!]
    \begin{center} 
    \subfigure[Densities $\rho_1, \rho_2, \rho$ and volume fraction $\alpha_1$.]{
    	\begin{tabular}{cc}        \includegraphics[width=0.35\textwidth]{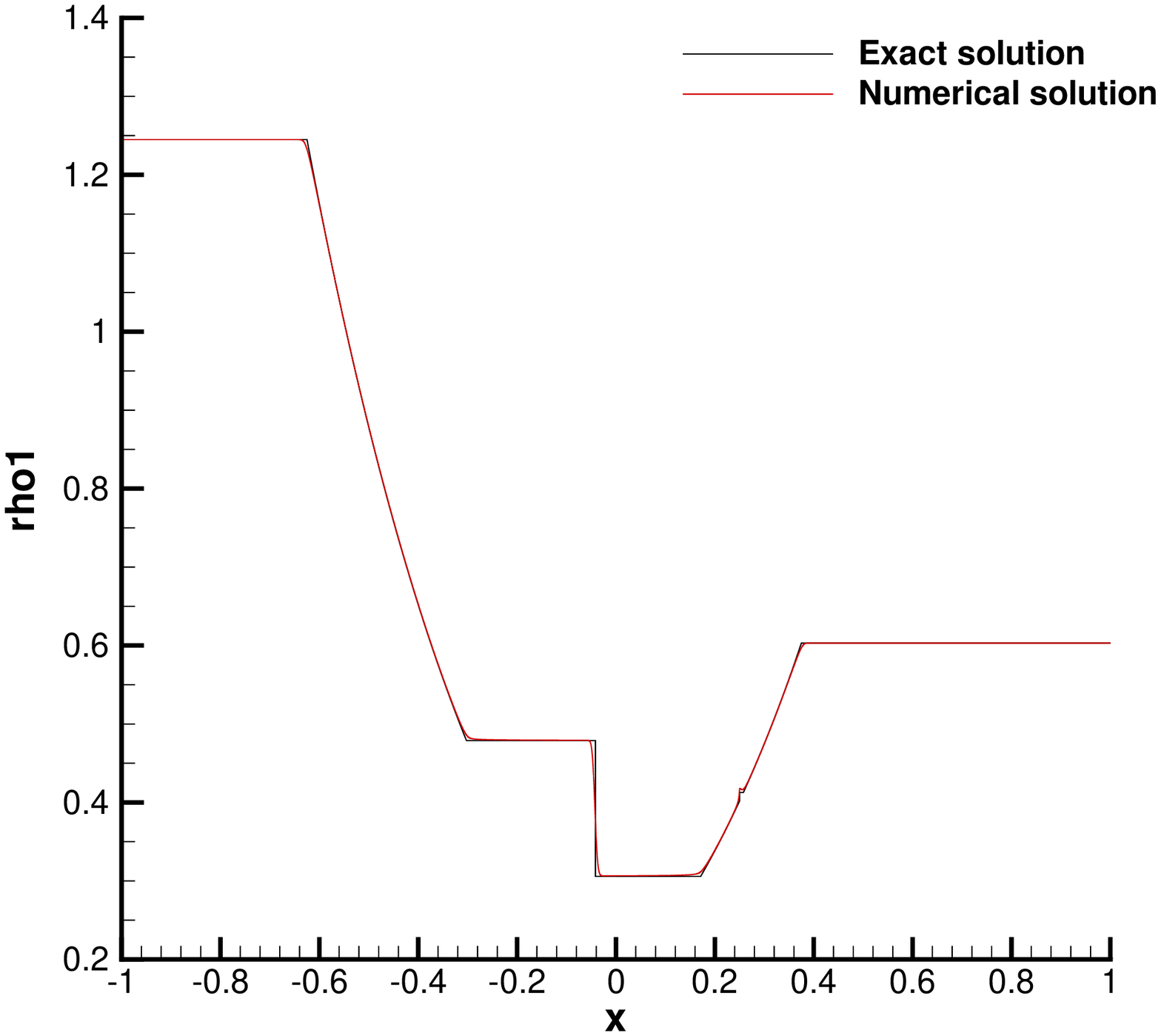}  &  
		\includegraphics[width=0.35\textwidth]{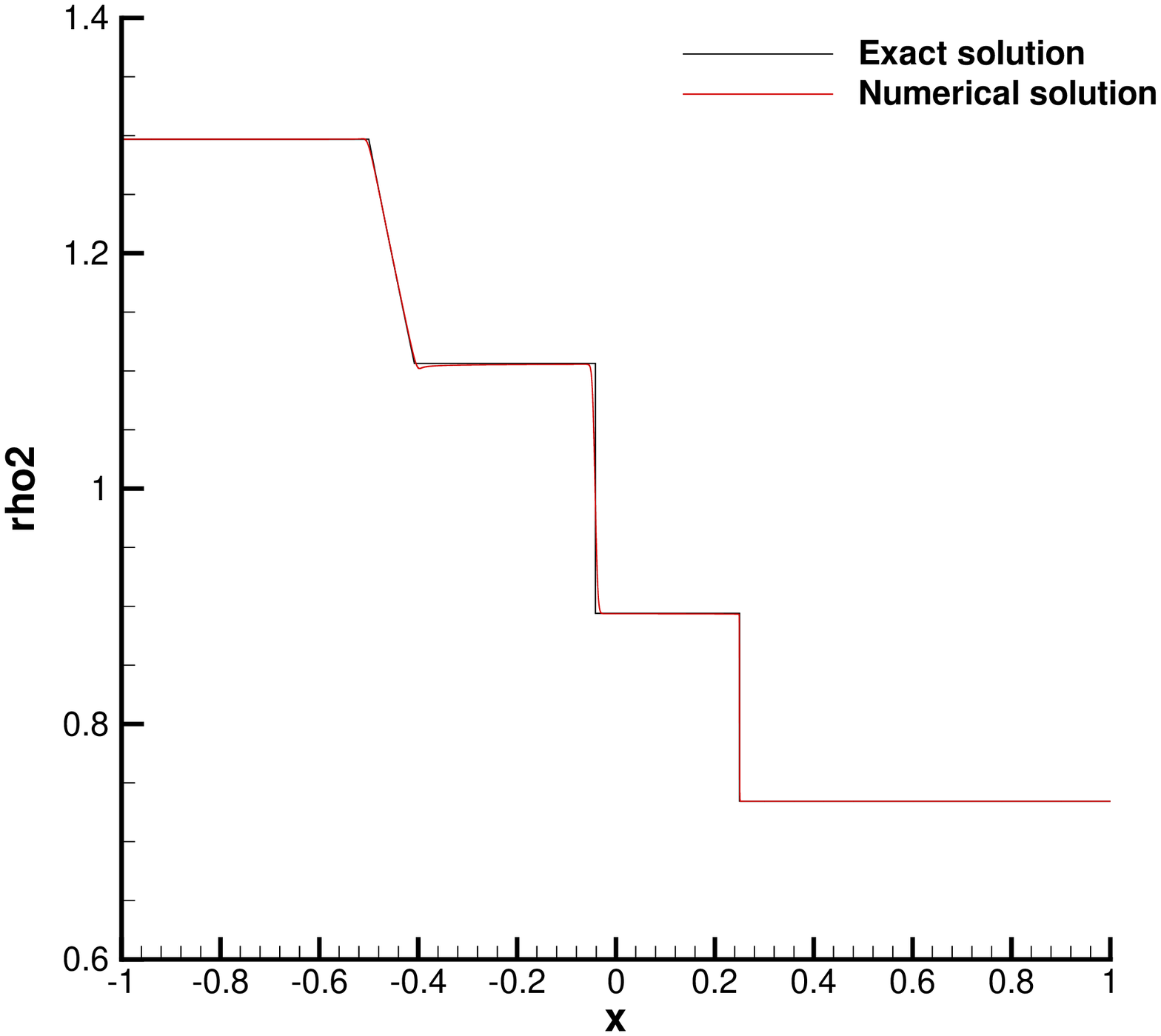}  \\ 	
		\includegraphics[width=0.35\textwidth]{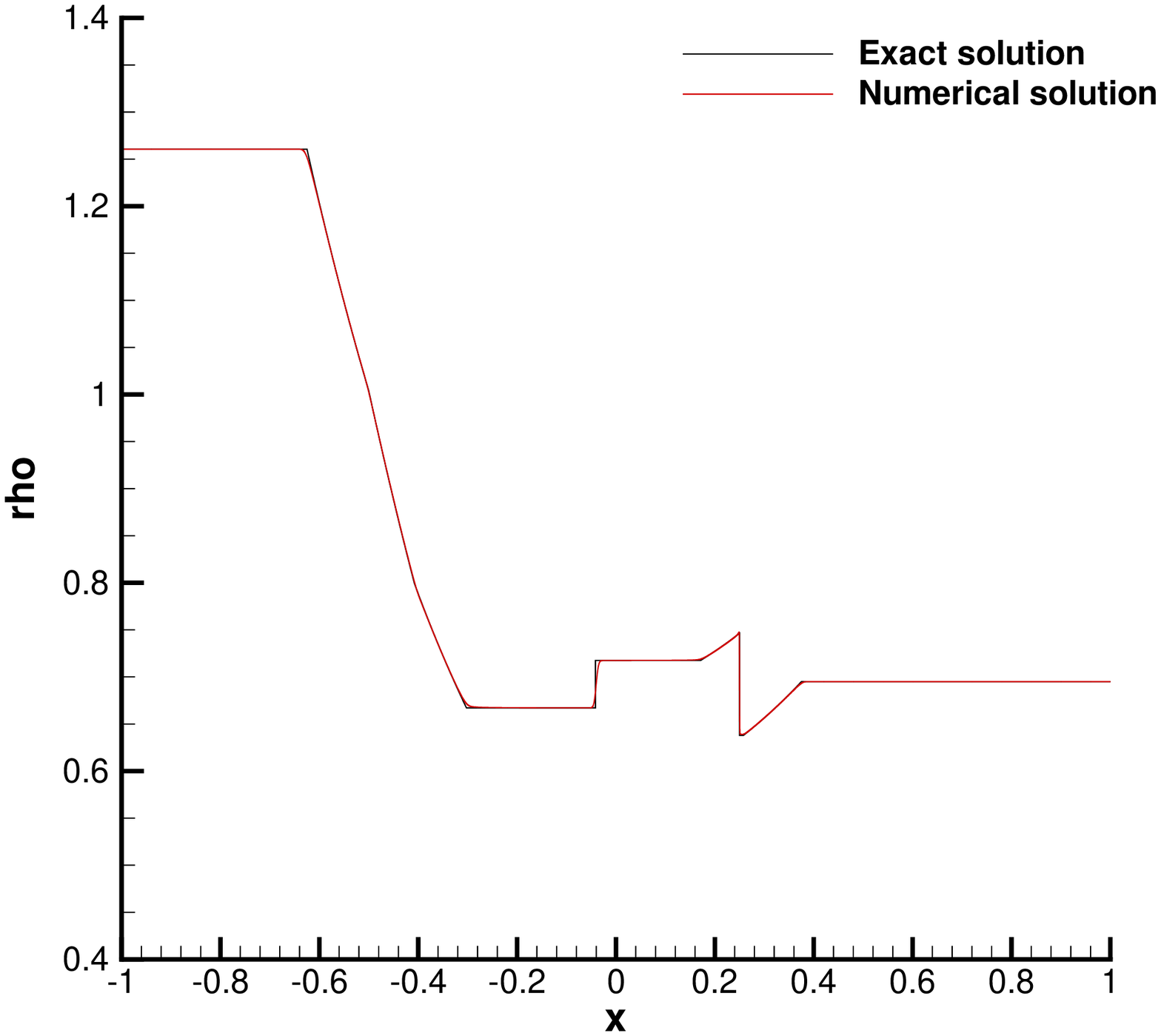}  &  
		\includegraphics[width=0.35\textwidth]{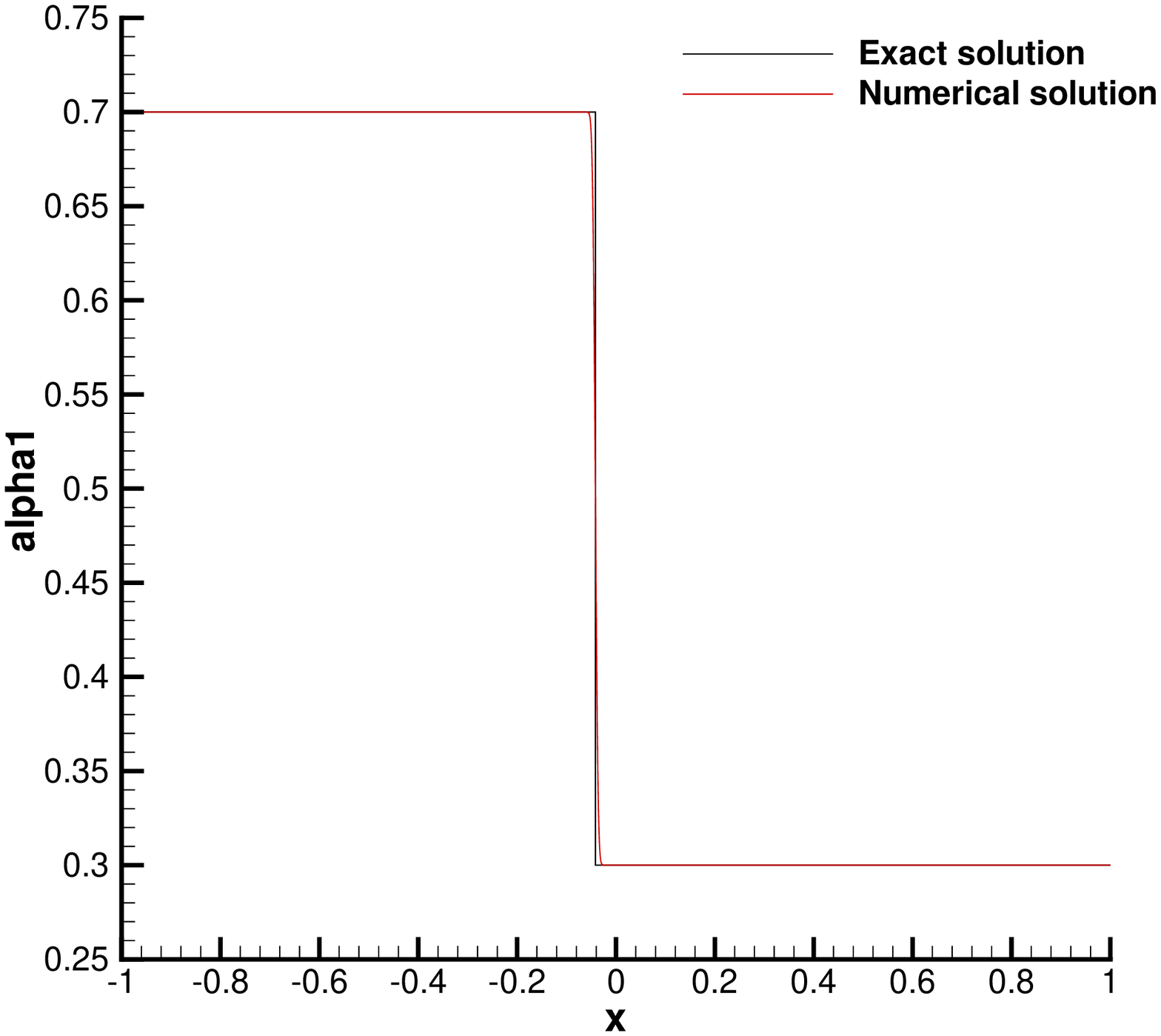}    	
		\end{tabular}
     }\\
    \subfigure[Velocities $u_1, u_2, u, w$.]{
    	\begin{tabular}{cc}        \includegraphics[width=0.35\textwidth]{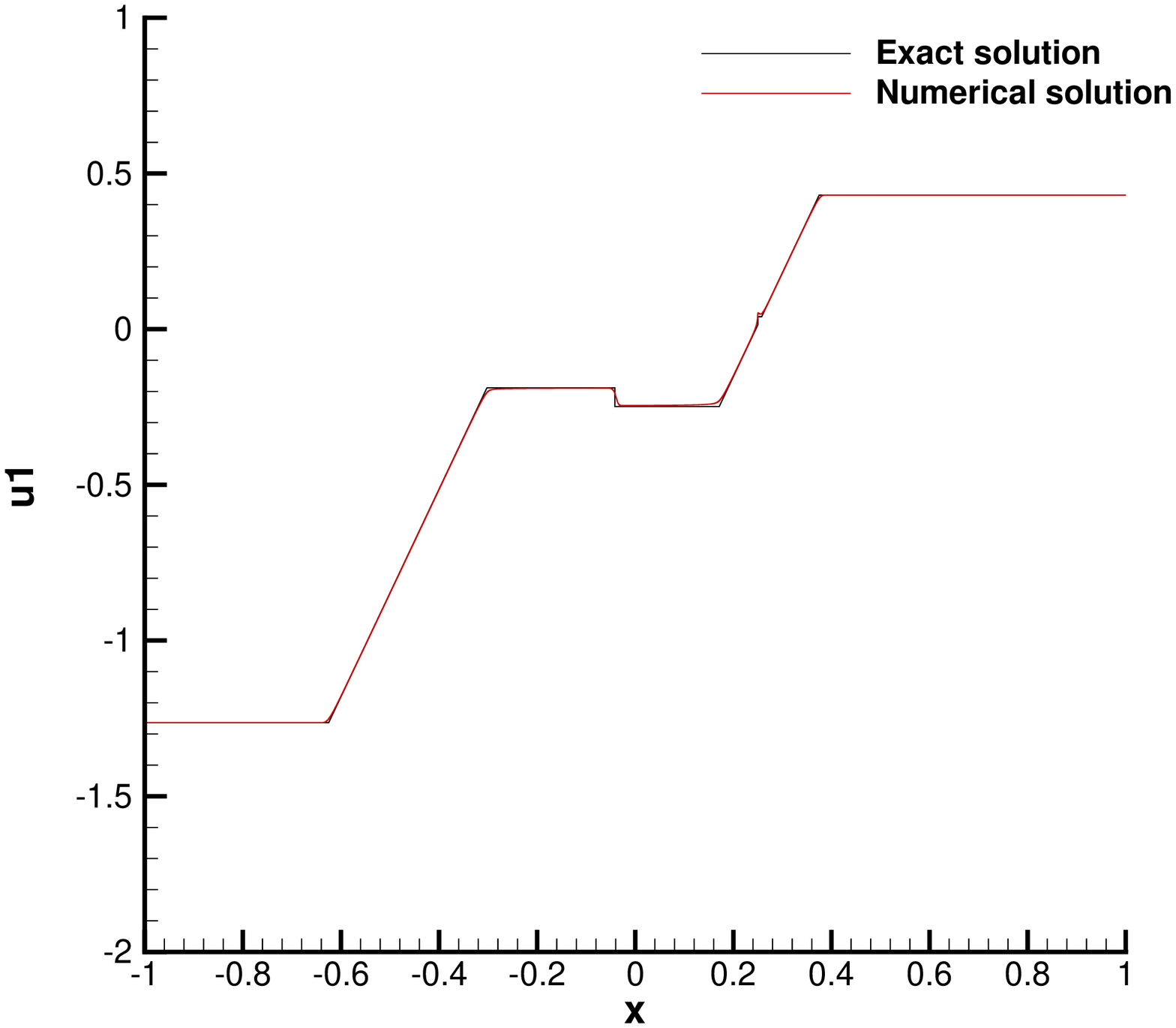}  &  
		\includegraphics[width=0.35\textwidth]{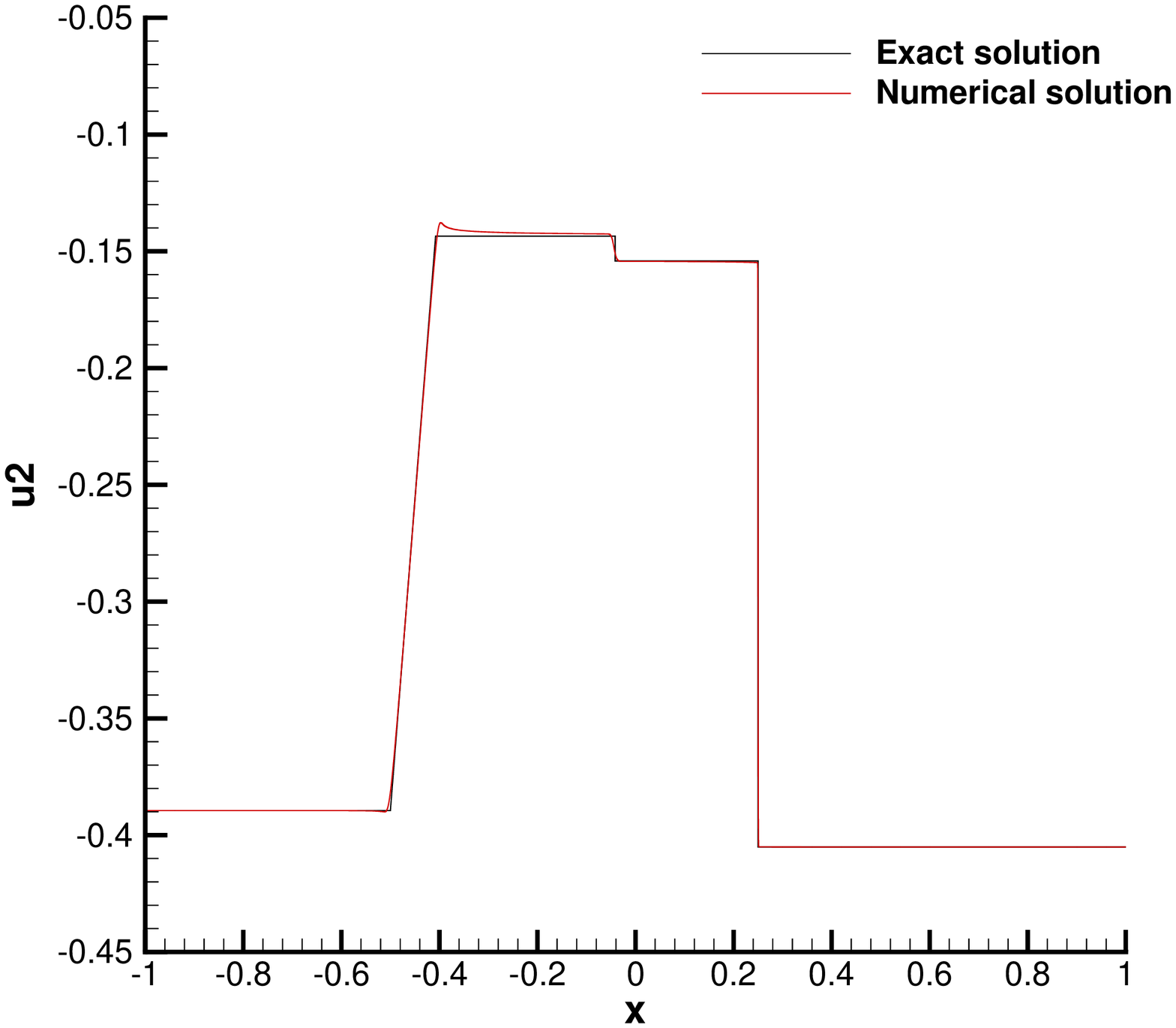}  \\ 	
		\includegraphics[width=0.35\textwidth]{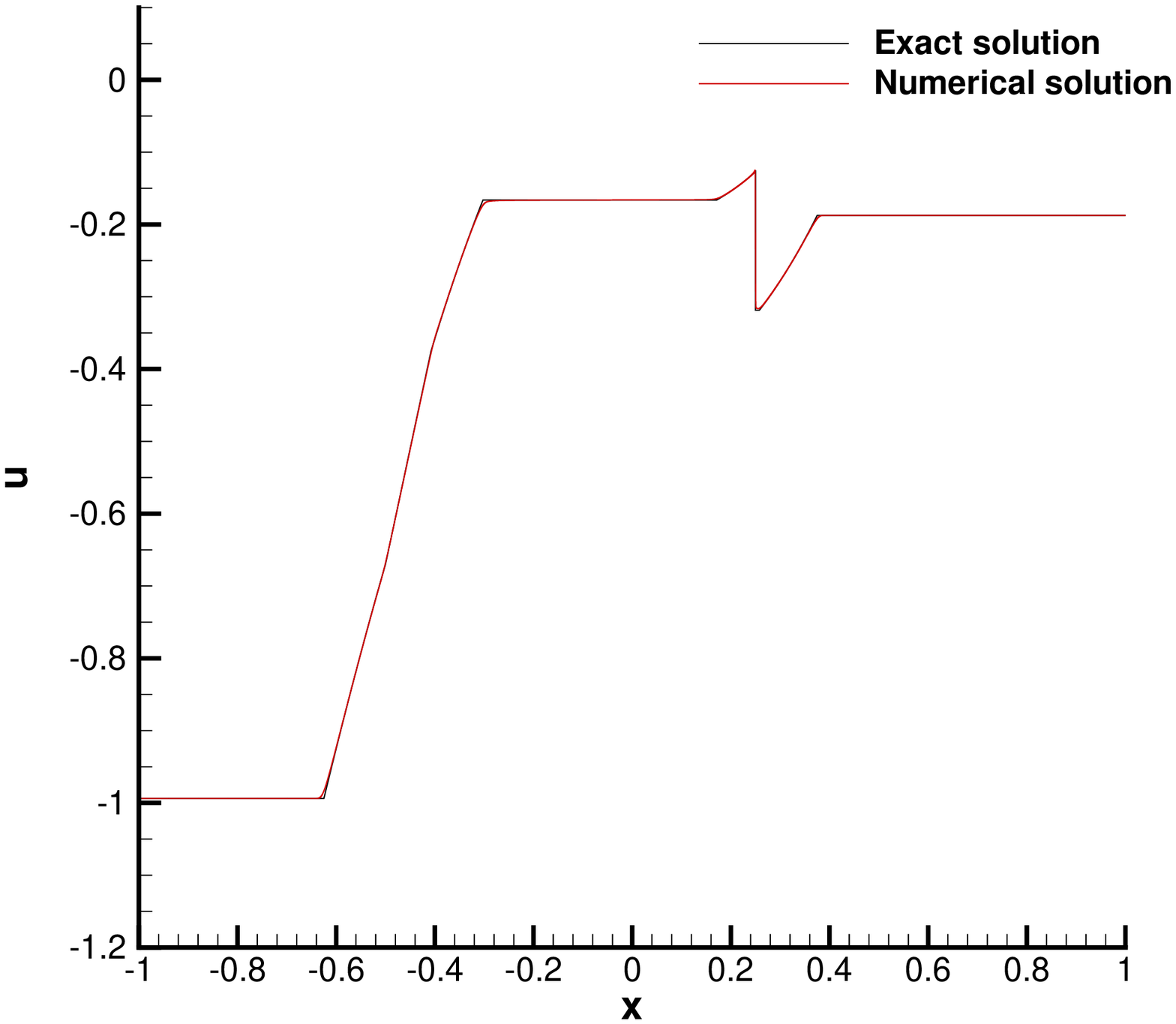}  &  
		\includegraphics[width=0.35\textwidth]{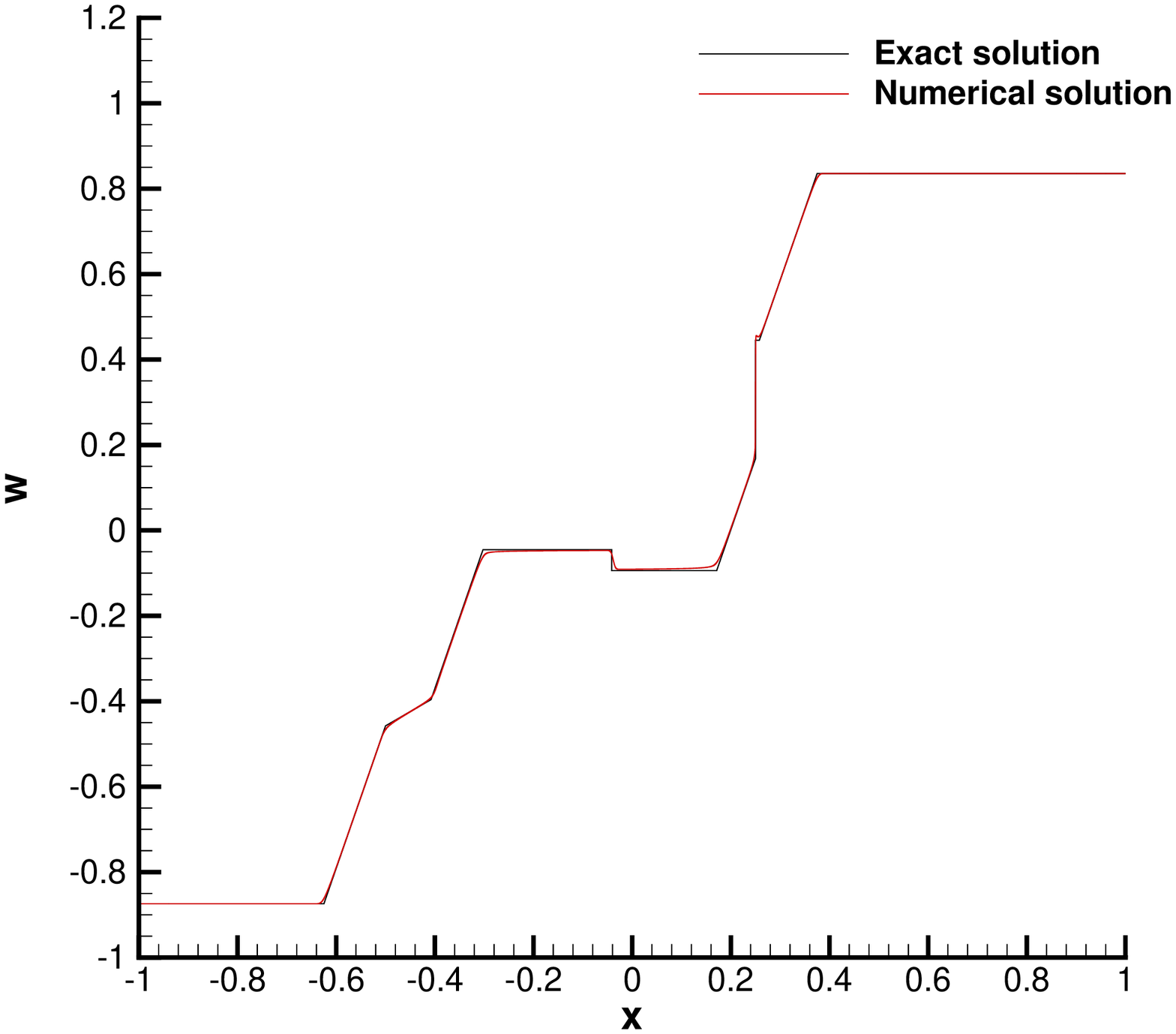}    	
		\end{tabular}
     }
	\caption{Exact solution (black) and numerical solution (red) of Riemann problem RP1.}
    \label{fig:ex1_p1}
    \end{center} 
\end{figure}

\begin{figure}
    \begin{center} 
    \includegraphics[width=0.8\textwidth]{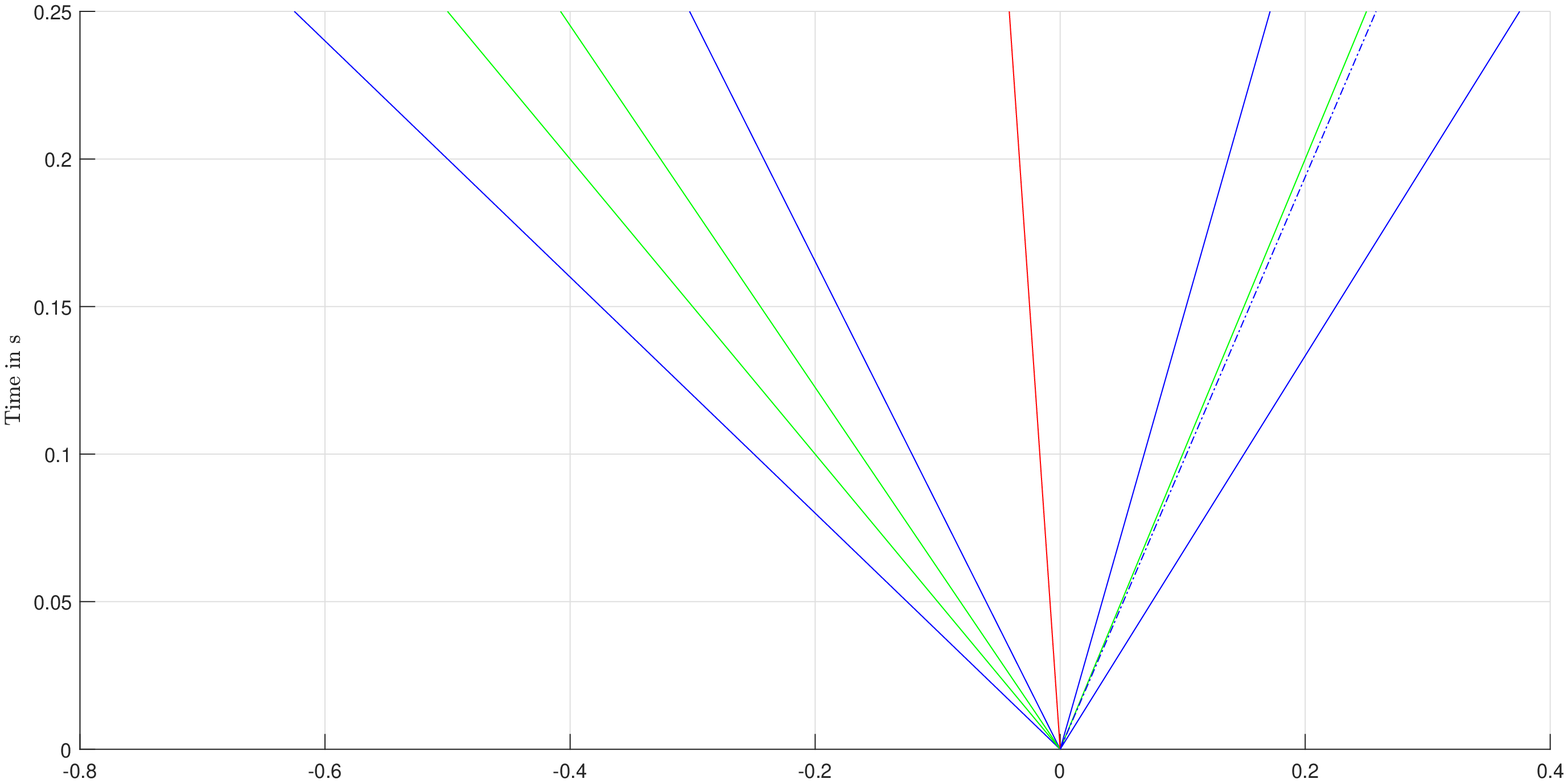} \\
    \includegraphics[width=0.8\textwidth]{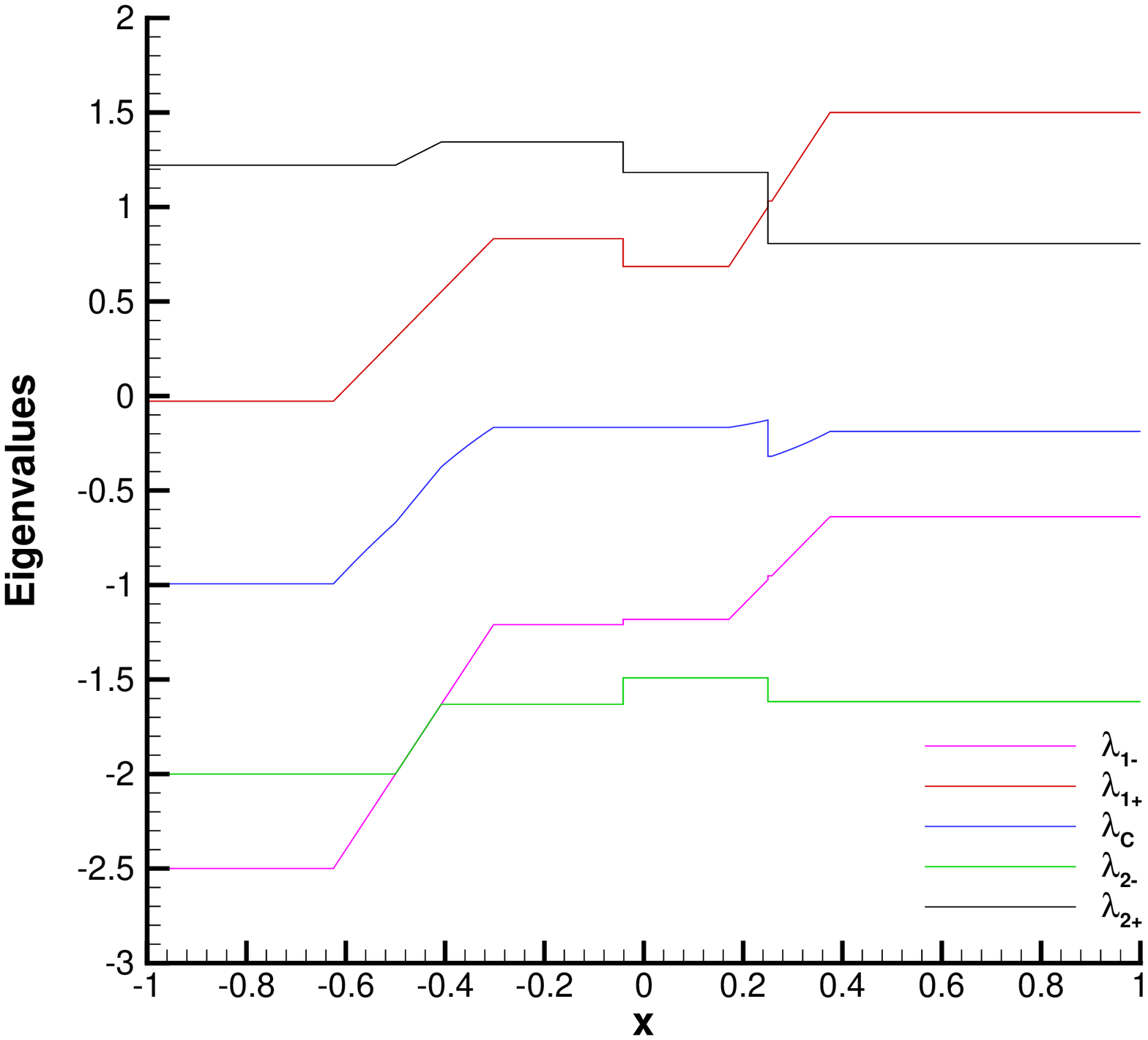}
    \caption{Wave structure of Riemann problem RP1 (top): phase one (blue), contact (red) and phase 2 (green). Eigenvalues of RP1 (bottom). }
    \label{fig:ex1_struct}
    \end{center}
\end{figure}

The exact solution is obtained by using a contact centered inverse construction of the solution as described above. The states of this first Riemann problem (RP1 )are given in Table \ref{tab:states_prob_1}.
Due to the shock inside the rarefaction we have an additional state behind the shock given by $\bar{U}$.
\begin{table}[h!]
    \centering
    \renewcommand*{\arraystretch}{1.25}
    \begin{tabular}{c|lllllll}
                   &  \multicolumn{1}{c}{$U_L$}              &  \multicolumn{1}{c}{$U^\ast_L$} &  \multicolumn{1}{c}{$U^{\ast\ast}_L$} &  \multicolumn{1}{c}{$U_R^{\ast\ast}$}
                   &  \multicolumn{1}{c}{$\bar{U}$}          &  \multicolumn{1}{c}{$U^\ast_R$} &  \multicolumn{1}{c}{$U_R$}            \\
        \hline
        $\alpha_1$ &  \phantom{-}0.7    &  \phantom{-}0.7     &  \phantom{-}0.7     &  \phantom{-}0.3     &  \phantom{-}0.3     &  \phantom{-}0.3      &  \phantom{-}0.3     \\
        $\rho_1$   &  \phantom{-}1.2449 &  \phantom{-}0.47883 &  \phantom{-}0.47883 &  \phantom{-}0.30577 &  \phantom{-}0.40186 &  \phantom{-}0.41275  &  \phantom{-}0.60312 \\
        $\rho_2$   &  \phantom{-}1.2969 &  \phantom{-}1.2969  &  \phantom{-}1.1064  &  \phantom{-}0.894   &  \phantom{-}0.894   &  \phantom{-}0.73436  &  \phantom{-}0.73436 \\
        $u_1$      &  -1.2638           &  -0.18865           &  -0.18865           &  -0.24825           &  \phantom{-}0.01399 &  \phantom{-}0.040001 &  \phantom{-}0.43059 \\
        $u_2$      &  -0.38947          &  -0.38947           &  -0.14351           &  -0.15416           &  -0.15416           &  -0.40507            &  -0.40507           
        %
    \end{tabular}
    \caption{Primitive states of Riemann problem RP1.}
    \label{tab:states_prob_1}
\end{table}
The following computation was performed using the ideal gas EOS
\begin{align*}
    p_i(\rho_i) = \rho_i^{\gamma_i},\, i\in\{1,2\}\quad\text{with}\quad \gamma_1 = 1.4,\,\gamma_2 = 2
\end{align*}
and the parameters
\begin{align*}
    \Delta x = 0.5\cdot 10^{-4}\,\si{m},\; C_{CFL} = 0.25,\, t_{end} = 0.25\,\si{s} \quad\text{and}\quad x \in [-1,1]\,\si{m}.
\end{align*}
In Figures \ref{fig:ex1_p1} and \ref{fig:ex1_struct} the numerical results together with the exact solution are shown.
The overall wave structure of the solution is depicted in Figure \ref{fig:ex1_struct}.
It consists of a $\lambda_{2-}$ - rarefaction (green) which is completely contained inside the $\lambda_{1-}$ - rarefaction (blue).
The contact (red) is (by construction) the middle wave. On the right we have a $\lambda_{1+}$-rarefaction (blue) which contains a $\lambda_{2+}$-shock (green) inside.
The dashed blue line marks the intermediate tail of the rarefaction, i.e.\ the position when the rarefaction starts again after the shock.
Looking at the densities and velocities in Figure \ref{fig:ex1_p1} we verify that a rarefaction only affects the related phase.
However, the interaction of the two rarefactions manifests itself in the mixture quantities $\rho, u, w$.
For the volume fraction we observe that it only jumps at the contact as shown before and takes the initial values on each side.
Special attention has to be paid to the density of the first phase.
Right of the contact the rarefaction starts (not affecting phase two) until the shock occurs.
There the density jumps according to the jump conditions. Continuing on the right we have a plateau of the right state of the shock and then the rarefaction continues again.

The eigenvalues are shown in the upper right panel of 
Figure \ref{fig:ex1_struct}. One clearly sees the overlapping eigenvalues in the region where the rarefaction waves coincide.
Moreover, one can see that in most of the states the system is strictly hyperbolic, but the ordering of the eigenvalues changes throughout the whole domain.
\FloatBarrier
\subsubsection{Solution without contact}
\begin{figure}[h!]
    \begin{center} 
    \subfigure[Densities $\rho_1, \rho_2, \rho$ and volume fraction $\alpha_1$]{
    	\begin{tabular}{cc}        \includegraphics[width=0.35\textwidth]{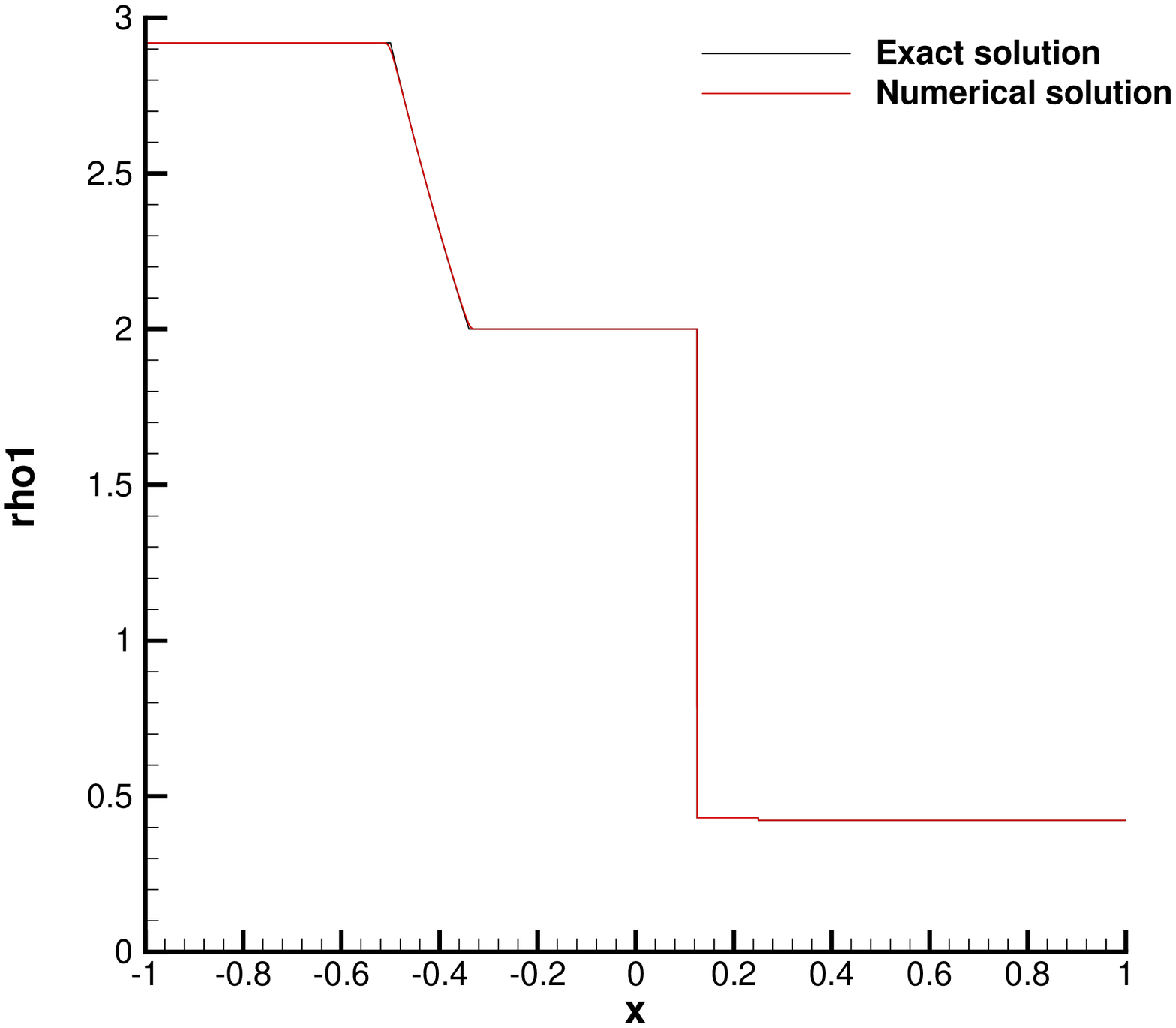}  &  
	\includegraphics[width=0.35\textwidth]{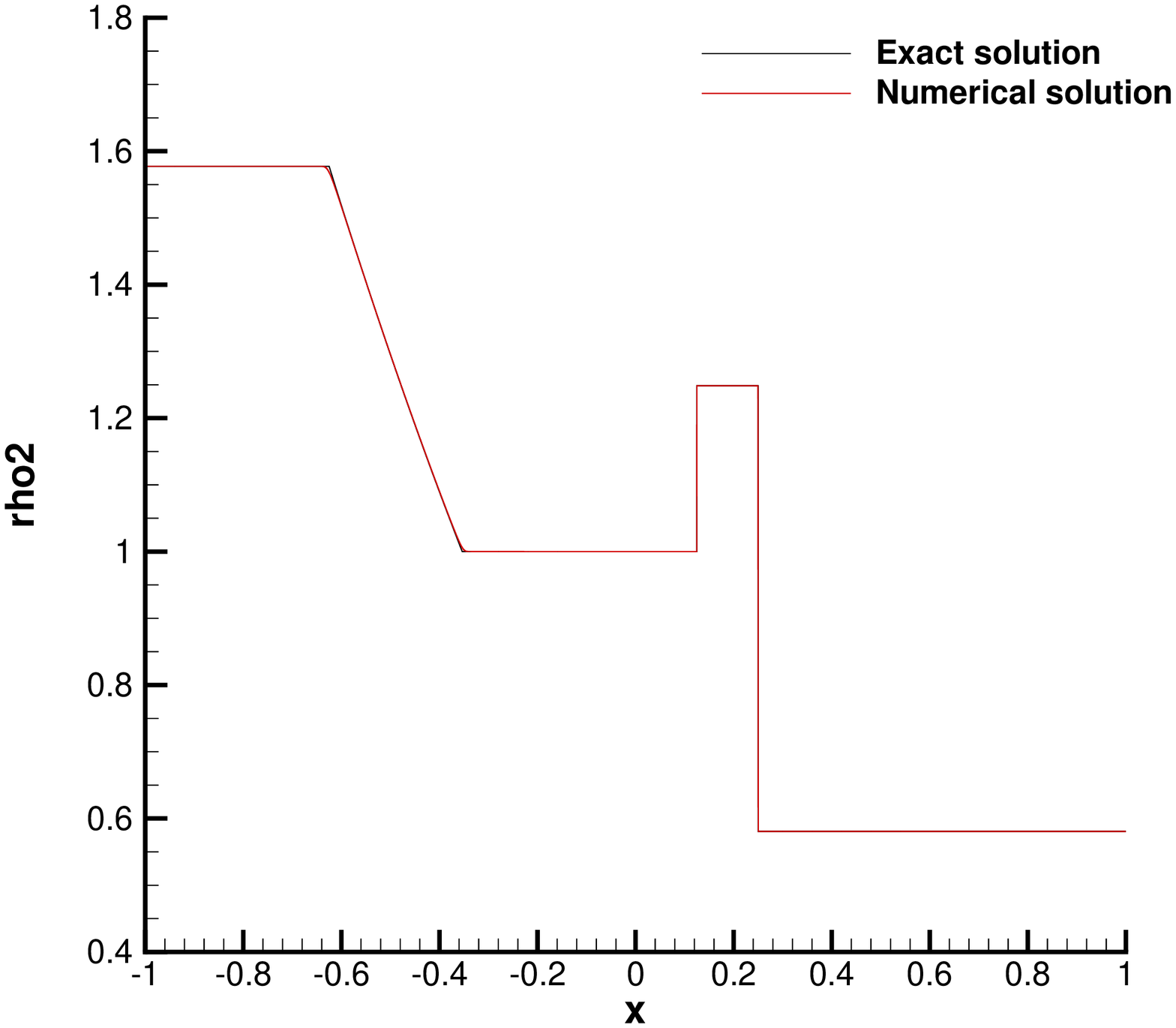}  \\ 	
	\includegraphics[width=0.35\textwidth]{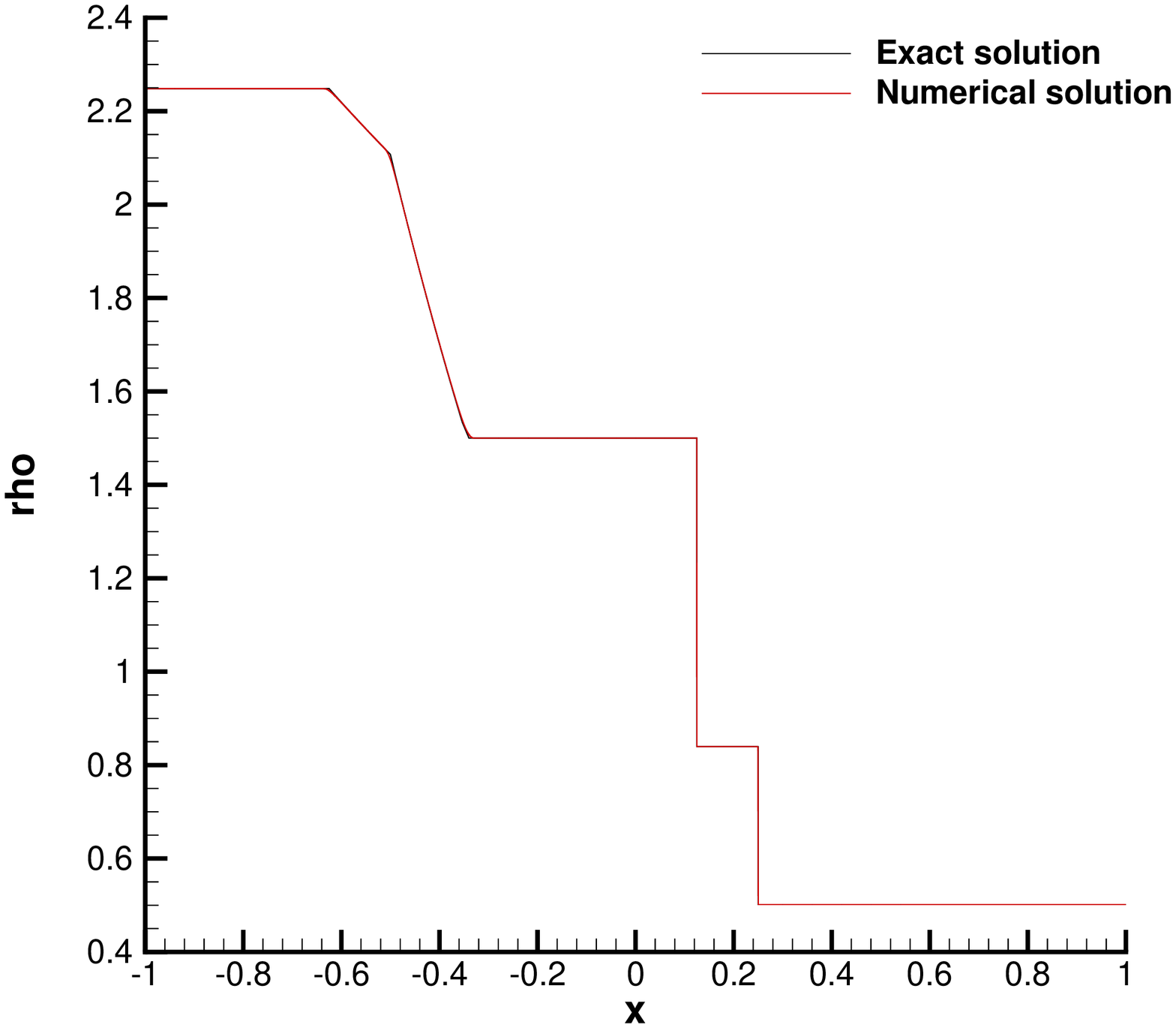}  &  
	\includegraphics[width=0.35\textwidth]{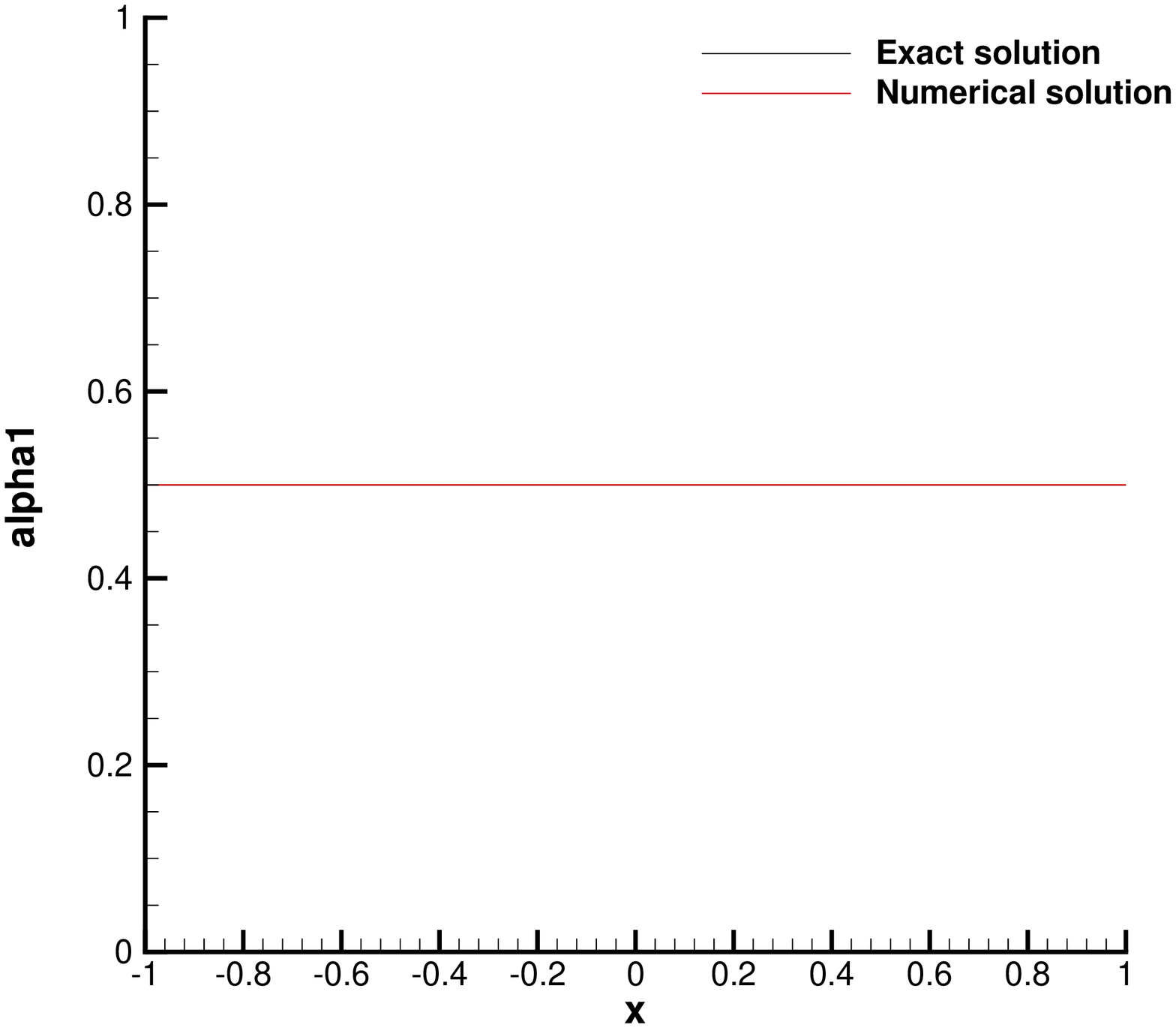}    	
\end{tabular}
	}
    \subfigure[Velocities $u_1, u_2, u ,w$.]{
    	\begin{tabular}{cc}        \includegraphics[width=0.35\textwidth]{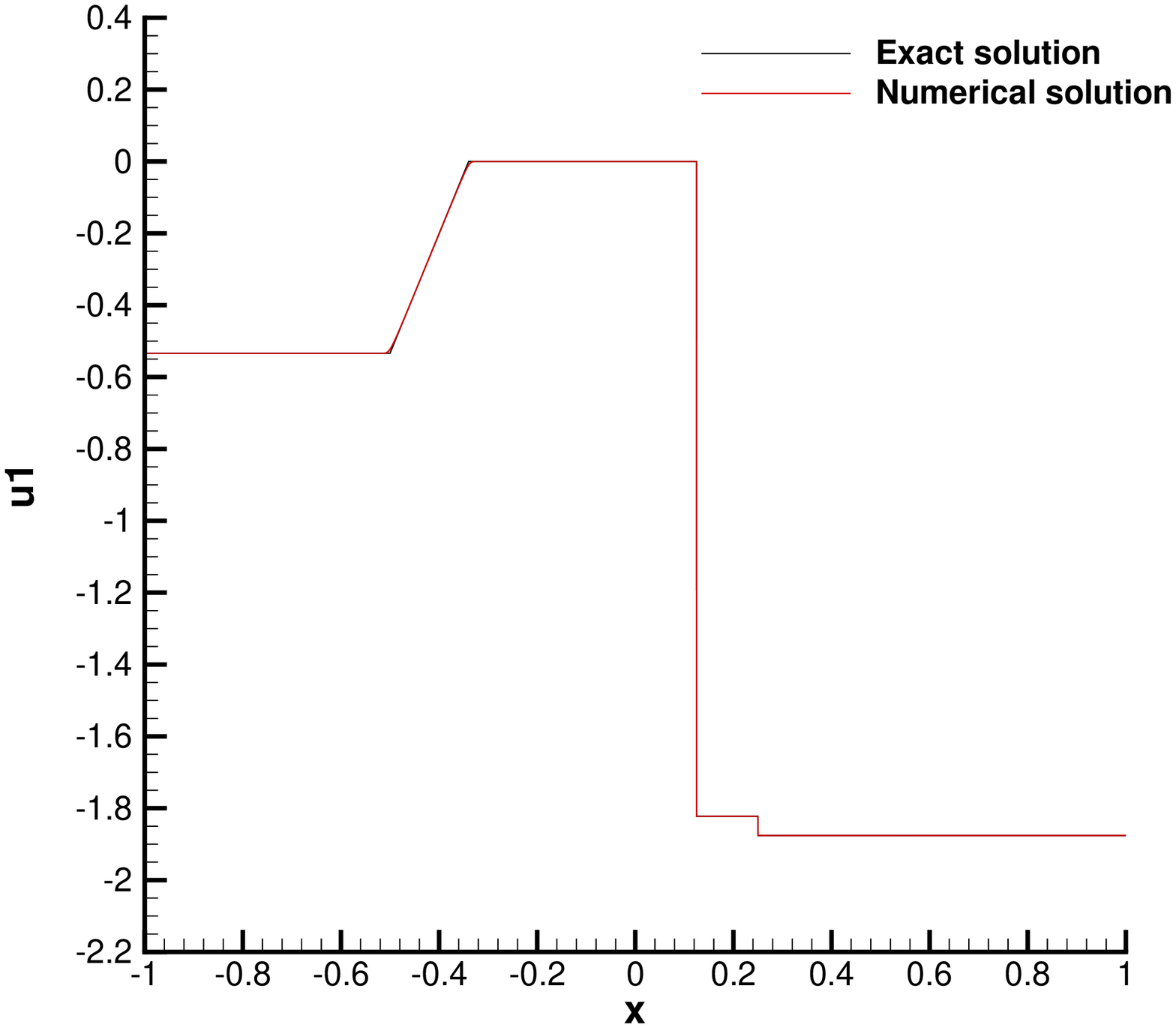}  &  
		\includegraphics[width=0.35\textwidth]{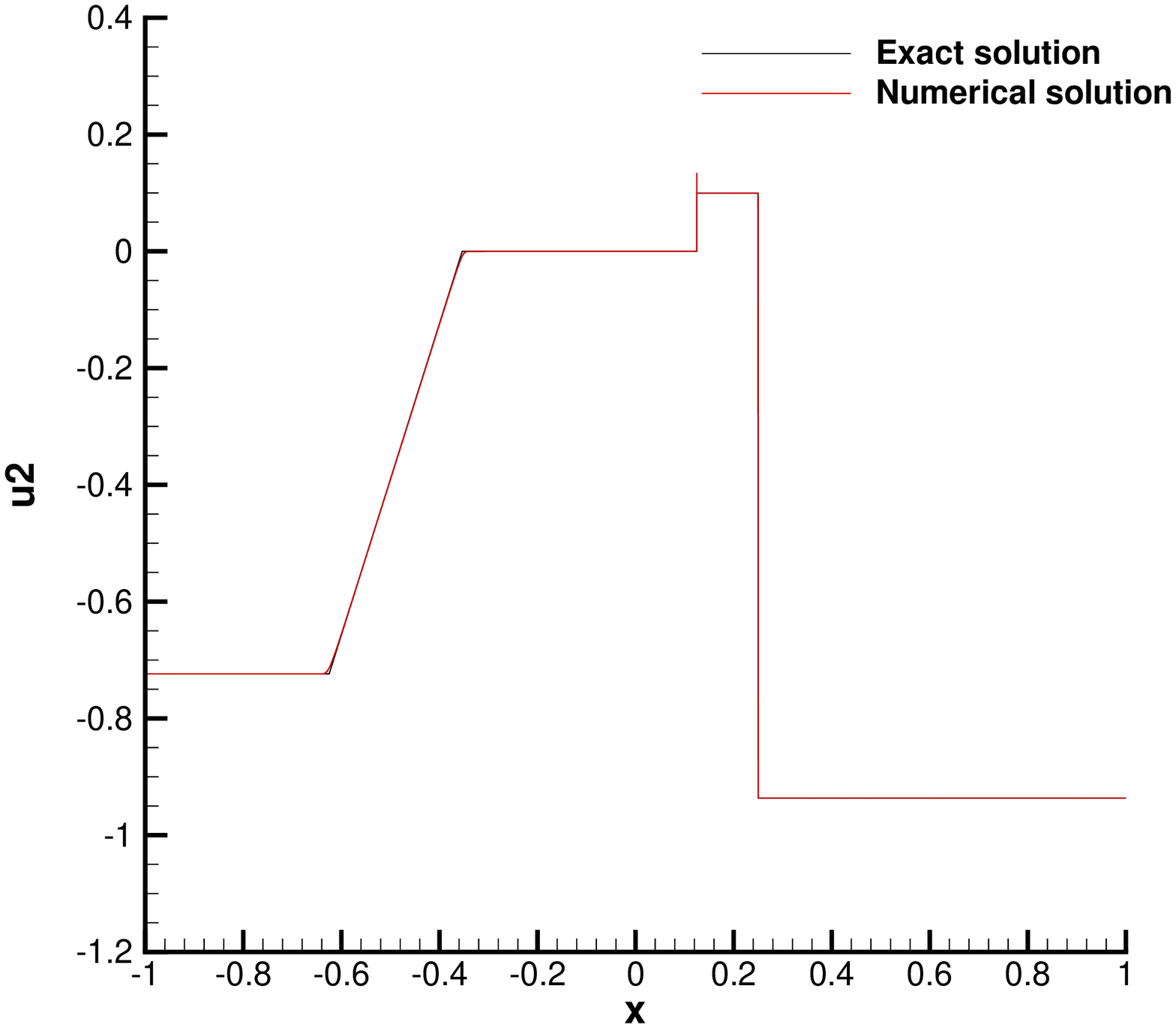}  \\ 	
		\includegraphics[width=0.35\textwidth]{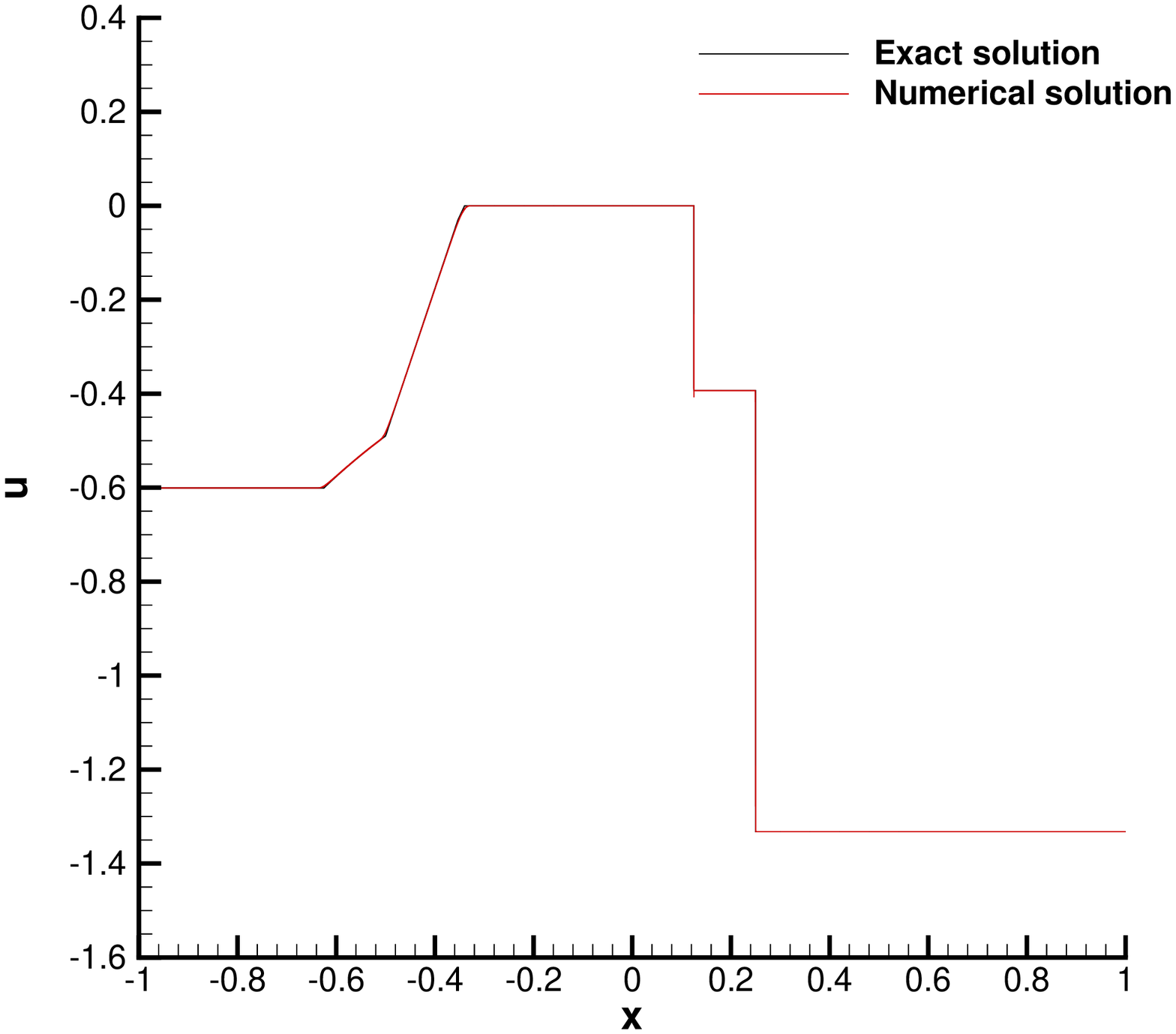}  &  
		\includegraphics[width=0.35\textwidth]{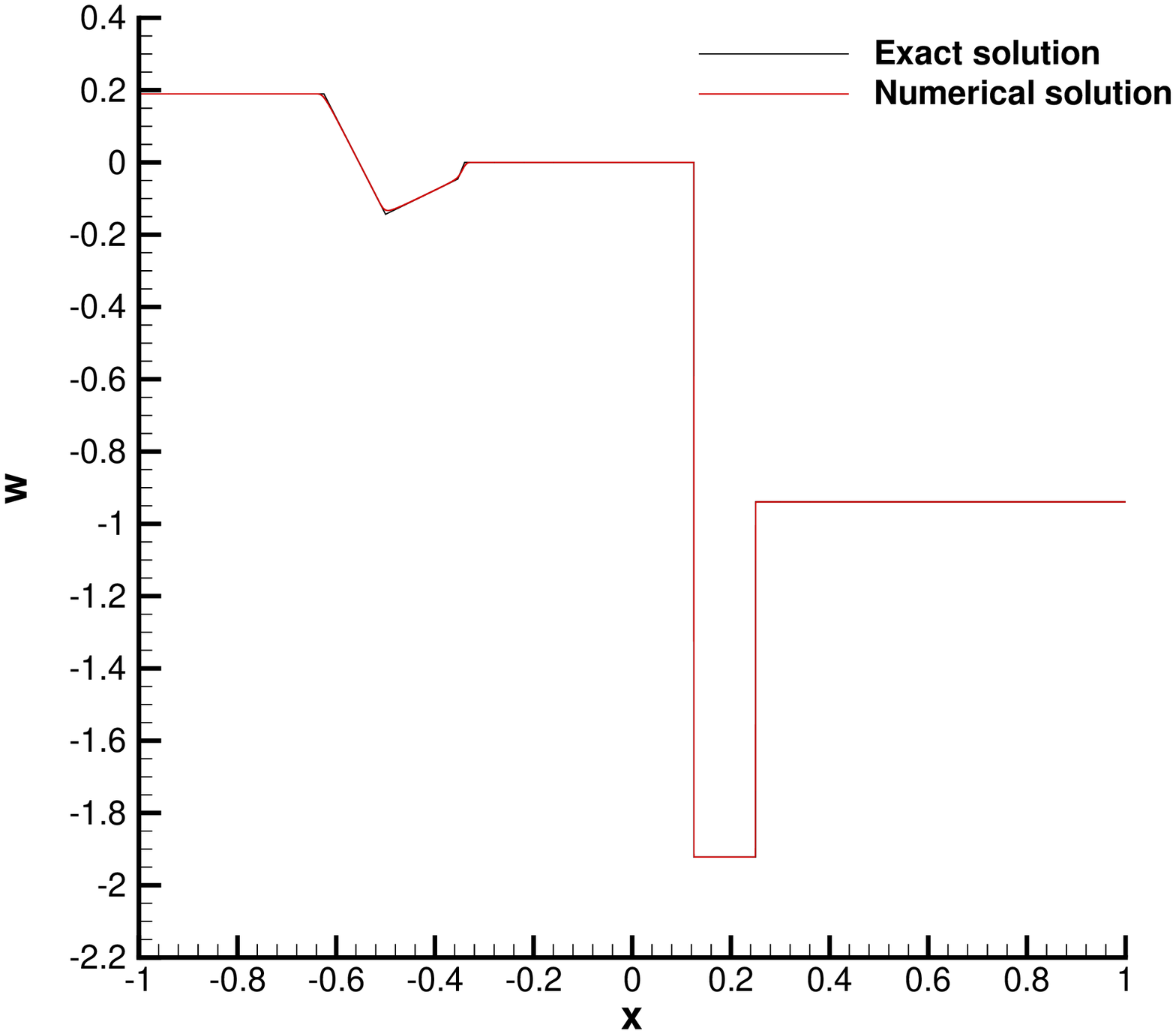}    	
 	\end{tabular}
    }
    \caption{Exact solution (black) and numerical solution (red) of Riemann problem RP2.}
    \label{fig:ex2_p1}
	\end{center} 
\end{figure}

\begin{figure}
	\begin{center} 
    \includegraphics[width=0.8\textwidth]{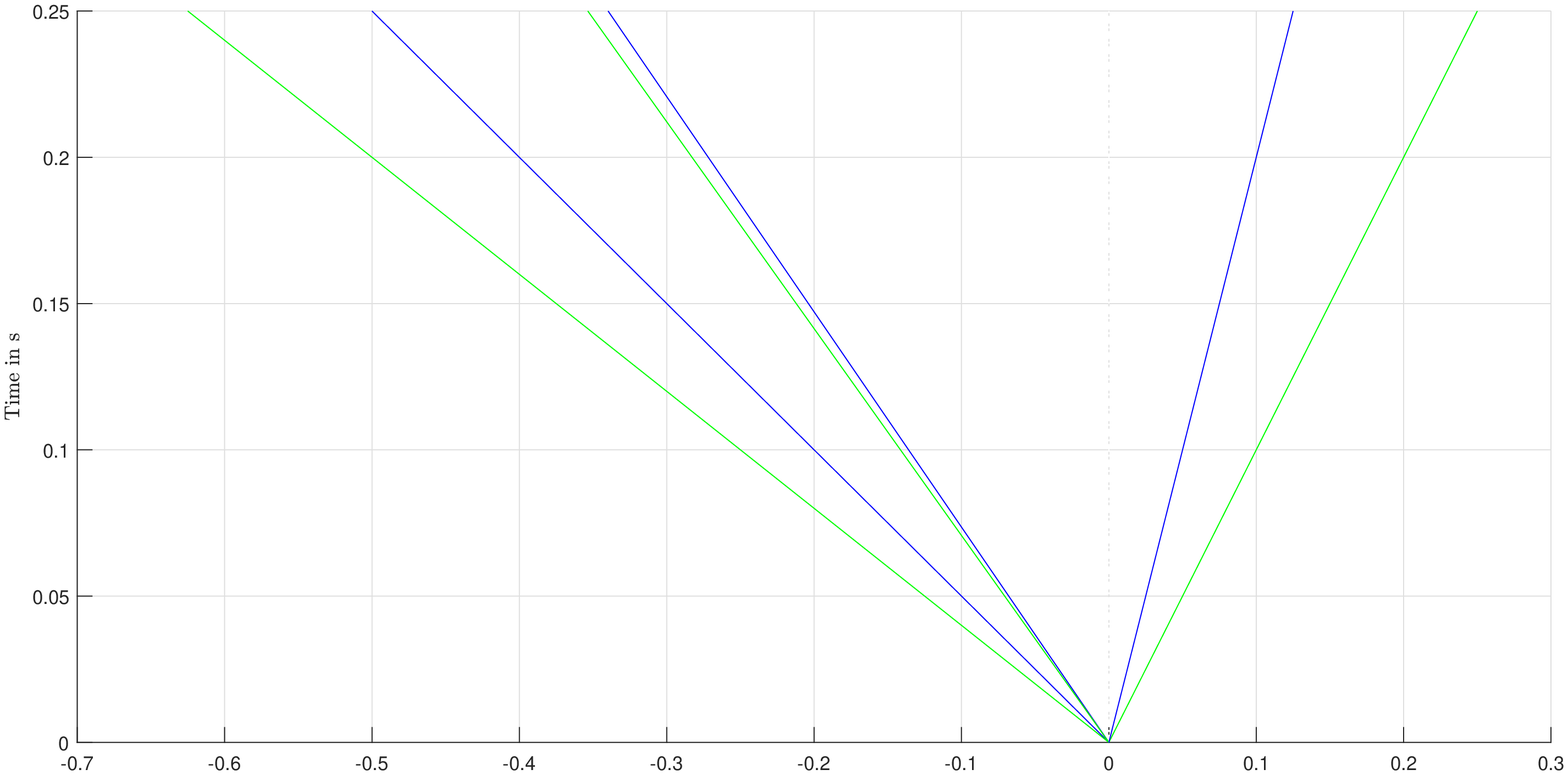} \\
	\includegraphics[width=0.8\textwidth]{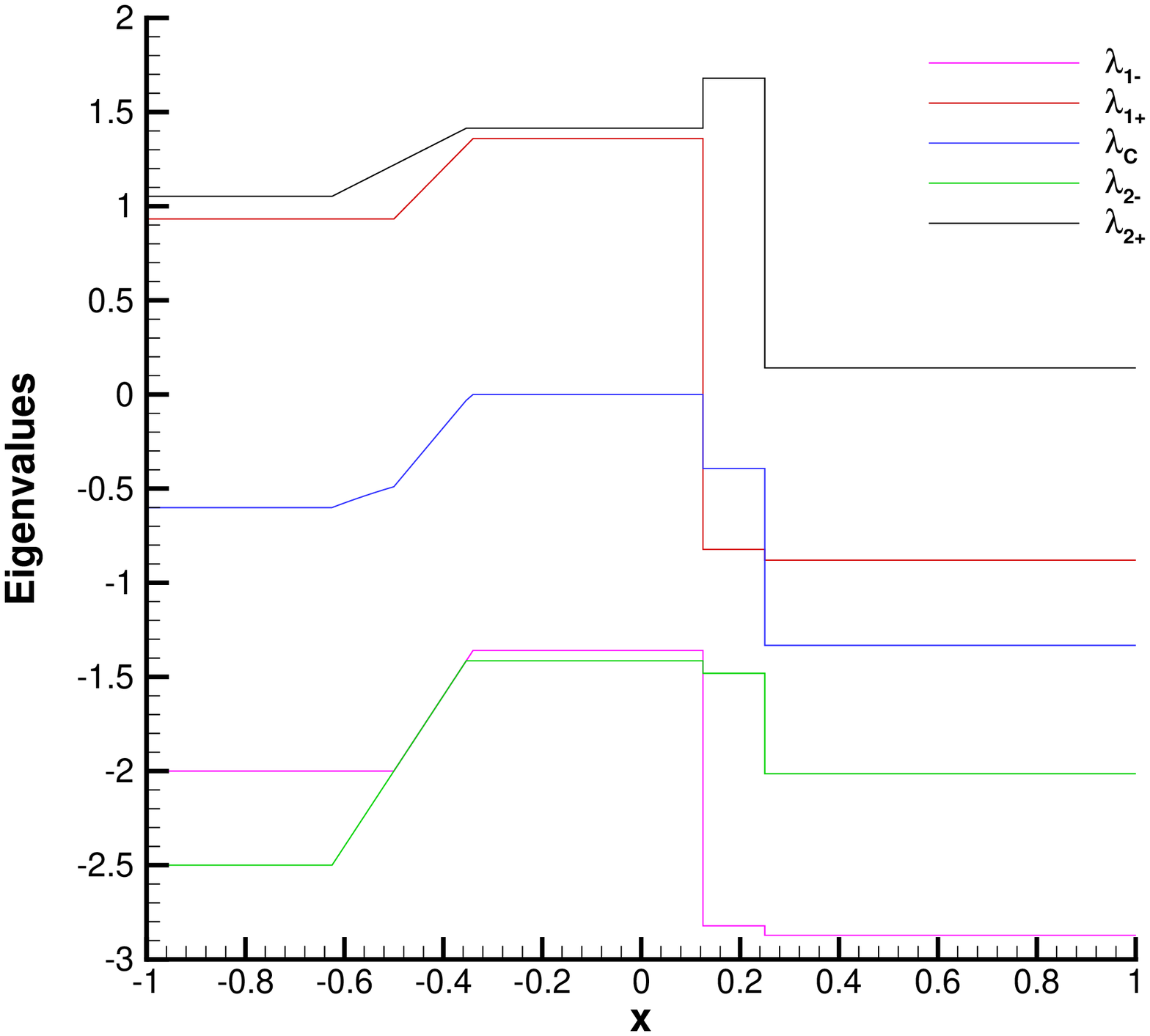}
    \caption{Wave structure of RP2 (top): phase one (blue), contact (red) and phase 2 (green). Eigenvalues of RP2 (bottom).)}
    \label{fig:ex2_p2}
    \end{center} 
\end{figure}

The states of this second Riemann problem RP2 are given in Table \ref{tab:states_prob_2} and we have $\alpha_1 = 0.5$ for all states.
\begin{table}[h!]
    \centering
    \renewcommand*{\arraystretch}{1.25}
    \begin{tabular}{c|llllll}
                 &  \multicolumn{1}{c}{$U_L$}              &  \multicolumn{1}{c}{$U^\ast_L$} &  \multicolumn{1}{c}{$U^{\ast\ast}_L$} &  \multicolumn{1}{c}{$U_R^{\ast\ast}$}
                 &  \multicolumn{1}{c}{$U^\ast_R$}         &  \multicolumn{1}{c}{$U_R$}            \\
        \hline
        $\rho_1$ &  \phantom{-}2.9194 &  \phantom{-}2.9194 &  2 &  2 &  \phantom{-}0.43057 &  \phantom{-}0.42256 \\
        $\rho_2$ &  \phantom{-}1.5773 &  \phantom{-}1      &  1 &  1 &  \phantom{-}1.2486  &  \phantom{-}0.58056 \\
        $u_1$    &  -0.53404          &  -0.53404          &  0 &  0 &  -1.8225            &  -1.876             \\
        $u_2$    &  -0.72386          &  \phantom{-}0      &  0 &  0 &  \phantom{-}0.09954 &  -0.93653
    \end{tabular}
    \caption{Primitive states of Riemann problem RP2.}
    \label{tab:states_prob_2}
\end{table}
The following computation was performed using the ideal gas EOS
\begin{align*}
    p_i(\rho_i) = \rho_i^{\gamma_i},\, i\in\{1,2\}\quad\text{with}\quad \gamma_1 = 1.4,\,\gamma_2 = 2
\end{align*}
and the parameters
\begin{align*}
    \Delta x = 0.5\cdot 10^{-4}\,\si{m},\; C_{CFL} = 0.25,\, t_{end} = 0.25\,\si{s} \quad\text{and}\quad x \in [-1,1]\,\si{m}.
\end{align*}
In Figures \ref{fig:ex2_p1} and \ref{fig:ex2_p2} the numerical results together with the exact solution are shown.
It consists of a $\lambda_{2-}$ - rarefaction (green) which overlaps with the $\lambda_{1-}$ - rarefaction (blue).
The contact is (by construction) not visible since all quantities stay constant across it.
On the right we have an isolated $\lambda_{1+}$-shock (blue) which is followed by an isolated $\lambda_{2+}$-shock (green).
Looking at the densities and velocities in Figure \ref{fig:ex2_p1} we again verify that a rarefaction only affects the related phase.
The interaction of the two rarefactions where they overlap can be observed by the changing slope in the mixture quantities $\rho, u, w$.
The eigenvalues are shown at the right in Figure \ref{fig:ex2_p2}. One clearly sees the overlapping eigenvalues in the region where the rarefaction waves coincide.
Moreover one can see that in most of the states the system is strictly hyperbolic, but the ordering of the eigenvalues changes throughout the whole domain three times.
\FloatBarrier
\subsubsection{Symmetric double rarefaction}
\begin{figure}[h!]
    \begin{center}
    \subfigure[Densities $\rho_1, \rho_2, \rho$ and volume fraction $\alpha_1$.]{
    	\begin{tabular}{cc}        \includegraphics[width=0.35\textwidth]{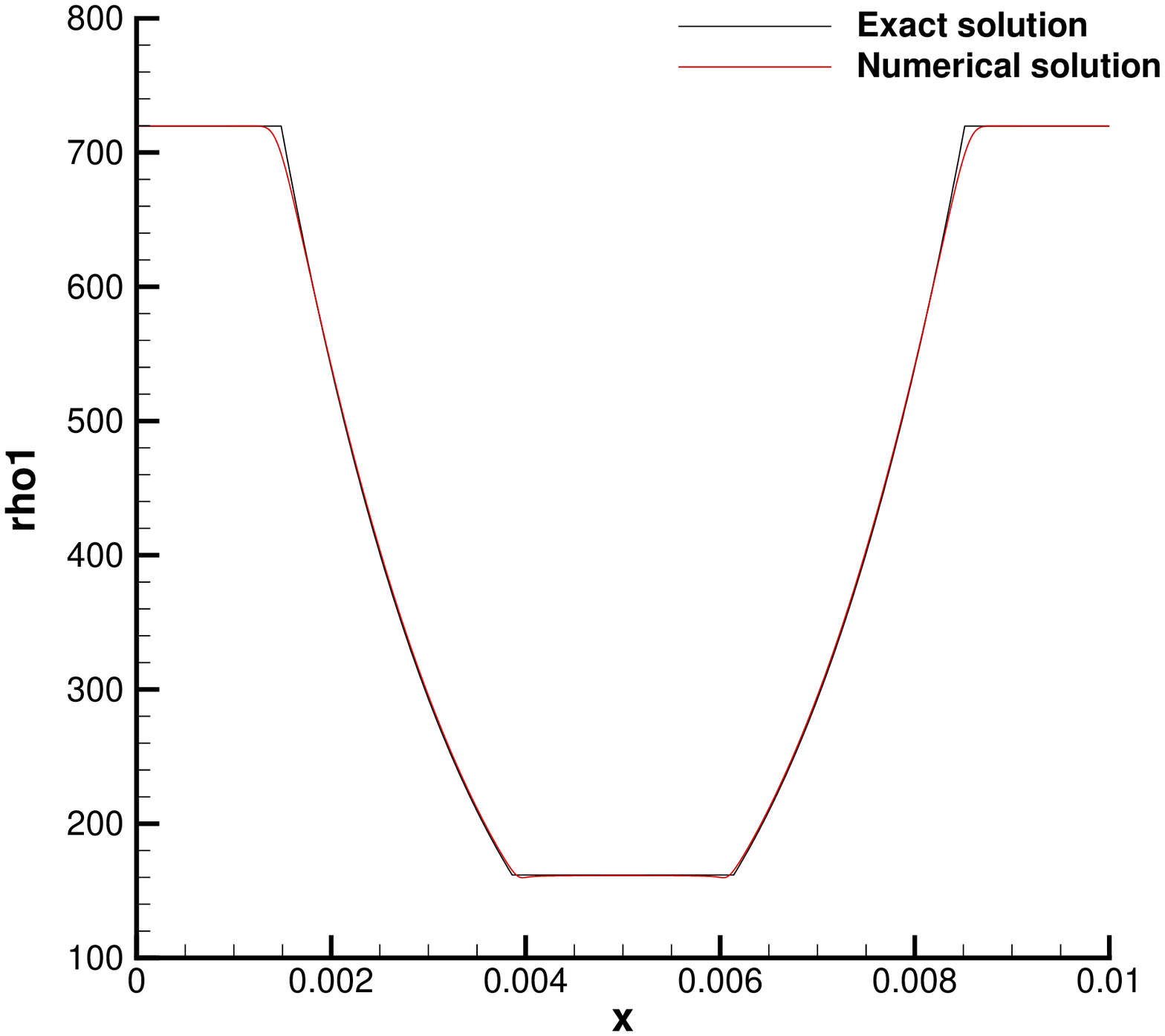}  &  
    	\includegraphics[width=0.35\textwidth]{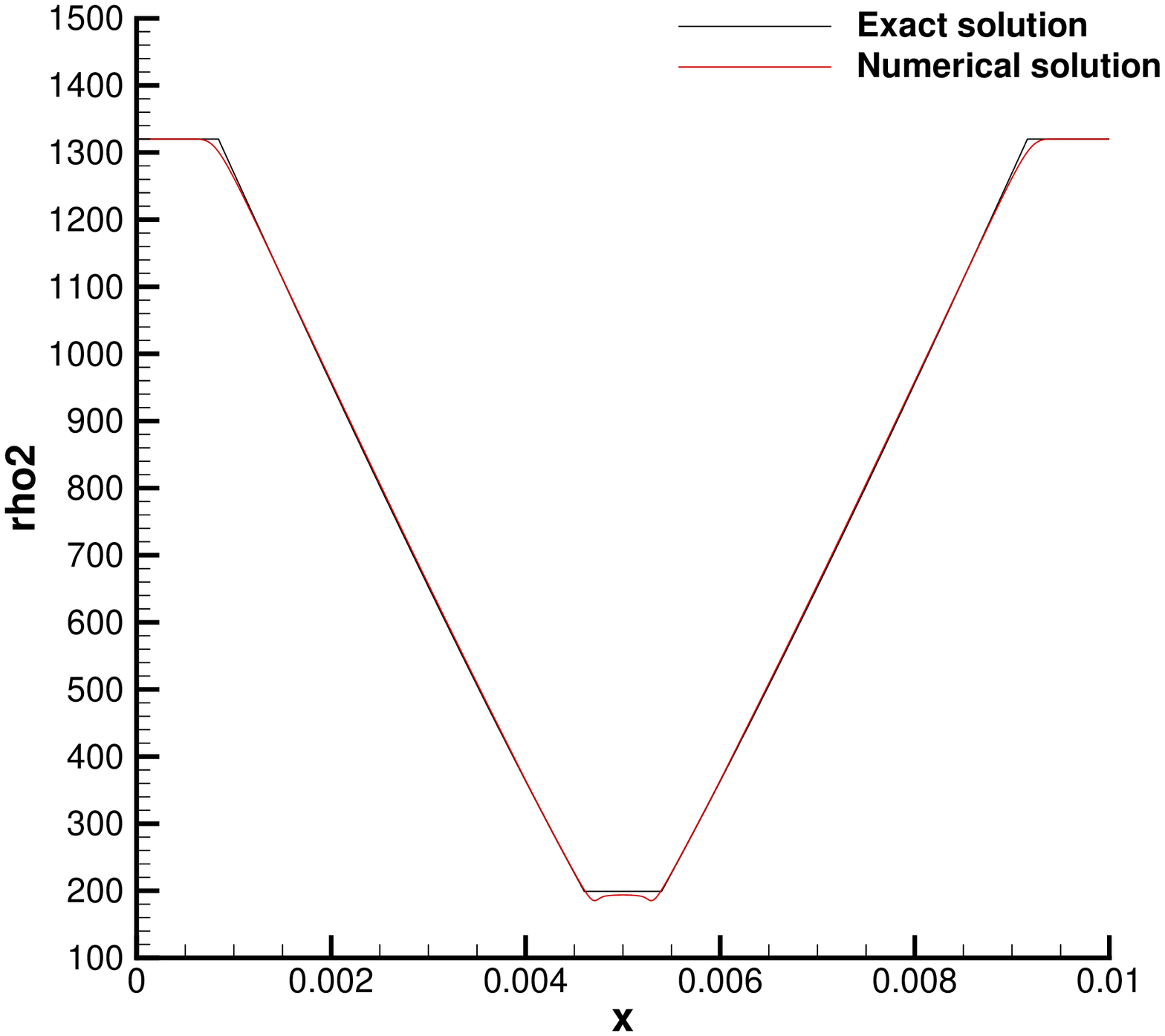}  \\ 	
        \includegraphics[width=0.35\textwidth]{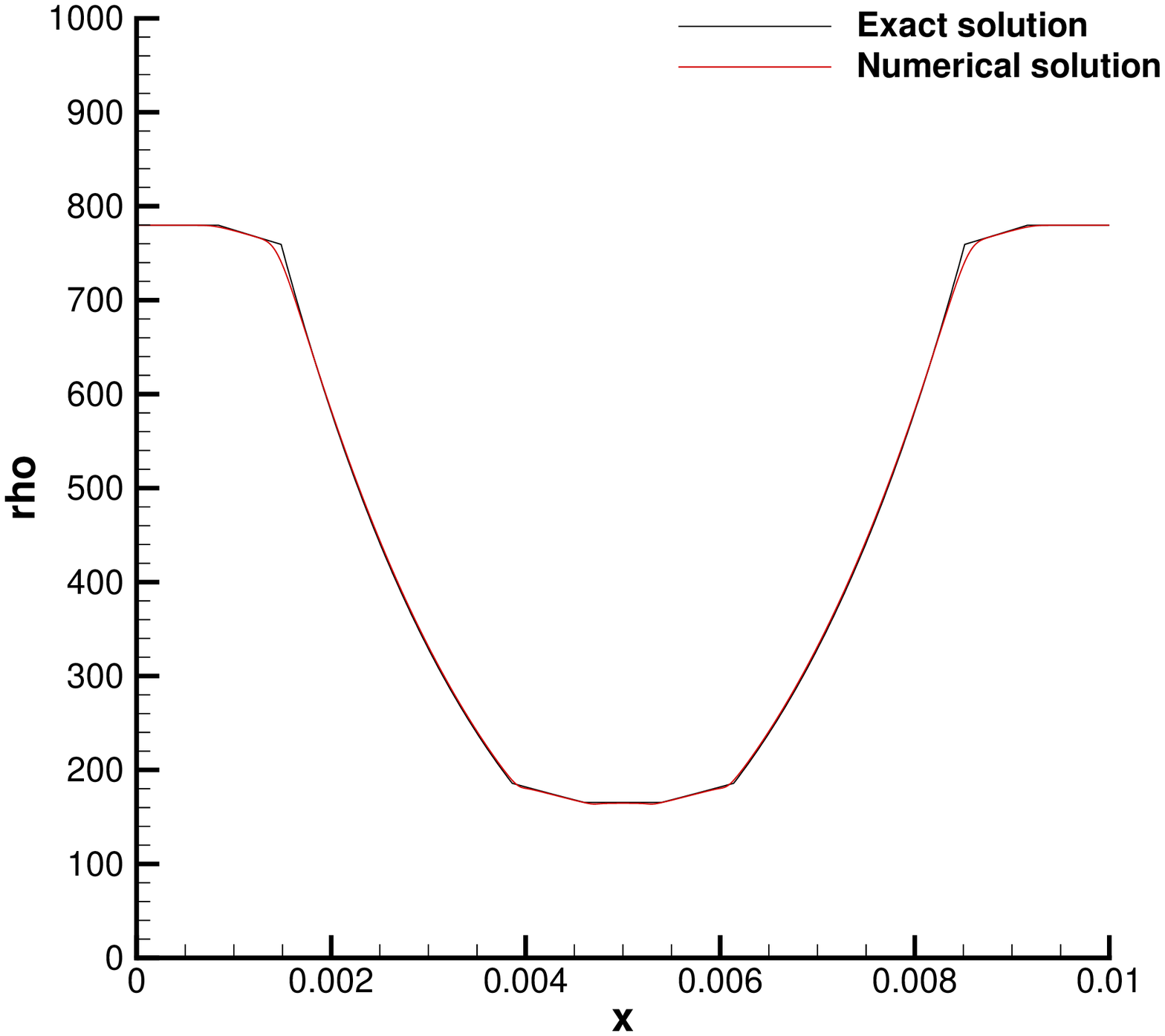}  &  
        \includegraphics[width=0.35\textwidth]{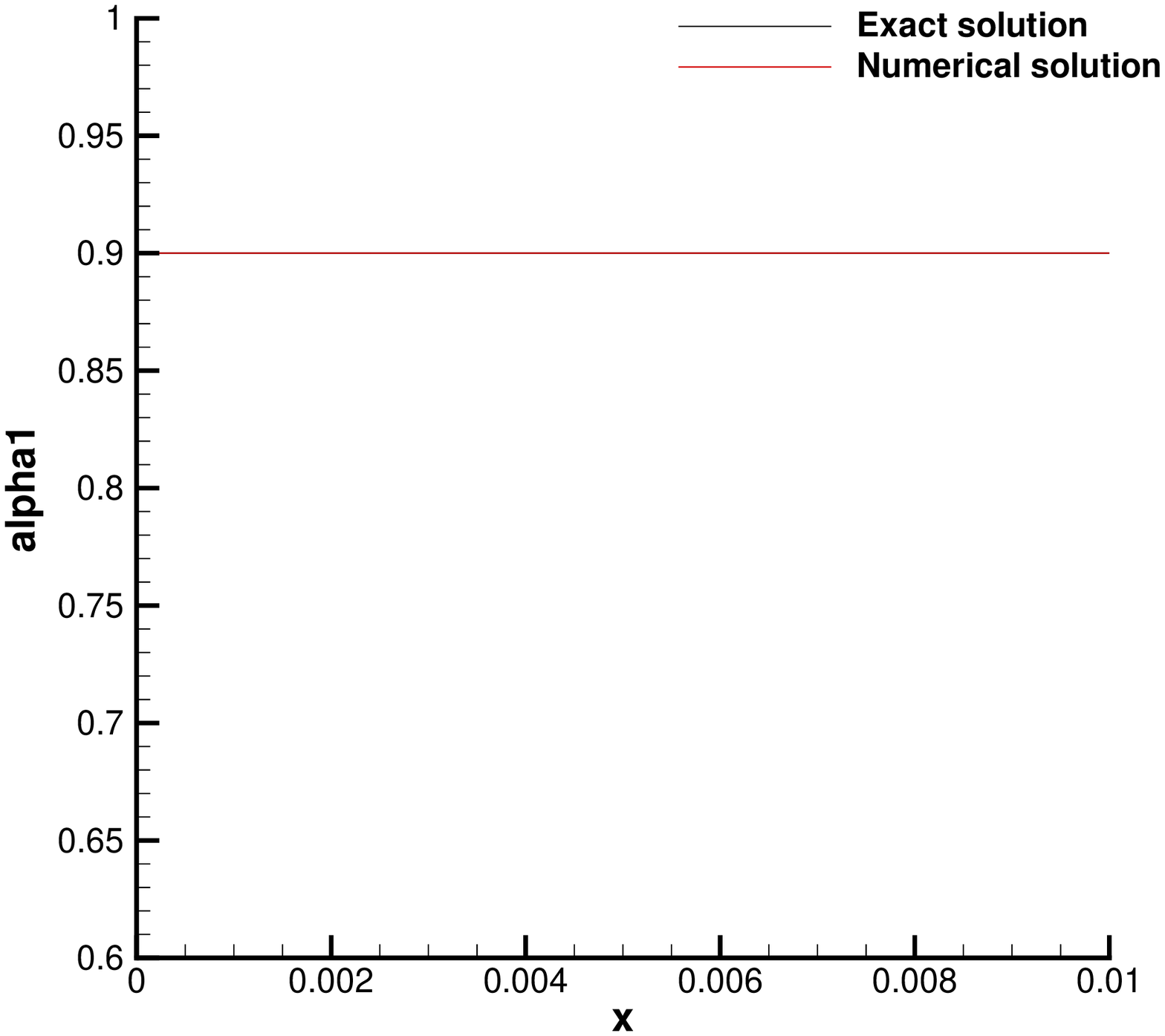}    	
    	\end{tabular}
     }\\
    \subfigure[Velocities $u_1, u_2, u, w$.]{
    	\begin{tabular}{cc}        \includegraphics[width=0.35\textwidth]{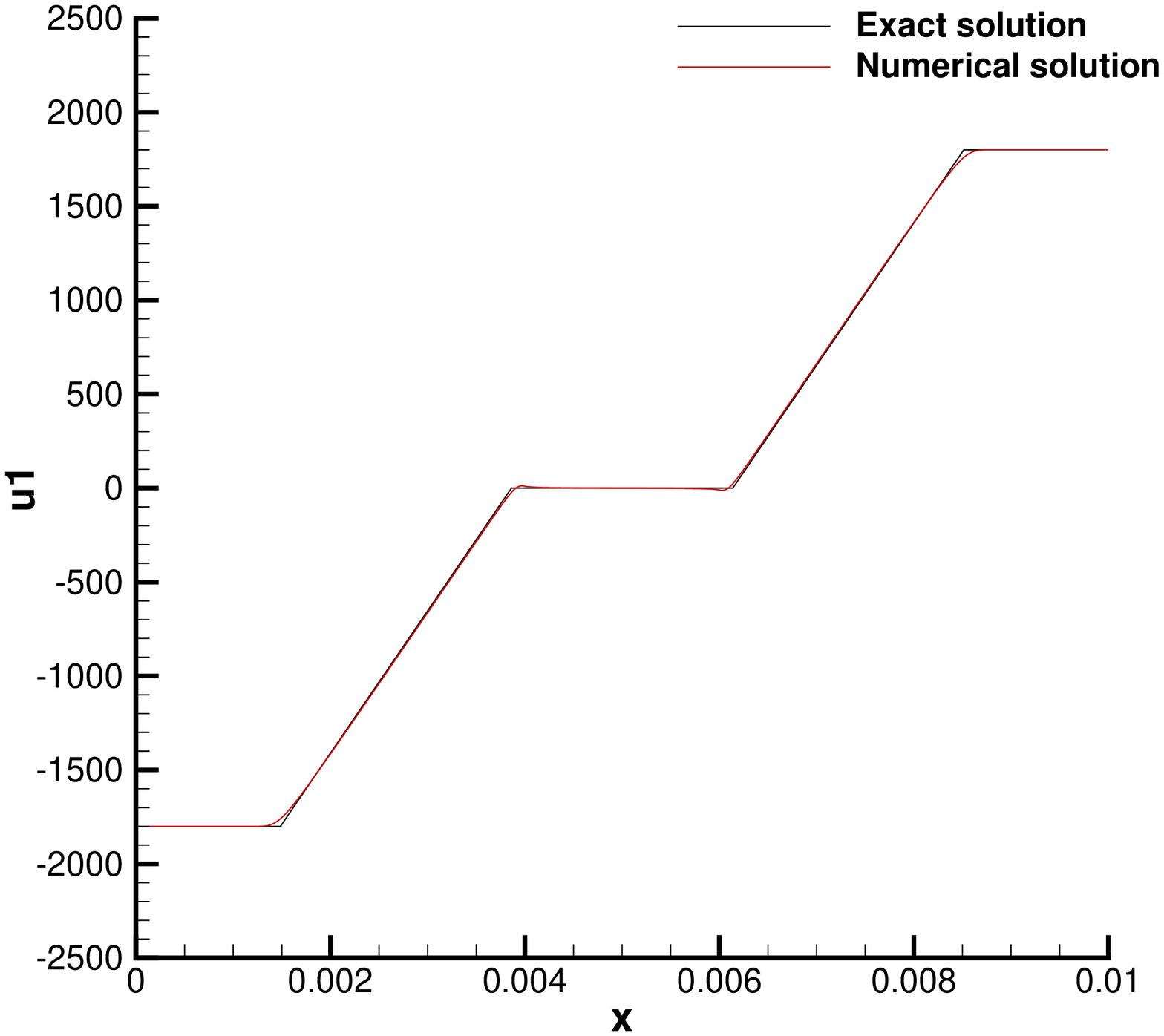}  &  
		\includegraphics[width=0.35\textwidth]{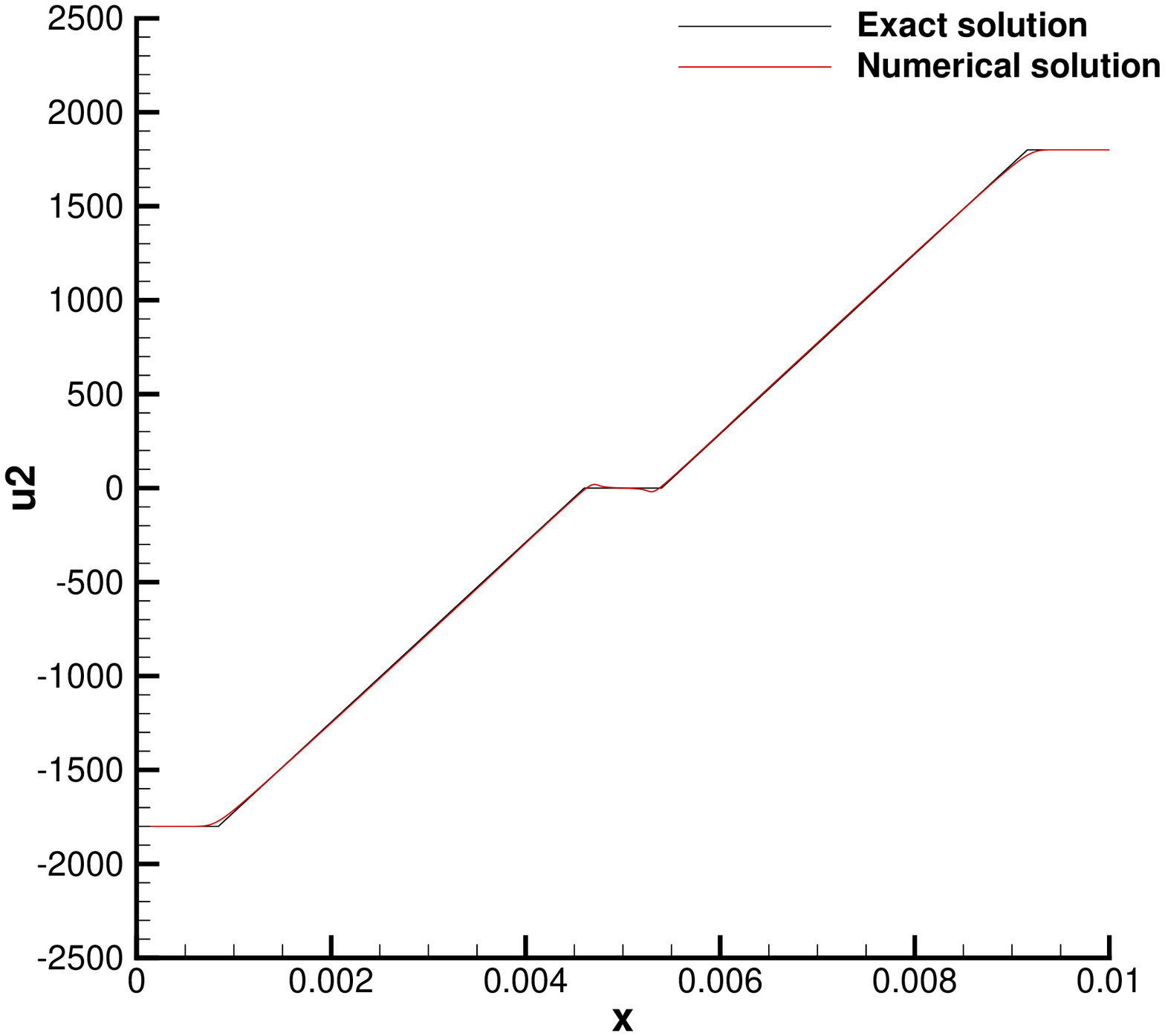}  \\ 	
		\includegraphics[width=0.35\textwidth]{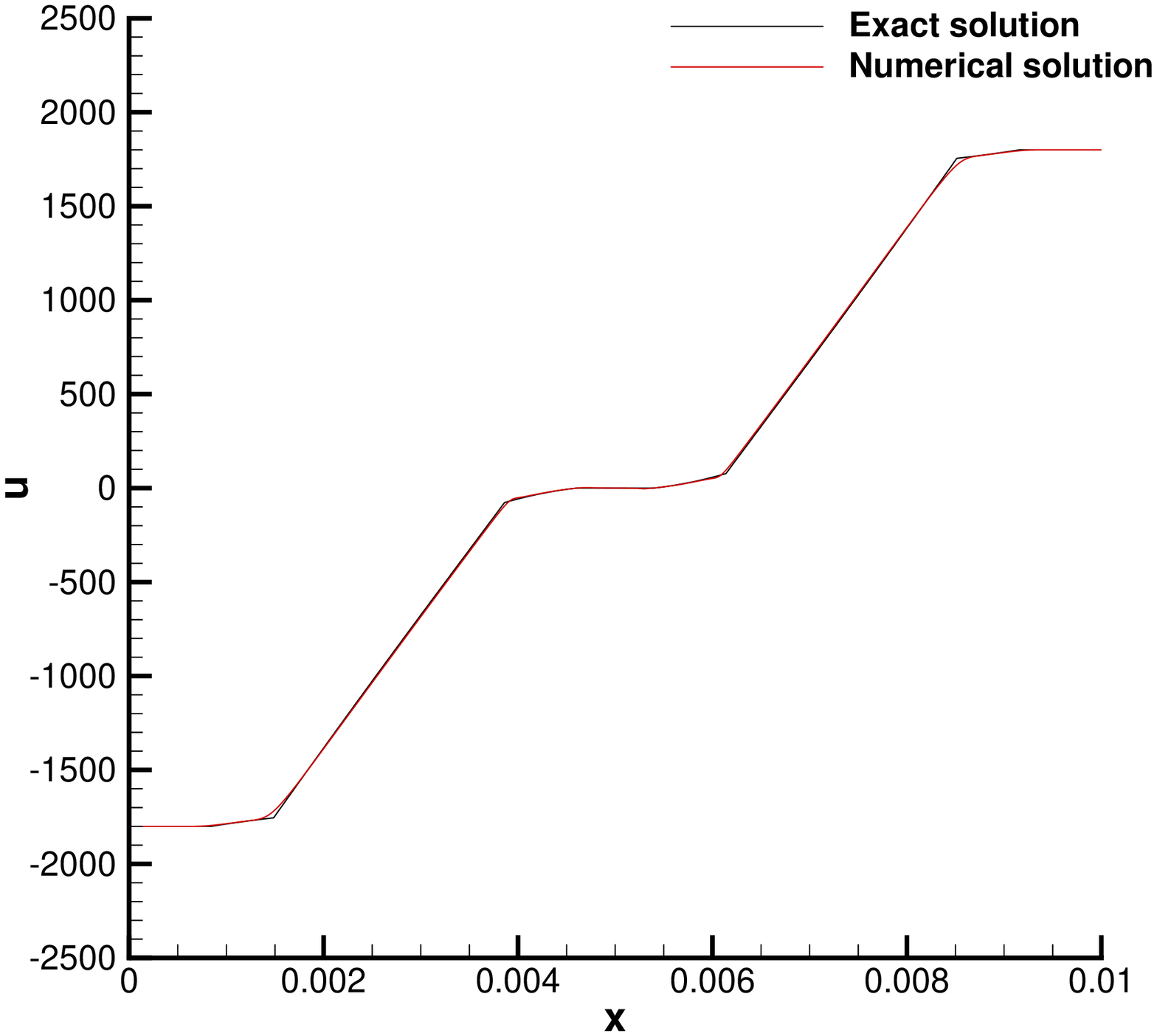}  &  
		\includegraphics[width=0.35\textwidth]{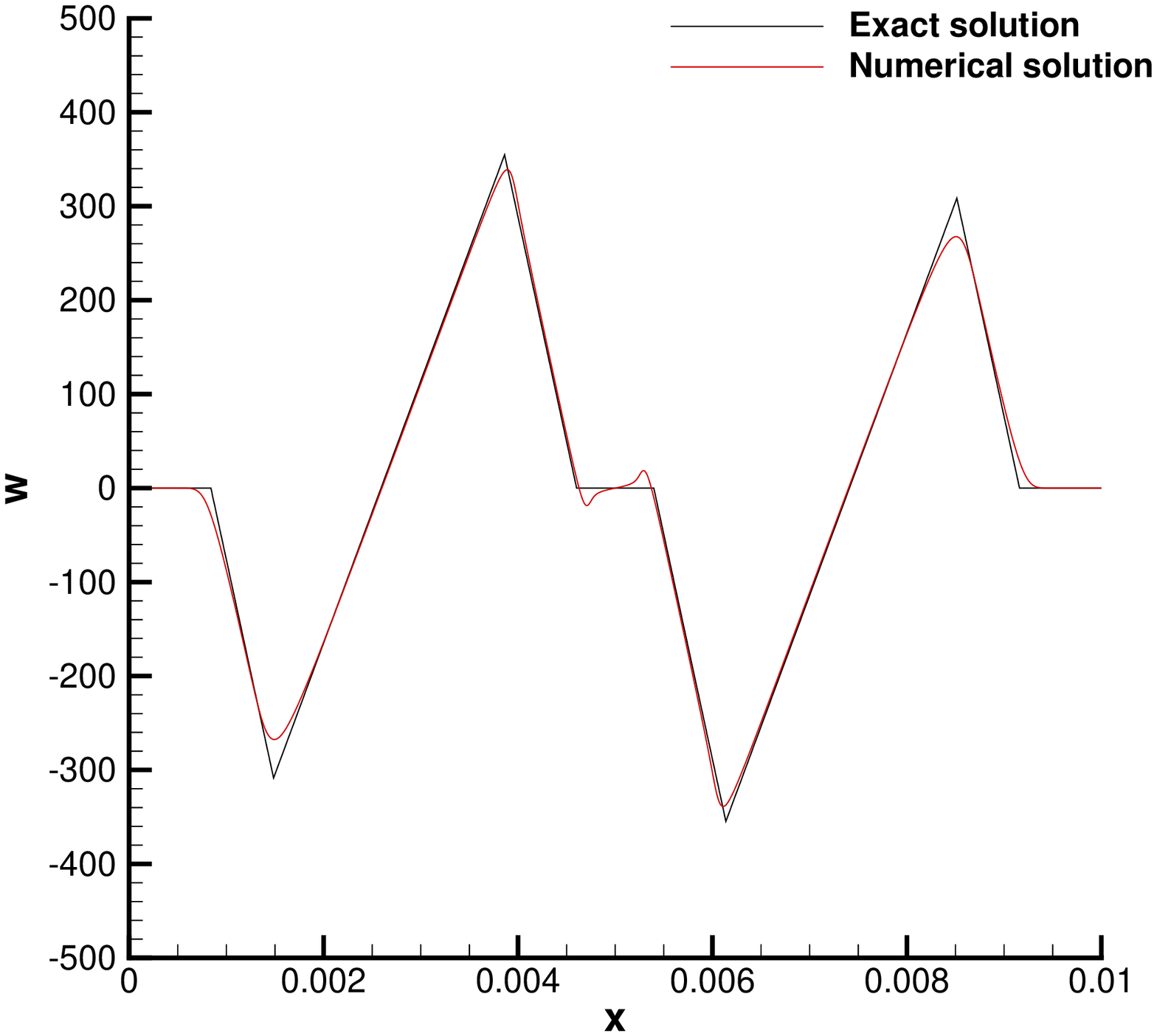}    	
		\end{tabular}
	}
    \caption{Exact solution (black) and numerical solution (red) of Riemann problem RP3.}
    \label{fig:ex3_p1}
    \end{center}
\end{figure}

\begin{figure}
	\begin{center} 
    	\includegraphics[width=0.8\textwidth]{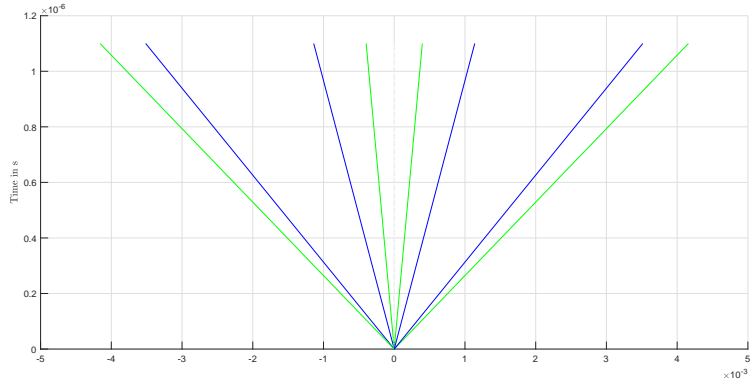} \\
    	\includegraphics[width=0.8\textwidth]{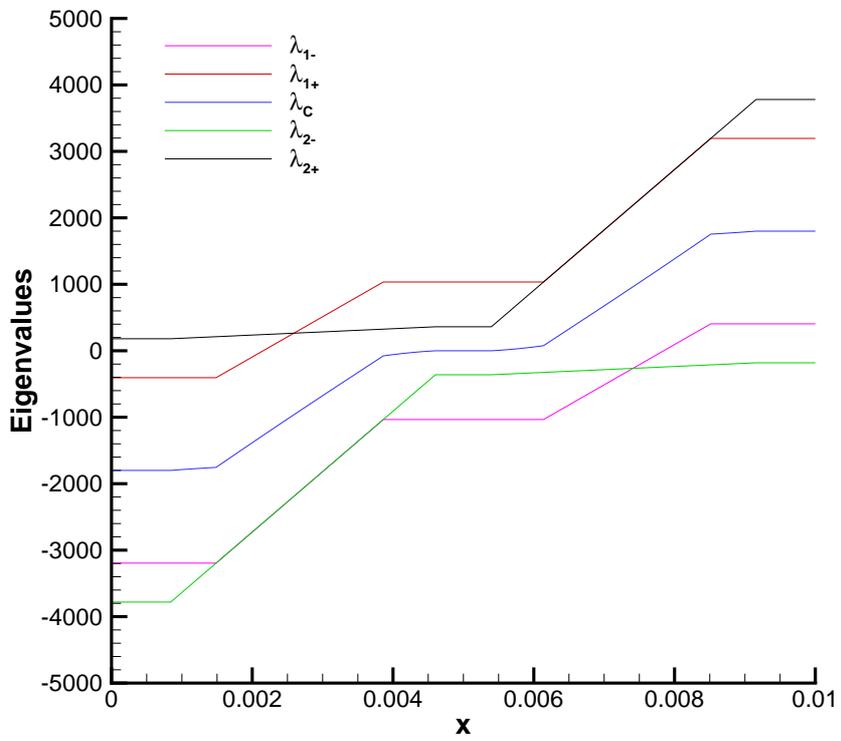}    	
    \caption{Wave Structure of RP3 (top): phase one (blue), contact (red) and phase 2 (green). Eigenvalues of RP3 (bottom). }
    \label{fig:ex3_p2}
    \end{center} 
\end{figure}

The following example (RP3) is taken from a paper of Romenski and Toro \cite{Romenski2004}.
The exact solution is obtained by using a contact centered inverse construction of the solution as described previously.
The states are given in Table \ref{tab:states_prob_3} and we have $\alpha_1 = 0.9$ for all states.
\begin{table}[h!]
    \centering
    \renewcommand*{\arraystretch}{1.25}
    \begin{tabular}{c|llllll}
                 &  \multicolumn{1}{c}{$U_L$}              &  \multicolumn{1}{c}{$U^\ast_L$} &  \multicolumn{1}{c}{$U^{\ast\ast}_L$} &  \multicolumn{1}{c}{$U_R^{\ast\ast}$}
                 &  \multicolumn{1}{c}{$U^\ast_R$}         &  \multicolumn{1}{c}{$U_R$}            \\
        \hline
        $\rho_1$ &  \phantom{-}789.79932 &   \phantom{-}160       & 160 &   160 &   160       & 789.79932 \\
        $\rho_2$ &  \phantom{-}1270.0579 &   \phantom{-}1270.0579 & 200 &   200 &   1270.0579 & 1270.0579 \\
        $u_1$    &   -1942.0873          &   \phantom{-}0         & 0   &   0   &   0         & 1942.0873 \\
        $u_2$    &   -1722.9353          &   -1722.9354           & 0   &   0   &   1722.9354 & 1722.9354
    \end{tabular}
    \caption{Primitive states for Riemann problem RP3.}
    \label{tab:states_prob_3}
\end{table}
The following computation was performed using the following EOS
\begin{align*}
    p_i(\rho_i) &= A_i\left(\frac{\rho_i}{\rho_{ref}^{(i)}}\right)^{\gamma_i} + B_i,\, i\in\{1,2\},\\
    \text{with}\quad A_1 &= 10^5\,\si{Pa},\,\gamma_1 = 1.4,\, \rho_{ref}^{(1)} = 1\,\si{kg m^{-3}},\, B_1 = 0\,\si{Pa},\\
    \text{and}\quad A_2 &= 8.5\cdot 10^8\,\si{Pa},\,\gamma_2 = 2.8,\, \rho_{ref}^{(2)} = 10^3\,\si{kg m^{-3}},\, B_2 = 8.4999\cdot 10^8\,\si{Pa}.
\end{align*}
and the parameters
\begin{align*}
    \Delta x = \frac{0.01}{5\cdot 10^3}\,\si{m},\; C_{CFL} = 0.25,\, t_{end} = 0.11\cdot 10^{-5}\,\si{s} \quad\text{and}\quad x \in [0,0.01]\,\si{m}.
\end{align*}
In Figures \ref{fig:ex3_p1} and \ref{fig:ex3_p2} the numerical results together with the exact solution are shown.
For the numerical solution we used the Force Flux together with Godunov's method as exemplary shown in \cite{Toro2009}.
This example again shows the behaviour of rarefaction waves quite nicely.
Each rarefaction is seen individually in the corresponding phase but the interaction of the overlapping rarefaction waves is observed in the mixture quantities.
In Figure \ref{fig:ex3_p2} on can see that the rarefaction waves of phase one are contained in the rarefaction waves of phase two and that the corresponding eigenvalues coincide where they overlap.
\FloatBarrier
\subsubsection{Symmetric Double Shock}
\begin{figure}[h!]
    \begin{center}
    \subfigure[Densities $\rho_1, \rho_2, \rho$ and volume fraction $\alpha_1$.]{
    \begin{tabular}{cc}     
		\includegraphics[width=0.35\textwidth]{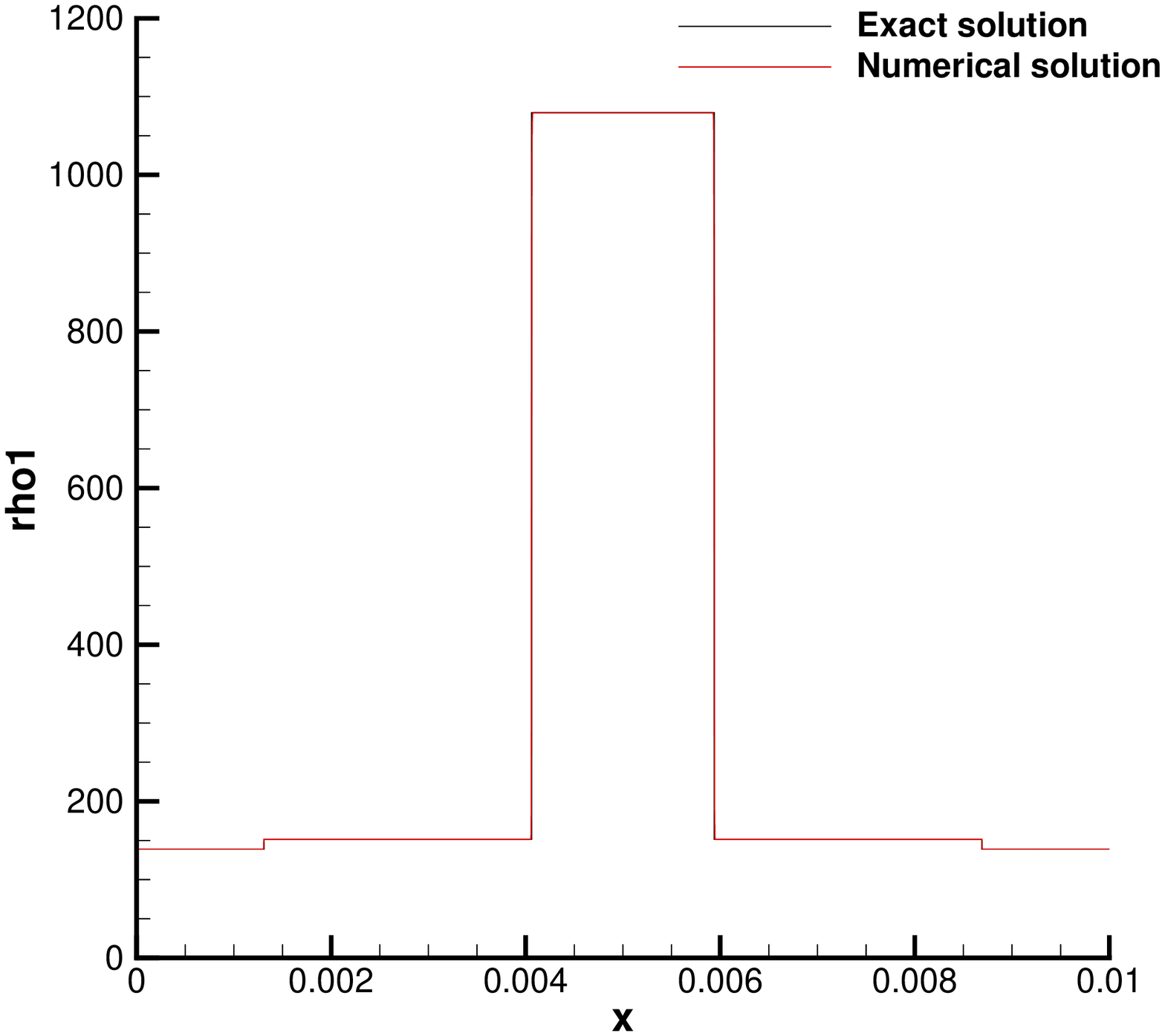}  &  
		\includegraphics[width=0.35\textwidth]{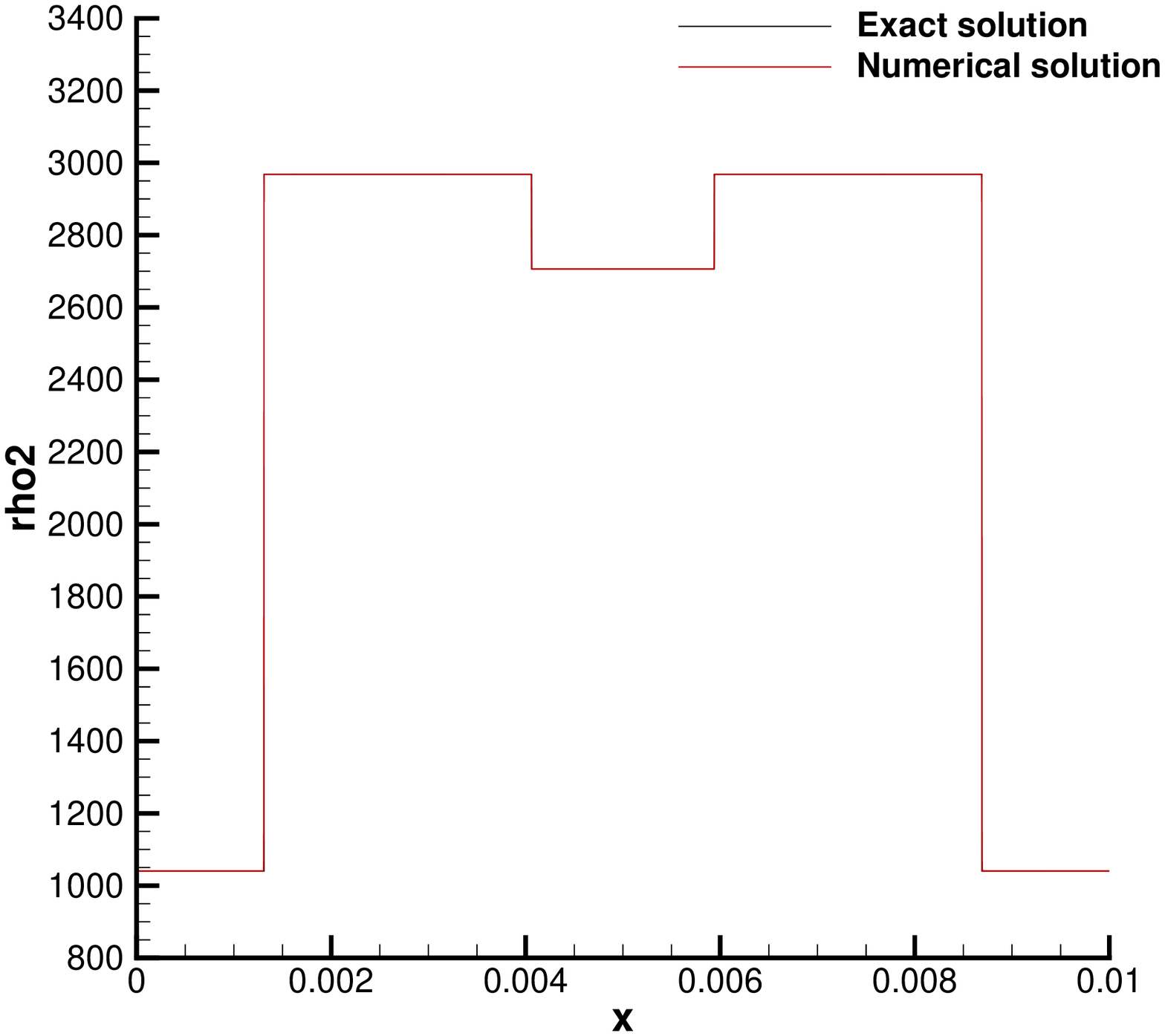}  \\ 	
		\includegraphics[width=0.35\textwidth]{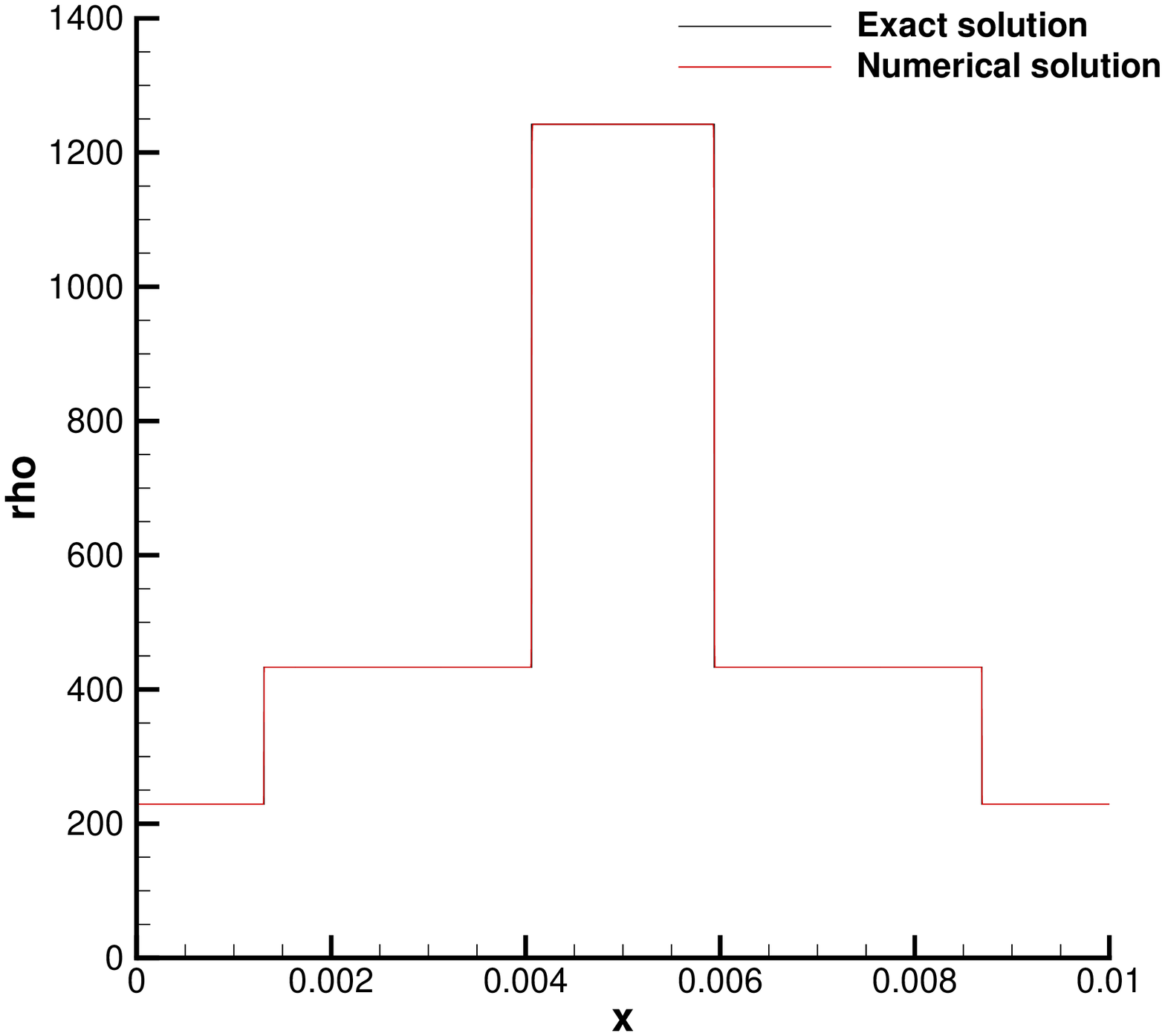}  &  
		\includegraphics[width=0.35\textwidth]{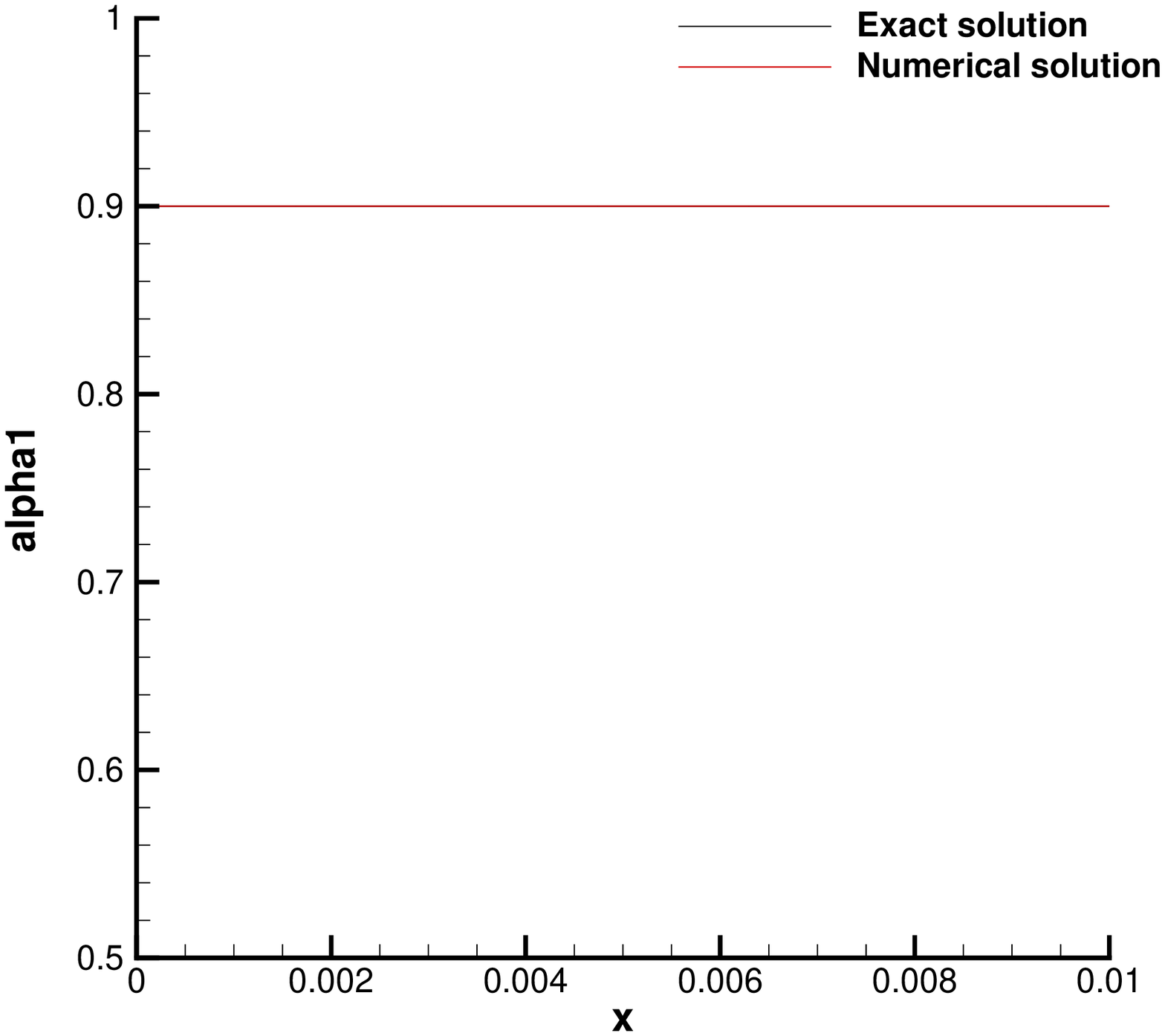}    	
	\end{tabular}     
	}\\
    \subfigure[Velocities $u_1, u_2, u, w$.]{
    	\begin{tabular}{cc}        \includegraphics[width=0.35\textwidth]{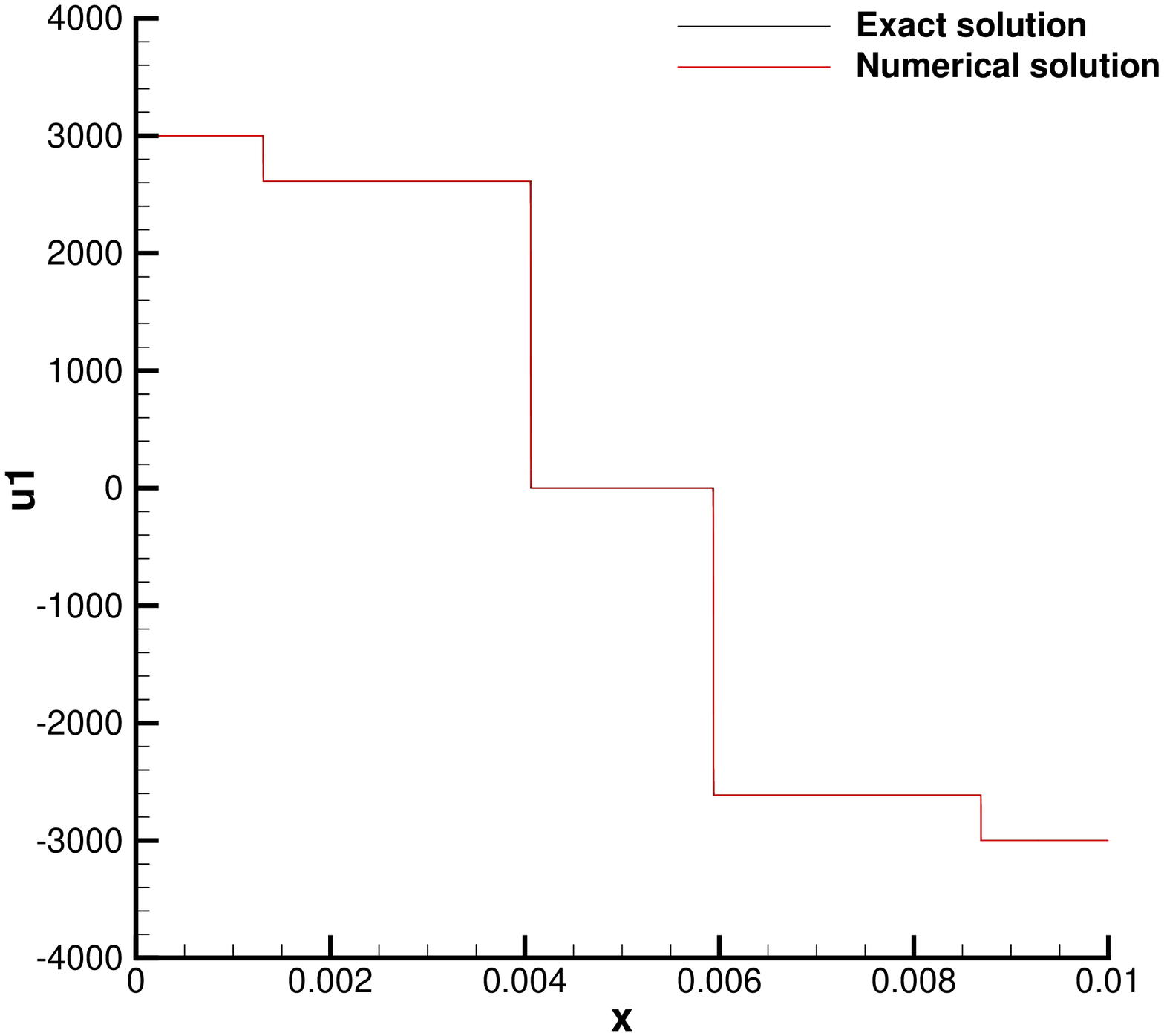}  &  
		\includegraphics[width=0.35\textwidth]{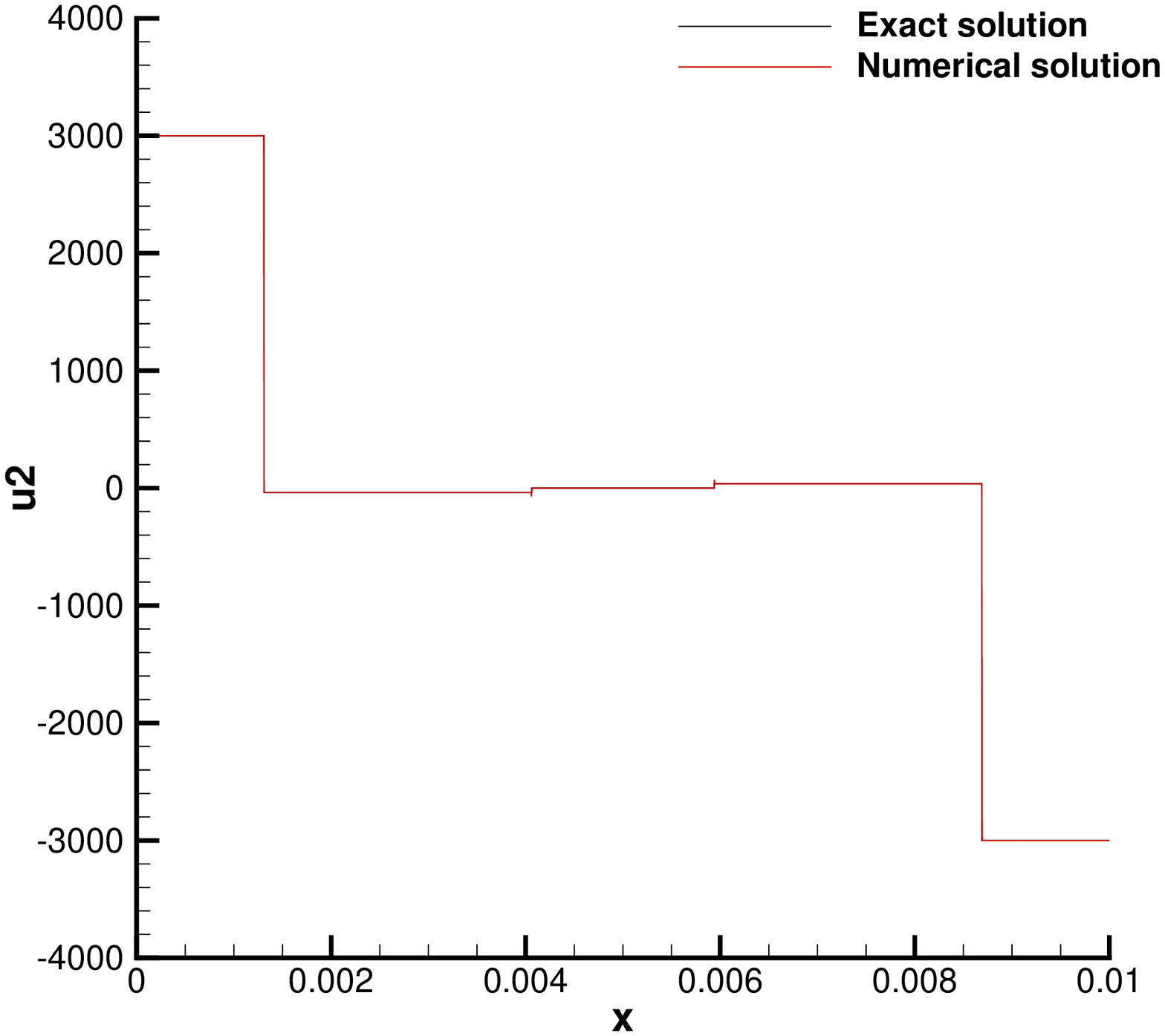}  \\ 	
		\includegraphics[width=0.35\textwidth]{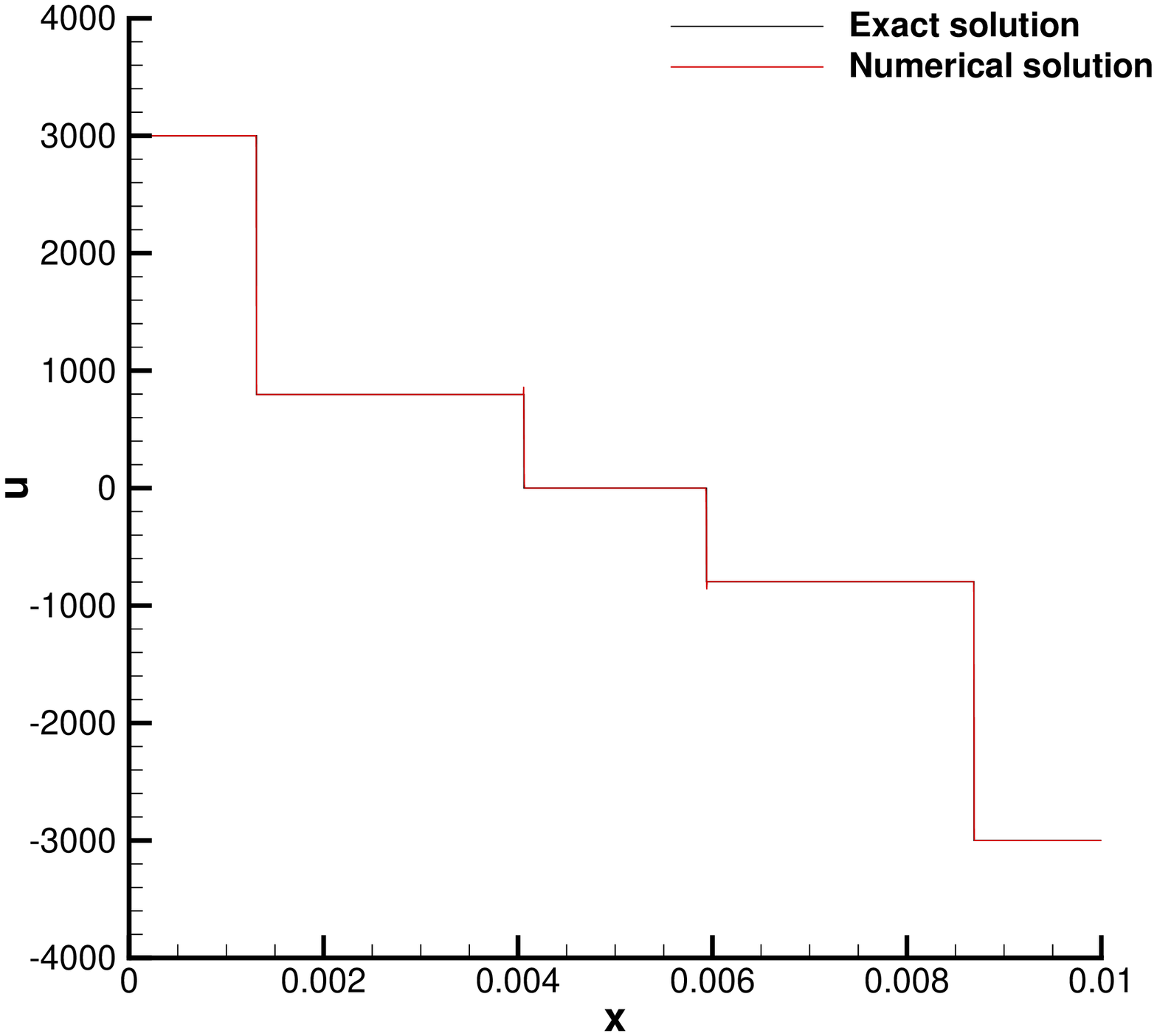}  &  
		\includegraphics[width=0.35\textwidth]{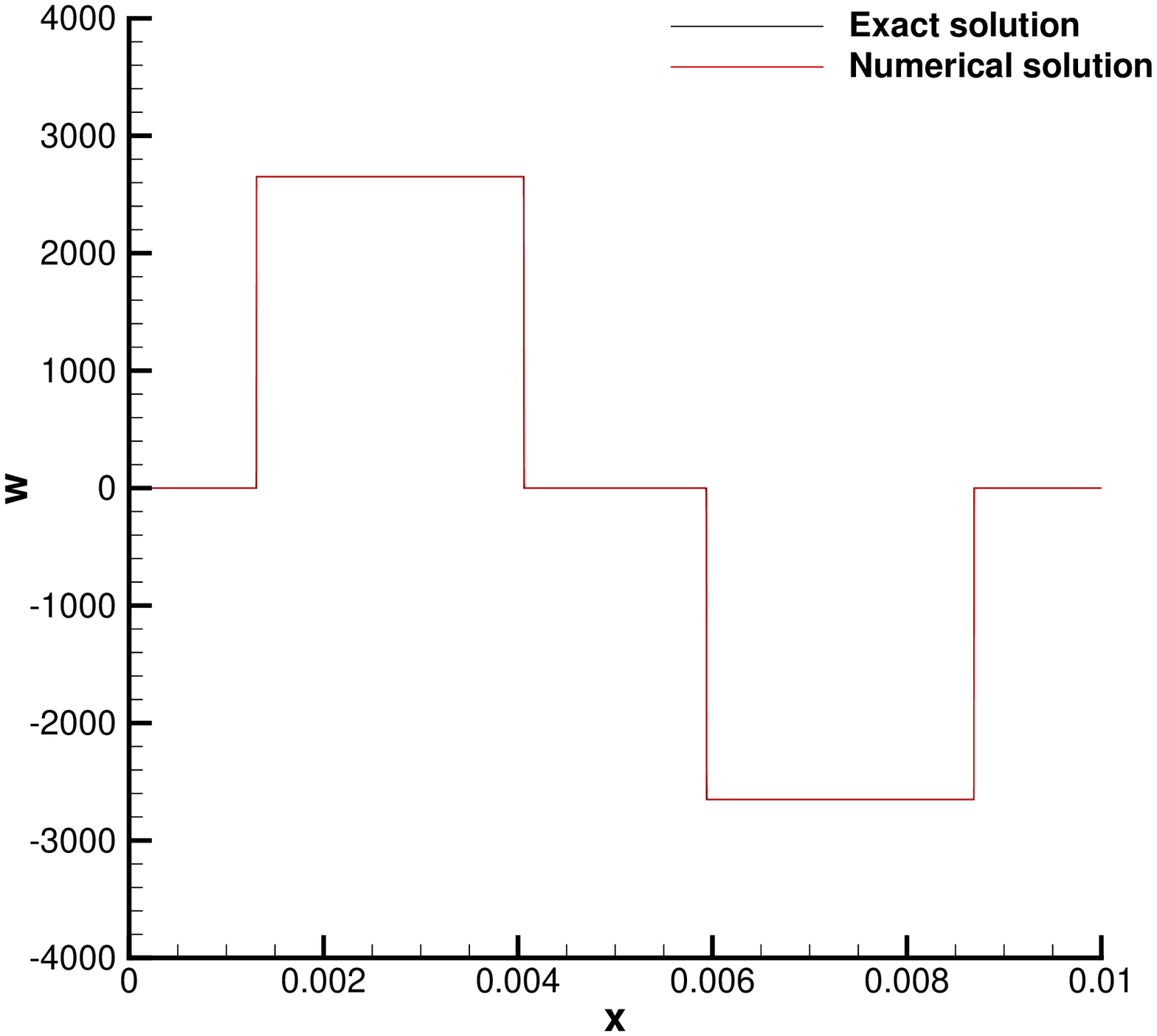}    	
	\end{tabular}
	}
    \caption{Exact solution (black) and numerical solution (red) of Riemann problem RP4.}    
    \label{fig:ex4_p1}
    \end{center}
\end{figure}

\begin{figure}
	\begin{center} 
	\includegraphics[width=0.8\textwidth]{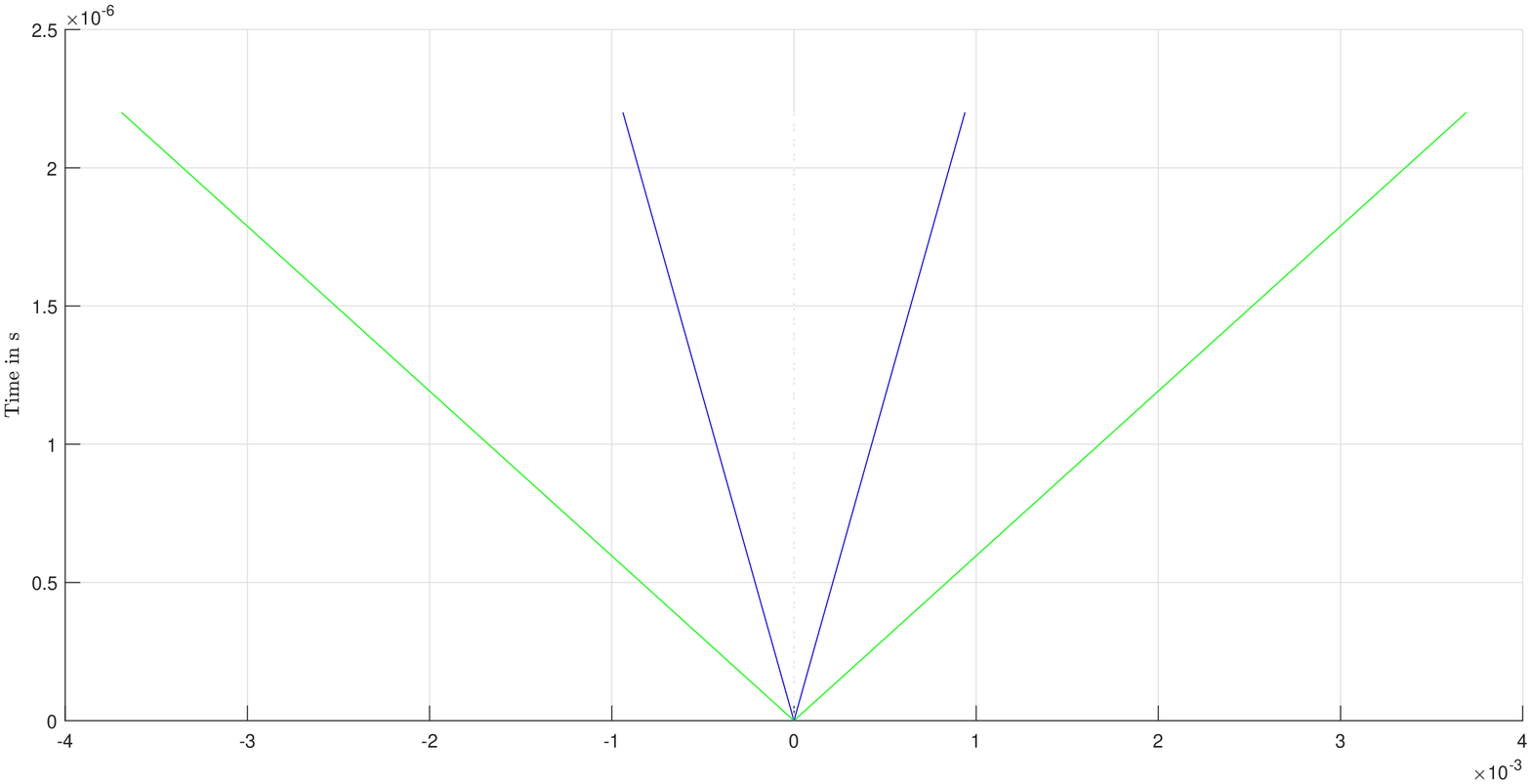} \\
	\includegraphics[width=0.8\textwidth]{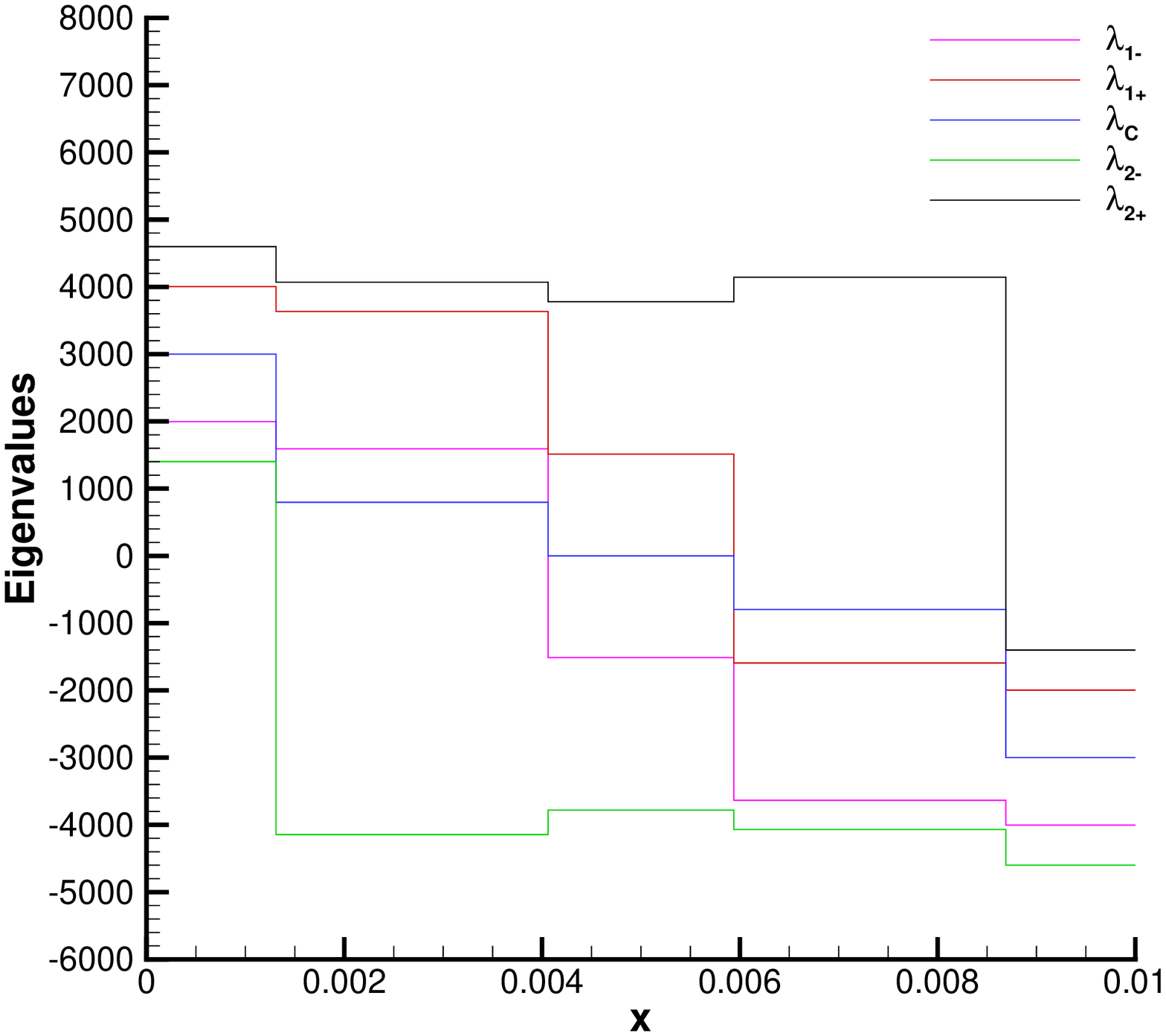}    	
    \caption{Wave structure of Riemann problem RP4 (top): phase one (blue), contact (red) and phase 2 (green). Eigenvalues of RP4 (bottom). }
    \label{fig:ex4_p2}
    \end{center} 
\end{figure}

The following example is taken from a paper of Romenski and Toro \cite{Romenski2004}.
The exact solution is obtained by using a contact centered inverse construction of the solution as described above.
The states are given in Table \ref{tab:states_prob_4} and we have $\alpha_1 = 0.9$ for all states.
\begin{table}[h!]
    \centering
    \renewcommand*{\arraystretch}{1.25}
    \begin{tabular}{c|llllll}
                 &  \multicolumn{1}{c}{$U_L$}              &  \multicolumn{1}{c}{$U^\ast_L$} &  \multicolumn{1}{c}{$U^{\ast\ast}_L$} &  \multicolumn{1}{c}{$U_R^{\ast\ast}$}
                 &  \multicolumn{1}{c}{$U^\ast_R$}         &  \multicolumn{1}{c}{$U_R$}            \\
        \hline
        $\rho_1$ &  131.01705 & \phantom{-}142.98406 &  1079 &  1079 &  \phantom{-}142.98406 &  \phantom{-}131.01705 \\
        $\rho_2$ &  1040.1358 & \phantom{-}2983.4101 &  2706 &  2706 &  \phantom{-}2983.4101 &  \phantom{-}1040.1358 \\
        $u_1$    &  3075.6226 & \phantom{-}2677.4348 &  0    &  0    &  -2677.4348           &  -3075.6226           \\
        $u_2$    &  3033.3793 & -38.030561           &  0    &  0    &  \phantom{-}38.030561 &  -3033.3793
    \end{tabular}
    \caption{Primitive states for Riemann problem RP4.}
    \label{tab:states_prob_4}
\end{table}
The following computation was performed using the following EOS
\begin{align*}
    p_i(\rho_i) &= A_i\left(\frac{\rho_i}{\rho_{ref}^{(i)}}\right)^{\gamma_i} + B_i,\, i\in\{1,2\},\\
    \text{with}\quad A_1 &= 10^5\,\si{Pa},\,\gamma_1 = 1.4,\, \rho_{ref}^{(1)} = 1\,\si{kg m^{-3}},\, B_1 = 0\,\si{Pa},\\
    \text{and}\quad A_2 &= 8.5\cdot 10^8\,\si{Pa},\,\gamma_2 = 2.8,\, \rho_{ref}^{(2)} = 10^3\,\si{kg m^{-3}},\, B_2 = 8.4999\cdot 10^8\,\si{Pa}.
\end{align*}
and the parameters
\begin{align*}
    \Delta x = 0.01\cdot 10^{-4}\,\si{m},\; C_{CFL} = 0.25,\, t_{end} = 0.22\cdot 10^{-5}\,\si{s} \quad\text{and}\quad x \in [0,0.01]\,\si{m}.
\end{align*}
In Figures \ref{fig:ex4_p1} and \ref{fig:ex4_p2} the numerical results together with the exact solution are shown.
For the numerical solution we used the Force Flux together with Godunov's method as exemplary shown in \cite{Toro2009}.
This is another good test as it shows four separated shock waves one in each phase. As noted before, in contrast to the Baer-Nunziato type systems, the shocks affect every phase.
In Figure \ref{fig:ex4_p2} we can see that the shock waves related to phase one are slower than the shocks related to phase two.
Moreover the order of the eigenvalues changes across every shock.
\FloatBarrier
\subsection{Comparison of the SHTC model with the Baer-Nunziato model}
In this last section we present a comparison of the numerical solutions obtained for the conservative SHTC model discussed in this paper, and the well-known Baer-Nunziato model. We show the solution of two Riemann problems, without and with stiff relaxation source terms. 

\begin{figure}[h!]
	\begin{center}
		\subfigure[Densities $\rho_1, \rho_2, \rho$ and volume fraction $\alpha_1$.]{
			\begin{tabular}{cc}        \includegraphics[width=0.33\textwidth]{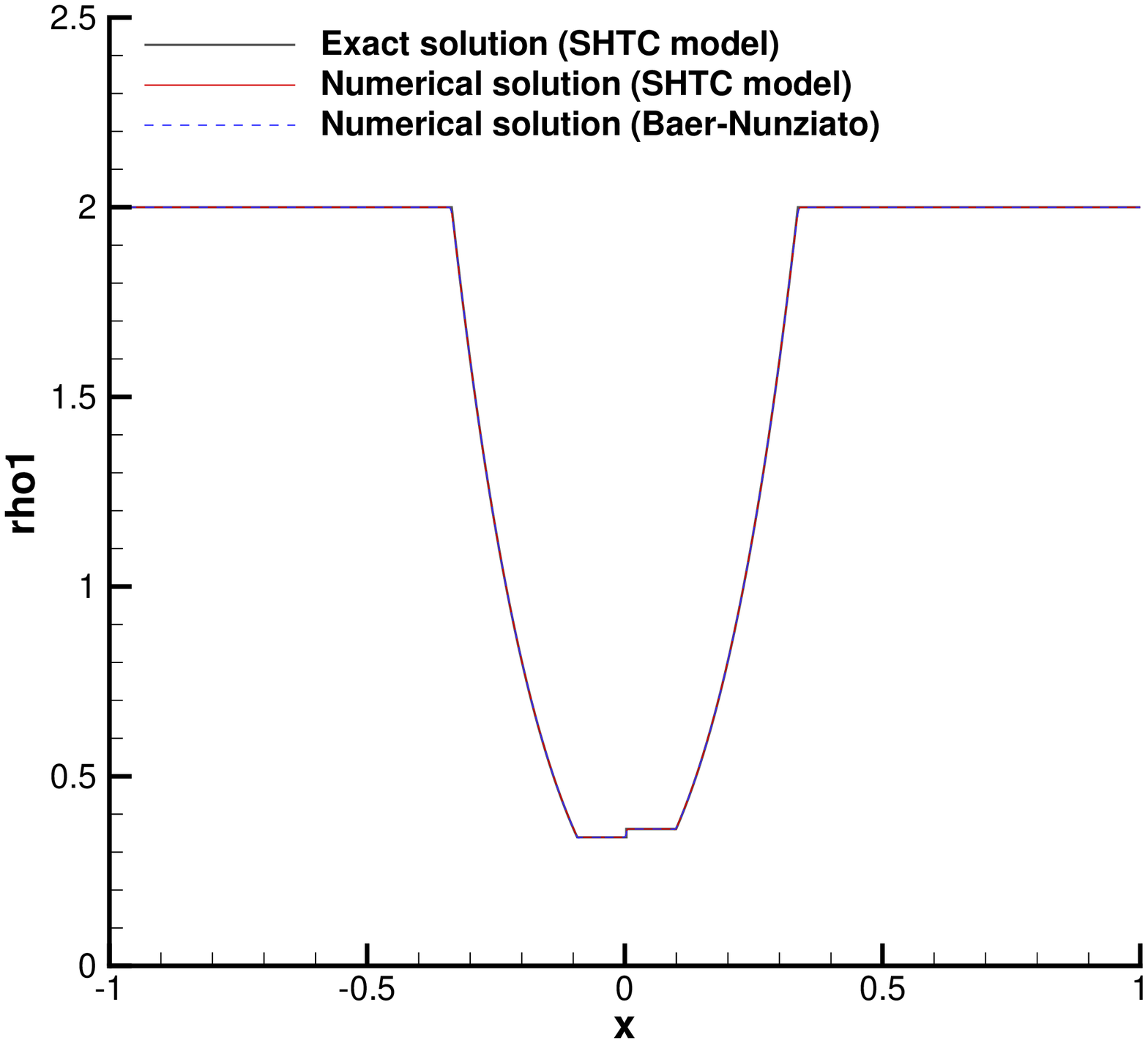}  &  
				\includegraphics[width=0.33\textwidth]{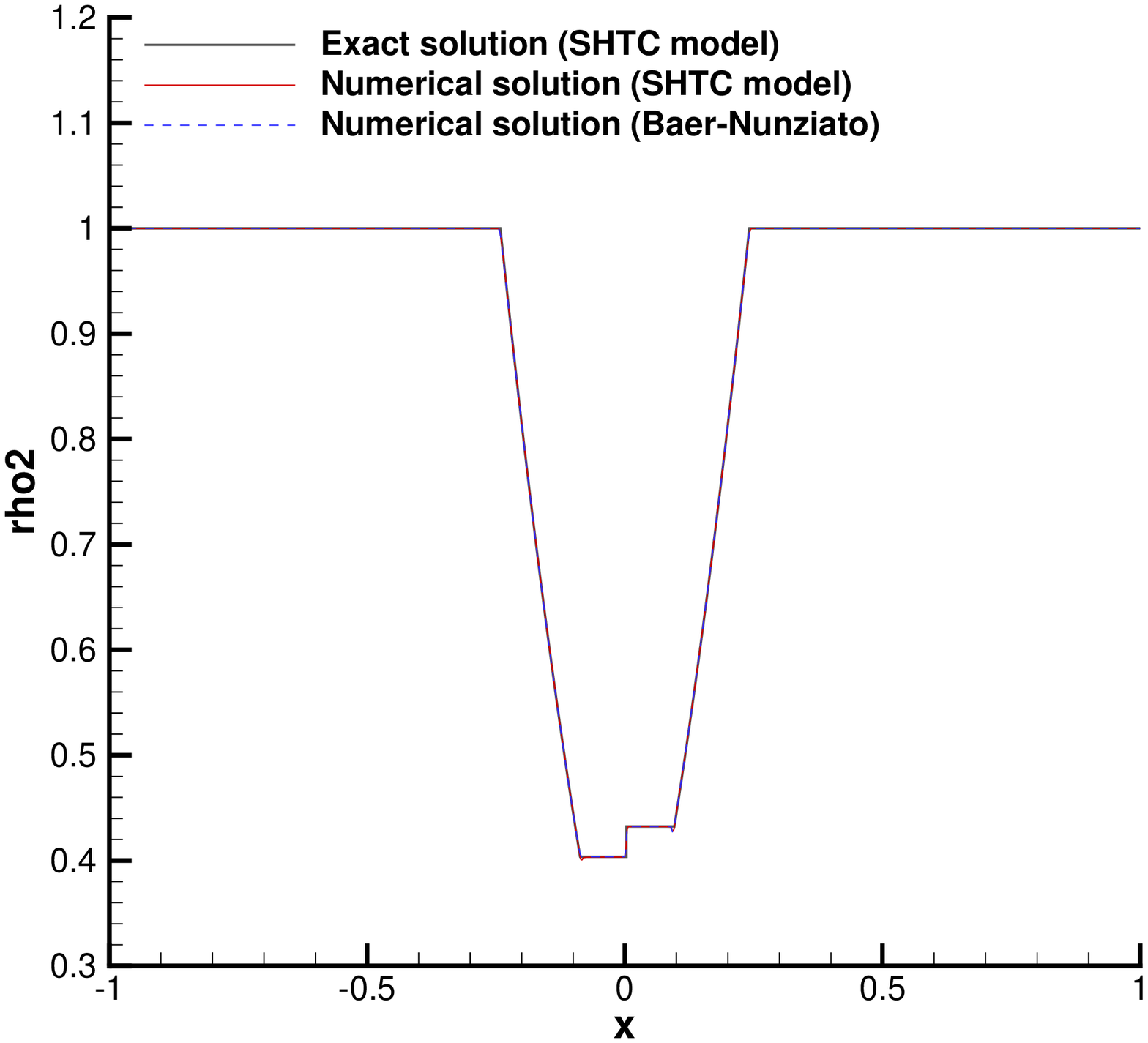}  \\ 	
				\includegraphics[width=0.33\textwidth]{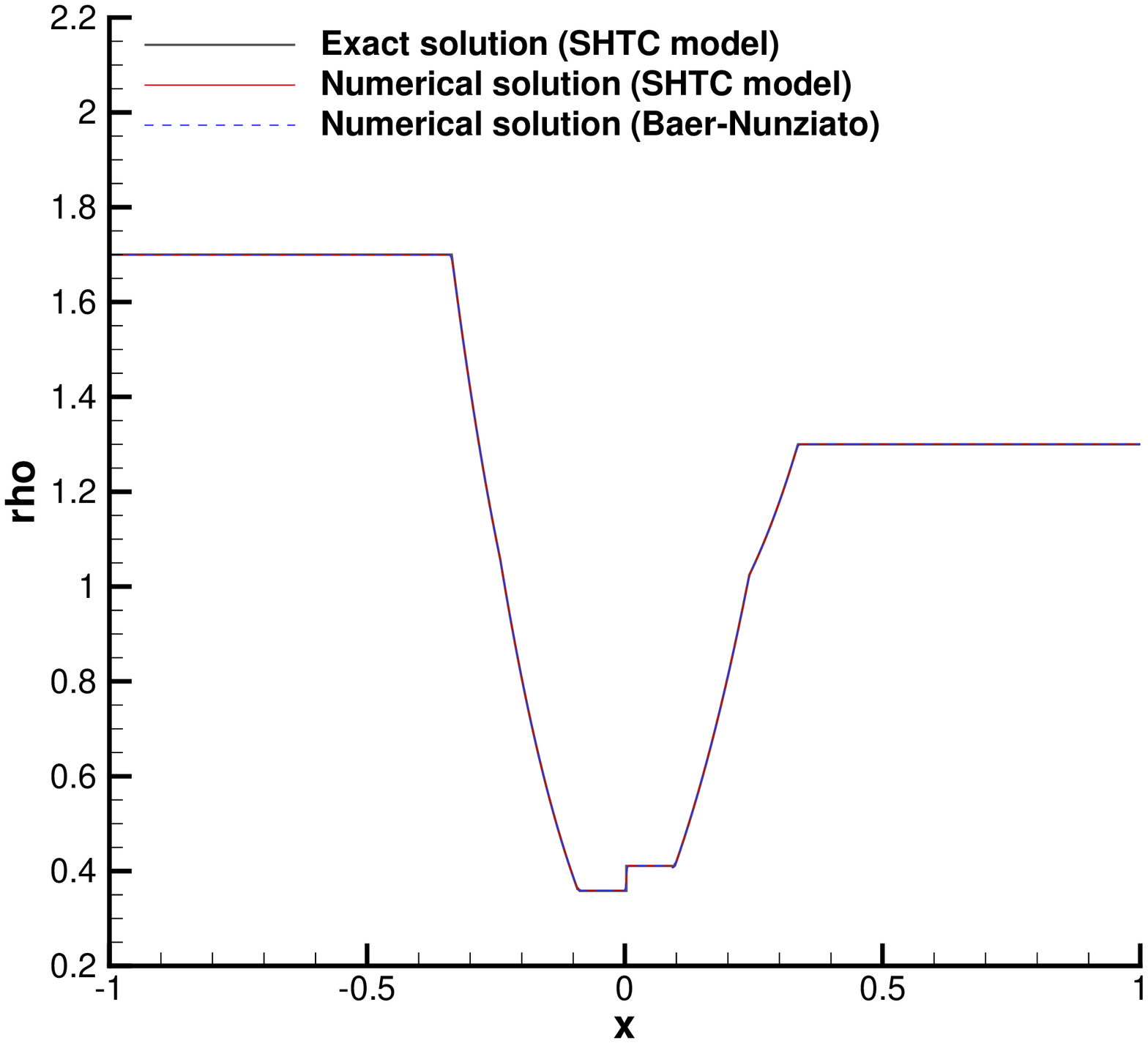}  &  
				\includegraphics[width=0.33\textwidth]{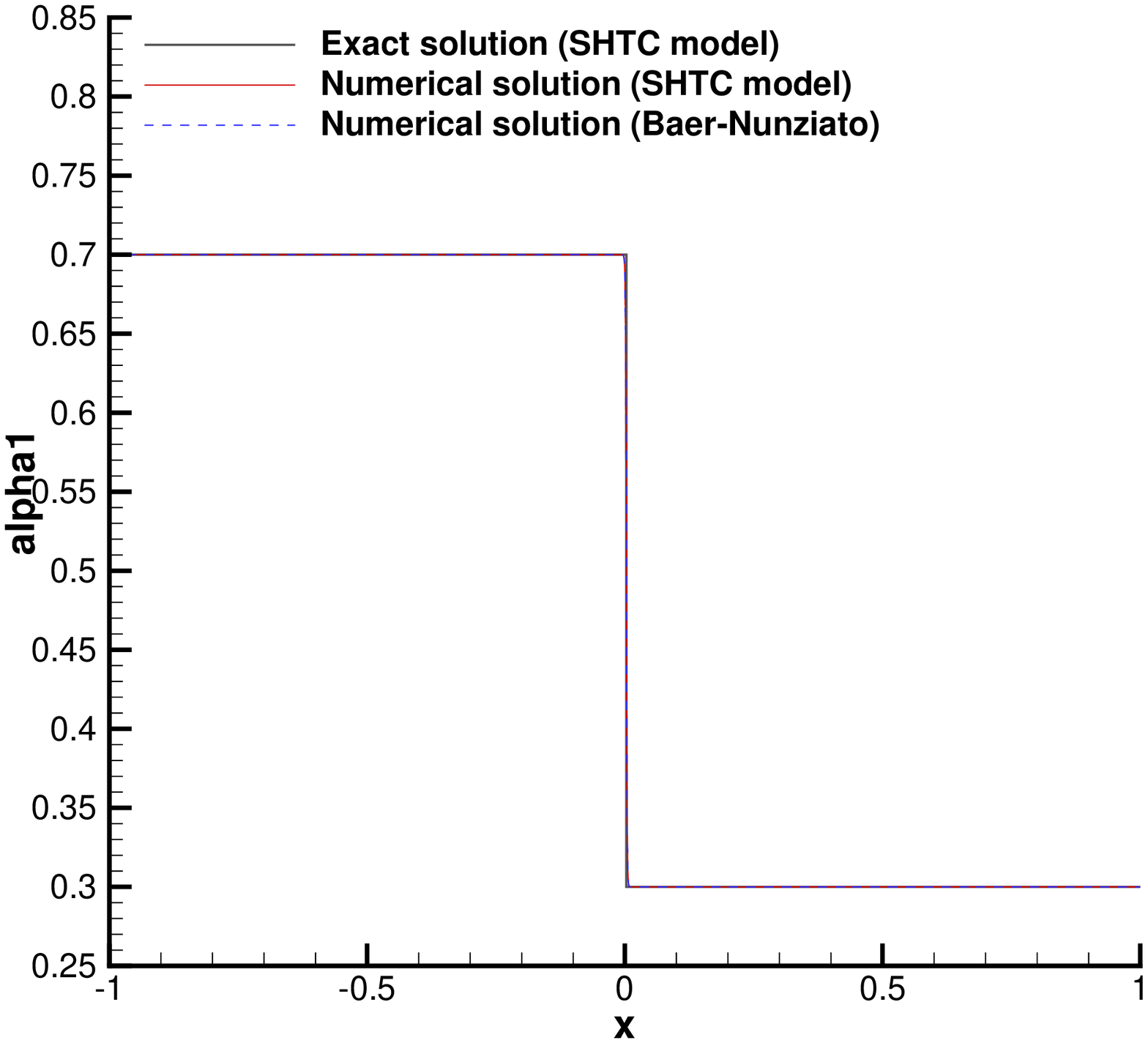}    	
			\end{tabular}
		}\\
		\subfigure[Velocities $u_1, u_2, u, w$.]{
			\begin{tabular}{cc}        \includegraphics[width=0.33\textwidth]{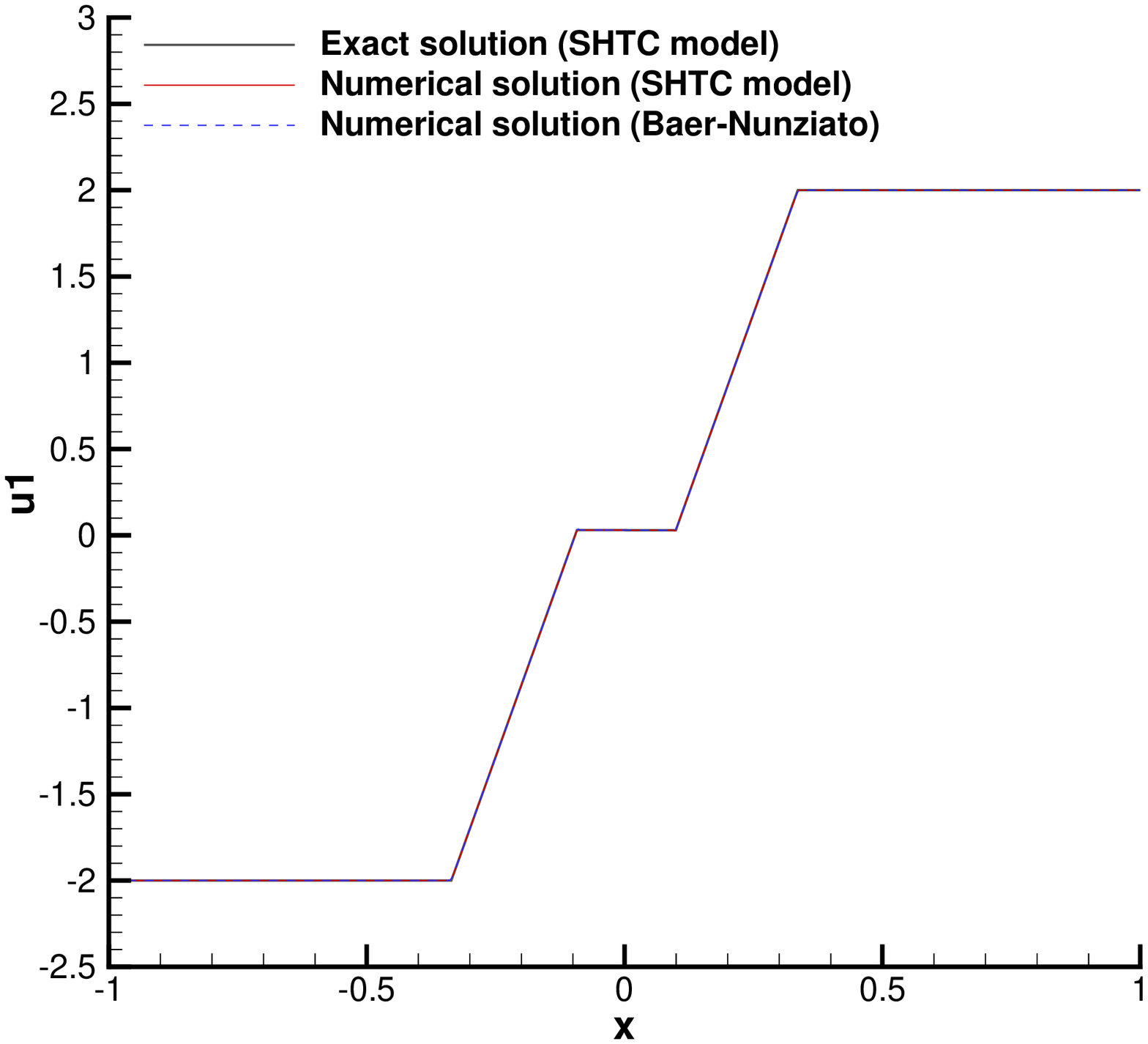}  &  
				\includegraphics[width=0.33\textwidth]{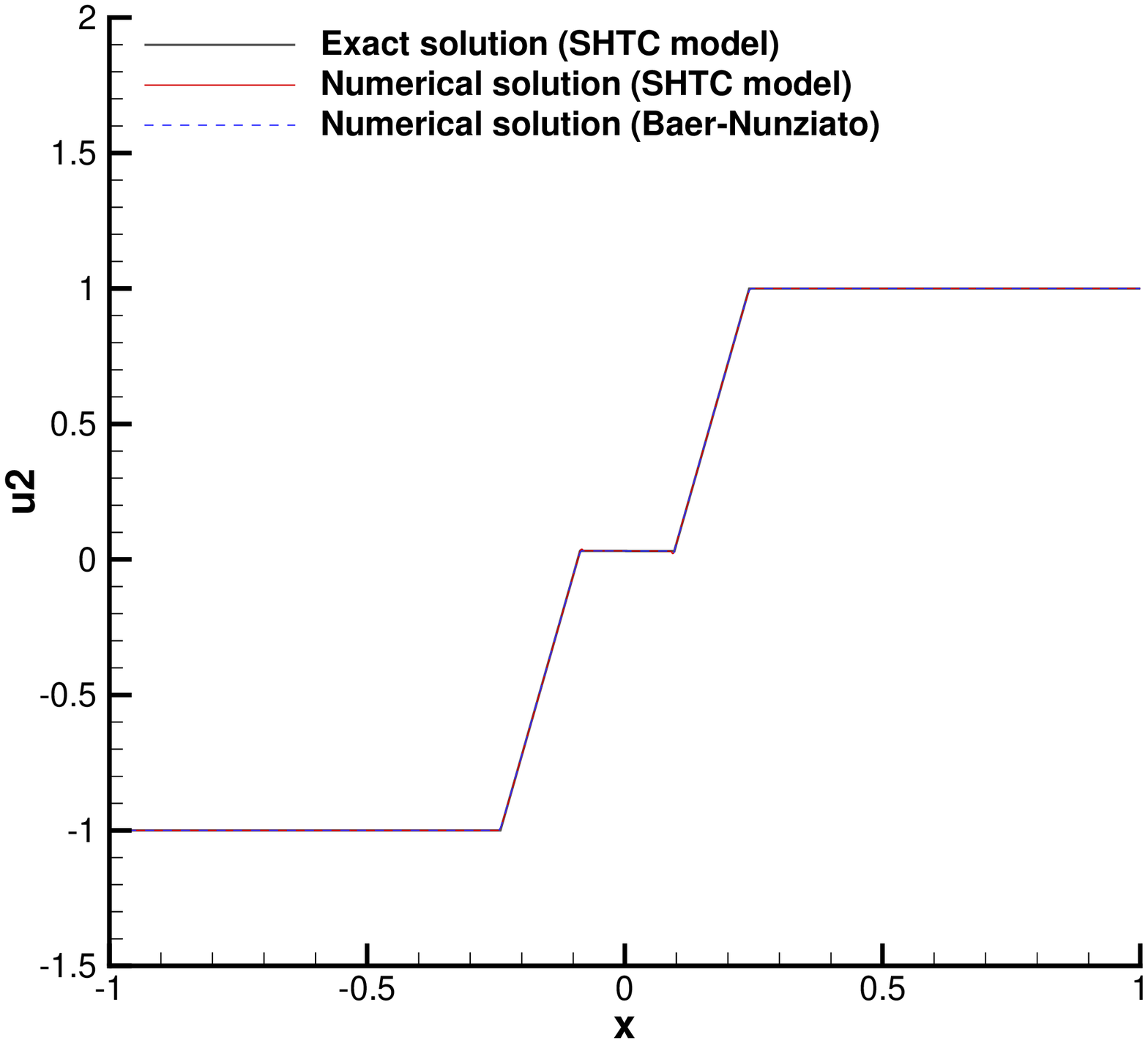}  \\ 	
				\includegraphics[width=0.33\textwidth]{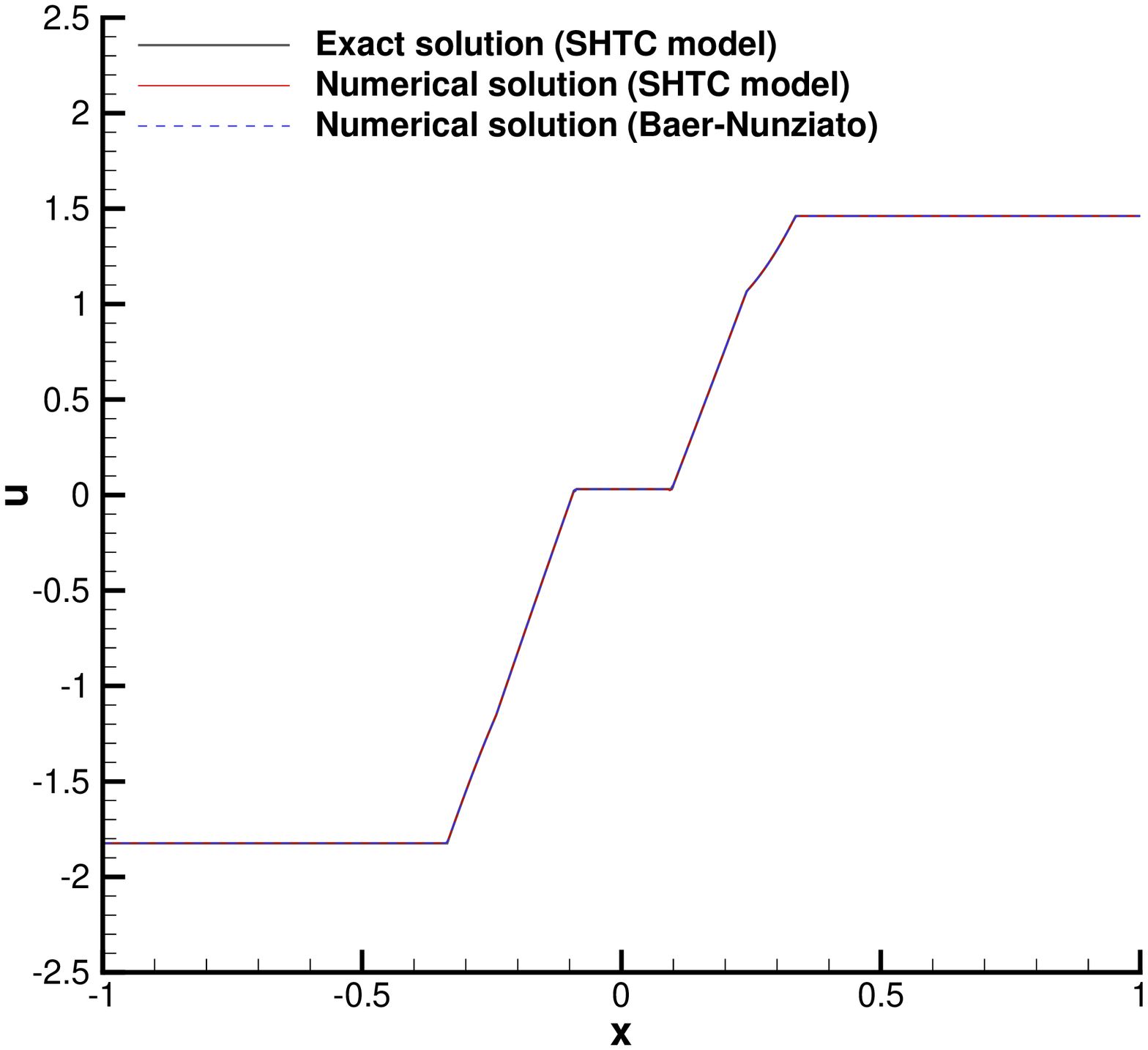}  &  
				\includegraphics[width=0.33\textwidth]{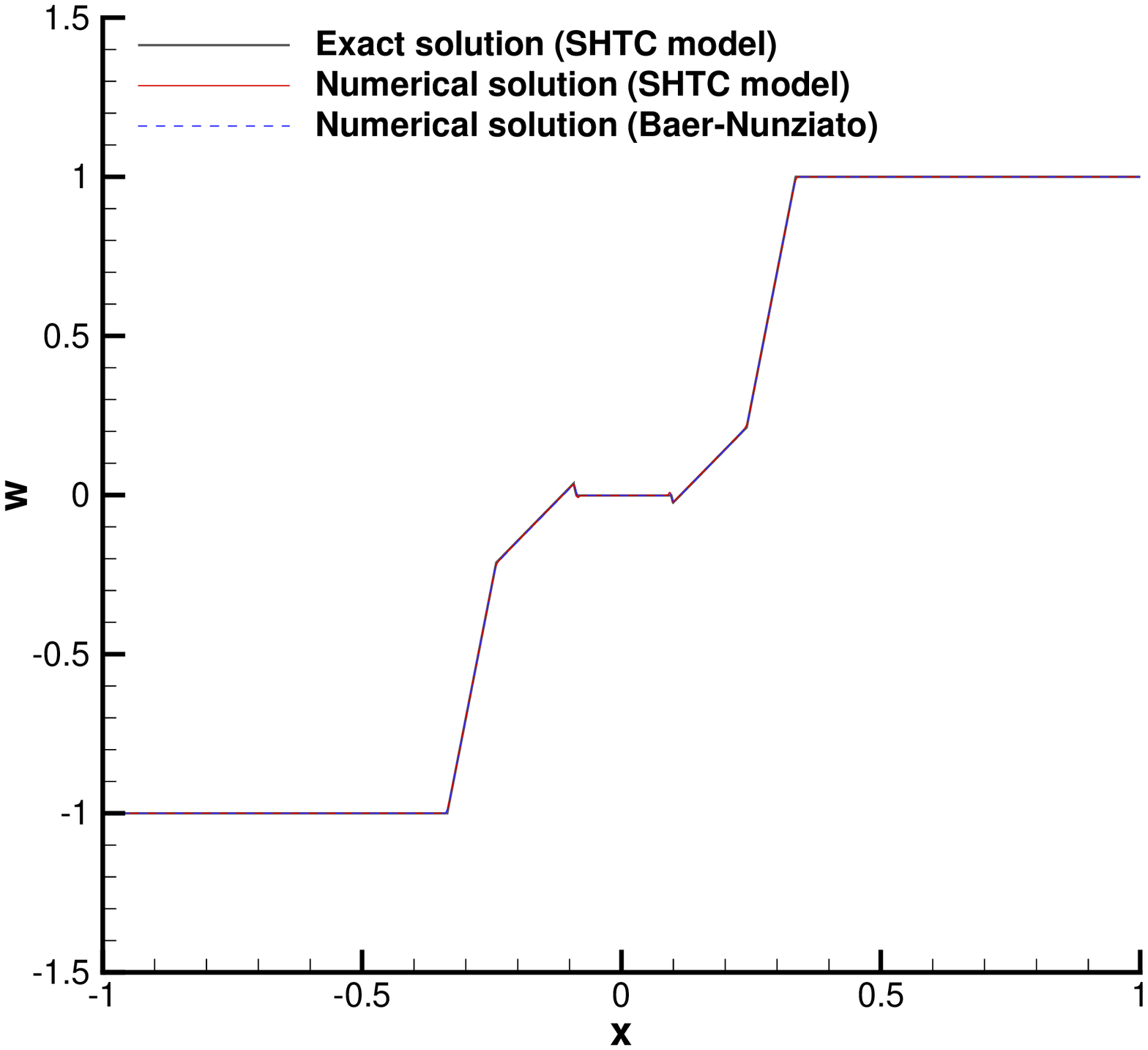}    	
			\end{tabular}
		}
		\caption{Exact solution of the SHTC system (black), numerical solution of the SHTC system (red) and numerical solution of the Baer-Nunziato model (blue) of Riemann problem RP5a without relaxation source terms.}
		\label{fig:rp5}
	\end{center}
\end{figure}

\begin{figure}[h!]
	\begin{center}
			\begin{tabular}{cc}        \includegraphics[width=0.4\textwidth]{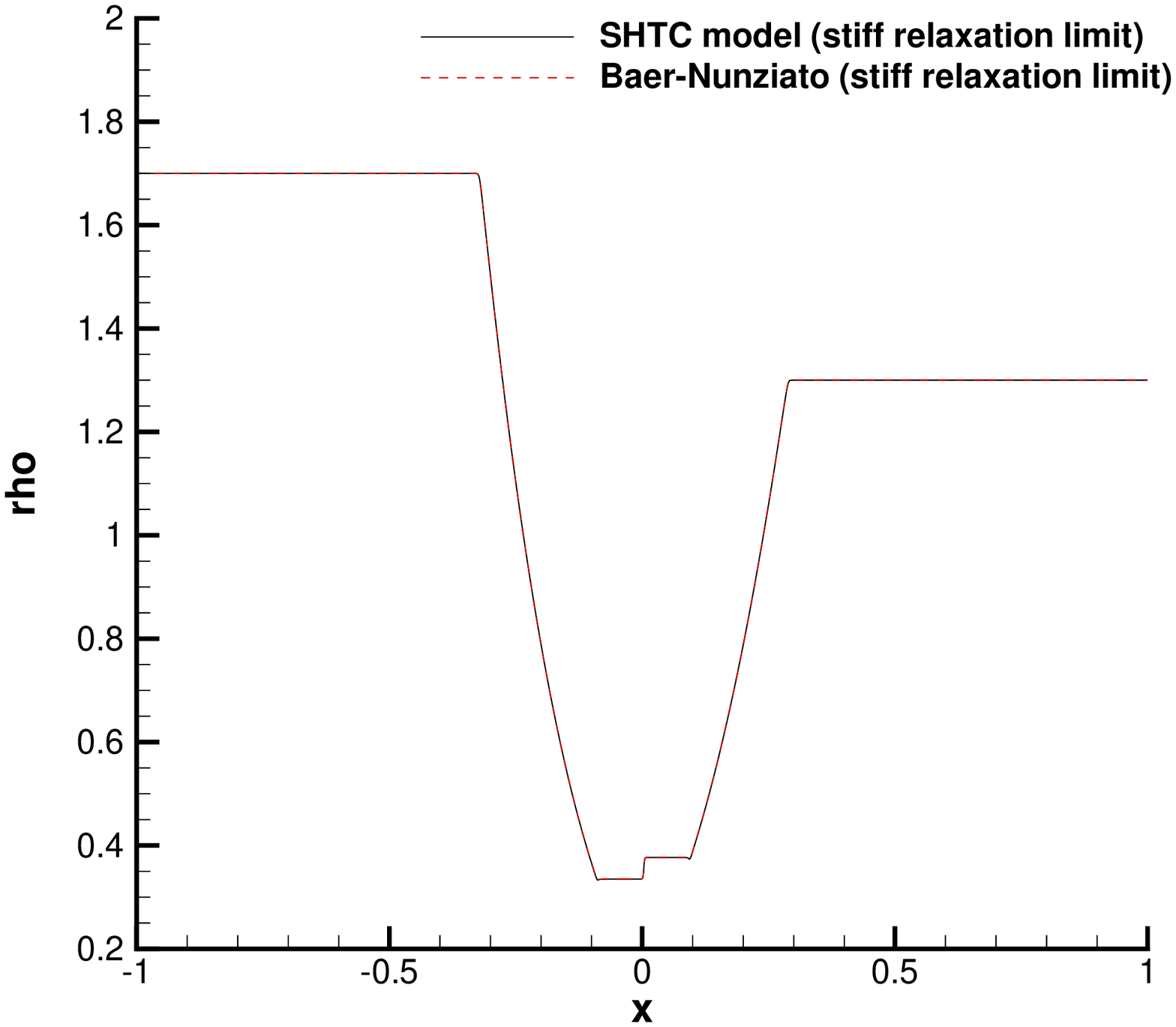}  & 
				\includegraphics[width=0.4\textwidth]{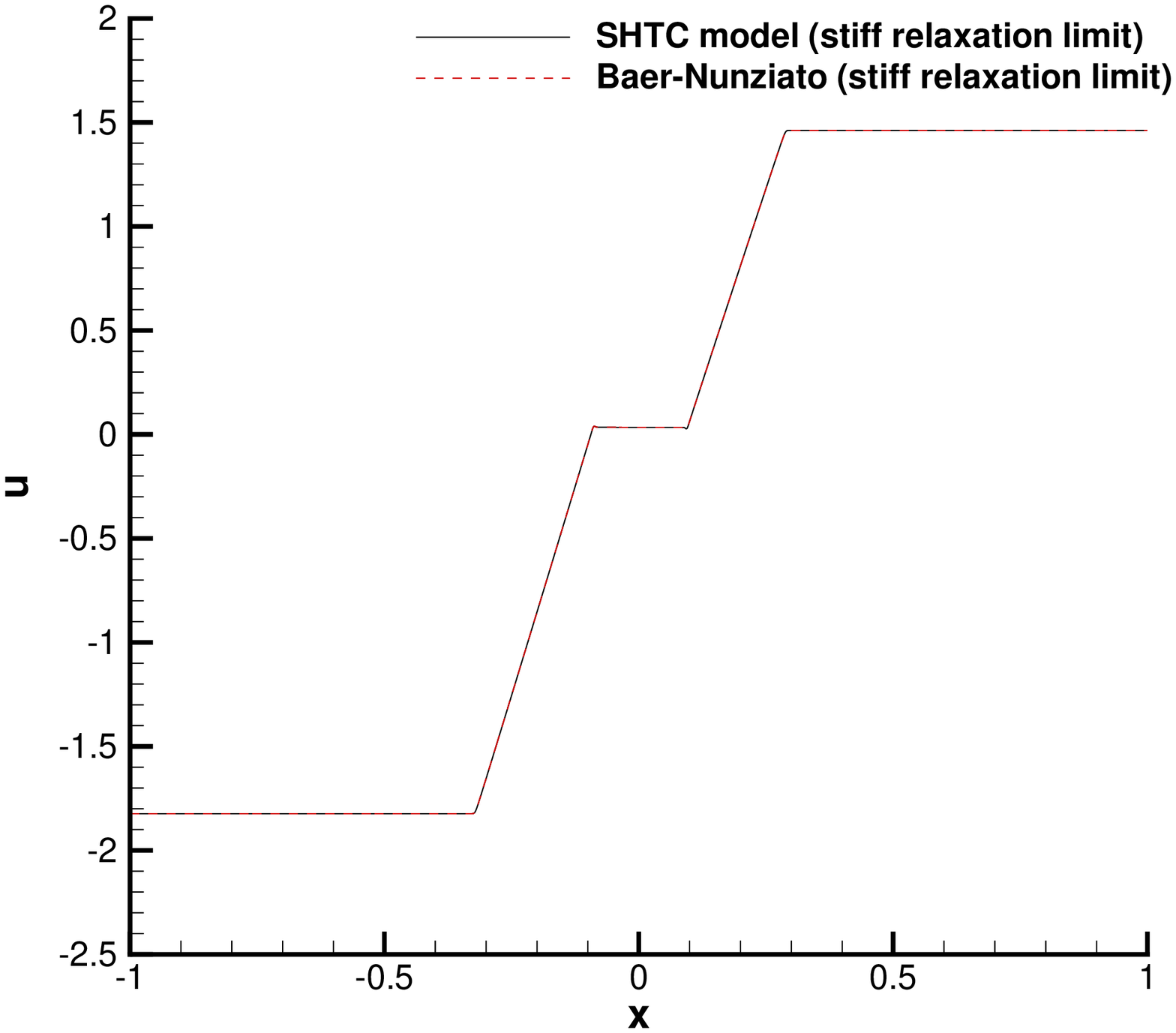}  \\ 	
				\includegraphics[width=0.4\textwidth]{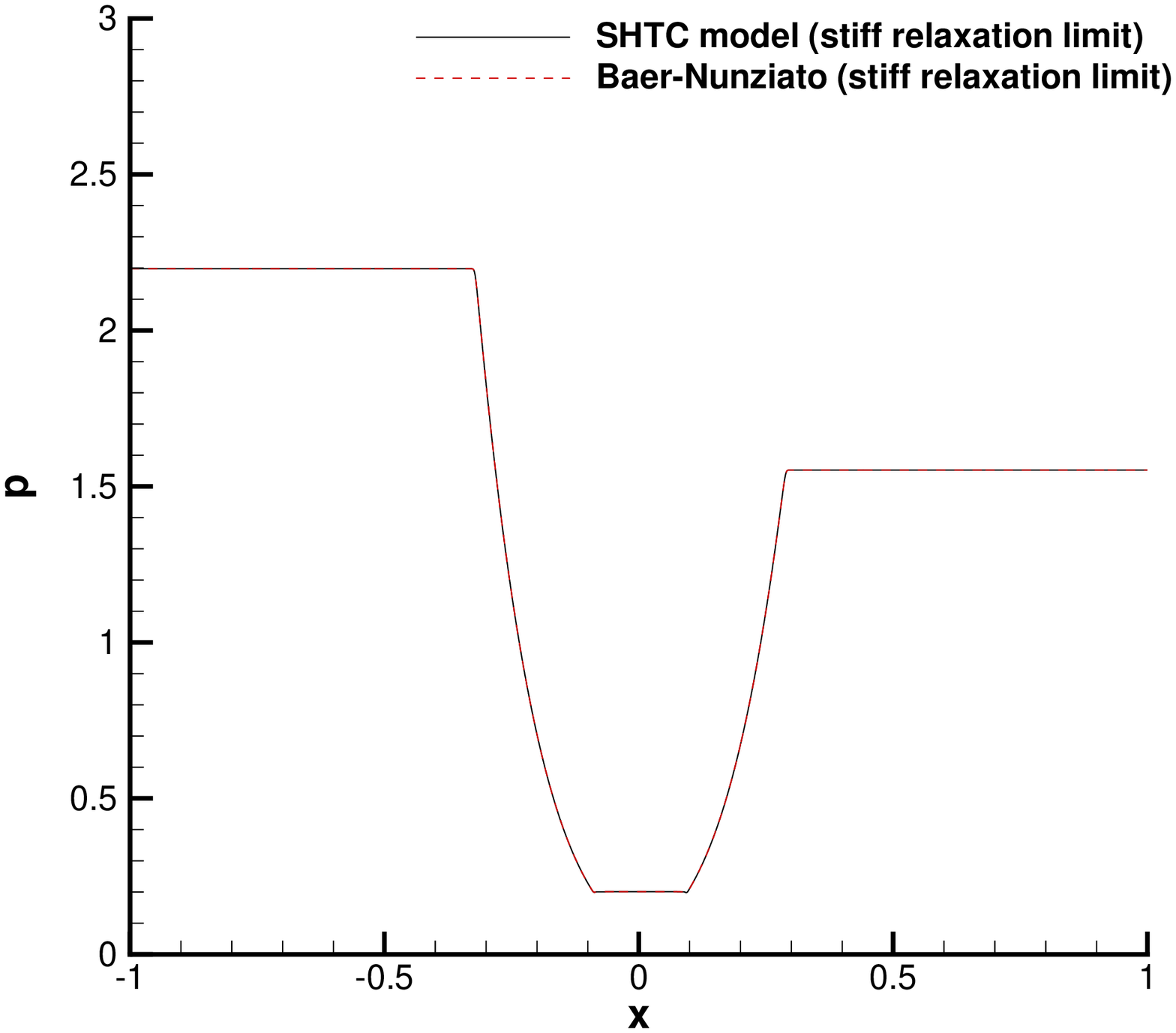} &
				\includegraphics[width=0.4\textwidth]{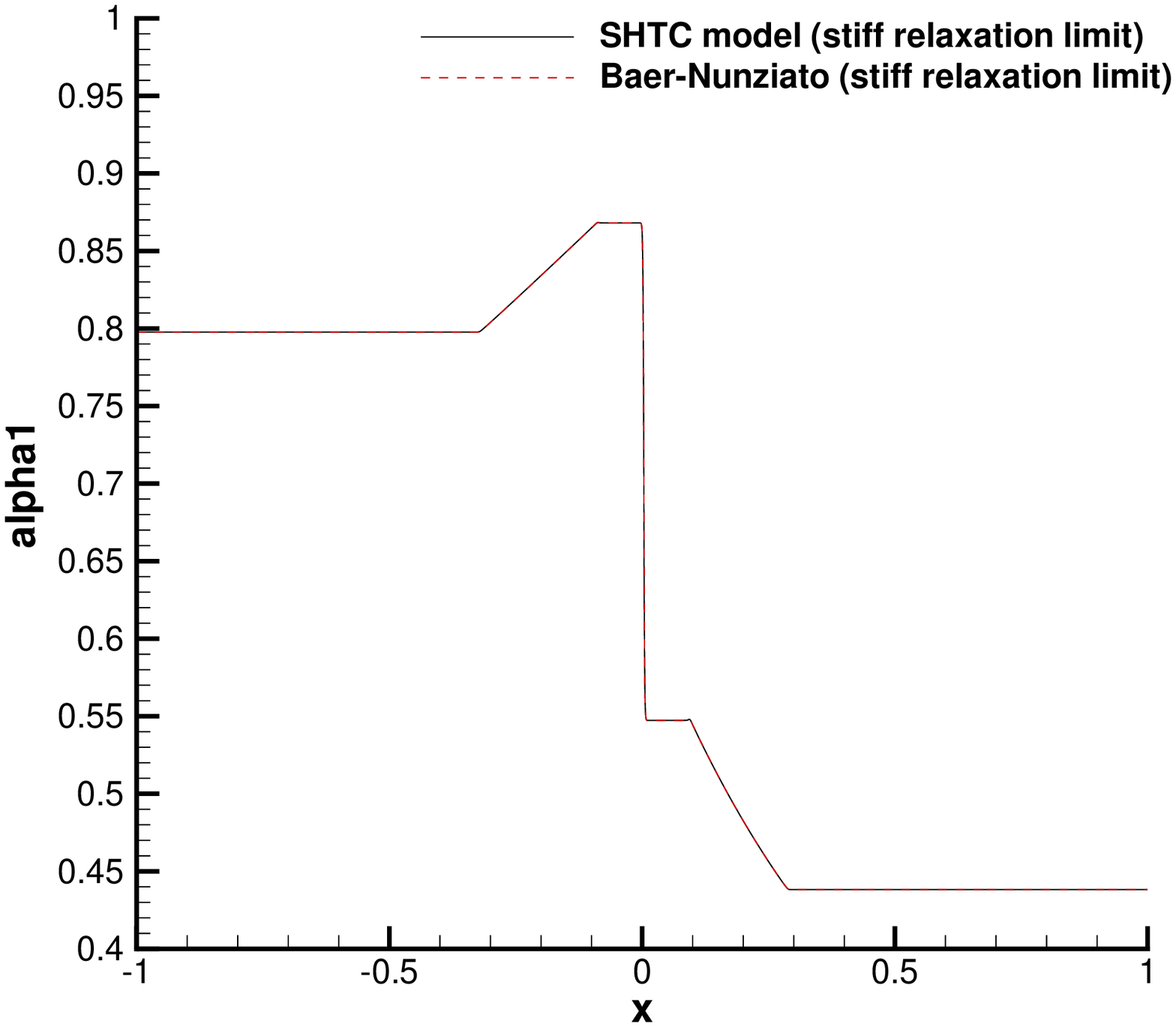}       
			\end{tabular}				
		\caption{Numerical solution of the SHTC system (black) and of the Baer-Nunziato model (red) of Riemann problem RP5b \textit{with} stiff relaxation source terms. From top left to bottom right: mixture density $\rho$, mixture velocity $u$, mixture pressure $p$ and volume fraction $\alpha_1$.}
		\label{fig:rp5b}
	\end{center}
\end{figure}

The test problem RP5 has the following initial data: $\rho_{1,L} = \rho_{1,R} = 2$,
$\rho_{2,L} = \rho_{2,R} = 1$, $u_{1,L}=-2$, $u_{1,R}=+2$, $u_{2,L}=-1$, $u_{2,R}=+1$, $\alpha_L = 0.7$, $\alpha_R=0.3$. The computational domain is the interval $[-1,1]$, which is discretized with 10000 uniform grid cells and the final time of the simulation is $t=0.1$. In Figure \ref{fig:rp5} we compare the exact solution of the Riemann problem for the SHTC system with the numerical solutions obtained for the SHTC system and the Baer-Nunziato model without any relaxation source terms. Since the wave structure consists of two rarefaction waves in both phases, the numerical results of the SHTC model and the numerical results obtained for the Baer-Nunziato system agree perfectly well with each other and with the exact solution of the Riemann problem, as expected. In Figure \ref{fig:rp5b} we show the results obtained for the same initial data, but with stiff pressure and velocity relaxation, choosing the relaxation parameters as $\theta_1 = 10^{-3}$ and $\theta_2=10^{-8}$, which corresponds to the Kapila limit of both systems. We find that the numerical solutions obtained for the conservative SHTC model and the Baer-Nunziato model agree perfectly well with each other.   

\begin{figure}[h!]
	\begin{center}
		\subfigure[Densities $\rho_1, \rho_2, \rho$ and volume fraction $\alpha_1$.]{
			\begin{tabular}{cc}        \includegraphics[width=0.33\textwidth]{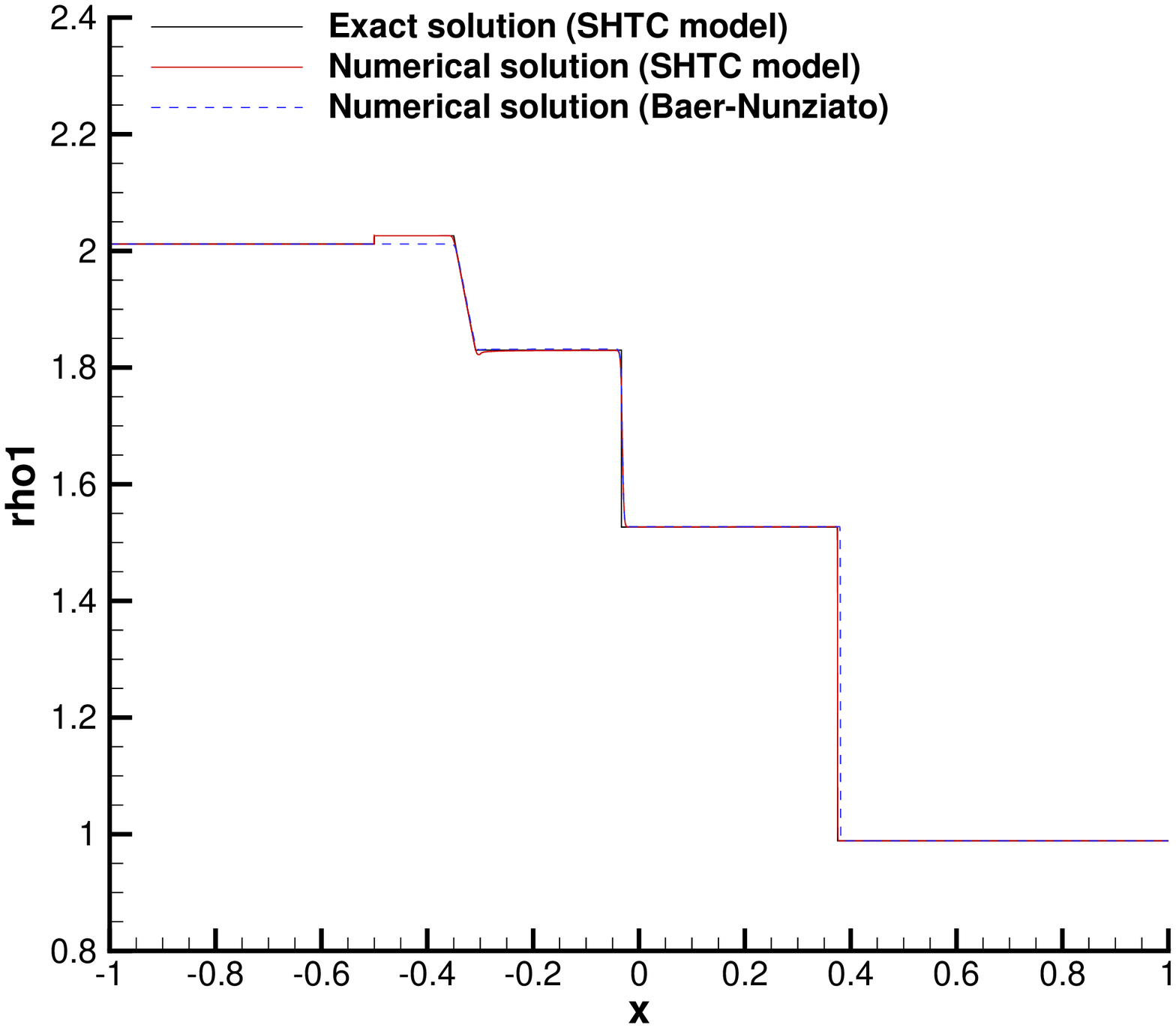}  &  
				\includegraphics[width=0.33\textwidth]{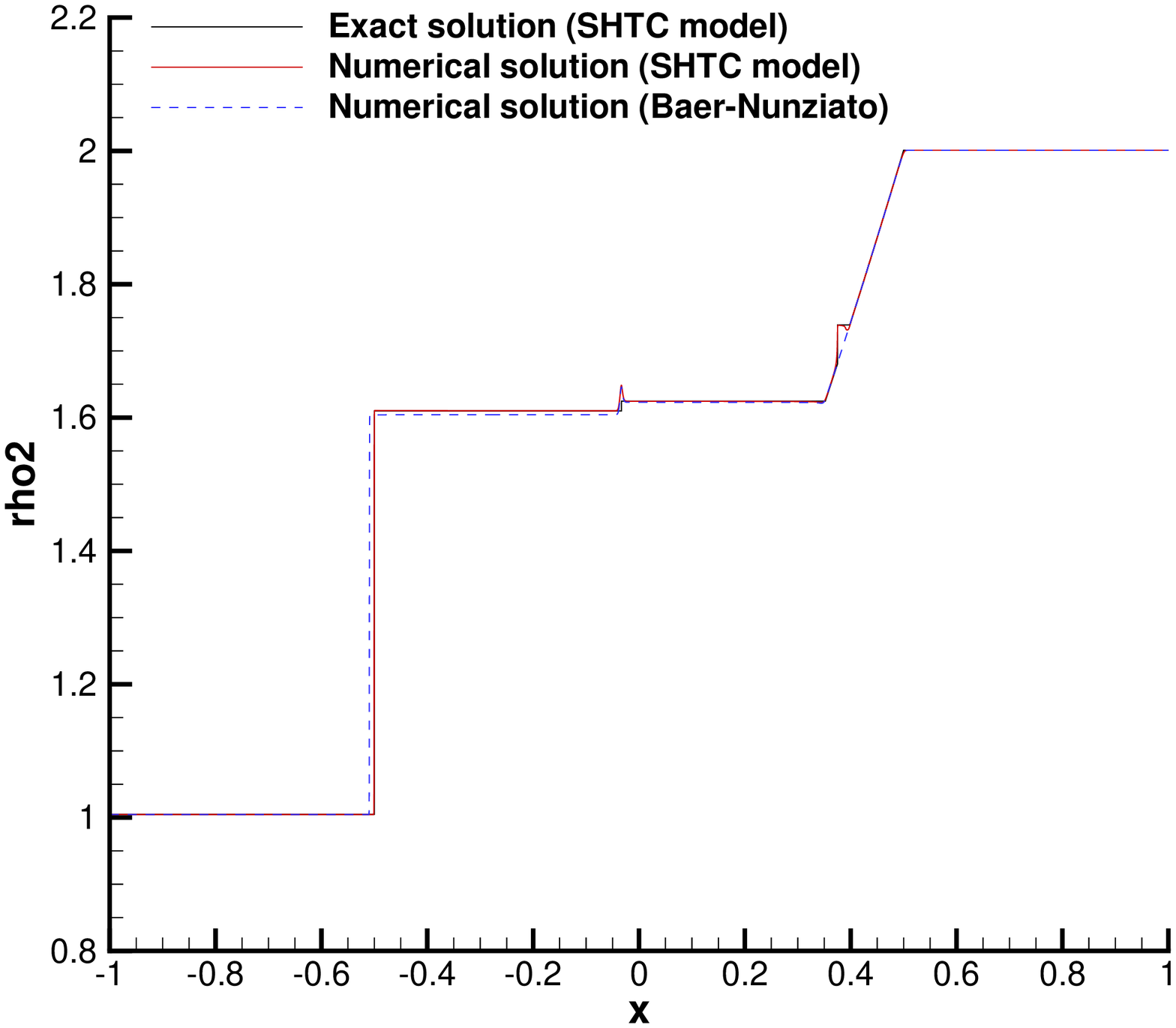}  \\ 	
				\includegraphics[width=0.33\textwidth]{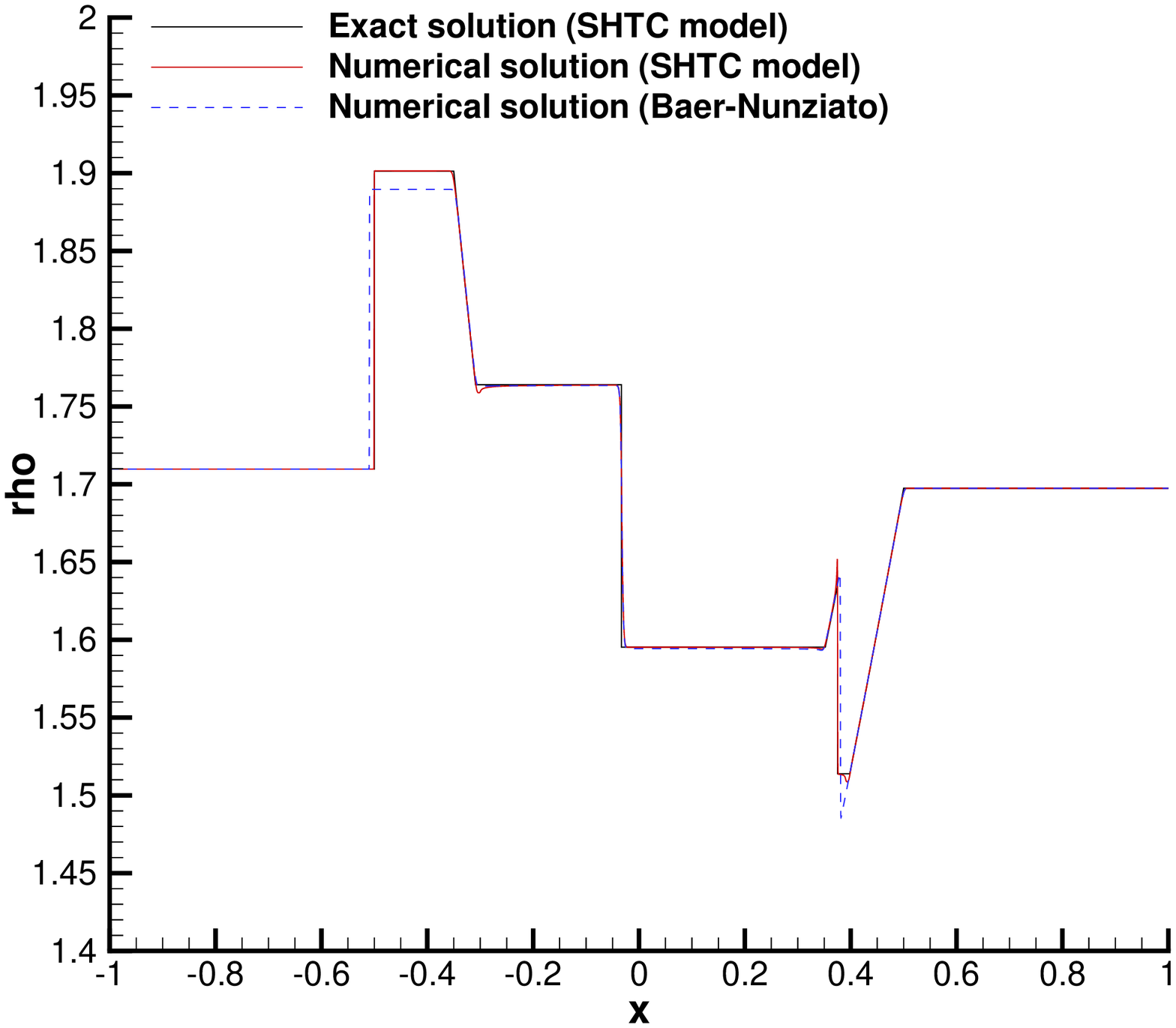}  &  
				\includegraphics[width=0.33\textwidth]{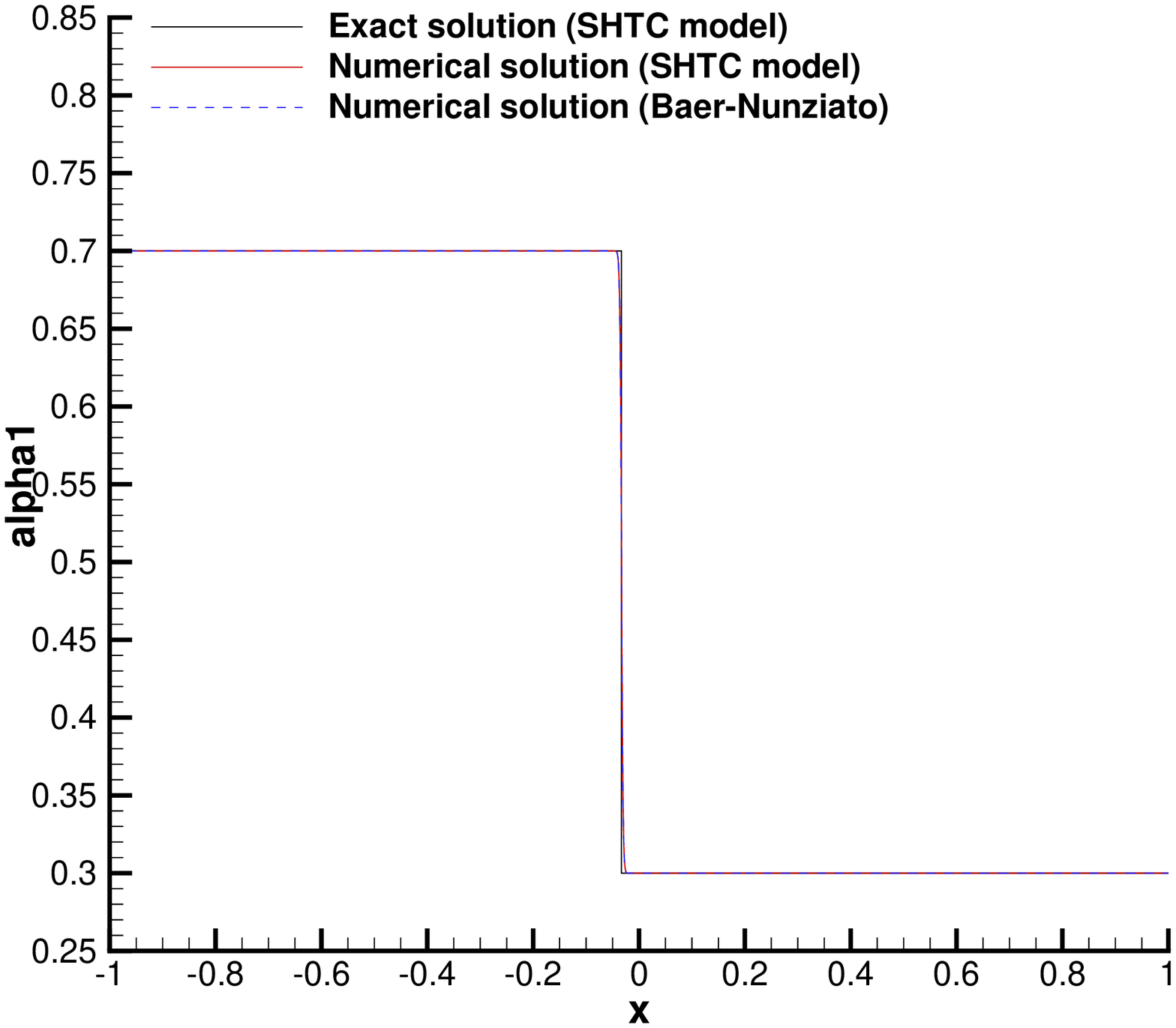}    	
			\end{tabular}
		}\\
		\subfigure[Velocities $u_1, u_2, u, w$.]{
			\begin{tabular}{cc}        \includegraphics[width=0.33\textwidth]{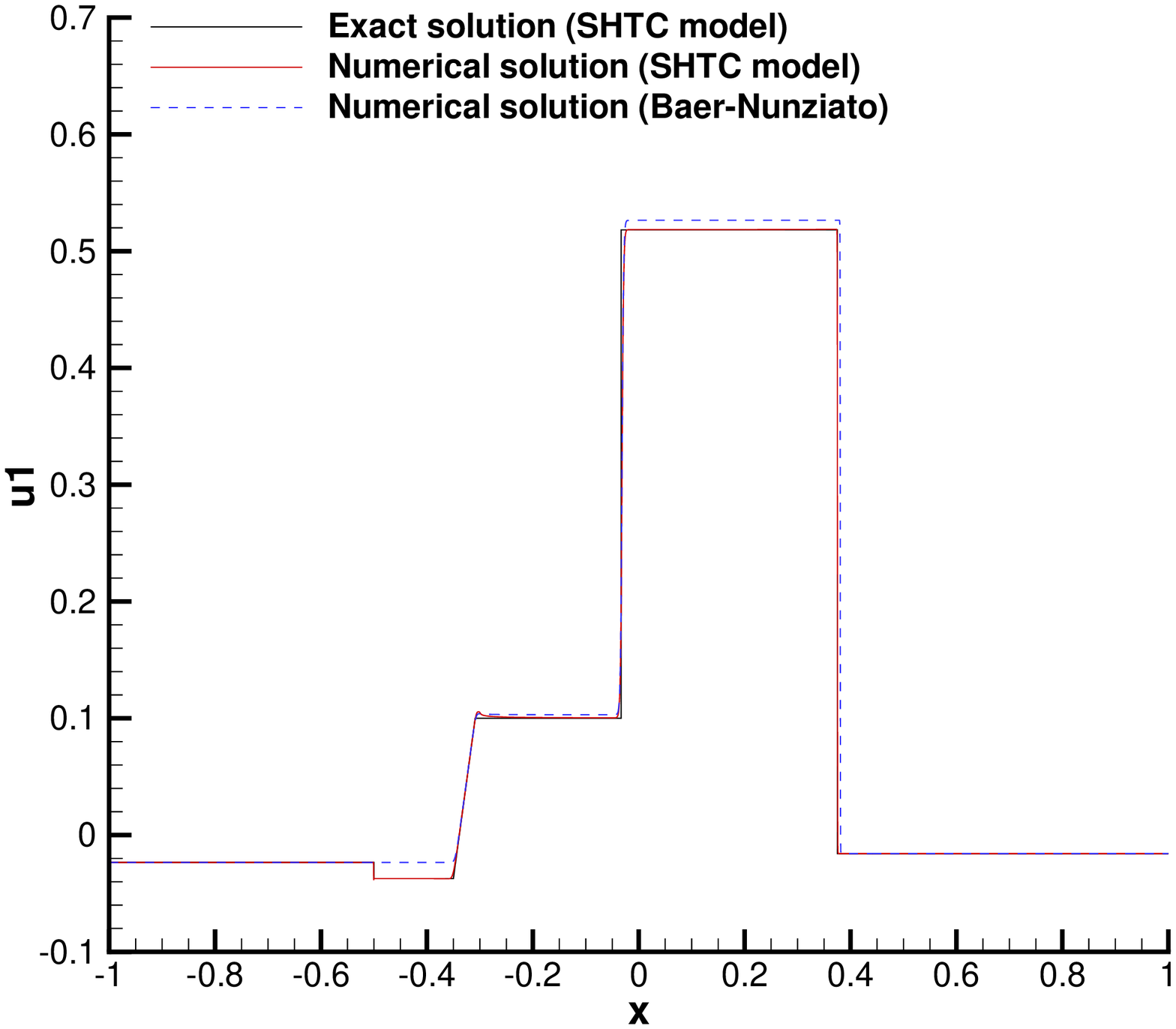}  &  
				\includegraphics[width=0.33\textwidth]{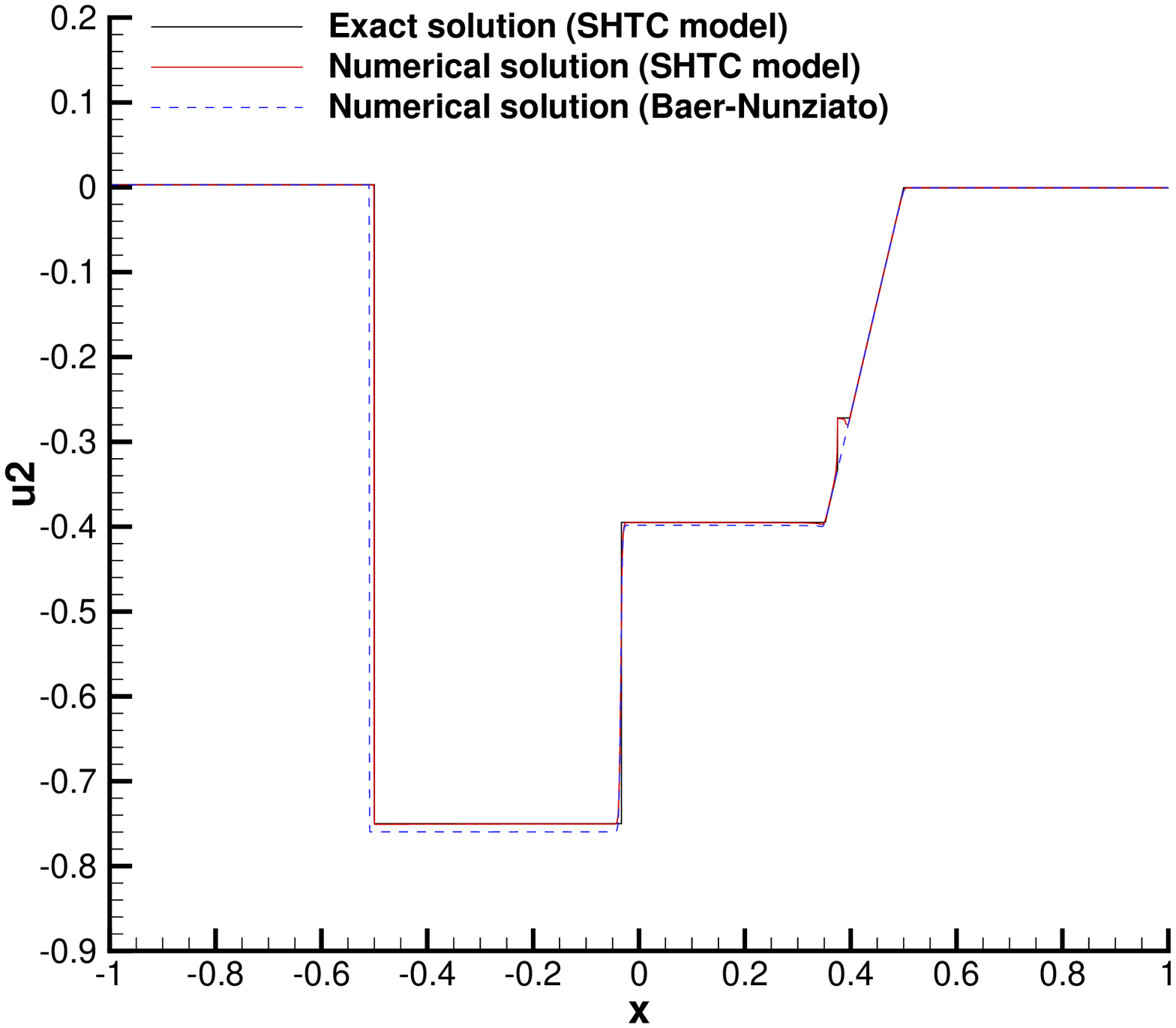}  \\ 	
				\includegraphics[width=0.33\textwidth]{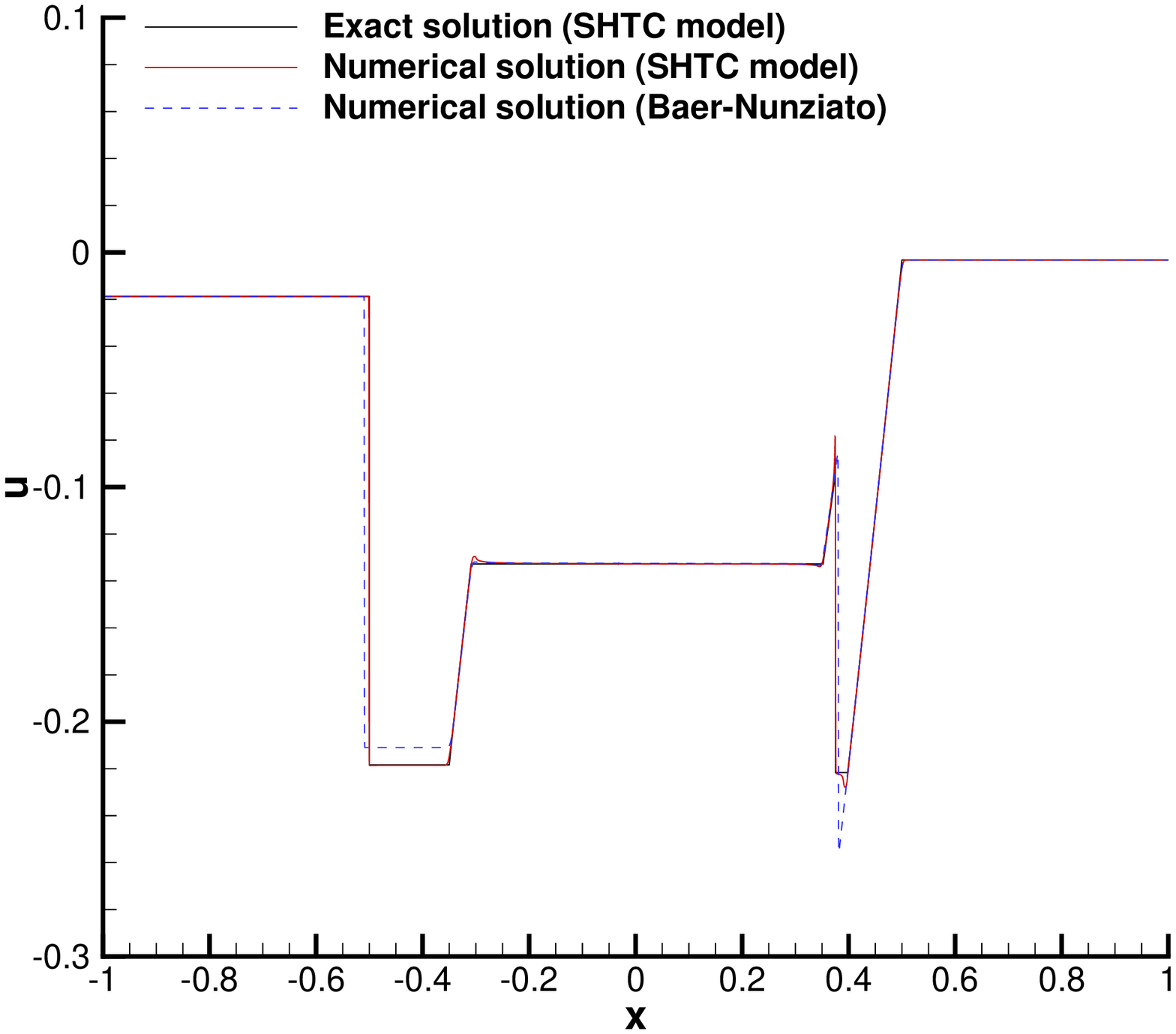}  &  
				\includegraphics[width=0.33\textwidth]{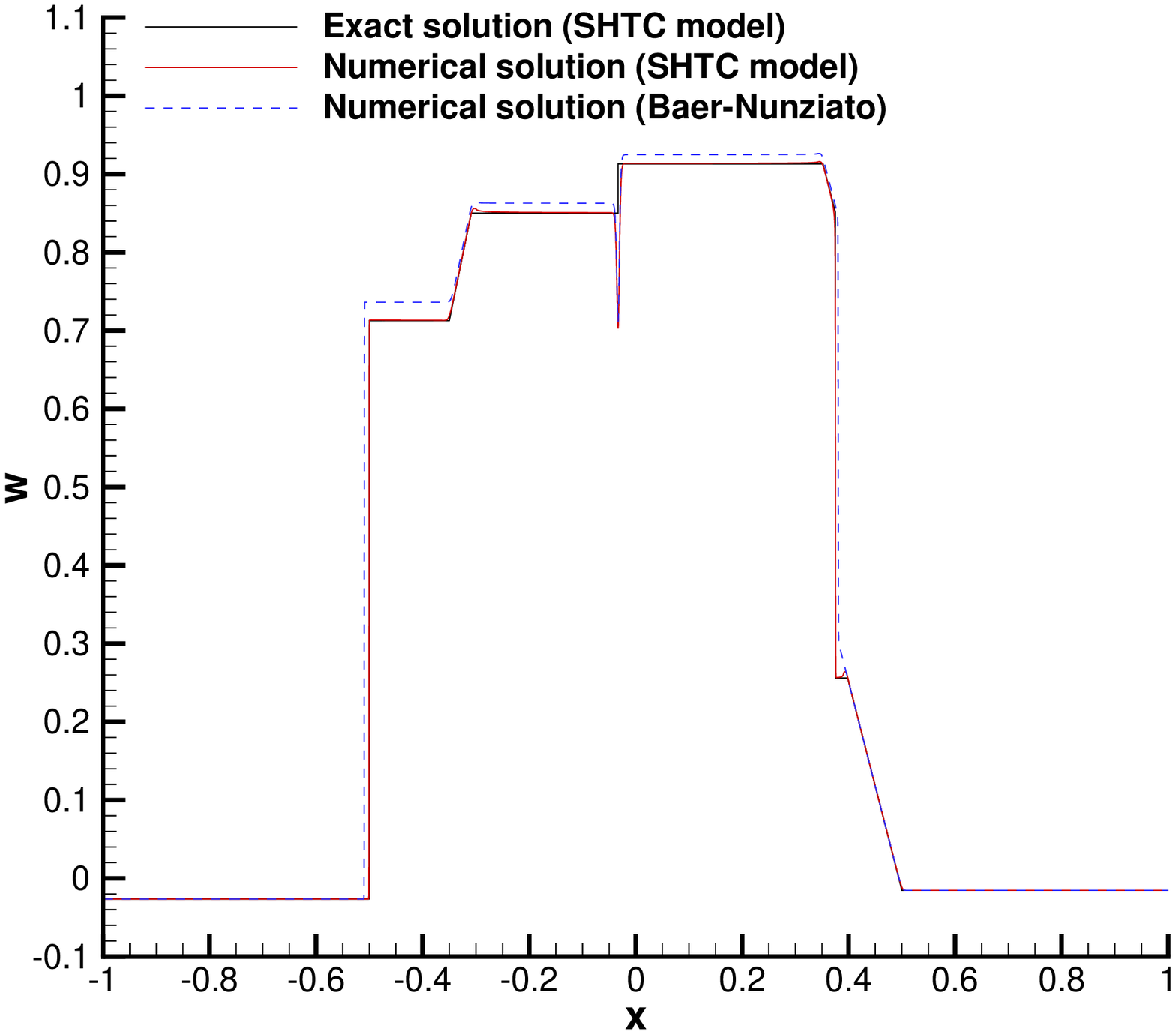}    	
			\end{tabular}
		}
		\caption{Exact solution of the SHTC system (black), numerical solution of the SHTC system (red) and numerical solution of the Baer-Nunziato model (blue) of Riemann problem RP6a without relaxation source terms.}
		\label{fig:rp6}
	\end{center}
\end{figure}

\begin{figure}[h!]
	\begin{center}
		\begin{tabular}{cc}        \includegraphics[width=0.4\textwidth]{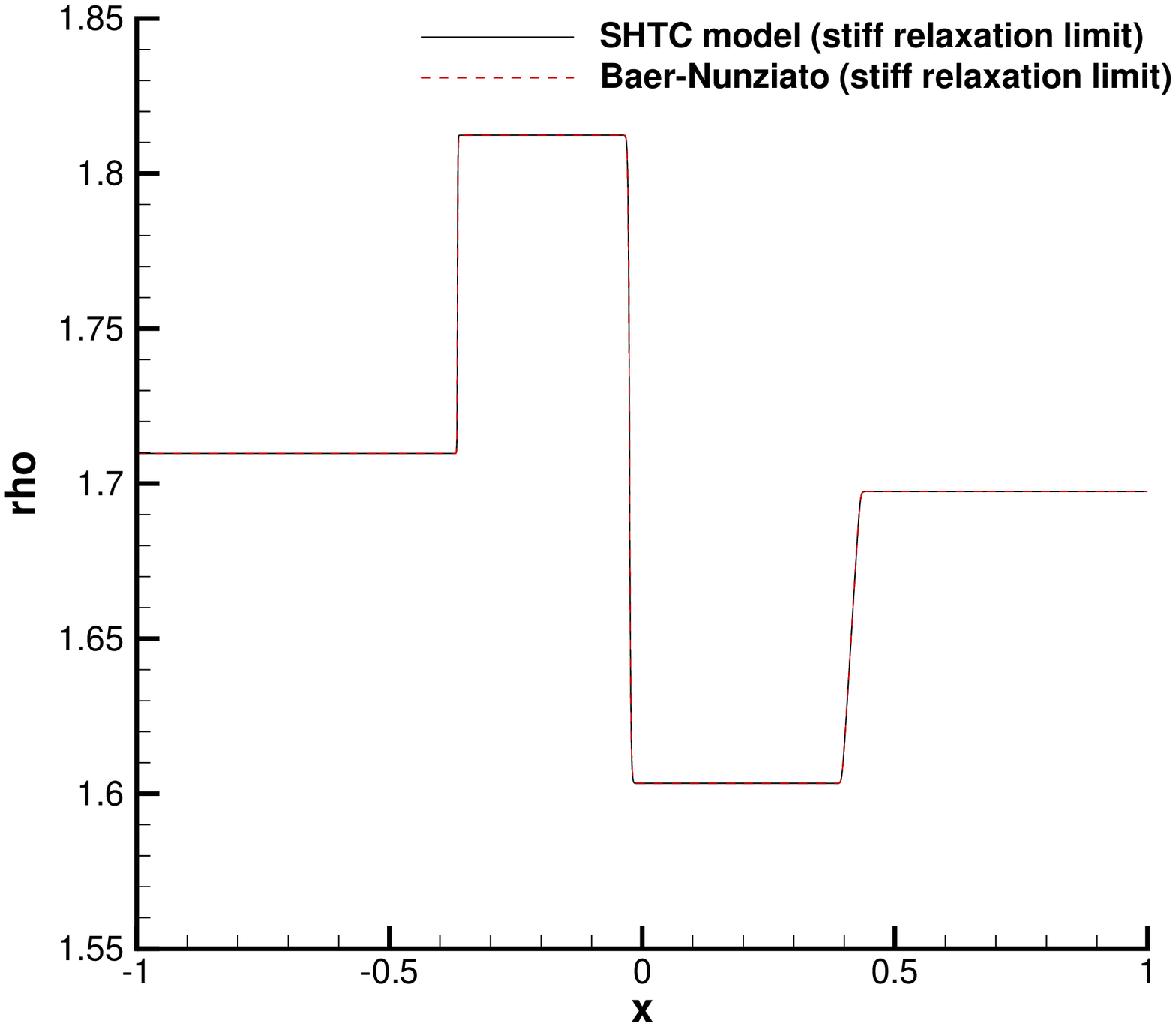}  & 
			\includegraphics[width=0.4\textwidth]{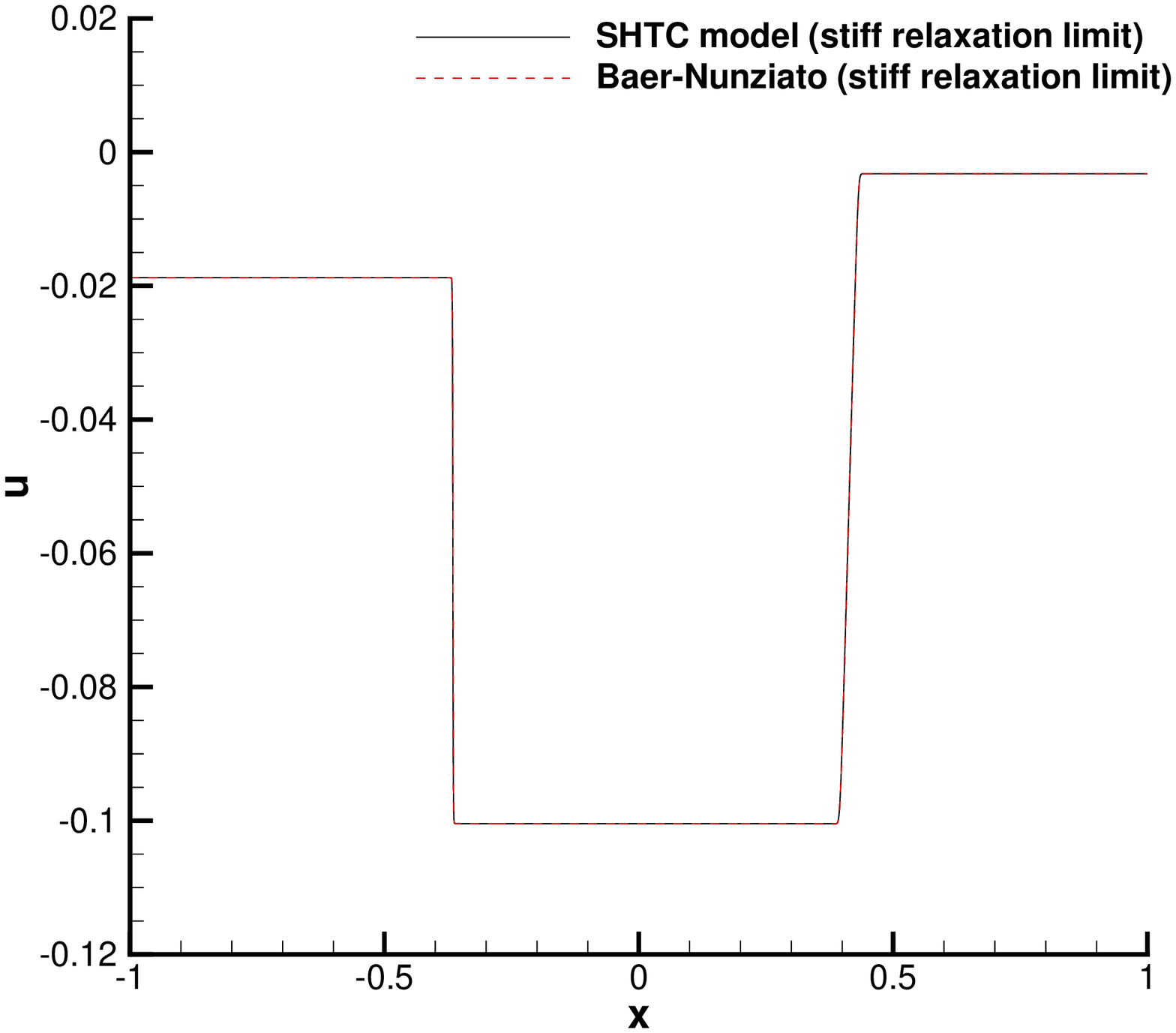}  \\ 	
			\includegraphics[width=0.4\textwidth]{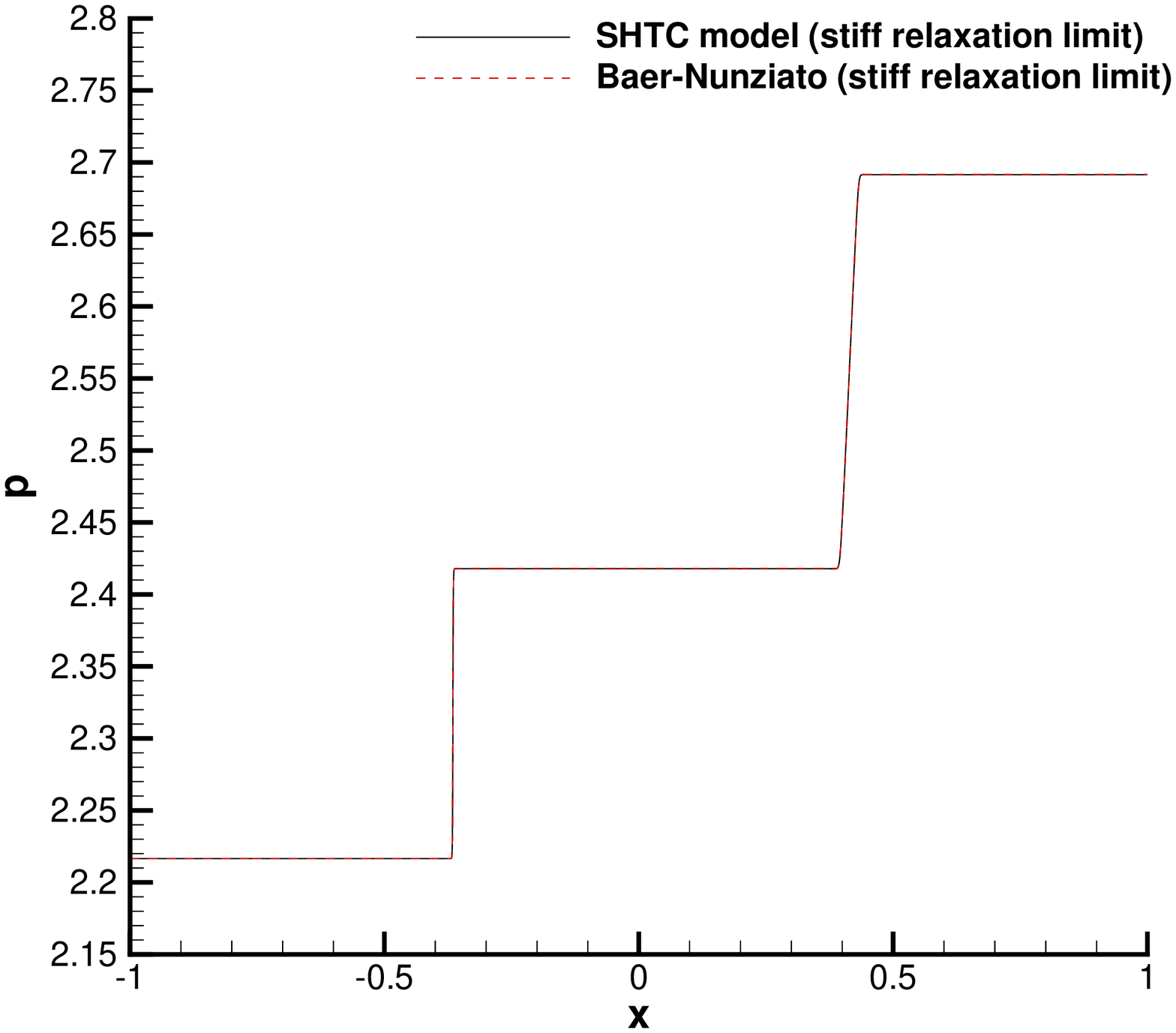} &
			\includegraphics[width=0.4\textwidth]{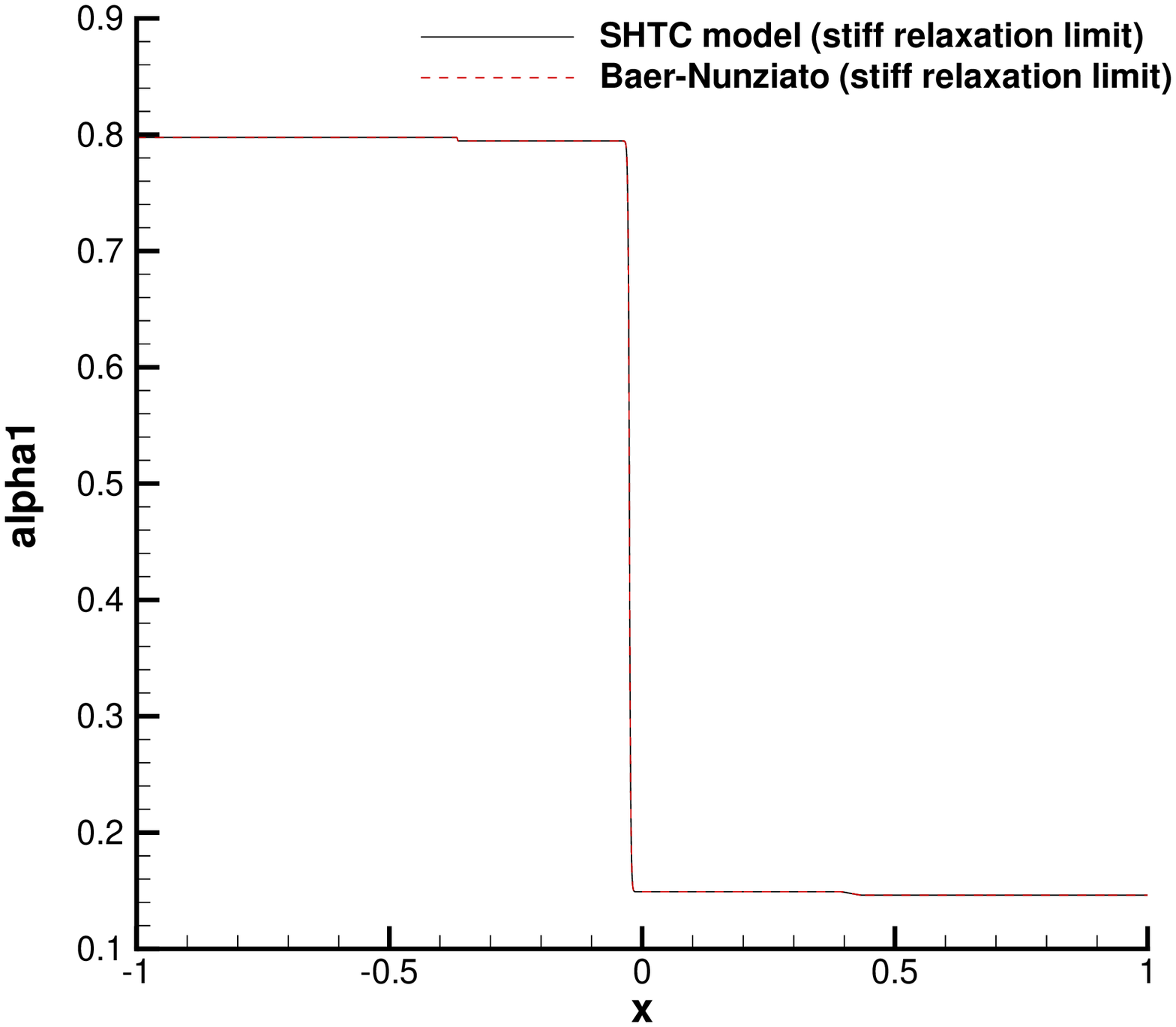}       
		\end{tabular}				
		\caption{Numerical solution of the SHTC system (black) and of the Baer-Nunziato model (red) of Riemann problem RP6b \textit{with} stiff relaxation source terms. From top left to bottom right: mixture density $\rho$, mixture velocity $u$, mixture pressure $p$ and volume fraction $\alpha_1$.}
		\label{fig:rp6b}
	\end{center}
\end{figure}

The last test problem RP6 has the following initial data: $\rho_{1,L} = 2$, $\rho_{1,R} = 1$, $\rho_{2,L} = 1$, $\rho_{2,R} = 2$, $u_{1,L}=u_{1,R}=0$, $u_{2,L}=u_{2,R}=0$, $\alpha_L = 0.7$, $\alpha_R=0.3$. The computational domain is again the interval $[-1,1]$, which is discretized with 10000 uniform grid cells and the final time of the simulation is $t=0.25$. In Figure \ref{fig:rp6} we compare the exact solution of the Riemann problem for the SHTC system with the numerical solutions obtained for the SHTC system and the Baer-Nunziato model without any relaxation source terms. Since the wave structure consists of a shock and a rarefaction wave in both phases, the numerical results of the SHTC model and the numerical results obtained for the Baer-Nunziato system do \textit{not} agree with each other any more, since the jump conditions of the non-conservative Baer-Nunziato system and the conservative SHTC model do not coincide. In particular, one can observe how the shock waves in one phase leave the other phase unchanged in the Baer-Nunziato model, while in the conservative SHTC system a shock in one phase always affects the other phase, which is a consequence of the jump conditions. In Figure \ref{fig:rp6b} we show the results obtained for the same initial data, but with stiff pressure and velocity relaxation, choosing the relaxation parameters as $\theta_1 = 10^{-3}$ and $\theta_2=10^{-8}$, which corresponds again to the Kapila limit of both systems. 
Despite the visible discrepancies observed in the homogeneous case, we find that the numerical solutions obtained in the stiff relaxation limit of the conservative SHTC model and the Baer-Nunziato model agree perfectly well with each other. This indicates that for numerical purposes it may be more beneficial to solve the SHTC system in the stiff relaxation limit, since standard numerical schemes for conservation laws can be applied, while in the Baer-Nunziato model special numerical techniques for the treatment of nonconservative terms are needed.  

\FloatBarrier

    \section{Conclusion}\label{sec:conclusion}
~
In this paper we have presented the exact solution of the Riemann problem for the  barotropic conservative two-phase model proposed in \cite{Romenski2004,Romenski2007,Romenski2009}, which belongs to the SHTC class of symmetric hyperbolic and thermodynamically compatible systems. We have discussed the characteristic fields of the system and the admissibility criteria of shock waves. The Riemann invariants as well as the Rankine-Hugoniot jump conditions have been provided. 
 A particular feature of the model under consideration is that some eigenvalues of the two phases may coincide, creating particular wave phenomena for the case where a shock in one phase interacts with a rarefaction wave in the other phase.  
 We have discussed possible wave patterns and based on the mathematical entropy inequality associated with the system and the admissibility criteria of shock waves we have ruled out non-admissible wave configurations.
 
Independent of the potential physical applications of the model discussed in this paper, its mathematical structure and properties as well as the possible wave configurations that may appear in the solution of the Riemann problem make it a very interesting object of study from a mathematical point of view.  
 
 At the end of the paper we show exact solutions of some example Riemann problems, providing also a detailed comparison with numerical results obtained with a second order TVD scheme. For two problems, we also provide a comparison with numerical solutions obtained for the barotropic Baer-Nunziato model, both, for the homogeneous case without source terms and for the case of stiff pressure and velocity relaxation. 
 In the case of smooth solutions, the results of the SHTC model and the Baer-Nunziato model agree perfectly well with each other, as expected, while in the presence of shock waves the numerical solutions disagree in the homogeneous case. On the contrary, in the presence of stiff velocity and pressure relaxation source terms, the numerical results obtained in the relaxation limit of the SHTC system and of the Baer-Nunziato model (Kapila limit) agree perfectly well with each other. 
 
 Future work will consider the extension of the Riemann solver presented in this paper to the non-barotropic case, solving the full model \cite{Romenski2007,Romenski2009} and as also shown in the appendix. 
 We furthermore plan to develop provably thermodynamically compatible finite volume schemes for the SHTC system, following the ideas outlined in \cite{HTCMHD,SWETurbulence,HTCGPR}.


    \appendix
    \section{Appendix}\label{sec:app}
\subsection{Derivation of the Barotropic Submodel}
%
A PDE system for compressible two-phase, two-temperature flow using a hyperbolic heat conduction model, was previously discussed in Romenski et al.\ \cite{Romenski2007,Romenski2009}.
Note, that there is another version of the SHTC model for hyperbolic heat transfer originally proposed in \cite{Romenski1998} and used in \cite{Dumbser2016,DPRZ2017,GRGPR,Boscheri2021}. 
The mentioned version seems to be more suitable from the point of view of the general theory of SHTC models, since its PDE system can be derived from the variational principle \cite{Peshkov2018}.
However, since this heat transfer model is not implemented in the SHTC two-phase flow model, we use the equations formulated in \cite{Romenski2007}.
Moreover, the conclusions drawn in this section should be the same for both models, since the difference in the equations for the thermal impulse does not affect the final result.
Written in terms of the generalized energy the system under consideration reads as
\begin{subequations}\label{eqn.HPRFF_app}
    \begin{align}
        \frac{\partial \rho \alpha_1}{\partial t} + \frac{\partial \rho \alpha_1 u^k }{\partial x_k} &= -\phi,\label{eqn.alphaFF_app}\\
        \frac{\partial \rho c_1}{\partial t} + \frac{\partial (\rho c_1 u^k+\rho E_{w_k})}{\partial x_k} &=  -\psi,\label{eqn.contiFF1_app}\\
        \frac{\partial \rho}{\partial t} + \frac{\partial \rho u^k}{\partial x_k} &= 0,\label{eqn.contiFF_app}\\
        \frac{\partial \rho u^i}{\partial t} + \frac{\partial (\rho u^i u^k + p \delta_{ik} + \rho w^iE_{w^k} )}{\partial x_k} & = 0,\label{eqn.momentumFF_app}\\
        \frac{\partial w^k}{\partial t} + \frac{\partial(w^lu^l+E_{c_1})}{\partial x_k} + u^l\left(\frac{\partial w^k}{\partial x_l} - \frac{\partial w^l}{\partial x_k}\right)
        &= -\frac{1}{\rho}\lambda_0.\label{eqn.relvel_app}\\
        \frac{\partial \rho j^i_1}{\partial t} + \frac{\partial (\rho j_1^i u^k + E_{S_1}\delta_{ik})}{\partial x_k} & = -\lambda_1^i,\label{eqn.thermal_impulse1}\\
        \frac{\partial \rho j^i_2}{\partial t} + \frac{\partial (\rho j_2^i u^k + E_{S_2}\delta_{ik})}{\partial x_k} & = -\lambda_2^i,\label{eqn.thermal_impulse2}\\
        \frac{\partial \rho S_1}{\partial t} + \frac{\partial (\rho S_1 u^k + E_{j_1^k})}{\partial x_k} & = \Pi_1 - \pi_1,\label{eqn.entropy1}\\
        \frac{\partial \rho S_2}{\partial t} + \frac{\partial (\rho S_2 u^k + E_{j_2^k})}{\partial x_k} & = \Pi_2 - \pi_2,\label{eqn.entropy2}
    \end{align}
\end{subequations}
Here, $\alpha_1$ is the volume fraction of the first phase which is connected with the volume fraction of the second phase $\alpha_2$ by the saturation law $\alpha_1+\alpha_2 = 1$,
$\rho$ is the mixture mass density which is connected with the phase mass densities $\rho_1,\rho_2$ by the relation $\rho = \alpha_1\rho_1 + \alpha_2\rho_2$.
The phase mass fractions are defined as $c_1=\alpha_1 \rho_1/\rho,\, c_2=\alpha_2 \rho_2/\rho$ and it is easy to see that $c_1 + c_2 = 1$.
Eventually, $u^i = c_1u_1^i + c_2u_2^i$ is the mixture velocity, $w^i=u_1^i - u_2^i$ is the phase relative velocity.
The quantities $j_i^k$ is the thermal impulse of phase $i$ and $S_i = c_is_i$ is the partial entropy of phase $i$.
The equations describe the balance law for the volume fraction, the balance law for the mass fraction, the conservation of total mass, the total momentum conservation law,
the balance for the relative velocity, the two equations for the partial thermal impulses and the tow equations for the partial entropies.
Moreover by construction the system is equipped with a conservation law for the total energy
\begin{align}
    \dfrac{\partial \rho\left(E + \frac{1}{2}u_lu^l\right)}{\partial t}
    + \dfrac{\partial \left(\rho u^k\left(E + \frac{1}{2}u_lu^l + \frac{p}{\rho} + w^lE_{w^l}\right) + \rho E_cE_{w_k} + E_{j_i^k}E_{S_i}\right)}{\partial x_k} = 0.\label{eqn:tot_energy}
\end{align}
In this work the source terms are of no special interest. The details can be found in \cite{Romenski2007,Romenski2009}.
We want to present a derivation of the barotropic submodel in order to emphasize the differences in the conservation law for the total energy since this will affect the mathematical entropy for the
system given above. By barotropic we understand that the internal energy of each phase solely depends on the corresponding phase density.
More precisely we consider either an isentropic (constant entropy) or an isothermal (constant temperature) process.
Let us first consider an isentropic process. This implies that the derivatives $E_{S_i}$ vanish which leads to the following conservation law for the total energy
\begin{align}
    \dfrac{\partial \rho\left(E + \frac{1}{2}u_lu^l\right)}{\partial t}
    + \dfrac{\partial \left(\rho u^k\left(E + \frac{1}{2}u_lu^l + \frac{p}{\rho} + w^lE_{w^l}\right) + \rho E_cE_{w_k}\right)}{\partial x_k} = 0.\label{eqn:tot_energy_isentropic}
\end{align}
The isothermal case is a bit more involved. The reason is, that for an isothermal process we necessarily need the thermal impulse.
We first rewrite the equation for the total energy (\ref{eqn:tot_energy})
\begin{align*}
    &\dfrac{\partial \rho\left(E + \frac{1}{2}u_lu^l\right)}{\partial t}
    + \dfrac{\partial \left(\rho u^k\left(E + \frac{1}{2}u_lu^l + \frac{p}{\rho} + w^lE_{w^l}\right) + \rho E_cE_{w_k}\right)}{\partial x_k} + \frac{\partial E_{j_i^k}E_{S_i}}{\partial x_k}\\
    = &\dfrac{\partial \rho\left(E + \frac{1}{2}u_lu^l\right)}{\partial t}
    + \dfrac{\partial \left(\rho u^k\left(E + \frac{1}{2}u_lu^l + \frac{p}{\rho} + w^lE_{w^l}\right) + \rho E_cE_{w_k}\right)}{\partial x_k}
    + E_{S_i}\frac{\partial E_{j_i^k}}{\partial x_k}\\
    + &E_{j_i^k}\frac{\partial E_{S_i}}{\partial x_k} = 0
\end{align*}
Since we have an isothermal process the partial space derivative of $E_{S_i} = T_i$ vanishes identically.
We further reformulate the partial entropy balances (\ref{eqn.entropy1}) and (\ref{eqn.entropy2}), i.e.
\begin{align*}
    \frac{\partial E_{j_1^k}}{\partial x_k} & = \Pi_1 - \pi_1 - \frac{\partial \rho S_1}{\partial t} - \frac{\partial (\rho S_1 u^k)}{\partial x_k},\\
    \frac{\partial E_{j_2^k}}{\partial x_k} & = \Pi_2 - \pi_2 - \frac{\partial \rho S_2}{\partial t} - \frac{\partial (\rho S_2 u^k)}{\partial x_k}
\end{align*}
and insert these two equations into the equation for the total energy. From now on we also assume that the flow has a single temperature $T = T_1 = T_2$ and thus the sources $\pi_i$ vanish.
This gives
\begin{align*}
    &\dfrac{\partial \rho\left(E + \frac{1}{2}u_lu^l\right)}{\partial t}
    + \dfrac{\partial \left(\rho u^k\left(E + \frac{1}{2}u_lu^l + \frac{p}{\rho} + w^lE_{w^l}\right) + \rho E_cE_{w_k}\right)}{\partial x_k}
    + E_{S_i}\frac{\partial E_{j_i^k}}{\partial x_k}\\
    = &\dfrac{\partial \rho\left(E + \frac{1}{2}u_lu^l\right)}{\partial t}
    + \dfrac{\partial \left(\rho u^k\left(E + \frac{1}{2}u_lu^l + \frac{p}{\rho} + w^lE_{w^l}\right) + \rho E_cE_{w_k}\right)}{\partial x_k}\\
    + &T\sum_{i=1}^2\left[\Pi_i - \frac{\partial \rho S_i}{\partial t} - \frac{\partial (\rho S_i u^k)}{\partial x_k}\right] = 0.\label{eqn:tot_energy_isothermal}
\end{align*}
Before we further simplify the energy equations for the isentropic and the isothermal case we need to precisely discuss the mixture equation of state in each case.
In the following we will discuss the isentropic and isothermal case separately.
\subsubsection{Isentropic Case}
Using relation (\ref{e_deriv:iso_S}) we obtain the following derivatives for the mixture internal energy (\ref{eq:mix_int_e})
\begin{align}
    \frac{\partial e}{\partial\alpha} &= c\left(\frac{\partial e_1}{\partial\rho_1}\right)_s\left(-\frac{c\rho}{\alpha^2}\right)
                                      + (1 - c)\left(\frac{\partial e_2}{\partial\rho_2}\right)_s\frac{(1 - c)\rho}{(1 - \alpha)^2}\notag\\
                                      &= \frac{(1 - c)^2\rho}{(1 - \alpha)^2}\frac{p_2}{\rho_2^2} - \frac{c^2\rho}{\alpha^2}\frac{p_1}{\rho_1^2} = \frac{p_2 - p_1}{\rho},\\
    \frac{\partial e}{\partial c} &= e_1 + c\left(\frac{\partial e_1}{\partial\rho_1}\right)_s\frac{\rho}{\alpha}
                                  - e_2 - (1 - c)\left(\frac{\partial e_2}{\partial\rho_2}\right)_s\frac{\rho}{1 - \alpha}\notag\\
                                  &= e_1 + \frac{c\rho}{\alpha}\frac{p_1}{\rho_1^2} - \left(e_2 + \frac{(1 - c)\rho}{1 - \alpha}\frac{p_2}{\rho_2^2}\right) = h_1(\rho_1) - h_2(\rho_2),\\
    \frac{\partial e}{\partial\rho} &= c\left(\frac{\partial e_1}{\partial\rho_1}\right)_s\frac{c}{\alpha}
                                    + (1 - c)\left(\frac{\partial e_2}{\partial\rho_2}\right)_s\frac{1 - c}{1 - \alpha}\notag\\
                                    &= \frac{c^2}{\alpha}\frac{p_1}{\rho_1^2} + \frac{(1 - c)^2}{1 - \alpha}\frac{p_2}{\rho_2^2} = \frac{\alpha p_1 + (1 - \alpha)p_2}{\rho^2} = \frac{p}{\rho^2}.
\end{align}
Here we have introduced the \emph{specific enthalpy} of phase $i$
\begin{align*}
    h_i(\rho_i) = e_i(\rho_i) + \frac{p_i}{\rho_i}
\end{align*}
and the \emph{mixture pressure}
\begin{align*}
    p = \alpha_1p_1 + \alpha_2p_2.
\end{align*}
Note that the derivatives with respect to $\alpha$ and $\rho$ are equal to the corresponding derivatives in the general case, see \cite{Romenski2007}.
However, the derivative with respect to $c$ is different (enthalpy vs.\ chemical potential). Using these results we obtain the following for the mixture EOS
\begin{align}
    \frac{\partial E}{\partial\alpha} &= \frac{\partial e}{\partial\alpha} = \frac{p_2 - p_1}{\rho},\label{isoS_E_deriv_alpha}\\
    \frac{\partial E}{\partial c} &= \frac{\partial e}{\partial c} + (1 - 2c)\frac{w_iw^i}{2} = h_1(\rho_1) - h_2(\rho_2) + (1 - 2c)\frac{w_iw^i}{2},\label{isoS_E_deriv_c}\\
    \frac{\partial E}{\partial\rho} &= \frac{\partial e}{\partial\rho} = \frac{p}{\rho^2},\label{isoS_E_deriv_rho}\\
    \frac{\partial E}{\partial w_i} &= c(1 - c)w_i.\label{isoS_E_deriv_w}
\end{align}
\subsubsection{Review Isothermal Thermodynamics}
Due to the extra term with the derivative of the entropy in (\ref{e_deriv:iso_T}) we need to revisit some thermodynamic relations before we calculate the derivatives.
First, we have the following well known Maxwell relation (omitting the phase index for the following considerations)
\begin{align}
    \left(\frac{\partial p}{\partial T}\right)_\rho = -\rho^2\left(\frac{\partial s}{\partial \rho}\right)_T\label{maxwell_rel:p_s}
\end{align}
which bases on the equality of second order derivatives of the internal energy depending on the canonical variables $\rho$ and $s$, exemplary see \cite{Landau1987}.
This relation then can be used to derive a more complicated Maxwell relation which is also known but not that obvious. More precisely we want to show the following relation
\begin{align}
    \left(\frac{\partial (\rho e)}{\partial \rho}\right)_T = -T^2\left(\dfrac{\partial\left(\dfrac{g}{T}\right)}{\partial T}\right)_\rho.\label{maxwell_rel:e_g}
\end{align}
This is also a well established relation, see \cite{Mueller1985, Landau1987}, but for the convenience of the reader we show how it can be derived.
Here $g$ denotes the \emph{specific Gibbs energy} which may be generalized to the \emph{chemical potential} $\mu$ when more substances are involved, see \cite{Landau1987}.
The Gibbs energy (or sometimes free enthalpy) is given by
\begin{align*}
    g = e - Ts + \frac{p}{\rho}.
\end{align*}
The differentials for the internal energy and the Gibbs energy depending on the variables $\rho$ and $T$ are given by
\begin{align}
    \dd e &= T\dd s + \frac{p}{\rho^2}\dd\rho = \left(T\left(\frac{\partial s}{\partial\rho}\right)_T + \frac{p}{\rho^2}\right)\dd \rho + T\left(\frac{\partial s}{\partial\rho}\right)_T\dd T,
    \label{differntial_e}\\
    \dd g &= \frac{1}{\rho}\dd p - s\dd T = \frac{1}{\rho}\left(\frac{\partial p}{\partial\rho}\right)_T\dd \rho + \frac{1}{\rho}\left(\left(\frac{\partial p}{\partial T}\right)_\rho - s\right)\dd T.
    \label{differential_g}
\end{align}
Now we can show the stated relation (\ref{maxwell_rel:e_g})
\begin{align*}
    \left(\frac{\partial (\rho e)}{\partial \rho}\right)_T &\stackrel{\phantom{(\ref{maxwell_rel:p_s})}}{=} e + \rho\left(\frac{\partial e}{\partial \rho}\right)_T
    \stackrel{(\ref{e_deriv:iso_T})}{=} e + \rho\left(\frac{p}{\rho^2} + T\left(\frac{\partial s}{\partial \rho}\right)_T\right)\\
    &\stackrel{(\ref{maxwell_rel:p_s})}{=} e + \frac{p}{\rho} - \frac{T}{\rho}\left(\frac{\partial p}{\partial T}\right)_\rho
    \stackrel{(\ref{def:gibbs_energy})}{=} g + Ts - \frac{T}{\rho}\left(\frac{\partial p}{\partial T}\right)_\rho \stackrel{(\ref{differential_g})}{=} g - T\left(\frac{\partial g}{\partial
    T}\right)_\rho\\
    &\stackrel{\phantom{(\ref{maxwell_rel:p_s})}}{=} -T^2\left(\dfrac{\partial\left(\dfrac{g}{T}\right)}{\partial T}\right)_\rho
\end{align*}
Since we are considering the isothermal case we obtain for this specific case
\begin{align*}
    -T^2\left(\dfrac{\partial\left(\dfrac{g}{T}\right)}{\partial T}\right)_\rho = g.
\end{align*}
\subsubsection{Isothermal Case}
Now we can calculate the derivatives of the mixture internal energy (\ref{eq:mix_int_e}) using (\ref{maxwell_rel:e_g})
\begin{align}
    \frac{\partial e}{\partial\alpha} &= c\left(\frac{\partial e_1}{\partial\rho_1}\right)_T\left(-\frac{c\rho}{\alpha^2}\right)
                                      + (1 - c)\left(\frac{\partial e_2}{\partial\rho_2}\right)_T\frac{(1 - c)\rho}{(1 - \alpha)^2}\notag\\
                                      &= -\frac{c}{\alpha}\rho_1\left(\frac{\partial e_1}{\partial\rho_1}\right)_T
                                      + \frac{1 - c}{1 - \alpha}\rho_2\left(\frac{\partial e_2}{\partial\rho_2}\right)_T\notag\\
                                      &= \frac{1 - c}{1 - \alpha}\left(\frac{\partial(\rho_2e_2)}{\partial\rho_2}\right)_T - \frac{c}{\alpha}\left(\frac{\partial(\rho_1e_1)}{\partial\rho_1}\right)_T
                                      + \frac{c}{\alpha}e_1 - \frac{1 - c}{1 - \alpha}e_2\notag\\
                                      &= \frac{1 - c}{1 - \alpha}(g_2 - e_2) - \frac{c}{\alpha}(g_1 - e_1) = \frac{p_2 - p_1}{\rho} + T(c_1s_1 - c_2s_2),\\
    \frac{\partial e}{\partial c} &= e_1 + c\left(\frac{\partial e_1}{\partial\rho_1}\right)_T\frac{\rho}{\alpha}
                                  - e_2 - (1 - c)\left(\frac{\partial e_2}{\partial\rho_2}\right)_T\frac{\rho}{1 - \alpha}\notag\\
                                  &= e_1 + \rho_1\left(\frac{\partial e_1}{\partial\rho_1}\right)_T - \left(e_2 + \rho_2\left(\frac{\partial e_2}{\partial\rho_2}\right)_T\right)\notag\\
                                  &= g_1(\rho_1) - g_2(\rho_2),\\
    \frac{\partial e}{\partial\rho} &= c\left(\frac{\partial e_1}{\partial\rho_1}\right)_T\frac{c}{\alpha}
                                    + (1 - c)\left(\frac{\partial e_2}{\partial\rho_2}\right)_T\frac{1 - c}{1 - \alpha}\notag\\
                                    &= \frac{c}{\rho}\rho_1\left(\frac{\partial e_1}{\partial\rho_1}\right)_T + \frac{1 - c}{\rho}\rho_2\left(\frac{\partial e_2}{\partial\rho_2}\right)_T\notag\\
                                    &= \frac{c}{\rho}\left(\frac{\partial(\rho_1 e_1)}{\partial\rho_1}\right)_T + \frac{1 - c}{\rho}\left(\frac{\partial(\rho_2 e_2)}{\partial\rho_2}\right)_T
                                    - \frac{ce_1 + (1 - c)e_2}{\rho}\notag\\
                                    &= \frac{cg_1 + (1 - c)g_2}{\rho} - \frac{ce_1 + (1 - c)e_2}{\rho}
                                    = \dfrac{c\left(\dfrac{p_1}{\rho_1} - Ts_1\right) + (1 - c)\left(\dfrac{p_2}{\rho_2} - Ts_2\right)}{\rho}\notag\\
                                    &= \frac{\alpha_1p_1 + \alpha_2p_2}{\rho^2} - \frac{T(c_1s_1 - c_2s_2)}{\rho} = \frac{p}{\rho^2} - \frac{T(c_1s_1 - c_2s_2)}{\rho}.
\end{align}
Note that in the general case we also obtain the difference of the Gibbs energies (or chemical potentials) for the derivative with respect to $c$, see \cite{Romenski2007}.
Using these results we obtain the following for the mixture EOS
\begin{align}
    \frac{\partial E}{\partial\alpha} &= \frac{\partial e}{\partial\alpha} = \frac{p_2 - p_1}{\rho} + T(c_1s_1 - c_2s_2),\label{isoT_E_deriv_alpha}\\
    \frac{\partial E}{\partial c} &= \frac{\partial e}{\partial c} + (1 - 2c)\frac{w_iw^i}{2} = g_1(\rho_1) - g_2(\rho_2) + (1 - 2c)\frac{w_iw^i}{2},\label{isoT_E_deriv_c}\\
    \frac{\partial E}{\partial\rho} &= \frac{\partial e}{\partial\rho} = \frac{p}{\rho^2} - \frac{T(c_1s_1 - c_2s_2)}{\rho},\label{isoT_E_deriv_rho}\\
    \frac{\partial E}{\partial w_i} &= c(1 - c)w_i.\label{isoT_E_deriv_w}
\end{align}
\subsubsection{Simplifying the Energy Equations}
Using the results obtained above we can further simplify the energy equations.
For the isentropic case we the relations \eqref{isoS_E_deriv_alpha} - \eqref{isoS_E_deriv_w} can be used to yield
\begin{align*}
    0 &= \dfrac{\partial \rho\left(E + \frac{1}{2}u_lu^l\right)}{\partial t}
    + \dfrac{\partial \left(\rho u^k\left(E + \frac{1}{2}u_lu^l + \frac{p}{\rho} + w^lE_{w^l}\right) + \rho E_cE_{w_k}\right)}{\partial x_k}\\
    &= \sum_{i=1}^2\frac{\partial\alpha_i\rho_i\left(e_i + \frac{1}{2}u_i^2\right)}{\partial t}
        + \frac{\partial\alpha_i\rho_i u^k_i\left(h_i + \frac{1}{2}u_i^2\right)}{\partial x_k}
\end{align*}
For the isothermal case we first rewrite the energy equality again using $S_i = c_is_i$ and $S = S_1 + S_2$, i.e.
\begin{align*}
    &\dfrac{\partial \rho\left(E + \frac{1}{2}u_lu^l\right)}{\partial t}
    + \dfrac{\partial \left(\rho u^k\left(E + \frac{1}{2}u_lu^l + \frac{p}{\rho} + w^lE_{w^l}\right) + \rho E_cE_{w_k}\right)}{\partial x_k}\\
    + &T\sum_{i=1}^2\left[\Pi_i - \frac{\partial \rho S_i}{\partial t} - \frac{\partial (\rho S_i u^k)}{\partial x_k}\right]\\
    \Leftrightarrow\;-&T\sum_{i=1}^2\Pi_i = \dfrac{\partial \rho\left(E + \frac{1}{2}u_lu^l\right)}{\partial t}
    + \dfrac{\partial \left(\rho u^k\left(E + \frac{1}{2}u_lu^l + \frac{p}{\rho} + w^lE_{w^l}\right) + \rho E_cE_{w_k}\right)}{\partial x_k}\\
    - &T\sum_{i=1}^2\left[\frac{\partial \rho S_i}{\partial t} + \frac{\partial (\rho S_i u^k)}{\partial x_k}\right]\\
    = &\dfrac{\partial \rho\left(E - TS + \frac{1}{2}u_lu^l\right)}{\partial t}
    + \dfrac{\partial \left(\rho u^k\left(E - TS + \frac{1}{2}u_lu^l + \frac{p}{\rho} + w^lE_{w^l}\right) + \rho E_cE_{w_k}\right)}{\partial x_k}.
\end{align*}
Here $F = E - TS$ is generalized mixture free energy and we can further write $\Pi = \Pi_1 + \Pi_2$ for the non-negative total entropy production.
Now we use relations \eqref{isoT_E_deriv_alpha} - \eqref{isoT_E_deriv_w} to obtain
\begin{align*}
    -&T\Pi = \dfrac{\partial \rho\left(E - TS + \frac{1}{2}u_lu^l\right)}{\partial t}
    + \dfrac{\partial \left(\rho u^k\left(E - TS + \frac{1}{2}u_lu^l + \frac{p}{\rho} + w^lE_{w^l}\right) + \rho E_cE_{w_k}\right)}{\partial x_k}\\
    = &\sum_{i=1}^2\frac{\partial\alpha_i\rho_i\left(e_i - Ts_i + \frac{1}{2}u_i^2\right)}{\partial t}
        + \frac{\partial\alpha_i\rho_i u^k_i\left(g_i + \frac{1}{2}u_i^2\right)}{\partial x_k}.
\end{align*}
%
In \cite{Romenski2009} a single temperature model is presented. The governing equations are the same as in the isothermal case, the energy equation reads
\begin{align*}
    \sum_{i=1}^2\left[\frac{\partial\alpha_i\rho_i\left(e_i + \frac{1}{2}u_i^2\right)}{\partial t} + \frac{\partial\alpha_i\rho_i u^k_i\left(h_i + \frac{1}{2}u_i^2\right)}{\partial x_k}\right]
    - \frac{\partial\rho c_1c_2(u^k_1 - u^k_2)(s_1 - s_2)T}{\partial x_k} = 0.
\end{align*}
However, after some manipulations this equation may also be written in the same form as it is used in this work. Thus the thermal impulses are again basically hidden in the energy equality.
\subsection{Calculation for Coinciding Eigenvectors}
We want to give the detailed calculation for the statement \eqref{cond:absent_contact}.
Let us first assume that the vector $\mathbf{R}_C$ is a multiple of a fixed vector $\mathbf{R}_{i\pm}$, see \eqref{sys_eigenvectors_prim}.
We can directly conclude that
\begin{align*}
    \varepsilon_i = 0\quad\Leftrightarrow\quad 0 = (u - u_i)^2 - a_i^2\quad\Leftrightarrow\quad u = u_i \pm a_i\quad\Leftrightarrow\quad \lambda_C = \lambda_{i\pm}.
\end{align*}
Note that the specific correct sign $\pm$ is determined by the non-zero components of the eigenvector.
To show the reverse direction we assume that $\lambda_C$ coincides with a fixed $\lambda_{i\pm}$, see \eqref{sys_eigenvalues}.
From this we get that the corresponding $\varepsilon_i$ vanishes. To simplify the notation let us assume w.l.o.g. that $i = 1$ and thus $\varepsilon_1 = 0$.
We further obtain
\begin{align*}
    \delta_1 = \frac{p_1 - p_2}{\rho} - \frac{(u - u_1)^2}{\alpha_1} = \frac{\alpha_1(p_1 - p_2) - \rho a_1^2}{\alpha_1\rho} = \rho_1\gamma_1.
\end{align*}
Thus we get
\begin{align*}
    \delta_1\varepsilon_2 = \rho_1\gamma_1\varepsilon_2 = \frac{\rho_1}{a_1}a_1\gamma_1\varepsilon_2 = \pm\frac{\rho_1}{a_1}(u - u_1)\gamma_1\varepsilon_2.
\end{align*}
This shows the desired relation.
\subsection{Regularity of the Mass Flux System Matrix}
In Section \ref{subsec:shock} we derived a linear system for the squares of the partial mass fluxes \eqref{lin_sys_Q_square}.
It remains to verify that the matrix is regular, i.e.\ the determinant
\[
  \det(\bb{M}) = -\frac{1}{2}\left(\alpha_1\dbl \dfrac{1}{\rho_1}\dbr\dbl\dfrac{1}{\rho_2^2}\dbr + \alpha_2\dbl\dfrac{1}{\rho_1^2}\dbr\dbl \dfrac{1}{\rho_2}\dbr\right)
\]
does not vanish. We can immediately verify that
\[
  \textup{sgn}\left(\dbl \dfrac{1}{\rho_1}\dbr\dbl\dfrac{1}{\rho_2^2}\dbr\right) = \textup{sgn}\left(\dbl\dfrac{1}{\rho_1^2}\dbr\dbl \dfrac{1}{\rho_2}\dbr\right)
\]
and thus the determinant will only vanish iff at least one of densities does not jump across the shock.
If the partial density would not jump across the shock we obtain from the continuity of the corresponding partial mass flux that the associated partial velocity also would not jump.
From the Lax condition of the eigenvalues we hence can conclude that for a shock in phase $\mu \in \{1,2\}$ the partial density of phase $\mu$ must jump.
Let us now consider the other phase $\nu \in \{1,2\},\;\nu \neq \mu$.
If we assume $\dbl\rho_\nu\dbr = 0$ we obtain for \eqref{jc:velocity}
\[
  \dbl u_\nu\dbr = -Q_\nu\dbl\frac{1}{\rho_\nu}\dbr = 0.
\]
Thus we have for the momentum jump condition \eqref{jc:momentum} and for the relative velocity jump condition \eqref{jc:w}
\begin{align}
    -Q_\mu\dbl u_\mu\dbr + \dbl p_\mu\dbr &= 0,\label{jc:momentum_deg}\\
        \dbl\frac{1}{2}(u_\mu - S)^2 + \Psi_\mu\dbr &= 0.\label{jc:w_deg}
\end{align}
Equation \eqref{jc:w_deg} can be rewritten using
\[
  \dbl\frac{1}{2}(u_\mu - S)^2\dbr = -\frac{1}{2}\dbl p_\mu\dbr\left(\frac{1}{\rho_\mu^+} + \frac{1}{\rho_\mu^-}\right)
\]
and thus
\begin{align}
    0 = \dbl\frac{1}{2}(u_\mu - S)^2 + \Psi_\mu\dbr = \dbl\Psi_\mu\dbr - \frac{1}{2}\dbl p_\mu\dbr\left(\frac{1}{\rho_\mu^+} + \frac{1}{\rho_\mu^-}\right).\label{jc:w_deg2}
\end{align}
Using the reformulated jump conditions \eqref{jc:momentum_v4} and \eqref{jc:w_v4} we have in the case $\dbl\rho_\nu\dbr = 0$
\[
  \frac{Q_\mu^2}{2}\dbl\frac{1}{\rho_\mu^2}\dbr - \dbl\Psi_\mu\dbr = 0,\quad Q_\mu^2\dbl\frac{1}{\rho_\mu}\dbr + \dbl p_\mu\dbr = 0.
\]
Since we know that $\rho_\mu$ jumps we can combine these two equations and obtain
\begin{align}
    0 &= 2\dfrac{\dbl\Psi_\mu\dbr}{\dbl\dfrac{1}{\rho_\mu^2}\dbr} + \dfrac{\dbl p_\mu\dbr}{\dbl\dfrac{1}{\rho_\mu}\dbr}\notag\\
    \Leftrightarrow\quad 0 &= 2\dbl\Psi_\mu\dbr\dbl\dfrac{1}{\rho_\mu}\dbr + \dbl p_\mu\dbr\dbl\dfrac{1}{\rho_\mu^2}\dbr\notag\\
    &= 2\left[\dbl\Psi_\mu\dbr + \frac{1}{2}\dbl p_\mu\dbr\left(\frac{1}{\rho_\mu^+} + \frac{1}{\rho_\mu^-}\right)\right]\dbl\dfrac{1}{\rho_\mu}\dbr\notag\\
    \Leftrightarrow\quad 0 &= \dbl\Psi_\mu\dbr + \frac{1}{2}\dbl p_\mu\dbr\left(\frac{1}{\rho_\mu^+} + \frac{1}{\rho_\mu^-}\right).\label{jc:w_deg3}
\end{align}
Summing up \eqref{jc:w_deg2} and \eqref{jc:w_deg3} gives
\begin{align*}
    %
    0 =  2\dbl \Psi_\mu\dbr\quad\Leftrightarrow\quad 0 = \dbl\rho_\mu\dbr.
\end{align*}
This is a contradiction since the partial density $\rho_\mu$ must jump as argued above.
Hence both partial densities jump across a shock and thus the determinant does not vanish.


    \section*{Conflict of Interest}
    The authors declare that they have no conflict of interest
    %
    %
	\section*{Acknowledgements}
	M.D. is member of the INdAM GNCS group and acknowledges the financial support received from the Italian Ministry of Education, University and Research (MIUR) in the frame of the PRIN 2017 project \textit{Innovative numerical methods for evolutionary partial differential equations and applications}. 
	E.R. was supported by the Mathematical Center in Akademgorodok under agreement No. 075-15-2019-1613 with the Ministry of Science and Higher Education of the Russian Federation.     
    
    \phantomsection
    \renewcommand{\refname}{References}
    \bibliographystyle{abbrv}
    \bibliography{gpr_literature}
\end{document}